\documentclass[11pt]{amsart}
\usepackage[foot]{amsaddr}
\usepackage[T1]{fontenc}
\usepackage[utf8]{inputenc}
\usepackage[english]{babel}
\usepackage[a4paper,margin=2.65cm]{geometry}
\usepackage{textcomp}
\usepackage{amsmath,amssymb,amsthm,mathtools,mathrsfs}
\usepackage{newtxtext,newtxmath}
\usepackage[bb=boondox,cal=boondoxo,scr=boondoxo]{mathalfa}
\usepackage[dvipsnames]{xcolor}
\usepackage[renew-dots,renew-matrix]{nicematrix}
\mathtoolsset{showonlyrefs}
\allowdisplaybreaks
\usepackage{enumitem}
\usepackage{microtype,booktabs,array,xurl,needspace}
\usepackage[numbers,sort&compress]{natbib}
\usepackage[hidelinks]{hyperref}
\usepackage{orcidlink}
\AtBeginDocument{\let\coloneq\coloneqq}
\theoremstyle{plain}
\newtheorem{theorem}{Theorem}[section]
\newtheorem{proposition}[theorem]{Proposition}
\newtheorem{corollary}[theorem]{Corollary}
\newtheorem{lemma}[theorem]{Lemma}
\newtheorem*{maintheorem}{Main construction}
\newtheoremstyle{definitionstyle}{6pt}{12pt}{\normalfont}{}{\bfseries}{.}{0.5em}{}
\theoremstyle{definitionstyle}

\newtheoremstyle{remarkstyle}{6pt}{12pt}{\normalfont}{}{\itshape}{.}{0.5em}{}
\theoremstyle{remarkstyle}
\newtheorem{remark}[theorem]{Remark}
\newcommand{\CC}{\mathbb C}
\newcommand{\NN}{\mathbb N}
\newcommand{\OO}{\mathcal O}

\newcommand{\II}{\mathcal I}
\newcommand{\QQ}{\mathcal Q}
\newcommand{\cA}{\mathcal A}
\newcommand{\cB}{\mathcal B}
\newcommand{\cF}{\mathcal F}
\newcommand{\T}{\mathsf T}
\newcommand{\dd}{\mathrm{d}}
\DeclareMathOperator{\Mat}{Mat}
\DeclareMathOperator{\rank}{rank}
\DeclareMathOperator{\im}{im}
\DeclareMathOperator{\coker}{coker}

\DeclareMathOperator{\Span}{span}
\DeclareMathOperator{\Tor}{Tor}
\DeclareMathOperator{\diag}{diag}
\numberwithin{equation}{section}
\hypersetup{
 pdftitle={Bilateral Christoffel--Geronimus--Uvarov formulas for mixed-type orthogonality on algebraic supports},
 pdfauthor={Manuel Manas},
 pdfsubject={Manuscript version 19, 31 October 2026}}
\title[Bilateral Christoffel formulas]{Bilateral Christoffel--Geronimus--Uvarov formulas for mixed-type orthogonality on algebraic supports}
\author[Manuel Ma\~nas]{Manuel Ma\~nas\,\orcidlink{0000-0003-3764-5737}}
\address{Department of Theoretical Physics, Faculty of Physical Sciences,
Complutense University of Madrid, 28040 Madrid, Spain}
\email{manuel.manas@ucm.es}
\date{Manuscript version 19 --- 31 October 2026}
\dedicatory{Dedicated to Professor Francisco Marcell\'an\\
on the occasion of his 75th birthday}

\begin{document}
\raggedbottom
\begin{abstract}
Christoffel-type determinant formulas construct both normalized polynomial
families and their reproducing kernel for mixed-type orthogonality in several
variables after independent rational matrix transformations on the two
arguments. The pairing may have a rectangular matrix of functionals, with
separate degree bounds on each component, and may be supported on an
algebraic variety. The reduced matrices are calculated from the original
biorthogonal polynomials, rational moments, polynomial division and the
masses left free by division; their nonvanishing determinants characterize
existence and uniqueness. For measure pairings the rational data are
expressed by Cauchy and multivariate Fantappi\`e transforms. Exact dimension
counts incorporate the four matrix factors and the relations of the
support. For normal degree-ordered original families and fixed factors, when the
component degrees on each side differ by a
bounded amount and the normalization dimension grows at most like one
degree layer, the resulting orders are bounded on curves and grow at
most linearly on surfaces. For pairings allowing scalar factors to move between arguments, called
balanced pairings, an exact obstruction determines whether a bilateral
mass decomposes into unilateral masses when the denominators act
injectively.
Residue duality provides coordinates for division and mass quotients on
integral curves, including singular points, and on the finite
surface quotients specified here. Explicit singular-curve families retain
nondegenerate residue pairings as denominator zeros coalesce. Positive
congruences admit a further reduction by orthogonal projection, while
balanced pairings give simultaneous intertwining identities for all
coordinate multiplication matrices.
\end{abstract}
\maketitle
\setcounter{tocdepth}{2}
\begingroup\small
\tableofcontents
\endgroup

\section{Introduction}

The objective is to calculate orthogonal polynomials after independent
rational matrix transformations on the two arguments of a pairing. The
calculation must preserve the prescribed degrees of each polynomial
component and the equations defining its algebraic support. The results
give the transformed polynomials, their reproducing kernel, the complete
freedom of added masses, and the sizes of the determinants required to
compute them.

In mixed-type multiple orthogonality, a polynomial row is required to be
orthogonal to prescribed polynomial columns. A typical pairing is
\[
 \beta_\mu(f,g)\coloneq\int_Z f(x)\,\dd\mu(x)\,g(x),
\]
where $Z$ is the specified integration support, $\mu$ is a $q$-by-$p$
matrix of measures, $f$ is a polynomial row with $q$ components, and $g$
is a polynomial column with $p$ components. The product under the integral
is a scalar. Each component of $f$ and $g$ is assigned its own degree
bound. To construct a row polynomial, the column space supplies the
orthogonality tests; to construct a column polynomial, the row space
supplies them. In either case, the unknown polynomial is required to pair
to zero with every element of the opposite test space, and its remaining
coefficients are fixed by a prescribed normalization. More generally, moments can be prescribed by linear
functionals instead of measures. The algebraic formulation by a matrix
of linear forms belongs to the vector orthogonality framework of
Sorokin and Van Iseghem~\cite{SorokinVanIseghem1997}; the mixed-type
orthogonality conditions extend types I and II as in
Daems and Kuijlaars~\cite{DaemsKuijlaars2007}.
This notation is bilinear, hence linear in both arguments. When a
sesquilinear pairing is used, the second argument is conjugated explicitly.

Polynomial coefficients and matrix entries may be complex even when
the integration variables are real. In the real-support examples,
the integration domain is specified separately from the complex
coordinate algebra used for the polynomial calculations. Added
evaluation or derivative functionals may instead be supported at
complex points; their support is stated when they are introduced.

The transformations of Christoffel and Geronimus multiply or divide
the measure by polynomial factors, and Uvarov modifications add masses
\cite{Christoffel1858,Geronimus1940,Uvarov1969}. The transformations considered here use independent rational matrix factors on
both sides of the pairing. The order of the matrix factors matters:
the coordinate algebra is commutative, but these matrices generally
do not commute. The row and column may have different
component counts and component degree bounds. The variables satisfy
the equations of an algebraic support of any dimension. Reducible supports,
algebraic multiplicities and pairings depending on two independent
support points are included.

We construct the perturbed mixed-type multiple orthogonal polynomials
and their \emph{reproducing kernel}, a finite sum of polynomial
column--row products in two sets of variables. Pairing in one set of
variables reproduces polynomials from the chosen space in the other.
Division leaves \emph{homogeneous masses}, terms vanishing after both
denominator actions. Our Christoffel--Geronimus--Uvarov problem is to
determine all these masses and compute the transformed polynomials
and kernels while retaining the prescribed component degree bounds and the
defining relations of the support.

More precisely, the data are two $q$-by-$q$ polynomial matrices
$L_0,L_1$ and two $p$-by-$p$ polynomial matrices $R_0,R_1$.
The new pairing $\widehat\beta$ is required to satisfy
\[
 \widehat\beta(fL_1,R_1g)=\beta_\mu(fL_0,R_0g)
 \qquad\text{for every polynomial row $f$ and column $g$}.
\]
The identity is formulated using multiplication. Under the injectivity
hypotheses stated below, the rational factors are
$S_L\coloneq L_1^{-1}L_0$ and $S_R\coloneq R_0R_1^{-1}$. Evaluating the original
pairing on $fS_L,S_Rg$ gives a particular solution; its rational moments
are additional input whenever they are not fixed by the original
polynomial moments. Every other solution adds a mass that vanishes on
$fL_1,R_1g$. Section~\ref{sec:algebra} proves these assertions and
also treats the consistency question when the matrices are not injective.

The calculation begins with finite polynomial spaces and an original
biorthogonal block on which the pairing matrix is invertible. Multiply
the target rows and columns by the rational factors and enlarge the
original spaces until all these products can be expressed in them.
Their coordinate matrices determine linear equations whose solutions
are exactly the admissible images. These equations impose both polynomial
divisibility and the required component degrees. Pair the added functions
with the original polynomials, insert the chosen mass, and impose the
normalization. Eliminating the coefficients already determined by the
original orthogonality gives the reduced matrix of the construction.
The scalar calculation in Subsection~\ref{sec:scalarresidual} carries
out every step before the general matrix formula is introduced.

The final step is block elimination. Its inputs and admissible
subspaces are supplied by the construction: the ordered rational
matrix actions, the complete degree and divisibility equations on the
algebraic support, and the mass coordinates. In particular, the
constraints can be computed by polynomial reduction followed by a nullspace
calculation. Values and derivatives, and residues for the finite quotients treated
below, give alternative coordinates for the divisibility conditions.
The degree or section-space conditions are imposed alongside them. Pairings with the added
rational functions play the role of functions of the second kind; for
measures, Subsubsection~\ref{sec:secondkind} gives their Cauchy and
Fantappi\`e integral representations.

The algebraic support affects this calculation through its polynomial
relations. Two expressions that differ by a defining relation represent
the same function, so they must not be counted as independent directions.
At singular points or on supports with multiplicities, these relations
also restrict which derivative data can occur. We retain them throughout
the construction. Subsubsection~\ref{sec:mixedcurveexample} calculates
all the blocks on a curve, with noncommuting matrix factors on
both sides and a mass that cannot be split into unilateral masses, and verifies the resulting
polynomial and kernel directly.

The main construction can be summarized as follows. A square pairing
is called normal when its matrix in bases of the chosen spaces is
invertible. A normalized orthogonal block consists of polynomial vectors
orthogonal to the opposite test space and satisfying specified linear
normalization conditions.

\begin{maintheorem}
Let $A=\CC[x_1,\ldots,x_M]/I$ be the coordinate algebra of the support.
Assume that the four polynomial matrices $L_0,L_1,R_0,R_1$ have regular
determinants, meaning that multiplication by each determinant is
injective on $A$. Fix a rational extension of the original pairing,
a homogeneous mass, and finite polynomial row and column spaces.
For a normalized problem, take $k>0$, a candidate space of dimension
$r+k$, an opposite test space of dimension $r$, and a prescribed
$\CC^k$-valued linear normalization. Choose an original normal block
and finite rational enlargements as in Section~\ref{part:finite}.

The original biorthogonal data determine reduced square matrices for
the right and left normalized problems. Each problem has a unique
solution if and only if its reduced determinant is nonzero. Every
entry of the solution is a bordered determinant divided by that reduced
determinant. For equal-dimensional target spaces, the same reduced
matrix that tests normality also gives their reproducing kernel.

The orders of these matrices are computed from dimensions of polynomial
images, complete constraint spaces and the restricted mass. For fixed
factors and original degree-ordered blocks normal at the levels considered,
take target degree-ordered flags or component bounds of bounded dispersion
on each side. If the normalization dimension is of order at most
$n^{s-1}$, the degree-based
construction has order $\mathrm{O}(n^{s-1})$ as $n\to\infty$, where
$s=\dim A\ge1$ and $n$ is the largest target component degree.
\end{maintheorem}

The determinant identities and solution formulas are proved in
Theorem~\ref{thm:taucriterion}, Corollary~\ref{cor:polynomialborders}
and the kernel construction \eqref{eq:kernelborder};
Theorem~\ref{thm:orders} and Proposition~\ref{prop:ordergrowth} give
the size statements. A tau determinant, denoted by $\tau$, is the
determinant of the reduced system. Bordering means adding one row and
one column to that matrix. The prescribed component degree bounds and
normalization are imposed directly, so the construction applies without
normality assumptions on intermediate transformations. On each finite
space, Proposition~\ref{prop:completejets} also supplies a separating
list of value and derivative conditions for division; their concrete
coordinates are computed in the examples and specializations.

The order of the determinant is part of the output.
Subsection~\ref{sec:orders} calculates the orders for specified
constructions, including all four matrix factors. Theorem~\ref{thm:orders}
gives an explicit degree count and a reduced count obtained from two
polynomial coefficient matrices. Corollary~\ref{cor:ordercurve} gives
the exact stable count on curves, including the row and column component counts, masses
and the position within the last degree block. The bounds of
Proposition~\ref{prop:ordergrowth} explain the change from bounded
orders on curves to at most linear growth on surfaces. The comparison with formulas using a fixed range of consecutive
polynomial degrees states which determinant representation is being
counted. It also explains how matrix division and component bounds
contribute information beyond the degree of a matrix determinant.

The mass equation supplies the remaining input to these matrices.
Theorem~\ref{thm:masses} describes every solution when numerator
and denominator actions are injective: no nonzero polynomial row or
column on the support is sent to zero.
Proposition~\ref{prop:compatibility} also decides existence when this
fails. Cancelling common rational factors can discard masses.
For pairings allowing polynomial factors to move between arguments,
Theorem~\ref{thm:tormasses} decides whether a bilateral mass splits
into left and right masses. Shared pieces of the denominator zero sets can prevent splitting.
At their intersections, a mass killed by each denominator separately
can be added to one summand and subtracted from the other without
changing the sum.

The geometric part of the article calculates particular inputs to the
reduced matrices.
A divisibility condition is a condition on a remainder class. A mass
also evaluates such classes, because it vanishes on the image of the
denominator action. Finite quotient coordinates therefore provide a common
method for both parts, with the quotient chosen for the matrix division
or denominator mass being calculated. These quotients can differ: a
Christoffel numerator imposes divisibility conditions even when the
denominators are identity matrices and there is no homogeneous mass.
Section~\ref{res:sec:residueconstruction} constructs dual coordinates
for the relevant quotients using residues. On an
integral curve, a geometric Cauchy kernel gives a matrix divided
difference that reproduces the remainder. At a singular point the
residue is interpreted using the module of dualizing differentials,
which retains the relations between the branches. On surfaces the
construction concerns finite quotients obtained by imposing a matrix
equation and an additional scalar equation. These coordinates supply
the linear tests and mass entries of the Christoffel-type determinant formulas.
The rational moments and the freely chosen mass parameters remain
part of the input; polynomial moments alone need not determine them.

These residue coordinates are then used in an explicit family.
The curve family in Section~\ref{sec:families} gives the factors,
residue functionals and finite Grams explicitly. Its normality is
verified on the stated component bounds by exploiting the triangular
matrix structure. The residue quotients remain perfectly paired when
denominator zeros reach a non-Gorenstein singular point. A constant
nonzero minor also proves that the matrix mass keeps its full rank,
so no mass direction is lost at the collision. Once this property and invertibility
of the finite Grams are established, the inverse Gram formulas give
regular limits of the normalized polynomials and kernels.

The same construction also yields the standard specializations.
On the line, Theorem~\ref{thm:univariate} and
Proposition~\ref{prop:jordanconstraints} express the division conditions
through principal parts, complete Jordan chains and forbidden
coefficients at infinity. The Cauchy data include shared poles and
their confluent limits. On Laurent supports, conjugation changes the
second argument and the degree bounds depend on the chosen embedding;
the finite construction retains both features. For symmetric
congruences, Proposition~\ref{prop:symmetricmass} describes the complete
mass quotient and classifies positive balanced masses on real curves.
The examples distinguish these evaluation masses from positive
bivariate pairings that couple distinct points.
Theorem~\ref{thm:positiveprojection} eliminates the common space by
orthogonal projection and gives the normalized polynomials and kernel
from two positive matrices of the same explicit order. Only the mass
on the overlap increases that order. Under the coprimality and degree
conditions of Corollary~\ref{cor:positiveorders}, every positive
homogeneous mass leaves the order unchanged.
Proposition~\ref{prop:positivenodalfamily} constructs a noncommutative
matrix family on a nodal curve, with a positive Sobolev mass involving
values and derivatives at the node and exact counts for unequal component
bounds. This derivative mass belongs to the positive, generally
nonbalanced theory. The coordinate identities below apply when both
pairings are balanced, as happens for the evaluation-mass specialization.

Once the polynomial bases have been constructed, they determine
coordinate multiplication matrices.
For balanced pairings on nested normal spaces that together contain
every polynomial vector, multiplication by a coordinate has a matrix
in the orthogonal polynomial bases. The four factors give connection
matrices between the original and transformed bases.
Theorem~\ref{thm:darboux} proves four-factor intertwining identities
simultaneously for all coordinates. The matrix identities preserve the equations
defining the support. For the stated positive congruences,
orthonormal bases make the multiplication matrices symmetric.

\subsection*{Relation to previous work and organization}

The scalar rational transformation problem has an established spectral
formulation in \cite{Zhedanov1997}. The algebraic theory of vector
orthogonality with a matrix of linear forms is developed by Sorokin and
Van Iseghem~\cite{SorokinVanIseghem1997}; the mixed-type formulation and
its Christoffel--Darboux formula are treated by Daems and
Kuijlaars~\cite{DaemsKuijlaars2007}. Matrix rational transformations
and bivariate pairings on the line are
developed in \cite{Matrix2016}, including noncommutative factors and
bilateral specializations; the matrix Christoffel case appears in
\cite{MatrixChristoffelIMRN}. Christoffel, Geronimus and Uvarov
transformations for mixed-type orthogonality on the step-line are treated in
\cite{MixedChristoffel,MixedGeronimus,MixedUvarov}. Kozhan and
Vaktn\"as~\cite[Theorems~3.1 and~3.3]{KV} give rational perturbation
formulas for types I and II on the line with general multi-indices and
componentwise scalar factors; their companion
paper~\cite{KVChristoffel} treats Christoffel transformations and their
relation to nearest-neighbour recurrences. Bueno and
Marcell\'an~\cite{BuenoMarcellan2006} formulate scalar polynomial
perturbations through separate factors in the two arguments of a
bilinear functional. Their detailed factorization results concern
one-sided modifications and the same modification in both arguments.
The construction here retains independent
matrix actions on both polynomial spaces over the coordinate algebra of
the support, with the prescribed component bounds.

For equal component bounds and a square positive pairing on the line,
Corollary~\ref{cor:matrixorthogonal} gives the monic matrix orthogonal
polynomials, norm matrices and Christoffel--Darboux kernel after a
rational congruence of the measure. This connects the construction
with the matrix weights of Dur\'an and Gr\"unbaum
\cite{DG2004,DGStructural,DG2007}. Subsection~\ref{sec:matrixorthogonal}
recovers their two Hermite matrix weights by polynomial triangular
factors, computes the reduced orders and displays a monic formula
in scalar Hermite polynomials. Their differential equations impose
additional conditions on the weights; the rational construction
supplies the polynomial and kernel transformation.

In several variables, Fern\'andez, P\'erez, Pi\~nar and
Xu~\cite{FPPX2010} give polynomial and reproducing-kernel formulas
for the addition of a Dirac mass to a positive definite functional.
The formulas of \cite{Granada2010} treat finitely many point-mass
additions. Christoffel modifications of total degree two
and Uvarov modifications are treated in \cite{Granada2016}; scalar
Geronimus transformations in two variables, their connection formulas
and block Jacobi factorizations appear in \cite{BivariateGeronimus2020}.
Multivariate rational formulas in \cite{AMRational} use scalar
polynomial factors. Related multivariate Christoffel constructions
appear in \cite{AMMultivariate,AMChristoffel}. Here the support
relations, the full bilateral mass equation and the finite component degree
bounds enter one determinant construction. Polynomial classes on
singular or nonreduced supports are included, subject to the stated
regularity and finite normality assumptions. Thus the extension is the
joint construction with independent matrix actions on the two polynomial
spaces over a quotient algebra, together with complete mass data and
determinant counts for those spaces.

Fasondini, Olver and Xu~\cite{FOXCubic,FOXGeneral} construct orthogonal
bases and computational recurrences on classes of planar algebraic
curves, including singular cubics. The present problem starts from an
original biorthogonal system and computes the systems produced by
independent rational matrix factors on its two arguments. The scalar
calculation in Subsection~\ref{sec:scalarresidual} and the full mixed-type
example in Subsubsection~\ref{sec:mixedcurveexample} illustrate the
construction of the finite matrices from the original data.

The geometric ingredients have their own established background.
Bertola~\cite{Bertola} relates scalar section biorthogonality on smooth
curves to matrix biorthogonality through a finite projection. His
Proposition~3.8 also gives rational left and right changes induced by
changes of geometric data. The curve residues below use the associated
Cauchy-kernel framework. On smooth complex curves, the geometric
Szeg\H{o} kernel has the classical theta-function description of
Fay~\cite{Fay} and the bundle formulation of
Ben-Zvi and Biswas~\cite{BenZviBiswas}. For smooth projective curves
over finite fields, Hu and Yu~\cite[Theorems~6.2 and~6.3]{HuYuResidue}
identify the finite residual Casimir tensor with the intrinsic
principal part of a normalized Szeg\H{o} kernel and prove its
reproduction property. This is the direct scalar, smooth antecedent
of the residue-quotient construction below. Here that construction
is carried out on integral curves with singular points, including
non-Gorenstein points, and for matrix quotients with the explicit
order of multiplication required by the bilateral formulas.
Rosenlicht's residue description of regular differentials,
presented in \cite[Chapter~IV]{SerreACF}, supplies the dualizing
module at those points; see also \cite{Lipman1984}.
The finite quotient
construction on surfaces is stated separately with its own hypotheses.
The recurrence interpretation is connected with the functional and
Darboux transformations of \cite{BuenoMarcellan2004} and with the
Hessenberg-matrix factorizations for polynomial perturbations of
bilinear functionals in \cite{BuenoMarcellan2006}. Here the same four
connection matrices intertwine every coordinate action and preserve
the polynomial relations of the support.

Section~\ref{part:masses} develops the functional equation and its
mass space. Section~\ref{part:finite} proves the determinant formulas
and explains their Gram realization. Section~\ref{res:sec:residueconstruction}
provides the curve and surface residue coordinates.
Section~\ref{sec:families} develops an explicit bilateral family
on singular curves and proves regularity of its normalized polynomials. Sections~\ref{part:classical}
and~\ref{sec:symmetric} treat standard, Laurent and positive
specializations. Section~\ref{sec:recurrence} proves the simultaneous
coordinate Darboux identities.

\section{Functional equation and homogeneous masses}\label{part:masses}

The functional equation prescribes the new pairing only on arguments
multiplied by the two denominator matrices. The first task is to decide
whether these prescribed values are consistent. The second is to describe
every value that remains free. These free values are the homogeneous
masses. On a curve their description reduces to finitely many polynomial
remainders, although a bivariate mass can still involve arbitrary
functionals in the opposite argument. The last part of the section
decides when a mass can be written as the sum of a left mass and a right
mass, each annihilated by just one denominator.

\subsection{The functional equation and all its masses}\label{sec:algebra}

We encode the support equations in the coordinate algebra $A\coloneq \CC[x_1,\ldots,x_M]/I$, where
$I$ is a proper ideal ($1\notin I$). Thus two polynomials represent the same element
of $A$ when their difference belongs to $I$; addition and multiplication
are performed subject to these relations. Write $Z^\circ$ for this
algebraic support, retaining multiplicities.\footnote{An ideal is closed
under addition and multiplication by arbitrary polynomials. It is prime
if $fg\in I$ implies $f\in I$ or $g\in I$. For prime $I$ the support
is integral: it cannot split into two proper algebraic subsets and its
algebra has no nonzero nilpotent classes. A class $u\ne0$ with $u^m=0$ is nilpotent; such classes retain
multiplicity data invisible in point values. See \cite[Section~0.1]{Eisenbud}.}
Here $\dim A$ is the dimension of the support, whereas vector-space
dimensions are over $\CC$. For instance, $\CC[x]$ has support dimension
one but infinitely many linearly independent powers of $x$.

The number $M$ counts coordinate variables. The row and column
sizes $q,p$ introduced below count polynomial components, not variables.

A \emph{scalar} means an element $a\in A$, possibly a nonconstant
polynomial, as opposed to a matrix over $A$. On a vector it multiplies
every component by $a$. It is \emph{regular}
if $ab=0$ forces $b=0$. Thus multiplication by a regular scalar is
injective, and cancellation by that scalar is valid.
Let $K$ consist of fractions with regular denominators. For integral $A$
this is the usual fraction field; on $A=\CC[x,y]/(xy)$, $x$ is not
regular because $xy=0$ with $y\ne0$, whereas $1+x+y$ is regular.

A linear functional is a linear rule assigning a complex number to each
vector; evaluations and moments are examples. We write
$W^*\coloneq \operatorname{Hom}_\CC(W,\CC)$ for the space of all such rules on $W$,
called its algebraic dual. Fix positive integers $q,p$ and write
$\Mat_{a\times b}(B)$ for $a$-by-$b$ matrices with entries in $B$,
with $\Mat_a(B)\coloneq \Mat_{a\times a}(B)$. Matrix entries have
comma-separated indices, as in $H=[H_{i,j}]$. The superscript $\T$ denotes
transpose; semicolons separate block rows in inline matrices.

We keep track of which argument a polynomial belongs to. For an arbitrary
bivariate pairing use $\mathscr A\coloneq A\otimes_\CC A$, the algebra of two
independent copies of the support. The maps $\iota_L$ and $\iota_R$
place a polynomial in the left or right copy:
$\iota_L(a)\coloneq a\otimes1$ and $\iota_R(a)\coloneq 1\otimes a$.
For $A=\CC[t]$, this means $\mathscr A=\CC[x,y]$,
$\iota_L(a)=a(x)$ and $\iota_R(a)=a(y)$, so $\iota_L(a)\iota_R(b)=a(x)b(y)$.
A \emph{balanced} pairing allows a polynomial factor to be moved from
one argument to the other. For this case one copy suffices:
$\mathscr A\coloneq A$ and $\iota_L(a)\coloneq \iota_R(a)\coloneq a$.
Both maps act entry by entry on matrices.

For a row $f\in A^{1\times q}$ and a column $g\in A^{p\times1}$,
the basic scalar tests are the products
$\iota_L(f_i)\iota_R(g_j)$. In the bivariate convention these are
$f_i(P)g_j(Q)$, with independent support points $P,Q$; in the balanced
convention they are products in the same algebra $A$.
Choose a scalar functional $\mu_{i,j}$ for each component pair and define
\[
 \beta_\mu(f,g)\coloneq
 \sum_{i=1}^{q}\sum_{j=1}^{p}
 \mu_{i,j}[\iota_L(f_i)\iota_R(g_j)].
\]
For a balanced matrix measure this is the familiar integral
$\int f\,\mathrm{d}\mu\,g$. The balanced identity
$\beta_\mu(af,g)=\beta_\mu(f,ag)$ follows because multiplication
of scalar entries in $A$ is commutative.
For $A=\CC[x]$, the matrix of linear forms in this balanced pairing
is the framework of Sorokin and Van Iseghem~\cite{SorokinVanIseghem1997}.
The coordinate algebra and the two-copy construction above allow
the same functional notation to express support relations and
pairings with independent arguments.

It is useful to collect the $qp$ tests in a matrix. Put
$V\coloneq\Mat_{q\times p}(\mathscr A)$ and write
\[
 H_{f,g}\coloneq
 [\iota_L(f_i)\iota_R(g_j)]_{i\in\{1,\ldots,q\},\,j\in\{1,\ldots,p\}}
 =\iota_L(f)^{\T}\iota_R(g)^{\T},\qquad
 \mu[H]\coloneq\sum_{i,j}\mu_{i,j}(H_{i,j}).
\]
Then $\beta_\mu(f,g)=\mu[H_{f,g}]$. A general test matrix $H\in V$
need not come from a single pair $f,g$. It is a finite sum of such
matrices: use matrix units for the component indices and, in the
bivariate case, express each tensor as a finite sum of elementary
tensors. Consequently, specifying the pairing specifies $\mu$ on all
of $V$. In the bivariate convention, any complex bilinear pairing is
represented in this way: take $\mu_{i,j}(a\otimes b)$ to be its value
on the row with component $a$ at $i$ and the column with component $b$
at $j$, and extend linearly. In the balanced convention the same rule
depends only on the product $ab$, because balance gives
$\beta(ae_i^{\T},be_j)=\beta(e_i^{\T},abe_j)$, where $e_i,e_j$
are coordinate columns of the appropriate sizes.
This explains precisely which pairings are represented by each algebra.

Fix $L_0,L_1\in\Mat_q(A)$ and $R_0,R_1\in\Mat_p(A)$.
We seek a new pairing whose value on $fL_1,R_1g$ equals the old
pairing on $fL_0,R_0g$, for every polynomial row $f$ and column $g$.
This prescription makes sense using multiplication alone, even when
some factors cannot be inverted. Subscripts $0$ and $1$ distinguish
the factors on the original and perturbed sides; under regularity
they become the numerator and denominator factors of rational ratios.

To express the prescription on arbitrary tests, first observe that
\[
 H_{fL_i,R_i g}
 =\iota_L(L_i)^{\T}H_{f,g}\iota_R(R_i)^{\T}.
\]
To verify the matrix order, fix a transformation index
$\epsilon\in\{0,1\}$. The $(a,b)$ entry of the transformed test is
\[
 \sum_{u,v}\iota_L((L_\epsilon)_{u,a})
       (H_{f,g})_{u,v}\iota_R((R_\epsilon)_{b,v}).
\]
The left matrix is therefore $\iota_L(L_\epsilon)^{\T}$ and the right
matrix is $\iota_R(R_\epsilon)^{\T}$. This expansion checks both
transposes directly from ordinary row and column multiplication. Extend this
transformation linearly to every $H\in V$ by setting
\[
 D_i(H)\coloneq\iota_L(L_i)^{\T}H\iota_R(R_i)^{\T},
 \qquad (L_i\mu R_i)[H]\coloneq\mu[D_i(H)].
\]
Thus multiplication of a functional means applying the corresponding
multiplication to its test before evaluation. Given $\mu$, the
functional equation for $\widehat\mu\in V^*$ is
\begin{equation}\label{eq:four}
                  L_1\widehat\mu R_1=L_0\mu R_0.
\end{equation}
Equivalently, the two pairings satisfy
\begin{equation}\label{eq:pairfour}
 \beta_{\widehat\mu}(fL_1,R_1g)=\beta_\mu(fL_0,R_0g)
 \quad\hbox{for every } f,g.
\end{equation}
\Needspace{5\baselineskip}
The next criterion allows factors that kill nonzero classes.
For a pairing $\beta$, the left radical consists of rows $f$ with
$\beta(f,g)=0$ for every column $g$; the right radical consists of
columns $g$ for which this holds for every row $f$.
Nondegeneracy means that these sets contain only the zero row and
zero column, respectively.

There is one possible obstruction to this prescription. Distinct tests
$H_1,H_2$ may give the same transformed test $D_1(H_1)=D_1(H_2)$.
They must then prescribe the same value of $\widehat\mu$. The next
proposition states that this consistency requirement is sufficient,
and identifies all the values that remain free.

\begin{proposition}[Existence for arbitrary factors]\label{prop:compatibility}
Equation~\eqref{eq:four} has a solution if and only if
\begin{equation}\label{eq:compatibility}
              \mu[D_0(H)]=0\quad\hbox{for every }H\in\ker D_1.
\end{equation}
All solutions are obtained from one by adding a functional on $V/D_1V$.
For $a,b\in A$ and a scalar balanced functional $u\in A^*$ with nondegenerate pairing,
the equation $bv=au$ has a solution if and only if
$\operatorname{Ann}(b)\subseteq\operatorname{Ann}(a)$,
where $\operatorname{Ann}(b)\coloneq \{f\in A:bf=0\}$.
The solution freedom is $(A/(b))^*$. If $v$ is also nondegenerate,
the two annihilators are equal.
\end{proposition}

\begin{proof}
Equation~\eqref{eq:four} asks for a linear functional $\widehat\mu$
satisfying $\widehat\mu[D_1(H)]=\mu[D_0(H)]$ for all $H\in V$.
If $D_1(H)=0$, the left-hand side is zero, proving necessity.

Conversely, assume \eqref{eq:compatibility}. If $D_1(H_1)=D_1(H_2)$,
then $H_1-H_2\in\ker D_1$, so $\mu[D_0(H_1)]=\mu[D_0(H_2)]$.
The prescribed value therefore depends only on $D_1(H)$. More
explicitly, the rule
\[
 \lambda(D_1(H))\coloneq\mu[D_0(H)]
\]
is well defined and linear on $D_1V$. Linearity follows by replacing
$H$ by $aH+bK$ for complex scalars $a,b$. Complete a basis of $D_1V$ to one of $V$
and choose arbitrary values on the added vectors. This extends the
functional to a solution $\widehat\mu$ on $V$.

The difference of two solutions vanishes on $D_1V$, so it gives the
same value to tests differing by an element of $D_1V$. It therefore
defines a functional on $V/D_1V$. Conversely, adding any such functional
preserves the equation.

In the scalar case, $bv=au$ means $v(bh)=u(ah)$ for every $h\in A$.
If $bf=0$, then $u((af)g)=v(bfg)=0$ for every $g\in A$.
Thus $af$ lies in the radical of the pairing of $u$; nondegeneracy
forces $af=0$. This proves $\operatorname{Ann}(b)\subseteq\operatorname{Ann}(a)$.
Conversely, this inclusion gives $bf=0$ implies $af=0$, hence
$u(af)=0$: precisely \eqref{eq:compatibility} for $D_1(h)=bh$
and $D_0(h)=ah$. The image of multiplication by $b$ is the ideal
$(b)$, so the freedom is $(A/(b))^*$.
Finally, if $v$ is also nondegenerate, $af=0$ implies
$v((bf)g)=u(afg)=0$ for every $g$, hence $bf=0$.
This gives the reverse inclusion and proves equality.
\end{proof}

A filtration is a nested sequence of finite-dimensional spaces; it is exhaustive
when their union is $A$. Invertibility of its Gram matrices at every
level (quasi-definiteness) implies this nondegeneracy. Indeed, a
nonzero element belongs to one finite level. The invertible Gram
matrix at that level supplies an opposite element with nonzero pairing.
It therefore cannot belong to either radical of the full pairing.
For pure division, $a=1$ gives $\operatorname{Ann}(a)=0$, so a solution
exists precisely when $b$ kills no nonzero class. Nonregular factors
are allowed if $\operatorname{Ann}(b)\subseteq\operatorname{Ann}(a)$;
equality is necessary if both pairings are nondegenerate.
Subsubsection~\ref{sec:polynomialelimination} uses fixed moments without
inverting the factors.

Assume each of the four determinants is regular: multiplication by it
kills no nonzero element of $A$. Adjugate formulas give inverses over
$K$. Put $S_L\coloneq L_1^{-1}L_0$ and $S_R\coloneq R_0R_1^{-1}$;
these send polynomial rows $f$ and columns $g$ to rational rows $fS_L$
and columns $S_Rg$.

Set $d\coloneq\iota_L(\det L_1)\iota_R(\det R_1)$. The notation
$\mathscr A[1/d]$ means the algebra of fractions $h/d^m$, with
$h\in\mathscr A$ and $m\in\{0,1,\ldots\}$; $m=0$ includes the original
elements.\footnote{This algebra is the localization at $d$. In the bivariate case,
multiplication by $d$ is the tensor product over $\CC$ of the two injective
determinant multiplications, hence is injective. Regularity makes the
original algebra inject into this localization. Inverting all regular elements of $A$ gives its total quotient ring $K$.
See \cite[Section~2.1]{Eisenbud}.}
For $\mathscr A=\CC[x]$ and $d=x-1$, this includes $(x+2)/(x-1)^2$.
Thus $\Mat_{q\times p}(\mathscr A[1/d])$ consists of matrices with these
fractional entries.

The functional $\mu$ acts only on polynomial tests. To evaluate rational
tests, extend a complex basis of $V$ to one of the larger matrix space
above, keeping the values prescribed by $\mu$ and choosing values on the
added basis elements. Linearity then determines an extension
$\widetilde\mu$ on the whole space. Denote its pairing by
$\widetilde\beta_\mu$, and write $\widehat\beta\coloneq\beta_{\widehat\mu}$
for the perturbed pairing. If $\mu$ is represented by a measure, integration
provides such an extension whenever every rational test is integrable.
The polynomial moments alone do not prescribe the new values: for
example, a value on $H/d$ must be supplied in addition unless that
fraction is already a polynomial test. A finite construction below only
uses finitely many such additional values.

The choice of extension fixes a particular solution, rather than the
set of solutions. If $\widetilde\mu_1,\widetilde\mu_2$ are two extensions,
their induced difference is
\[
 \nu_{12}(f,g)\coloneq
 (\widetilde\beta_{\mu,1}-\widetilde\beta_{\mu,2})(fS_L,S_Rg).
\]
On $fL_1,R_1g$ it becomes the difference of the two extensions on
the polynomial arguments $fL_0,R_0g$, hence is zero. Thus changing
the extension only changes which homogeneous mass is called zero.

\Needspace{4\baselineskip}
A homogeneous mass assigns arbitrary values to the classes of
$\QQ\coloneq V/D_1V$.
As a functional $\nu$ on $V$, it satisfies $\nu[D_1(H)]=0$;
we also write $\nu(f,g)$ as a pairing.\footnote{Finite generation over $\mathscr A$ means that every element is an
$\mathscr A$-linear combination of finitely many fixed generators.
For example, $\CC[x,y]/(x)$ is generated by $1$ over $\CC[x,y]$,
although the classes $1,y,y^2,\ldots$ are complex linearly independent.
See \cite[Section~0.3]{Eisenbud}.}
For a balanced scalar pairing on $\CC[t]$, take $q=p=1$ and
$L_1=R_1=t$. Then $D_1(h)=t^2h$, so the remainders are $a+bt$.
A mass is therefore a linear combination of evaluation and first derivative at zero.

When $\dim A=1$, put
$d_L\coloneq \dim_\CC\bigl(A^{1\times q}/(A^{1\times q}L_1)\bigr)$ and
$d_R\coloneq \dim_\CC\bigl(A^{p\times1}/(R_1A^{p\times1})\bigr)$.
Curves may be reducible or nonreduced.

The next theorem separates the construction into a particular solution
and the complete freedom left by division. Its first part gives the
particular solution by rational substitution. The second identifies
the masses without choosing coordinates. The last two parts describe
their size on curves: balanced masses have finitely many scalar
parameters, whereas a bivariate mass still has finite rank but can
depend on arbitrary functionals in the opposite argument.
Here finite rank means a representation by finitely many products of
one row functional and one column functional, as in
\eqref{eq:massfactor}.

\begin{theorem}[All rational bilateral perturbations]\label{thm:masses}
\textup{(i)} A functional $\widehat\mu$ solves \eqref{eq:four} if and only if, for some $\nu$,
\begin{equation}\label{eq:pairing}
 \widehat\beta(f,g)=\widetilde\beta_\mu(fS_L,S_Rg)+\nu(f,g),
 \qquad L_1\nu R_1=0.
\end{equation}

\textup{(ii)} The full mass space is $\QQ^*$; $\QQ$ is finitely generated
over $\mathscr A$ and annihilated by $d$.

\textup{(iii)} For balanced pairings with $\dim A=1$, $\QQ$ is finite
dimensional over $\CC$ and satisfies
$\dim\QQ=p\,d_L+q\,d_R$.

\textup{(iv)} For balanced or bivariate pairings with $\dim A=1$, each
mass has a representation with scalar linear functionals
$\ell_a:A^{1\times q}\to\CC$ and $r_a:A^{p\times1}\to\CC$:
\begin{equation}\label{eq:massfactor}
       \nu(f,g)=\sum_{a=1}^h\ell_a(f)r_a(g),
       \qquad h\le d_L+d_R.
\end{equation}
The bivariate mass space can be infinite dimensional even on a curve.
\end{theorem}

\begin{proof}
\textup{(i)} Apply the rational term in \eqref{eq:pairing} to
$fL_1,R_1g$. Since $L_1S_L=L_0$ and $S_RR_1=R_0$, its value is
$\widetilde\beta_\mu(fL_0,R_0g)=\beta_\mu(fL_0,R_0g)$: both arguments
are now polynomial, where the extension agrees with the original pairing.
Thus it solves \eqref{eq:pairfour}. Subtracting it from any solution
leaves a homogeneous mass, and adding any homogeneous mass gives another solution.

\textup{(ii)} Since $\nu[D_1(H)]=0$, tests differing by an element of
$D_1V$ have the same value. Thus $\nu$ defines a functional on $\QQ$.
Conversely, a functional on $\QQ$ gives a mass by evaluating quotient
classes. Hence the mass space is $\QQ^*$.
Since $V$ has $qp$ matrix-unit generators over $\mathscr A$, its quotient
$\QQ$ is finitely generated. Applying $D_1$ to
$\iota_L(\operatorname{adj}L_1)^{\T}H\iota_R(\operatorname{adj}R_1)^{\T}$
gives $dH$ by the adjugate identities. Hence $dV\subset D_1V$, so
multiplication by $d$ kills every quotient class.

This also locates the freedom geometrically: away from the zeros
of the denominator determinants, $d$ is invertible and the quotient
is zero. All homogeneous mass data therefore belong to the part
where at least one denominator ceases to be invertible.

\textup{(iii)} A regular element belongs to no minimal prime of $A$: the
minimal primes of a Noetherian ring are associated primes, and
an element in an associated prime is a zero divisor. Thus a
regular determinant does not vanish identically on any irreducible
curve component. Its zero set is zero dimensional. The quotient
of the finitely generated complex algebra $A$ by that determinant
is consequently Artinian and finite dimensional over $\CC$. Here
Artinian means that descending chains of ideals terminate; for an
affine complex algebra of dimension zero it is equivalent to having
only finitely many independent complex coefficients. Isolated
zero-dimensional components cause no exception: on a finite-dimensional
component, injective multiplication is surjective, so a regular scalar
is a unit there and its quotient is zero. The argument also retains
multiplicities and embedded point contributions.
It is the precise reason that only finitely many remainder
coefficients survive on a curve. Adjugates show
that the row and column remainder spaces are annihilated by their
respective determinants $\det L_1$ and $\det R_1$. These finitely generated remainder spaces therefore have finite dimensions $d_L,d_R$.
For balanced pairings, $\mathscr A=A$. Use the chain
$D_1V=L_1^{\T}VR_1^{\T}\subset L_1^{\T}V\subset V$.
The quotient $V/(L_1^{\T}V)$ has $p$ columns, each contributing $d_L$
remainders: transposition gives the row quotient defining $d_L$.
Multiplying $L_1^{\T}H=0$ by $\operatorname{adj}(L_1)^{\T}$ gives $(\det L_1)H=0$;
regularity forces $H=0$. This injective multiplication identifies
$(L_1^{\T}V)/D_1V$ with $V/(VR_1^{\T})$, which has $q$ row remainders
of dimension $d_R$. The two successive quotients fit into
\[
 0\longrightarrow V/(VR_1^{\T})
 \xrightarrow{\ [H]\mapsto[L_1^{\T}H]\ }\QQ
 \longrightarrow V/(L_1^{\T}V)\longrightarrow0.
\]
The first map is injective because $L_1^{\T}$ is injective. Its image
is exactly the kernel of the last quotient map. Choosing bases of that
kernel and of the final quotient gives a basis of $\QQ$, and therefore
$\dim_\CC\QQ=qd_R+pd_L$.

\textup{(iv)} For either convention on a curve, put $\mathscr F\coloneq A^{1\times q}$,
$\mathscr G\coloneq A^{p\times1}$, $\mathscr F_0\coloneq\mathscr FL_1$
and $\mathscr G_0\coloneq R_1\mathscr G$. Choose complements of dimensions $d_L,d_R$ to
$\mathscr F_0,\mathscr G_0$, with bases $\phi_a,\psi_b$ and remainder
coordinates $\ell_a,\eta_b$.
Write $f=f_0+r_L$, $g=g_0+r_R$, where
$f_0\in\mathscr F_0$, $g_0\in\mathscr G_0$,
$r_L\coloneq\sum_a\ell_a(f)\phi_a$ and $r_R\coloneq\sum_b\eta_b(g)\psi_b$;
put $\pi_{\mathscr F}f\coloneq r_L$.
Homogeneity gives $\nu(f_0,g_0)=0$, hence
\[
 \nu(f,g)=\nu(r_L,g)+\nu(f_0,r_R).
\]
The first term already includes the pairing of both remainders.
Define $T_a(g)\coloneq\nu(\phi_a,g)$ and
$S_b(f_0)\coloneq\nu(f_0,\psi_b)$. Expanding the two remainders gives
\begin{equation}\label{eq:bivariantmasses}
 \nu(f,g)=\sum_{a=1}^{d_L}\ell_a(f)T_a(g)
          +\sum_{b=1}^{d_R}S_b(f-\pi_{\mathscr F}f)\eta_b(g).
\end{equation}
The data are recovered uniquely by evaluating the mass on
$(\phi_a,g)$ and $(f_0,\psi_b)$, respectively. Arbitrary
$T_a\in\mathscr G^*$ and $S_b\in\mathscr F_0^*$ give all bivariate masses:
when $f\in\mathscr F_0$ and $g\in\mathscr G_0$, every term vanishes.
For balanced masses the data must also satisfy $\nu(af,g)=\nu(f,ag)$.
The same $d_L+d_R$ products prove \eqref{eq:massfactor} for either convention.
For $A=\CC[z]$ and $d_L+d_R>0$, pairing one nonzero remainder coordinate
with arbitrary opposite functionals gives an infinite-dimensional bivariate mass space.
\end{proof}

For balanced masses on a curve, at a smooth point with parameter $t$,
use local Smith normal forms.\footnote{A local parameter vanishes to order one. Invertible local row and
column operations give diagonal powers of $t$; their exponents
are the partial multiplicities. See \cite[Section~11.1]{Eisenbud}
and \cite{DopicoNoferini}.}
Write $L_1=U_L\diag(t^{\alpha_i})V_L$ and
$R_1=U_R\diag(t^{\beta_j})V_R$, where $U_L,V_L,U_R,V_R$ and their inverses have no poles at the point. These changes of basis reduce
the denominator action on $\eta\coloneq V_L\nu U_R$ to multiplication of
$\eta_{i,j}$ by $t^{\alpha_i+\beta_j}$. Thus that entry depends only on
Taylor coefficients $\mathsf c_k(a)\coloneq [t^k]a$ with
$k\in\{0,\ldots,\alpha_i+\beta_j-1\}$, or equivalently derivatives of
these orders. At a singular point the remainder module $\QQ$ supplies
the tests; the following cusp shows why a smooth parameter can miss one.

Let $e_i$ be the $i$th coordinate column and $E_{i,j}$ the matrix
unit with entry $1$ at $(i,j)$, of the required size.
For $A=\CC[t^2,t^3]$, write $x=t^2$ and take $L_1=\diag(x,1)$,
$R_1=\diag(x,1,1)$ and $\nu=E_{1,1}\mathsf c_5$.
The algebra contains $1,t^2,t^3,t^4,\ldots$, but not $t$.
Consequently $\mathsf c_5(x^2h)=0$ for every $h\in A$: a term of
degree five would require a degree-one term in $h$. On the other
hand, $\mathsf c_5(xh)=\mathsf c_3(h)$ can be nonzero.
Thus $x^2\mathsf c_5=0$ but $x\mathsf c_5=\mathsf c_3\ne0$.
The quotient $A/(x^2)$ has basis $1,t^2,t^3,t^5$. Pairing two basis
elements by extracting the coefficient of $t^5$ gives the matrix
\[
 \begin{bmatrix}0&0&0&1\\0&0&1&0\\0&1&0&0\\1&0&0&0\end{bmatrix}.
\]
It is invertible, so this mass attains rank $d_L+d_R=4$.
Normalization adjoins the missing parameter $t$, giving $\CC[t]$. There, masses annihilated by $t^4$ use only derivatives
of orders zero through three. Since $t^5=t^4t$, they vanish at $t^5$,
whereas $\mathsf c_5(t^5)=1$; none restricts to this cusp mass.

For the curve defined by the two circle equations,
$A=\CC[x,a,b]/(x^2+a^2-1,x^2+b^2-1)$ is free over $\CC[x]$
with basis $1,a,b,ab$. Hence $1-x$ is regular and
$A/(1-x)=\CC[a,b]/(a^2,b^2)$ retains four classes,
including the one detected by $\partial_a\partial_b$.

Indeed, after setting $x=1$, every class has a unique representative
$c_{0,0}+c_{1,0}a+c_{0,1}b+c_{1,1}ab$. Evaluation and the three
derivatives $\partial_a$, $\partial_b$, $\partial_a\partial_b$ at
$a=b=0$ recover these four coefficients. Point evaluation alone
would retain only $c_{0,0}$ and would miss three allowed mass data.

\begin{remark}[Cancellations preserve a ratio, not an equation]
The rational term depends on $S_L,S_R$ and the extension; $\QQ^*$
depends on the denominator factors. Cancelling common factors in
\eqref{eq:four} loses their mass freedom unless it is retained separately.
For example, take $q=p=1$, $A=\CC[t]$ and
$L_0=L_1=R_0=R_1=t$. Both rational ratios equal one, but the
equation is $t^2\widehat\mu=t^2\mu$. Its solutions are
$\widehat\mu=\mu+c_0\mathsf c_0+c_1\mathsf c_1$.
Cancelling the factors and replacing the equation by
$\widehat\mu=\mu$ would remove both free constants.
\end{remark}

\subsubsection{Finite normality and reproducing kernels}\label{sec:kernelsetup}

We record the finite Gram construction used throughout the paper.
For the comparison in this subsubsection, both target spaces will have
the same dimension as the original spaces. The original
row space $E$ and column space $H$ have the same dimension $N$.
Choose bases $e_1,\ldots,e_N$ and $h_1,\ldots,h_N$, and set
\[
 G_0\coloneq[\beta_\mu(e_i,h_j)]_{i,j\in\{1,\ldots,N\}}.
\]
The pairing is \emph{normal on these spaces} when $G_0$ is invertible.
This is a finite matrix condition; no inverse of a semi-infinite
moment matrix is assumed. Set
$e(Q)\coloneq[e_i(Q)]_{i\in\{1,\ldots,N\}}$ (rows stacked) and
$h(P)\coloneq[h_1(P)\ \cdots\ h_N(P)]$. The reproducing kernel is
\begin{equation}\label{eq:kerneldefinition}
 K_E(P,Q)\coloneq h(P)G_0^{-1}e(Q).
\end{equation}
To check reproduction, write $g=\sum_j h_j c_j$. The column of values
$[\beta_\mu(e_i,g)]_i$ is $G_0c$; multiplying by $h(P)G_0^{-1}$
returns $h(P)c=g(P)$. For a row $f=ae$, the row of its pairings with the column basis
is $aG_0$. Multiplication by $G_0^{-1}e(Q)$ gives $ae(Q)=f(Q)$,
which is reproduction on $E$.
For independence of bases, replace $e$ by $Ce$ and $h$ by $hD$, with
$C,D$ invertible. The Gram becomes $CG_0D$, and the new kernel is
\[
 h(P)D(CG_0D)^{-1}Ce(Q)=h(P)G_0^{-1}e(Q).
\]
Thus the kernel depends on the spaces and their pairing, not on the
chosen bases.

Choose target spaces $X\subset A^{1\times q}$ and
$Y\subset A^{p\times1}$ of dimension $N$, normal for
$\widehat\beta$, and let $\widehat K$ be their kernel. Multiplication by the rational
factors gives the following kernel and its difference from $K_E$:
\begin{equation}\label{eq:transportedkernel}
 \cF(P,Q)\coloneq S_R(P)\widehat K(P,Q)S_L(Q),\qquad
 \Delta\coloneq\cF-K_E.
\end{equation}
They belong to the tensor product of
$W_R\coloneq H+S_RY$ and $W_L\coloneq E+XS_L$, over $\CC$.
Here $P,Q$ label independent arguments, even for balanced pairings.
The \emph{separation rank} of $\Delta$ is the smallest $r$ for which
$\Delta=\sum_{j=1}^r u_j(P)v_j(Q)$ with $u_j\in W_R,v_j\in W_L$.
In bases of these two spaces it is the ordinary rank of the coefficient
matrix: each separated product gives an outer product of coordinate
vectors, and a rank factorization of that matrix gives the converse
decomposition. We also call it tensor rank. On a nonreduced support this
is an identity of polynomial classes; point values alone need not
distinguish the classes.

We shall also use finite spaces of sections of vector bundles. Their
coordinates are polynomial or rational vectors in each local basis,
and the same Gram-matrix construction defines their kernels. A balanced
mass acts on the tensor product over $A$ of the two section modules.
All finite-dimensional arguments above use complex bases of section
spaces and are therefore unchanged by this extension.

\subsection{When a bilateral mass splits into unilateral terms}\label{sec:masssplitting}

Can a bilateral mass split into two masses, each killed by one denominator?
For arbitrary bivariate masses, formula~\eqref{eq:bivariantmasses}
already supplies such a decomposition. Its first sum is killed by
$L_1$, because all its left remainder coordinates vanish on
$\mathscr F_0$. Its second sum is killed by $R_1$, because the right
remainder coordinates vanish on $\mathscr G_0$. Requiring each summand
to be balanced introduces an additional condition: these two sums need
not separately allow scalar factors to move between arguments.
Put $\mathscr P_L\coloneq A^{1\times q}$, $\mathscr P_R\coloneq A^{p\times1}$, $\mathscr D_L\coloneq \mathscr P_L L_1$, $\mathscr D_R\coloneq R_1\mathscr P_R$,
and $\mathscr C_L\coloneq \mathscr P_L/\mathscr D_L$, $\mathscr C_R\coloneq \mathscr P_R/\mathscr D_R$. Tensor products in this subsection are
over $A$, and stars denote full complex algebraic duals.\footnote{The
tensor product over $A$ imposes $(af)\otimes g=f\otimes(ag)$ as well
as additivity. A complex linear functional on this tensor product
is therefore exactly a balanced bilinear pairing.
See \cite[Section~2.2]{Eisenbud}.} Identify $V$
with $\mathscr P_L\otimes \mathscr P_R$ by $f\otimes g\mapsto f^{\T}g^{\T}$.
Let $\mathscr M$ be the space of balanced homogeneous bilateral masses and let
$\mathscr M_L,\mathscr M_R$ be its subspaces satisfying $L_1\nu=0$
and $\nu R_1=0$, respectively.
The concrete obstruction is the quotient
\[
 \frac{(\mathscr D_L\otimes\mathscr P_R)
       \cap(\mathscr P_L\otimes\mathscr D_R)}
      {\mathscr D_L\otimes\mathscr D_R}.
\]
The intersection consists of tests divisible on each side separately;
the denominator discards those divisible on both sides at once.
The remaining classes form the obstruction, denoted here by
$\Tor_1^A(\mathscr C_L,\mathscr C_R)$.\footnote{For the free
presentation $0\to\mathscr D_L\to\mathscr P_L\to\mathscr C_L\to0$,
this Tor is the kernel of $\mathscr D_L\otimes_A\mathscr C_R\to
\mathscr P_L\otimes_A\mathscr C_R$, namely the injectivity lost after
tensoring; see \cite[Section~6.2]{Eisenbud}.}

For scalar denominators $a,b\ne0$, the obstruction is
$((a)\cap(b))/(ab)$: separate divisibility need not imply divisibility
by the product. For $a=b=t$ in $\CC[t]$, the extra mass direction
is represented by $\nu[f]=f'(0)$. It annihilates $t^2A$ but not $tA$, whereas
unilateral masses are multiples of evaluation at zero.
The
determinant divisors are the codimension-one zero sets of
$\det L_1$ and $\det R_1$, with their vanishing orders.
The local formula uses normality of the support.\footnote{A domain is normal if elements of its fraction field satisfying
monic equations over it already belong to it. At codimension one
this gives a discrete valuation ring, independently of Gram
normality. See \cite[Sections~4.2 and~11.1]{Eisenbud}.}

The exact sequence below records two separate issues. Its leftmost
term measures bilateral masses that cannot be written as sums of
unilateral masses. Its rightmost term measures the ambiguity of a
decomposition when one exists. Thus failure of existence and
nonuniqueness of the two summands have different algebraic sources.

\begin{theorem}[The obstruction to unilateral decomposition]\label{thm:tormasses}
For balanced masses, with regular denominator determinants, there is
a canonical exact sequence
\begin{equation}\label{eq:tormodule}
0\to\Tor_1^A(\mathscr C_L,\mathscr C_R)\to\QQ
 \to(\mathscr C_L\otimes\mathscr P_R)\oplus(\mathscr P_L\otimes\mathscr C_R)
 \to\mathscr C_L\otimes\mathscr C_R\to0.
\end{equation}
The middle maps take $[H]$ to $([H],[H])$ and
$([H_1],[H_2])$ to $[H_1-H_2]$ in the indicated quotients.
Consequently
$\mathscr M/(\mathscr M_L+\mathscr M_R)
 \simeq\Tor_1^A(\mathscr C_L,\mathscr C_R)^*$ and
$\mathscr M_L\cap\mathscr M_R\simeq(\mathscr C_L\otimes \mathscr C_R)^*$.
If $Z^\circ$ is smooth and the two determinant divisors have no common
irreducible component, every balanced bilateral mass is a sum of balanced unilateral masses.
For the following local formula, assume that $Z^\circ$ is normal at
the generic point $\eta$ of the common component. Its local ring
$B\coloneq \OO_{Z^\circ,\eta}$ inverts functions not identically zero on that
component. It is a discrete valuation ring: every nonzero element is
a unit times a power of a local parameter $t$.
With denominator Smith orders $\alpha_i,\beta_j$,
$\Tor_1^A(\mathscr C_L,\mathscr C_R)\otimes_A B
 \simeq\bigoplus_{i,j} B/(t^{\min(\alpha_i,\beta_j)})$.
\end{theorem}

\begin{proof}
Inside $\mathscr P_L\otimes \mathscr P_R$ write $I\coloneq \mathscr D_L\otimes \mathscr P_R$, $J\coloneq \mathscr P_L\otimes \mathscr D_R$,
$K_0\coloneq \mathscr D_L\otimes \mathscr D_R$. All these embeddings are injective, since
$\mathscr P_L,\mathscr D_L,\mathscr P_R,\mathscr D_R$ are free. A class in $\QQ=V/K_0$ maps to its two classes in $V/I$ and
$V/J$. Both images vanish exactly when its representative lies
in $I\cap J$, giving kernel $(I\cap J)/K_0$.
Next, the difference of the two images is taken modulo $I+J$.
To compute its kernel, suppose $H_1-H_2=i+j$, with $i\in I$ and
$j\in J$. The element $H\coloneq H_1-i=H_2+j$ has class $[H_1]$
modulo $I$ and class $[H_2]$ modulo $J$. Thus every pair in the
kernel is represented by one common element of $V$, which is exactly
the image of $\QQ$. The difference map
is surjective, with quotient
$V/(I+J)=\mathscr C_L\otimes\mathscr C_R$.
This verifies exactness directly at the three quotient terms.
For the leftmost term, tensor the free resolution
$0\to\mathscr D_L\to\mathscr P_L\to\mathscr C_L\to0$
with $\mathscr C_R$. Both modules in the resolution are free:
multiplication by the regular matrix $L_1$ identifies
$\mathscr P_L$ with its image $\mathscr D_L$. By the definition
of the first Tor group, the new kernel is
\[
 \ker(\mathscr D_L\otimes\mathscr C_R
               \longrightarrow\mathscr P_L\otimes\mathscr C_R).
\]
The source is $I/K_0$ and the target is $V/J$, because tensoring
the presentation of $\mathscr C_R$ with either free module
preserves exactness. This map sends a representative $i\in I$
to its class modulo $J$. Its kernel consists exactly of
$i\in I\cap J$, modulo $K_0$. This identifies the kernel with
$(I\cap J)/K_0$ and proves the left end of \eqref{eq:tormodule}
\cite{Eisenbud}.

The obstruction can be evaluated without further homological notation.
Given a mass $\nu$, restrict it to $I\cap J$:
\[
 \operatorname{obs}(\nu)([H])\coloneq\nu[H],
 \qquad [H]\in(I\cap J)/K_0.
\]
It is well defined because every mass vanishes on $K_0$. A sum of
unilateral masses vanishes on $I\cap J$, so its obstruction is zero.
Conversely, if this restriction is zero, $\nu$ factors through the
image of $\QQ$ in $(V/I)\oplus(V/J)$. Extend that functional to the
whole direct sum by completing a complex basis. Its two restrictions
are the required unilateral masses. This proves the existence criterion
directly. Any functional on $(I\cap J)/K_0$ extends to $\QQ$ by the
same basis argument, so every obstruction value occurs.

Equivalently, duality over $\CC$ is exact and reverses arrows. Since masses are
functionals on quotient tests, dualizing gives the addition map
$\mathscr M_L\oplus\mathscr M_R\to\mathscr M$, whose kernel consists of
pairs $(\omega,-\omega)$ from $(\mathscr C_L\otimes\mathscr C_R)^*$.
Its cokernel is $\mathscr M/(\mathscr M_L+\mathscr M_R)$, proving both assertions.

It remains to relate this obstruction to the denominator zero sets.
We first prove vanishing when the support is smooth and the two
determinant divisors have no common component.

There is a direct divisibility proof. Smoothness makes the local
ring regular, and a regular local ring is a unique factorization
domain by \cite[Theorem~19.19]{Eisenbud}. Thus a nonzero function
factors into prime elements, with multiplicities, up to a unit.
These prime factors describe the local codimension-one components
of its zero set. Put $s\coloneq\det L_1$ and $r\coloneq\det R_1$ in that
ring. The hypothesis on components says that $s,r$ have no common
prime factor. Suppose that a column class $[v]\in\mathscr C_R$ is
killed by $s$. This means $sv=R_1w$ for some column $w$. Multiplying
by the adjugate gives
\[
 s\,\operatorname{adj}(R_1)v=rw.
\]
If a prime element occurs to exponent $m$ in $r$, it does not divide
$s$. The displayed equality therefore forces that prime to occur to
at least exponent $m$ in every entry of $\operatorname{adj}(R_1)v$.
Doing this for all prime factors of $r$ proves that $r$ divides each
entry. Thus
$z\coloneq r^{-1}\operatorname{adj}(R_1)v$ is a column over the local
ring, and $R_1z=v$. Hence $[v]=0$. We have proved that multiplication
by $s$ is injective on $\mathscr C_R$. The adjugate of $L_1$ now shows
that the map defining the Tor kernel is injective. The kernel vanishes
at every point and is therefore zero.

At the generic point of a common component, the discrete valuation
ring $B$ admits Smith bases. They split the two cokernels into
cyclic modules $B/(t^{\alpha_i})$ and $B/(t^{\beta_j})$.
For a pair of exponents $a,b$, use the free resolution
\[
 0\longrightarrow B\xrightarrow{\ t^a\ }B
   \longrightarrow B/(t^a)\longrightarrow0.
\]
After tensoring with $B/(t^b)$, the Tor module is the kernel of
multiplication by $t^a$ on $B/(t^b)$. If $a\ge b$, that multiplication
is zero, so its kernel is $B/(t^b)$. If $a<b$, a class $[w]$ belongs
to the kernel exactly when $t^aw\in(t^b)$. A nonzero element of the
discrete valuation ring is a unit times a power of $t$; hence this
condition means $w\in(t^{b-a})$. The kernel is therefore
$t^{b-a}B/(t^b)$, and the map
\[
 B/(t^a)\longrightarrow t^{b-a}B/(t^b),\qquad
 [w]\longmapsto[t^{b-a}w]
\]
is onto and has zero kernel. Both cases give
$B/(t^{\min(a,b)})$; if either exponent is zero, the module is zero.
Localization preserves the free presentations used here, and tensor
products and their kernels distribute over these finite direct sums.
Applying the calculation to every pair of Smith components yields
\[
 \Tor_1^A(\mathscr C_L,\mathscr C_R)\otimes_A B
 \simeq\bigoplus_{i\in\{1,\ldots,q\}}
             \bigoplus_{j\in\{1,\ldots,p\}}
             B/(t^{\min(\alpha_i,\beta_j)}).
\]
This proves the local formula.
\end{proof}

The no-common-component assertion also has a homological interpretation.
The vanishing argument above uses divisibility; the following interpretation
records how smoothness controls the associated primes. The support of a finitely
generated module is the set of prime ideals where its localization
is nonzero. For the cokernel of a square matrix
with regular determinant, it lies on the determinant zero set: where
the determinant is invertible the adjugate gives an inverse.
An associated prime is the annihilator of a nonzero module element.
The zero divisors on a finitely generated module over a Noetherian
ring form the union of its associated primes
\cite[Chapter~3]{Eisenbud}. Thus, to prove that a scalar acts
injectively, it suffices to show that it belongs to none of those
primes. We apply this criterion to $\det L_1$ on $\mathscr C_R$.

Depth counts the successive parameters that act injectively after
quotienting by the preceding ones. Over a local ring, projective
dimension is the length of a shortest finite free resolution, when
such a resolution exists. The relation needed here is the Auslander--Buchsbaum formula for
a finite module of finite projective dimension over a regular local
ring. Here the injective square presentation gives projective
dimension at most one. At a point of the support the cokernel is a
nonzero torsion module, so it is not free and its projective dimension
is exactly one. The formula therefore subtracts precisely one from
the depth of the local ring. These are the only homological inputs
in the no-common-component assertion.

Smoothness makes each local ring $A_{\mathfrak p}$ regular: its maximal
ideal is generated by $\dim A_{\mathfrak p}$ parameters. This ring
property differs from a regular scalar (a nonzerodivisor).
For nonzero $\mathscr C_R$, its injective square presentation and the depth formula give
$\operatorname{depth}((\mathscr C_R)_{\mathfrak p})=\dim A_{\mathfrak p}-1$
at primes in its support.
At an associated prime, depth is zero. The displayed formula forces
the local ring to have dimension one, so the prime has height one:
it is the generic point of a codimension-one component.
These primes belong to the determinant divisor of $R_1$; the
disjoint-component hypothesis therefore makes $\det L_1$ a
nonzerodivisor on $\mathscr C_R$. By the adjugate identity,
$\mathscr D_L\otimes\mathscr C_R\to\mathscr P_L\otimes\mathscr C_R$
is injective, so its Tor kernel vanishes.

For $A=\CC[x,y]$, scalar denominators $x,y$ give
$\QQ=A/(xy)$ and zero Tor. A polynomial divisible by both $x$ and $y$
is divisible by $xy$, so $(x)\cap(y)=(xy)$. The final quotient in
\eqref{eq:tormodule} is $A/(x,y)=\CC$. A mass splits into functionals on the two
axes; transferring a multiple of evaluation at their intersection
between the two summands is the only ambiguity. For denominators $x,x$,
$\QQ=A/(x^2)$ and the Tor submodule is $xA/(x^2)$. The missing
unilateral classes have the form
$a\mapsto\lambda(\partial_xa(0,y))$, with arbitrary
$\lambda\in\CC[y]^*$. For $n\ge0$, a basis of the classes in $\QQ$ represented by
polynomials of total degree at most $n$ is
\[
 1,y,\ldots,y^n,\quad x,xy,\ldots,xy^{n-1},
\]
where the second list is empty for $n=0$. The submodule $xA/(x^2)$
has the second list as its basis. Its dimension is $n$, while the
whole quotient at this degree has dimension $2n+1$.
Thus a proper crossing and a common component have different effects.
The full-dual statement concerns algebraic functionals; imposing
positivity or a topology requires additional hypotheses.

\section{Finite formulas for polynomials and kernels}\label{part:finite}

We now calculate the transformed polynomials and the reproducing
kernel. The inputs are finite lists of original dual polynomials,
moments of rational functions, the four transformation matrices,
and the prescribed mass values. The rational moments are values of
the extension fixed in Theorem~\ref{thm:masses}; a finite list of them
is supplied together with the polynomial moments. They are evaluated
by integration when a measure and integrable rational tests are
given, or specified as additional functional data. The unknown polynomial must satisfy
three requirements: its components have the allowed bounds, it is
orthogonal to the specified lower space, and its leading coefficients
have the chosen normalization.

The method writes the rational images of these polynomials in finite
bases. Linear equations select exactly those combinations that have
allowed polynomial preimages. Original duality then gives an identity
block in the pairing matrix; eliminating that block leaves the reduced
system used in all subsequent determinant formulas, with the original
coefficients eliminated. Its order is calculated from the dimensions
of the added spaces and the constraints; comparison with the direct
normalized Gram order is made below. We first carry out
the elimination for a scalar example, then prove it for arbitrary
component sizes and polynomial bounds. The final two subsections count
the determinant orders and give direct finite-rank update formulas.

\subsection{A scalar calculation with a three-by-three determinant}
\label{sec:scalarresidual}

We begin with a scalar problem in one variable. The calculation
shows how original polynomials and Cauchy moments become the entries
of a small system, and how its solution gives the transformed
polynomial. Every step will then be used in the general mixed-type orthogonality problem.
Take $q=p=1$ and choose
\[
 L_0=z-a,\qquad L_1=z-\lambda,\qquad
 R_0=z-b,\qquad R_1=z-\xi,
 \qquad a\ne\lambda,\quad b\ne\xi.
\]
Let $\mu$ be a balanced scalar functional, and choose its extension
$\widetilde\mu$ to the rational tests with these denominator poles.
With zero homogeneous mass, the transformed functional is
\[
 \widehat\mu[h]=\widetilde\mu\left[
  \frac{(z-a)(z-b)}{(z-\lambda)(z-\xi)}h(z)\right].
\]
For $n\ge1$, the problem is to find the monic polynomial
$\widehat B_n$ of degree $n$, with leading coefficient one, such that
\[
 \widehat\mu[f\widehat B_n]=0
 \qquad\text{for every }f\in\CC[z]_{<n}.
\]
Thus $X=\CC[z]_{<n}$ and $Y=\CC[z]_{\le n}$, with one more
unknown coefficient than orthogonality conditions.

\medskip\noindent\textbf{Original data and the rational basis.}
Assume that $\mu$ has monic orthogonal polynomials $B_0,\ldots,B_n$
with $h_i\coloneq\mu[B_i^2]\ne0$. Define $A_i\coloneq B_i/h_i$;
then $\mu[B_iA_j]=\delta_{i,j}$. Choose
$E=H=\CC[z]_{\le n}$ as the original spaces, of dimension
$N=n+1$, and arrange their dual bases as
\[
 \cB(z)\coloneq[B_0(z)\ \cdots\ B_n(z)]^{\T},\qquad
 \cA(z)\coloneq[A_0(z)\ \cdots\ A_n(z)].
\]
Append $1/(z-\lambda)$ on the left and $1/(z-\xi)$ on the
right. Pairing the two enlarged lists gives
\[
 \begin{bmatrix}I_N&C\\D&E_0\end{bmatrix},
 \qquad
 \begin{aligned}
  C_i&\coloneq\widetilde\mu[B_i/(z-\xi)],\\
  D_i&\coloneq\widetilde\mu[A_i/(z-\lambda)],\\
  E_0&\coloneq\widetilde\mu[((z-\lambda)(z-\xi))^{-1}].
 \end{aligned}
\]
Here $C$ is the column with entries $C_i$ and $D$ the row with
entries $D_i$, for $i\in\{0,\ldots,n\}$. They are Cauchy values of the
original polynomials. The scalar $E_0$ pairs the two added pole
functions. If the poles coincide, it is the moment
$\widetilde\mu[(z-\lambda)^{-2}]$; the construction is unchanged.

\medskip\noindent\textbf{Image conditions and normalization.}
Let $e_N$ be the last coordinate column in $\CC^N$ and put
\[
 v\coloneq\cB(a),\qquad w\coloneq\cA(b)^{\T},\qquad
 \eta\coloneq h_n^{-1}e_N^{\T}.
\]
A general left ambient function is
$\phi(z)=u^{\T}\cB(z)+u_p/(z-\lambda)$. Its only possible
polynomial preimage is
\[
 \phi(z)S_L^{-1}(z)=\frac{z-\lambda}{z-a}\phi(z).
\]
Multiplication by $z-\lambda$ removes the original pole. The result
is divisible by $z-a$ precisely when $\phi(a)=0$, since
$a\ne\lambda$. After this division, the degree is at most $n$.
The coefficient of degree $n$ is $e_N^{\T}u$, because $B_n$ is
monic and $S_L^{-1}$ tends to one at infinity. Requiring this
coefficient to vanish gives exactly the bound defining $X$.
Hence the two left constraints are
\[
 \Gamma_L\begin{bmatrix}u\\u_p\end{bmatrix}=0,
 \qquad
 \Gamma_L=
 \begin{bmatrix}v^{\T}&(a-\lambda)^{-1}\\e_N^{\T}&0\end{bmatrix}.
\]
On the right, a general ambient function has the form
$\psi(z)=\cA(z)z_0+z_p/(z-\xi)$. Its preimage is
$(z-\xi)\psi(z)/(z-b)$. Vanishing at $b$ is the only image
condition, since $Y$ permits degree $n$. Its leading coefficient
is $\eta z_0$. Thus the right image condition and monic normalization
are, respectively,
\[
 w^{\T}z_0+(b-\xi)^{-1}z_p=0,\qquad \eta z_0=1.
\]

\medskip\noindent\textbf{Eliminating the original coefficients.}
The vector of pairings of the enlarged left basis with $\psi$ is
\[
 \begin{bmatrix}z_0+Cz_p\\Dz_0+E_0z_p\end{bmatrix}.
\]
Orthogonality requires its scalar product with every vector in
$\ker\Gamma_L$ to be zero.
Because $\Gamma_L$ has full row rank, its transpose has image
equal to this space of vectors orthogonal to $\ker\Gamma_L$.
Therefore there are unique scalars
$\alpha_1,\alpha_2$ such that
\[
 z_0+Cz_p=v\alpha_1+e_N\alpha_2,\qquad
 Dz_0+E_0z_p=(a-\lambda)^{-1}\alpha_1.
\]
The first equation gives
$z_0=-Cz_p+v\alpha_1+e_N\alpha_2$. Substituting this into the
second equation, the right image condition and the normalization
leaves three equations for $z_p,\alpha_1,\alpha_2$:
\begin{equation}\label{eq:explicitthree}
 T\begin{bmatrix}z_p\\\alpha_1\\\alpha_2\end{bmatrix}=e_3,
 \qquad
 T\coloneq\begin{bmatrix}
 E_0-DC&Dv-(a-\lambda)^{-1}&De_N\\
 (b-\xi)^{-1}-w^{\T}C&w^{\T}v&A_n(b)\\
 -\eta C&\eta v&\eta e_N
 \end{bmatrix}.
\end{equation}
Here $e_3=[0\ 0\ 1]^{\T}\in\CC^3$. The first row is the remaining
orthogonality equation, the second enforces polynomiality on the
right, and the third fixes the leading coefficient. Its order is
three for every $n$: the original polynomial coefficients have been
eliminated. The general matrix \eqref{eq:small} will implement the
same elimination with vector-valued constraints and masses.

\medskip\noindent\textbf{Recovering the polynomial.}
For $\det T\ne0$, solve the preceding system and multiply $\psi$
by $S_R^{-1}=(z-\xi)/(z-b)$. With
\[
 r(z)\coloneq
 [1-(z-\xi)\cA(z)C\quad(z-\xi)\cA(z)v\quad(z-\xi)A_n(z)],
\]
the result is $\widehat B_n(z)=r(z)T^{-1}e_3/(z-b)$.
Write $T[r]$ for $T$ with its last row replaced by $r$.
Expansion of its determinant along that row gives
$\det T[r]/\det T=rT^{-1}e_3$, and hence
\begin{equation}\label{eq:explicitthreepoly}
 \widehat B_n(z)=\frac{\det T[r(z)]}{(z-b)\det T}.
\end{equation}
All objects in this formula are original polynomials, their values
at $a,b$, and the rational moments $C,D,E_0$.

The apparent pole at $b$ cancels directly: $r(b)$ is $(b-\xi)$
times the second row of $T$, so the numerator determinant vanishes
there. More explicitly, its three entries at $b$ are
\[
 1-(b-\xi)\cA(b)C,\qquad
 (b-\xi)\cA(b)v,\qquad (b-\xi)A_n(b),
\]
which are precisely $(b-\xi)$ times the second row of $T$.
The numerator has degree at most $n+1$, and division by $z-b$
gives degree at most $n$. Since the degree-$n$ coefficient of
$\cA$ is $\eta$, the degree-$n+1$ coefficient of $r$ is
$[-\eta C\ \eta v\ \eta e_N]$, the last row of $T$.
The leading numerator coefficient is therefore $\det T$, proving
that the quotient has leading coefficient one. The pairing equations above give
orthogonality of $S_R\widehat B_n$ to every product $fS_L$, with
$f\in X$, for the chosen rational pairing. Therefore
$\widehat\mu[f\widehat B_n]=0$ for every $f\in X$. They also show
uniqueness, since all elimination steps are reversible and $T$ is
invertible.

\subsection{The finite construction for prescribed component degree bounds}\label{sec:finite}

We now extend the scalar calculation to polynomial rows and columns
with arbitrary prescribed finite bounds. The objective is to determine a polynomial column from
prescribed orthogonality and normalization conditions, using data of
the original pairing. We first choose the polynomial spaces, then
express the rational products in enlarged bases, and finally eliminate
the coefficients belonging to the original dual polynomials.
Throughout this subsection, the four determinants are regular and the
rational extension in Theorem~\ref{thm:masses} is fixed.
Let $X\subset A^{1\times q}$ and $Y\subset A^{p\times1}$ have dimensions
$r$ and $r+k$, respectively, where $k\ge0$. A square pairing is called
\emph{normal} when its Gram matrix is invertible.\footnote{The Gram matrix $[\beta(f_i,g_j)]$ in chosen bases is
invertible exactly when dual bases exist.} For $k>0$, prescribe a
linear map $\lambda:Y\to\CC^k$ and seek a $p\times k$ matrix $\widehat a$
with columns in $Y$ satisfying
$\widehat\beta(X,\widehat a)=0$ and $\lambda(\widehat a)=I_k$,
where $\lambda$ acts columnwise.
Each column has $r+k$ unknown coefficients, $r$ orthogonality
equations and $k$ normalization equations. In one variable, $X,Y$
impose the component degree bounds and $\lambda$ selects leading
coefficients. On a variety, we count independent polynomial classes
with the prescribed bounds. The restriction $\nu|_{X\times Y}$ is represented by a finite
matrix, so it has a rank factorization with $h$ terms of the form
\eqref{eq:massfactor}. This statement uses only the finite target
spaces. The bound $h\le d_L+d_R$ is specific to the curve case
of Theorem~\ref{thm:masses}; no such bound is being assumed in
higher dimension.

\subsubsection{Polynomial spaces and additional rational vectors}

Choose original polynomial spaces $E\subset A^{1\times q}$ and
$H\subset A^{p\times1}$ of dimension $N$, normal for $\beta_\mu$.
Their dual bases form matrices $\cB\coloneq [b_i]_{i\in\{1,\ldots,N\}}$ of size $N\times q$ and
$\cA\coloneq [a_1\ \cdots\ a_N]$ of size $p\times N$, with
$\beta_\mu(b_i,a_j)=\delta_{i,j}$.
We use $P,Q$ for independent copies of the support variables.
The spaces $E,H$ provide the original data, while $X,Y$ specify the
perturbed problem. They need not coincide. For a chosen basis of $X$,
multiply every row by $S_L$; for a chosen basis of $Y$, multiply every
column by $S_R$. Since these maps are linear, the resulting lists
span all the rational products that can occur in the problem.

On the left, start with the original basis of $E$ and add independent
vectors from the first list until it spans $E+XS_L$. On the right,
do the same with $H$ and $S_RY$. Independence is tested over $\CC$
and subject to the equations of the support. For example, in
$A=\CC[x,y]/(y^2-x^3)$, the fractions $y^2/d$ and $x^3/d$ are
the same element whenever the regular denominator $d$ is inverted.
They cannot contribute two independent directions.

Denote the added rows by $\phi_\alpha$ and the added columns by
$\psi_\beta$, and collect them in $\Phi_L,\Phi_R$. The minimal
choice of ambient spaces and the corresponding numbers of added
vectors are
\[
 \begin{aligned}
 W_L^{\min}&\coloneq E+XS_L,& W_R^{\min}&\coloneq H+S_RY,\\
 s_L^{\min}&\coloneq\dim(W_L^{\min}/E),&
 s_R^{\min}&\coloneq\dim(W_R^{\min}/H).
 \end{aligned}
\]
These spaces and dimensions are determined by the factors and the
chosen polynomial spaces. The individual added vectors are choices
of basis. Some may be polynomial, for example when multiplication
raises a degree beyond the original space.

For a scalar example, take $S_L(z)=(z-a)/(z-\xi)$ with $a\ne\xi$.
Polynomial division gives
\[
 f(z)S_L(z)=f(z)+(\xi-a)\frac{f(z)-f(\xi)}{z-\xi}
                  +\frac{(\xi-a)f(\xi)}{z-\xi}.
\]
The first two terms are polynomial. If their span is already in $E$,
the single added function $1/(z-\xi)$ suffices for the rational part.
With matrix factors, the component directions also matter; on an
algebraic support, the defining relations determine which of the
resulting fractions are independent.

Here an \emph{ambient space} means a finite space containing both the
original basis and the rational image to be calculated. Its extra
directions need not themselves be rational images of allowed polynomials.
The theorem permits larger ambient spaces when they are more convenient
for computation. Choose $W_L\supset E+XS_L$ and $W_R\supset H+S_RY$,
with dimensions $N+s_L,N+s_R$, and bases
\[
 \cB^+\coloneq\begin{bmatrix}\cB\\\Phi_L\end{bmatrix},
 \qquad \cA^+\coloneq[\cA\ \Phi_R].
\]
All products used in the pairings must lie in the chosen rational
extension domain. The minimal ambient spaces lie componentwise in
$A[1/(\det L_1)]$ on the left and $A[1/(\det R_1)]$ on the right,
as follows from the adjugate expressions for the inverses.

The image equations can now be constructed from polynomial coefficients.
Fix ordered bases $f_1,\ldots,f_r$ of $X$ and
$g_1,\ldots,g_{r+k}$ of $Y$. Choose common regular denominators
for the finite lists consisting of the ambient basis vectors and the
rational images $f_iS_L,S_Rg_j$. Multiplication by these denominators
turns equality of rational vectors into equality of polynomial classes
in $A$. Reduce all entries modulo $I$, using a Gr\"obner basis, and
collect the finitely many standard monomials that occur. Coefficients
in the component copies of these monomials give ordinary complex
column vectors.

On the left, let $\mathsf A_L$ have as columns the coefficient vectors
of the cleared ambient basis, and let $\mathsf B_L$ have as columns
those of the cleared images $f_iS_L$. Use the same scalar denominator
and the same polynomial coordinates in both matrices. The ambient
basis is independent and denominator multiplication is injective,
so $\mathsf A_L$ has full column rank $N+s_L$. Each column of
$\mathsf B_L$ belongs to its image. Solving the coefficient equations
therefore determines a unique matrix $Z_L$ by
\[
 \mathsf A_LZ_L=\mathsf B_L,
 \qquad Z_L\in\Mat_{(N+s_L)\times r}(\CC).
\]
Its columns are the coordinates of $f_iS_L$ in $\cB^+$.
Multiplication by $S_L$ is injective, so these columns have rank $r$.
The identical calculation on columns gives
\[
 \mathsf A_RZ_R=\mathsf B_R,
 \qquad Z_R\in\Mat_{(N+s_R)\times(r+k)}(\CC),
 \qquad \rank Z_R=r+k,
\]
where the ambient columns and the images $S_Rg_j$ are cleared
by the same regular denominator on the right.

Row reduction of $Z_L^{\T}$ gives a basis of its kernel.
Transposing those basis vectors produces independent equations on
the ambient coordinates; their common kernel is $\im Z_L=XS_L$.
Indeed they annihilate every column of $Z_L$, and their number is
$N+s_L-r$, leaving a kernel of dimension $r$. The right calculation
leaves a kernel of dimension $r+k$. Splitting the columns according
to the original and added basis vectors gives the constraint matrices
\begin{equation}\label{eq:constraints}
 \Gamma_L\coloneq [X_0\ X_p],\qquad\Gamma_R\coloneq [Y_0\ Y_p],
 \quad\ker\Gamma_L=XS_L,\quad\ker\Gamma_R=S_RY,
\end{equation}
where kernels are understood in the chosen coordinates. The rows are
independent, so their numbers are
$c_L\coloneq N+s_L-r$ and $c_R\coloneq N+s_R-r-k$.
In particular, $X_0,X_p$ have sizes $c_L\times N,c_L\times s_L$,
and $Y_0,Y_p$ have sizes $c_R\times N,c_R\times s_R$.
Thus the factors, the chosen finite spaces and polynomial coefficient
reduction determine every entry of the constraints. The division
construction below separates their polynomial meaning into remainder
conditions and bounds on the resulting quotients.
On their kernels,
$\pi_L(b)\coloneq bS_L^{-1}$ and $\pi_R(c)\coloneq S_R^{-1}c$ recover the polynomial
rows in $X$ and columns in $Y$, respectively. Extend these maps linearly to $W_L,W_R$, with values in
$X,Y$; only their values on the kernels enter the answers.

Such extensions exist because the spaces are finite dimensional:
complete a basis of each kernel and prescribe any values on the
added basis vectors, for example zero. A single enlarged-basis vector
need not itself correspond to an allowed polynomial. The maps
$\pi_L,\pi_R$ provide polynomial values for completing the tables
below; the constraint equations ensure that the final combinations
lie in the kernels, where the maps are the actual inverse
multiplications. This is why the final answers do not depend on the
chosen extensions.

\subsubsection{Constructing the divisibility conditions}
\label{sec:divisiondata}

For a rational row $b$, solving $b=fS_L$ requires two checks:
$f=bS_L^{-1}$ must be polynomial on the support, and it must belong
to the prescribed space $X$. A polynomial quotient can still violate
a component degree bound. For columns the order is reversed.
Choose regular scalars $a_L,a_R\in A$ clearing both the ambient spaces
and the rational matrices, and put
\begin{equation}\label{eq:polynomiallifts}
 P_L\coloneq a_LS_L,\qquad P_R\coloneq a_RS_R,
 \qquad a_LW_L\subset A^{1\times q},\quad a_RW_R\subset A^{p\times1}.
\end{equation}
For the minimal ambient spaces, the explicit choices
$a_L=\det L_1$, $P_L=\operatorname{adj}(L_1)L_0$,
$a_R=\det R_1$ and $P_R=R_0\operatorname{adj}(R_1)$ suffice; an enlarged ambient space may require additional scalar factors.
These lifts represent the rational images; the functional equation and
its mass module retain the original four factors.

\Needspace{4\baselineskip}
There are two distinct sources of constraints. A remainder detects
failure of polynomial divisibility. Once that remainder is zero, the
quotient may still lie outside the prescribed degree or section space.
The proposition counts these two groups separately. This distinction
prevents a divisible rational expression from being accepted as a
candidate with the wrong polynomial bounds.

\begin{proposition}[Polynomial division and the target quotient]
\label{prop:lifteddivision}
For $b\in W_L$ and $c\in W_R$,
\begin{equation}\label{eq:liftedmembership}
 \begin{aligned}
 b\in XS_L&\quad\text{if and only if}\quad a_Lb\in XP_L,\\
 c\in S_RY&\quad\text{if and only if}\quad a_Rc\in P_RY.
 \end{aligned}
\end{equation}
Each condition consists of polynomial matrix divisibility, followed by
membership of the unique quotient in $X$ or $Y$. Define $C(P_L)\coloneq A^{1\times q}/A^{1\times q}P_L$ and
$\rho_L\coloneq \dim\im[a_LW_L\longrightarrow C(P_L)]$,
and define $C(P_R)\coloneq A^{p\times1}/P_RA^{p\times1}$ and $\rho_R$
on the right. Exactly $\rho_L$ independent scalar divisibility tests
and $c_L-\rho_L$ quotient-space tests give $\Gamma_L$;
the corresponding numbers for $\Gamma_R$ are $\rho_R,c_R-\rho_R$.
\end{proposition}

\begin{proof}
Since $a_L$ is regular, multiplying an equality by $a_L$ neither loses
nor introduces solutions. Thus $b=fS_L$ with $f\in X$ is equivalent
to $a_Lb=fP_L$ with the same $f$. The right identity follows by
the same argument, keeping the matrix order $a_Rc=P_Rg$.

Consider first the map from $a_LW_L$ to the quotient $C(P_L)$.
Its image has dimension $\rho_L$. Choosing a basis of the dual of
this image gives exactly $\rho_L$ independent scalar conditions
for its remainder to vanish. Their common kernel therefore has
dimension $N+s_L-\rho_L$.

On that kernel, division by $P_L$ is well defined and unique,
because $P_L$ is invertible over $K$. It identifies the kernel with
a polynomial row space of dimension $N+s_L-\rho_L$ containing $X$.
Since $\dim X=r$, a further $N+s_L-\rho_L-r=c_L-\rho_L$
independent conditions select precisely $X$. Choose these conditions
on the divisible subspace and extend them linearly to the ambient
space. Together with the remainder conditions, they form $\Gamma_L$.
To verify independence, restrict a putative relation among all these
tests to the divisible subspace. Independence of the quotient-space
tests makes their coefficients zero. What remains is a relation among
the independent remainder tests, whose coefficients are also zero.
The common kernel is $XS_L$, and the number of independent tests is
therefore exactly $c_L$.

For columns, the remainder map has image dimension $\rho_R$ and
kernel dimension $N+s_R-\rho_R$. Division identifies this kernel
with a space containing $Y$, of dimension $r+k$. Thus the number
of further conditions is $N+s_R-\rho_R-(r+k)=c_R-\rho_R$.
\end{proof}

Remainders can be tested by values and derivatives. For example, modulo
$t^2$ the remainder is $a+bt$: its value at zero detects $a$, and its
first derivative detects $b$. On an algebraic support the derivatives
must also respect its equations. The following construction imposes
this compatibility, including at singular points. For a closed point $\xi$
let $\mathfrak m_\xi$ be its maximal ideal.\footnote{Here $\mathfrak m_\xi\coloneq \{f:f(\xi)=0\}$ and
$A/\mathfrak m_\xi=\CC$. The local ring $A_{\mathfrak m_\xi}$
permits denominators nonzero at $\xi$; quotienting by
$\mathfrak m_\xi^h$ retains finite-order data.
See \cite[Sections~2.1 and~4.5]{Eisenbud}.} A functional on
$C(P_L)/\mathfrak m_\xi^{h_\xi} C(P_L)$ acts on polynomial row representatives
as a finite jet of order below $h_\xi\ge1$. For $\alpha\in\NN^M$ put
$|\alpha|\coloneq \sum_j\alpha_j$, $\alpha!\coloneq \prod_j\alpha_j!$ and
$\partial^\alpha\coloneq \prod_j\partial_{x_j}^{\alpha_j}$; then
\begin{equation}\label{eq:compatiblejets}
 \gamma(p)\coloneq \sum_{i=1}^q\sum_{|\alpha|<h_\xi}
 c_{i,\alpha}\frac{\partial^\alpha p_i(\xi)}{\alpha!}.
\end{equation}
If $A=\CC[x_1,\ldots,x_M]/I$, the coefficients must annihilate
$I A_{\mathrm{amb}}^{1\times q}$ and every polynomial row multiple of a
representative of $P_L$, where $A_{\mathrm{amb}}\coloneq \CC[x_1,\ldots,x_M]$.
These relations make the jet independent of representatives on the
support. To see why these are precisely the necessary compatibility equations,
expand each polynomial about $\xi$ and discard all monomials of total
degree at least $h_\xi$. Its remaining coefficients are
$\partial^\alpha p_i(\xi)/\alpha!$. Every linear functional on this
truncated space has the displayed form. It descends first to the
support and then to the remainder module exactly when it vanishes on
the truncated support relations and on the truncated polynomial
matrix multiples. Only finitely many such truncated relations are
independent. The column definition uses column multiples and the
same Taylor coefficients.

\Needspace{4\baselineskip}
The value-only criterion uses radical ideals.\footnote{An ideal $J$ is radical if $f^m\in J$ implies $f\in J$.
Over $\CC$, point values then distinguish its quotient classes.
For $(t^2)$ they miss the class of $t$.
See \cite[Section~4.5]{Eisenbud}.}

The quotient description is intrinsic, but to fill a constraint
matrix we need scalar tests. The next result shows that finitely
many values and derivatives suffice on each prescribed ambient space.
It does not require the entire remainder module to be finite
dimensional: only its image from $a_LW_L$ or $a_RW_R$ is finite
dimensional. The compatibility conditions in
\eqref{eq:compatiblejets} ensure that the tests respect both the
support equations and multiplication by the polynomial matrix.

\Needspace{6\baselineskip}

The next result proves that a finite separating list of jets exists
for the chosen spaces. The actual points, orders and independent
rows are part of the constraint data; no uniform jet-order bound
is asserted. The coefficient construction above already computes the
constraints from any supplied polynomial presentation of the finite
spaces. The jet description identifies those equations with local
values and derivatives. Subsubsection~\ref{sec:mixedcurveexample}
computes the finite image conditions directly from its rational basis
vectors.

\begin{proposition}[Complete finite jets on a prescribed finite-dimensional space]
\label{prop:completejets}
The divisibility conditions in Proposition~\ref{prop:lifteddivision}
can be chosen as $\rho_L$ and $\rho_R$ independent compatible jets.
Directional values $p(\xi)v$, $v\in\ker P_L(\xi)$, suffice on the
left in either of the following cases: $(\det P_L)$ is radical in $A$;
or the support is a smooth affine curve and every local Smith exponent
of $P_L$ at a determinant zero is zero or one. The right version uses
left null vectors of $P_R(\xi)$ under the corresponding hypothesis.
The second criterion allows repeated zeros of the determinant.
\end{proposition}

\begin{proof}
Let $u$ be a nonzero remainder and choose a maximal ideal
$\mathfrak m$ containing its proper annihilator. Its localization
$u_{\mathfrak m}$ is nonzero: otherwise some $a\notin\mathfrak m$
would satisfy $au=0$, contradicting containment of the annihilator
in $\mathfrak m$. The localized remainder module is finite over
the Noetherian local ring $A_{\mathfrak m}$.
Krull's intersection theorem gives
\[
 \bigcap_{h\ge1}\mathfrak m^h C(P_L)_{\mathfrak m}=0.
\]
Hence $u_{\mathfrak m}$ has a nonzero image modulo some power
$\mathfrak m^h$. The quotient by this power is supported at the
single closed point, so localizing it does not change it. It is
a finite-dimensional vector space, and a linear functional on it
can be chosen nonzero on the image of $u$. Thus a finite Taylor
truncation detects $u$ \cite{Eisenbud}.
Its linear functionals are compatible jets \eqref{eq:compatiblejets}.
Restrict these jets to the image of $a_LW_L$ in $C(P_L)$.
Their common kernel is zero: every nonzero class is detected by at
least one of them. In a vector space of dimension $\rho_L$, a family that detects
every nonzero vector spans the dual: if its span had smaller dimension,
the common kernel would have positive dimension by rank-nullity. We can
therefore select $\rho_L$ independent jets, whose simultaneous
vanishing is equivalent to zero remainder. This proves finiteness
for the prescribed space, although the required points and orders
may change when the space is enlarged. The right construction is
obtained by transposition.

Where $\det P_L=0$, the columns of $\operatorname{adj}(P_L)$ lie in
$\ker P_L$, so the directional conditions make $p\operatorname{adj}(P_L)$
vanish there. Since $(\det P_L)$ is radical, the Nullstellensatz makes
every entry divisible by $\det P_L$. Thus
$pP_L^{-1}=p\operatorname{adj}(P_L)/\det P_L$ is polynomial.
Conversely, if $p=hP_L$ is divisible, then
$p(\xi)v=h(\xi)P_L(\xi)v=0$ for every
$v\in\ker P_L(\xi)$. Thus the directional tests are both necessary
and sufficient. Transposition proves the column assertion.

On a smooth curve write locally $P_L=U\diag(t^{\theta_i})V$, with
$U,V$ invertible. Here \emph{regular at a point} means without a pole;
regular scalars in Subsection~\ref{sec:algebra} are instead elements
that are not zero divisors.
If $\theta_i\in\{0,1\}$, the quotient
$pP_L^{-1}=pV^{-1}\diag(t^{-\theta_i})U^{-1}$ has no pole precisely
when $p(\xi)V(\xi)^{-1}e_i=0$ for each $\theta_i=1$.
These vectors form a basis of $\ker P_L(\xi)$. A rational function
without poles on the affine curve belongs to $A$, so the quotient is
polynomial. This proves the second criterion and its transposed version.
Differential descriptions of module membership in greater generality
are developed in \cite{CidRuizSturmfels}.
\end{proof}

For a rational matrix with local Smith exponents $(-1,0,\ldots,0)$,
a minimal local scalar lift has Smith form $\diag(1,t,\ldots,t)$. Its determinant vanishes to order $q-1$,
but the new criterion still uses values only. This illustrates why determinant multiplicity alone does not specify
the derivative order needed by the constraint tests. A larger scalar lift may
introduce redundant derivative tests.

For example, over $\CC[x,y]$, take
$P= \bigl[\begin{smallmatrix}x^2&0\\-xy&1\end{smallmatrix}\bigr]$.
Its inverse over $\CC(x,y)$ is
\[
 P^{-1}=\begin{bmatrix}x^{-2}&0\\y/x&1\end{bmatrix},\qquad
 pP^{-1}=[(p_1+xy p_2)/x^2\quad p_2].
\]
Write $p_1=A_0(y)+xA_1(y)+x^2R_1(x,y)$ and
$p_2=B_0(y)+xR_2(x,y)$. The first numerator becomes
$A_0+x(A_1+yB_0)+x^2(R_1+yR_2)$. It is divisible by $x^2$
exactly when $A_0=0$ and $A_1+yB_0=0$, or equivalently
\[
 p_1(0,y)=0,\qquad \partial_xp_1(0,y)+yp_2(0,y)=0.
\]
For total degree at most $n$, these two polynomials in $y$ have
degrees at most $n$ and $n+1$, respectively. Testing them at $n+2$
distinct values is complete, since a nonzero polynomial of degree
at most $n+1$ has at most $n+1$ roots; the tests may be redundant. Here the divisor is the line $x=0$: $\partial_x$ differentiates away
from it, while the coefficient $y$ varies along it.

The same constraints can be calculated by polynomial division, without
first choosing jet points. Fix a monomial order and a module
Gr\"obner basis for the submodule generated by the rows of $P_L$
and the support relations $Ie_j$. Reduce each polynomial row
$a_Lb_i^+$ by this basis. Collect the coefficients of its remainder
as column $i$ of a matrix $R_L$. Linearity of the normal form gives
\[
 R_Lu=0
 \quad\text{if and only if}\quad
 a_L\sum_i u_i b_i^+\text{ is divisible by }P_L\text{ in }A.
\]
Thus the kernel gives all divisible combinations of the ambient rows.

Take a basis of this kernel and divide each associated polynomial row
by $P_L$. The quotients are unique in $A$ under the regularity
hypothesis. Express them in polynomial coordinates modulo the support
ideal $I$. One more linear calculation now selects those combinations
whose quotients belong to the prescribed space $X$. For component
degree bounds this means setting every forbidden quotient coefficient
to zero; in a projective description it includes the conditions at
infinity. These are the quotient conditions in
Proposition~\ref{prop:lifteddivision}, in addition to divisibility.
Finally take a basis of the linear equations annihilating the selected
coefficient subspace. It gives exactly the $c_L$ independent rows of
$\Gamma_L$. Transposing the entire calculation gives the column
algorithm for $\Gamma_R$. Both polynomial reduction and compatible
jets therefore compute the same image constraints.

If $\gamma_a$ denotes the combined tests on $a_LW_L$, the actual entries
are $(\Gamma_L)_{a,i}=\gamma_a(a_Lb_i^+)$; similarly
$(\Gamma_R)_{b,j}=\delta_b(a_Ra_j^+)$.

In particular, the first $N$ columns of $\Gamma_L$ contain these
tests on the original polynomials $b_i$, and the remaining $s_L$
columns contain the same tests on the added rational rows after
clearing denominators. This is the meaning of the split
$\Gamma_L=[X_0\ X_p]$; $\Gamma_R=[Y_0\ Y_p]$ is constructed
in exactly the same order on the right. All these entries depend
on the factors and the polynomial spaces, before any perturbed
orthogonal polynomial has been found.

Section~\ref{res:sec:residueconstruction} realizes such remainder
functionals by residues on the curve and surface quotients specified
there. Evaluating those functionals on the cleared ambient basis
provides the divisibility rows of the constraint matrix. The remaining
rows impose the prescribed quotient space, as in
Proposition~\ref{prop:lifteddivision}. Residue coordinates also evaluate
a homogeneous mass once its values on a quotient basis have been
chosen. In this way the geometric formulas supply concrete entries
for the finite construction. The rational moments in the next
subsubsection are the values of the chosen extension, and the mass
coefficients specify the chosen solution of the functional equation.

\subsubsection{Second-kind data on the algebraic support}
\label{sec:secondkind}

The extra vectors $\Phi_L,\Phi_R$ require moments with denominators.
These are the \emph{second-kind data}. For measure pairings they are
expressed through Cauchy transforms in one variable and Fantappi\`e
transforms in several variables, as made explicit below. Once the
finite ambient bases are fixed, only the three finite tables below
are needed from the rational extension. Thus their entries are known
inputs to the elimination, whether evaluated from a measure or supplied
as additional moments consistently with that extension.
Write their rows and columns as
$\phi_\alpha\coloneq \varphi_\alpha/a_L$ and $\psi_\beta\coloneq \chi_\beta/a_R$, with
polynomial numerators. The moment blocks are
\begin{equation}\label{eq:secondkindblocks}
\begin{aligned}
 (C_0)_{i,\beta}&\coloneq \sum_{u,v}\widetilde\mu_{u,v}
                     [\iota_L(b_{i,u})\iota_R(\chi_{\beta,v}/a_R)],\\
 (D_0)_{\alpha,j}&\coloneq \sum_{u,v}\widetilde\mu_{u,v}
                     [\iota_L(\varphi_{\alpha,u}/a_L)\iota_R(a_{j,v})],\\
 (E_0)_{\alpha,\beta}&\coloneq \sum_{u,v}\widetilde\mu_{u,v}
                     [\iota_L(\varphi_{\alpha,u}/a_L)\iota_R(\chi_{\beta,v}/a_R)].
\end{aligned}
\end{equation}
Thus $C_0$ pairs original rows with extra columns, $D_0$ pairs extra
rows with original columns, and $E_0$ pairs both extra lists. Their sizes
are $N\times s_R$, $s_L\times N$ and $s_L\times s_R$. Balanced pairings
use one point in these moments; bivariate pairings keep two arguments.

Equivalently, the three tables are
\[
 (C_0)_{i,\beta}=\widetilde\beta_\mu(b_i,\psi_\beta),\qquad
 (D_0)_{\alpha,j}=\widetilde\beta_\mu(\phi_\alpha,a_j),\qquad
 (E_0)_{\alpha,\beta}=\widetilde\beta_\mu(\phi_\alpha,\psi_\beta).
\]
For a scalar balanced measure and added functions
$\phi(z)=1/(z-\xi_L)$, $\psi(z)=1/(z-\xi_R)$, these become
\[
 C_i=\int\frac{b_i(z)}{z-\xi_R}\,\mathrm{d}\mu(z),\qquad
 D_j=\int\frac{a_j(z)}{z-\xi_L}\,\mathrm{d}\mu(z),\qquad
 E=\int\frac{\mathrm{d}\mu(z)}{(z-\xi_L)(z-\xi_R)}.
\]
Thus $C_0,D_0$ are the familiar Cauchy data of the original
polynomials; $E_0$ pairs the additional functions with one another.
Repeated poles give derivatives with respect to pole parameters. For
example,
\[
 \frac{1}{(m-1)!}\partial_\xi^{m-1}\frac1{z-\xi}
 =\frac1{(z-\xi)^m},\qquad m\ge1,
\]
so the factorial and sign are fixed by differentiating the denominator
with respect to $\xi$, rather than with respect to $z$.
For an abstract functional, these rational moments are part of the
chosen extension $\widetilde\mu$. Polynomial moments alone need
not fix them. When a measure supplies convergent integrals, those
integrals specify the extension values used here.

For a finite matrix measure on a compact support, the rational
moments depend holomorphically on denominator coefficients as long
as the denominators have no zeros on that support. Indeed, near any
such choice their absolute values have a common positive lower
bound, which justifies differentiation under the integral against
the total variation of the measure. For a fixed polynomial $h$,
$q_{\boldsymbol t}\coloneq\sum_a t_ag_a$ and a multi-index $\alpha$, this gives
\[
 \partial_{\boldsymbol t}^{\alpha}
 \int\frac{h}{q_{\boldsymbol t}}\,\dd\mu_{u,v}
 =(-1)^{|\alpha|}|\alpha|!
 \int\frac{h\prod_a g_a^{\alpha_a}}
 {q_{\boldsymbol t}^{|\alpha|+1}}\,\dd\mu_{u,v}.
\]
The same argument applies to both denominators in the last block of
\eqref{eq:secondkindblocks}, also when they belong to independent
bivariate arguments. Thus repeated poles are obtained by parameter
derivatives. To exhibit the multivariate transform explicitly, fix a
denominator $q_0$ with no zeros on the support and write a nearby one as
$q_0-\sum_a t_ag_a$. Then
\[
 \int\frac{h}{q_0-\sum_a t_ag_a}\,\dd\mu_{u,v}
 =\int\frac{h/q_0}{1-\sum_a t_a(g_a/q_0)}\,\dd\mu_{u,v}.
\]
The right-hand side is a Fantappi\`e transform: the denominator is one
minus a linear combination of the functions $g_a/q_0$, integrated
against the weighted measure $(h/q_0)\mu_{u,v}$
\cite{McCarthyPutinar}. Thus the multivariate blocks use Fantappi\`e
transforms. The two-denominator block uses the corresponding iterated
transform, with one parameter family for each denominator. In the
univariate specialization, partial fractions recover the Cauchy values
and their derivatives displayed above.

\subsubsection{Eliminating the original polynomial coefficients}

We now eliminate the coefficients of the original dual polynomials.
The coefficients of the extra rational vectors remain, together with
the mass and division constraints. The required moment blocks are
\begin{equation}\label{eq:blocks}
 G_{\mathrm{rat}}\coloneq \begin{bmatrix}I_N&C_0\\D_0&E_0\end{bmatrix},\qquad
 F\coloneq \begin{bmatrix}F_0\\F_p\end{bmatrix},\qquad
 H_\nu\coloneq [H_0\ H_p],\qquad G\coloneq G_{\mathrm{rat}}+FH_\nu.
\end{equation}
Here $F_{i,a}\coloneq \ell_a(\pi_L b_i^+)$ and
$(H_\nu)_{a,j}\coloneq r_a(\pi_R a_j^+)$, using \eqref{eq:massfactor}.
Thus $F$ has size $(N+s_L)\times h$, $H_\nu$ has size
$h\times(N+s_R)$, and $G$ agrees with the target pairing on the
two rational images. Indeed, for admissible coordinates $v,z$, the
added entry is
\[
 v^{\T}FH_\nu z
 =\sum_{a=1}^h\ell_a\bigl(\pi_L(v^{\T}\cB^+)\bigr)
                  r_a\bigl(\pi_R(\cA^+z)\bigr)
 =\nu\bigl(\pi_L(v^{\T}\cB^+),\pi_R(\cA^+z)\bigr).
\]
This identity explains the factorization of the mass table: its
entries are sums of products of a left scalar value and a right
scalar value. The rational part supplies the other term of the
transformed pairing. Write $\eta\coloneq [\eta_0\ \eta_p]=\lambda\pi_R$,
a $k\times(N+s_R)$ matrix in column coordinates.

The split of the mass tables can be written entry by entry as
\[
 \begin{aligned}
 (F_0)_{i,a}&=\ell_a(\pi_Lb_i),&
 (F_p)_{\alpha,a}&=\ell_a(\pi_L\phi_\alpha),\\
 (H_0)_{a,j}&=r_a(\pi_Ra_j),&
 (H_p)_{a,\beta}&=r_a(\pi_R\psi_\beta).
 \end{aligned}
\]
For point and derivative masses, these entries are the corresponding
values and derivatives of the recovered polynomial rows and columns.
If there is no mass, $h=0$ and these blocks are absent. The entries
of $\eta_0,\eta_p$ are obtained by applying the prescribed
normalization to $\pi_Ra_j,\pi_R\psi_\beta$. For a monic scalar
problem, this means extracting the coefficient of the prescribed
highest power.

The identity block in $G_{\mathrm{rat}}$ comes from the original
duality $\beta_\mu(b_i,a_j)=\delta_{i,j}$. Every other block is
therefore already known from the original polynomials, their
second-kind data, the factors, the masses and normalization.

For the normalized problem $k>0$, write the rational product of a candidate as
$S_R\widehat a=\cA^+z$, where
$z\coloneq [z_0^{\T}\ z_p^{\T}]^{\T}$ has $k$ columns. Its conditions are $Gz=\Gamma_L^{\T}\alpha$,
$\Gamma_Rz=0$ and $\eta z=I_k$.
Here $\alpha$ is a $c_L\times k$ matrix of multipliers. The first
equation expresses orthogonality: a linear functional vanishing on
$\ker\Gamma_L$ is a linear combination of the rows of $\Gamma_L$.
The second enforces division and the target bounds; the third fixes
the normalization. Introduce $u\coloneq H_\nu z$, the $h$ mass
coordinates. In blocks, the full system is
\[
 \begin{aligned}
 z_0+C_0z_p+F_0u&=X_0^{\T}\alpha,\\
 D_0z_0+E_0z_p+F_pu&=X_p^{\T}\alpha,\\
 H_0z_0+H_pz_p-u&=0,\\
 Y_0z_0+Y_pz_p&=0,\\
 \eta_0z_0+\eta_pz_p&=I_k.
 \end{aligned}
\]
The first line can be solved directly because its coefficient of
$z_0$ is the identity. It gives
\[
 z_0=-C_0z_p-F_0u+X_0^{\T}\alpha.
\]
For example, substitution into the second line gives
\[
 (E_0-D_0C_0)z_p+(F_p-D_0F_0)u
       +(D_0X_0^{\T}-X_p^{\T})\alpha=0.
\]
For completeness, the other three substituted equations are
\[
 \begin{aligned}
 (H_p-H_0C_0)z_p+(-I_h-H_0F_0)u+H_0X_0^{\T}\alpha&=0,\\
 (Y_p-Y_0C_0)z_p-Y_0F_0u+Y_0X_0^{\T}\alpha&=0,\\
 (\eta_p-\eta_0C_0)z_p-\eta_0F_0u+\eta_0X_0^{\T}\alpha&=I_k.
 \end{aligned}
\]
These are, respectively, the mass-coordinate identity, the right
image constraints and the normalization after elimination. Thus its four row groups
come, in order, from the remaining pairing equations, the definition
of the mass coordinates, the right divisibility and degree conditions,
and normalization. No perturbed polynomial is needed to form them.
For instance,
\[
 (E_0-D_0C_0)_{\alpha,\beta}
 =\widetilde\beta_\mu(\phi_\alpha,\psi_\beta)
 -\sum_{j=1}^N\widetilde\beta_\mu(\phi_\alpha,a_j)
                    \widetilde\beta_\mu(b_j,\psi_\beta),
\]
which subtracts the contribution of the original polynomial block
from the pairing of two added rational vectors.

Retain the ordered bases $f_1,\ldots,f_r$ of $X$ and
$g_1,\ldots,g_{r+k}$ of $Y$ used to construct the image coordinates.
The target coefficient matrix is
\begin{equation}\label{eq:normalizedGram}
 \widehat M_{R,k}\coloneq 
 \begin{bmatrix}
 [\widehat\beta(f_i,g_j)]_{i\in\{1,\ldots,r\},\,j\in\{1,\ldots,r+k\}}\\
 [\lambda(g_1)\ \cdots\ \lambda(g_{r+k})]
 \end{bmatrix}.
\end{equation}
\Needspace{4\baselineskip}

The remaining coefficient matrix is the square matrix
\begin{equation}\label{eq:small}
\mathbb T_k\coloneq 
\begin{bmatrix}
E_0-D_0C_0&F_p-D_0F_0&D_0X_0^{\T}-X_p^{\T}\\
H_p-H_0C_0&-I_h-H_0F_0&H_0X_0^{\T}\\
Y_p-Y_0C_0&-Y_0F_0&Y_0X_0^{\T}\\
\eta_p-\eta_0C_0&-\eta_0F_0&\eta_0X_0^{\T}
\end{bmatrix}.
\end{equation}
Its column sizes are $s_R,h,c_L$, and its row sizes are
$s_L,h,c_R,k$. These give the same total:
\[
 s_R+h+c_L=s_L+h+c_R+k=N-r+s_L+s_R+h.
\]
Define this order as $t\coloneq N-r+s_L+s_R+h$ and set
$\tau_{R,k}\coloneq\det\mathbb T_k$. The term \emph{tau determinant}
means the ordinary determinant of the equations left after elimination.
The theorem below relates its nonvanishing to the prescribed
orthogonality and normalization.
The normalized Gram matrix \eqref{eq:normalizedGram} has order
$r+k$, whereas $\mathbb T_k$ has order $t$. Which one is smaller
depends on the original spaces and the transformation.

The dimensions in this count are computed using polynomial classes in
$A=\CC[x_1,\ldots,x_M]/I$. For component degree bounds, one counts the
independent allowed classes in each of the $q$ row components and $p$
column components; the relations in $I$ identify dependent monomials.
The rational factors then determine the extra dimensions $s_L,s_R$
relative to the chosen original polynomial spaces, while $h$ is the
number of terms in the chosen mass representation on these finite
spaces. Thus the number $M$ of ambient variables and the support
dimension $\dim A$ enter through the polynomial-space dimensions.
The rectangular size and these two dimensions alone do not fix the
order: the component bounds, the support relations, the rational
factors and the mass data also enter the count.
Subsection~\ref{sec:orders} makes this count numerical: it specifies
the ambient spaces by degrees, gives their exact determinant order,
and calculates the reduction obtained from polynomial coefficient ranks.

The passage from the full equations to $\mathbb T_k$ is block
elimination with the known invertible block $I_N$, or equivalently
a Schur-complement calculation. Its application here uses the image
equations already constructed, the supplied rational moment tables,
and the mass and normalization tables. The following theorem records
the resulting determinant relation and reconstructs the normalized
polynomial from the remaining coefficients. Since the eliminated
block belongs to the original pairing, the calculation applies
directly to the final bilateral problem without solving an
intermediate transformed orthogonality problem.

\begin{theorem}[Determinant criterion and normalized orthogonality]\label{thm:finite}\label{thm:taucriterion}

There is a nonzero factor $c$, determined by the ordered bases,
constraint complements and block orders, such that
\begin{equation}\label{eq:tauGram}
 \det\widehat M_{R,k}=c\,\tau_{R,k}.
\end{equation}
For $k>0$, the following conditions are equivalent:
\begin{equation}\label{eq:taunormality}
 \begin{gathered}
 \text{a unique normalized orthogonal block exists}\\
 \text{if and only if}\quad\det\widehat M_{R,k}\ne0
 \quad\text{if and only if}\quad\tau_{R,k}\ne0.
 \end{gathered}
\end{equation}
The left-normalized problem has the same criterion for the determinant
$\tau_{L,k}$ of its transposed reduced construction.
Consequently, for any prescribed collection $\Omega$ of spaces and
normalizations to which this construction applies,
\begin{equation}\label{eq:familytau}
 \begin{gathered}
 \text{the entire family is normal if and only if}\\
 \tau_\omega\ne0\quad\text{for every }\omega\in\Omega.
 \end{gathered}
\end{equation}
When $\tau_{R,k}\ne0$, put
\[
 V_*\coloneq 
 \begin{bmatrix}-C_0&-F_0&X_0^{\T}\\I_{s_R}&0&0\end{bmatrix},
 \qquad J_k\coloneq \begin{bmatrix}0_{(t-k)\times k}\\I_k\end{bmatrix}.
\]
The solution is
\begin{equation}\label{eq:normalizedanswer}
       \widehat a=\pi_R\bigl(\cA^+V_*\mathbb T_k^{-1}J_k\bigr).
\end{equation}
For $k=0$, omit the last block row of \eqref{eq:small}. The target square
pairing is normal exactly when $\det\mathbb T_0\ne0$.
\end{theorem}

\begin{proof}
We first compare solutions, then compare determinants. The preceding
substitution gives the reduced system
\[
 \mathbb T_k
 \begin{bmatrix}z_p\\u\\\alpha\end{bmatrix}=J_k.
\]
Given a solution of the original normalized problem, its coordinates
$z$ determine $u=H_\nu z$. Orthogonality determines $\alpha$
uniquely because $\Gamma_L$ has independent rows. These coordinates
solve the reduced system. Conversely, a solution of the reduced
system determines $z_0$ by the displayed elimination formula; the
four remaining equations then recover all the original conditions.
In particular $\Gamma_Rz=0$, so inverse multiplication gives a
polynomial in $Y$ with the required normalization.

The two systems therefore have the same solutions and the same
homogeneous freedom. For $k>0$, a normalized block is equivalent to a solution of
\[
 \widehat M_{R,k}C=\begin{bmatrix}0_{r\times k}\\I_k\end{bmatrix}.
\]
If the square coefficient matrix is invertible, this equation has
exactly one solution. Conversely, suppose a solution is unique.
A nonzero vector in its homogeneous kernel could be added to any
one of the $k$ columns of that solution, giving another solution.
Thus the homogeneous kernel is zero, and the matrix is invertible.
The correspondence just proved transfers the same conclusion to
$\mathbb T_k$. This proves the uniqueness criterion. For $k=0$,
the same correspondence identifies the right radicals of the two
square pairings and proves the normality criterion. Reconstructing
$z$ from the reduced unknowns is exactly multiplication by $V_*$;
applying $\pi_R$ gives \eqref{eq:normalizedanswer}.

To compare the determinants, use the coordinate matrices $Z_L,Z_R$
constructed from $f_iS_L,S_Rg_j$. Their columns are bases
of the respective constraint kernels. Choose matrices $Q_L,Q_R$
with $\Gamma_LQ_L=I_{c_L}$ and $\Gamma_RQ_R=I_{c_R}$.
Then
\[
 U_L\coloneq[Z_L\ Q_L],\qquad U_R\coloneq[Z_R\ Q_R]
\]
are invertible changes of coordinates: they separate the admissible
directions from the directions detected by the constraints.
The coefficient matrix before eliminating $z_0$ is
\[
 Q\coloneq\begin{bmatrix}
 G_{\mathrm{rat}}&F&-\Gamma_L^{\T}\\
 H_\nu&-I_h&0\\
 \Gamma_R&0&0\\
 \eta&0&0
 \end{bmatrix}.
\]
Elimination with the leading $I_N$ block uses determinant-preserving
row operations and leaves $\mathbb T_k$. Thus $\det Q=\tau_{R,k}$.
Alternatively, eliminate the mass variables using $-I_h$.
The pairing block becomes $G_{\mathrm{rat}}+FH_\nu=G$.
Change its row and column coordinates by $U_L^{\T}$ and $U_R$.
In these coordinates, $\Gamma_LU_L=[0\ I_{c_L}]$ and
$\Gamma_RU_R=[0\ I_{c_R}]$. After the mass elimination the
remaining matrix is
\[
 \begin{bmatrix}
 Z_L^{\T}GZ_R&Z_L^{\T}GQ_R&0\\
 Q_L^{\T}GZ_R&Q_L^{\T}GQ_R&-I_{c_L}\\
 0&I_{c_R}&0\\
 \eta Z_R&\eta Q_R&0
 \end{bmatrix}.
\]
Expand first along the multiplier columns, whose only nonzero
block is $-I_{c_L}$, removing the complementary left rows.
Then expand along the right-constraint rows, whose only nonzero
block is $I_{c_R}$, removing the complementary right columns.
The remaining matrix is
\[
 \begin{bmatrix}Z_L^{\T}GZ_R\\\eta Z_R\end{bmatrix}
 =\widehat M_{R,k}.
\]
We can record the sign as well. Moving the $h$ mass rows and columns
to the end of $Q$ gives the sign $(-1)^{h(c_R+k+c_L)}$;
eliminating their block $-I_h$ gives another factor $(-1)^h$.
After the changes of coordinates, reorder the remaining row groups
from $(r,c_L,c_R,k)$ to $(r,k,c_R,c_L)$. This contributes
$(-1)^{k(c_L+c_R)+c_Lc_R}$. The final two diagonal blocks are
$I_{c_R}$ and $-I_{c_L}$, with determinant $(-1)^{c_L}$.
Consequently, for precisely the block orders displayed above,
\[
 \begin{aligned}
 c&=(-1)^\epsilon\det U_L\det U_R,\\
 \epsilon&\coloneq
 h(1+c_L+c_R+k)+k(c_L+c_R)+c_Lc_R+c_L.
 \end{aligned}
\]
Both coordinate matrices are invertible, so $c\ne0$. This proves
\eqref{eq:tauGram}, including its determinant normalization.
Transposing exchanges the two arguments and proves the left
criterion. Applying the criterion at each prescribed index proves
\eqref{eq:familytau}.
\end{proof}

By \eqref{eq:tauGram}, nonvanishing is independent of bases and complete constraints.

\Needspace{4\baselineskip}
The inverse in \eqref{eq:normalizedanswer} can be replaced by
scalar determinants, one for each polynomial entry. The following
corollary records the required rows and columns, including the
left-normalized problem. Its last assertion also verifies
polynomial divisibility of the determinant numerators, so that the
formulas produce polynomials with the prescribed bounds.

\begin{corollary}[Bordered determinants for both normalized blocks]
\label{cor:polynomialborders}
Put $\mathcal W_R\coloneq \pi_R(\cA^+V_*)$. For the left problem
$\dim X=r+k$, $\dim Y=r$, prescribe $\lambda_L:X\to\CC^k$,
set $\zeta\coloneq \lambda_L\pi_L$ and form
$\mathbb T^L_k$ from \eqref{eq:small} by
\begin{equation}\label{eq:leftsubstitutions}
\begin{gathered}
 (C_0,D_0,E_0)\mapsto(D_0^{\T},C_0^{\T},E_0^{\T}),\quad
 (F_0,F_p,H_0,H_p)\mapsto(H_0^{\T},H_p^{\T},F_0^{\T},F_p^{\T}),\\
 (X_0,X_p,Y_0,Y_p,\eta)\mapsto(Y_0,Y_p,X_0,X_p,\zeta).
\end{gathered}
\end{equation}
Set $V^L_*\coloneq [-D_0^{\T}\ {-H_0^{\T}}\ Y_0^{\T};I_{s_L}\ 0\ 0]$
and $\mathcal W_L\coloneq \pi_L((V^L_*)^{\T}\cB^+)$.
Writing $\mathbb T^R_k\coloneq \mathbb T_k$, $\tau_{\sigma,k}\coloneq \det\mathbb T^\sigma_k$,
for a row $u$ of length $t$ and $j\in\{1,\ldots,k\}$, define the scalar border
\[
 \mathcal B_\sigma(u,j)\coloneq 
 -\tau_{\sigma,k}^{-1}
 \det\begin{bmatrix}\mathbb T^\sigma_k&e_{t-k+j}\\u&0\end{bmatrix}.
\]
Then the right and left normalized blocks are
\begin{equation}\label{eq:rightborder}
 \widehat a_{i,j}(P)=\mathcal B_R((\mathcal W_R(P))_{i,:},j),\qquad
 \widehat b_{i,j}(Q)=\mathcal B_L((\mathcal W_L(Q))_{:,j}^{\T},i).
\end{equation}
Here $i,:$ selects row $i$ and $:,j$ selects column $j$;
$\widehat a$ and $\widehat b$ have sizes $p\times k$ and $k\times q$.
Their denominators are nonzero exactly when the respective normalized
problems have unique solutions.

Divisibility is also explicit inside the borders. Replace $\mathcal W_R$
by $a_R\cA^+V_*$ and $\mathcal W_L$ by
$a_L(V^L_*)^{\T}\cB^+$ in \eqref{eq:rightborder}, and multiply the
results by $\tau_{R,k}$ and $\tau_{L,k}$, respectively. The resulting
polynomial matrices satisfy
\begin{equation}\label{eq:numeratordivision}
 N_R=\tau_{R,k}P_R\widehat a,\qquad
 N_L=\tau_{L,k}\widehat bP_L.
\end{equation}
Thus their columns and rows lie in $P_RY$ and $XP_L$.
\end{corollary}

\begin{proof}
For an invertible square matrix $T$, a column $v$ and a row $u$,
subtracting $uT^{-1}$ times the upper block row from the last row
gives
\[
 \det\begin{bmatrix}T&v\\u&0\end{bmatrix}
 =-\det(T)\,uT^{-1}v.
\]
The $j$th column of $J_k$ is $e_{t-k+j}$. Apply this identity
with $T=\mathbb T_k$, this column as $v$, and the $i$th row
of $\mathcal W_R(P)$ as $u$. The resulting scalar is exactly
the $(i,j)$ entry of \eqref{eq:normalizedanswer}, which proves
the right formula. Exchanging rows and columns transposes the
pairing and all its data. The substitutions in
\eqref{eq:leftsubstitutions} list these exchanges explicitly,
so the same calculation proves the left formula.

For the divisibility assertion, apply the border identity before
recovering the polynomial candidate, and multiply by the common
scalar denominator. The numerator is then
$N_R=\tau_{R,k}a_RS_R\widehat a$ on the right and
$N_L=\tau_{L,k}a_L\widehat bS_L$ on the left.
Since $a_RS_R=P_R$ and $a_LS_L=P_L$, these are
\eqref{eq:numeratordivision}. Membership of $\widehat a$ in $Y$
and of $\widehat b$ in $X$ also proves the stated quotient bounds.
\end{proof}

\subsubsection{The kernel from the same reduced matrix}

When the target pairing is square, write $\mathbb T\coloneq \mathbb T_0$ and put
\begin{equation}\label{eq:kernelmatrices}
 U_*\coloneq 
 \begin{bmatrix}-D_0&I_{s_L}\\-H_0&0\\-Y_0&0\end{bmatrix},
 \qquad J\coloneq \begin{bmatrix}I_N&0\\0&0\end{bmatrix}.
\end{equation}
Here $J$ has size $(N+s_R)\times(N+s_L)$. Let $K_E(P,Q)\coloneq \cA(P)\cB(Q)$, and let $\widehat K$ be the kernel of the
target pairing, with its column factor at $P$ and its row factor at $Q$.
This is the dual-basis form of \eqref{eq:kerneldefinition}.
The reproduction identities and their meaning on a nonreduced
support are given in Subsubsection~\ref{sec:kernelsetup}.

\Needspace{4\baselineskip}
The separation rank from Subsubsection~\ref{sec:kernelsetup} equals
the rank of the coefficient matrix in independent ambient bases.
The present theorem constructs that matrix.

For the kernel we solve the same orthogonality equations with an
arbitrary right-hand side, instead of a fixed normalization.
The resulting inverse pairing reproduces all members of the target
space simultaneously. This is why the same reduced matrix suffices
for both the polynomial and kernel formulas.

\begin{theorem}[Simultaneous kernel formula]\label{thm:kernel}
If the target pairing is normal, then
\begin{align}
 \cF(P,Q)&\coloneq S_R(P)\widehat K(P,Q)S_L(Q)
       =\cA^+(P)\bigl[J+V_*\mathbb T^{-1}U_*\bigr]\cB^+(Q),\label{eq:kernel}\\
 \cF-K_E&=\cA^+V_*\mathbb T^{-1}U_*\cB^+ .\label{eq:correction}
\end{align}
The correction has tensor rank $\rank(V_*\mathbb T^{-1}U_*)\le t$.
The full tensor $\cF$ has rank $r$ and factor spaces $S_RY$ and $XS_L$.
\end{theorem}

\begin{proof}
Choose coordinate inclusion matrices $Z_L,Z_R$ for bases of
$XS_L,S_RY$ in the enlarged bases. Their pairing matrix is
$Z_L^{\T}GZ_R$, which is invertible by target normality.
The coefficient matrix of
$S_R(P)\widehat K(P,Q)S_L(Q)$ in the enlarged bases is therefore
\[
 \mathcal P\coloneq
 Z_R(Z_L^{\T}GZ_R)^{-1}Z_L^{\T}.
\]
Indeed, write an admissible coefficient column as $z=Z_Rw$.
Its pairing values against the target row basis form
$Z_L^{\T}GZ_Rw$. To reproduce prescribed values, multiply by
the inverse of this square Gram matrix and then by $Z_R$.
The factor $Z_L^{\T}$ converts an ambient right-hand side to
the prescribed row values. This gives the displayed coefficient
matrix with the basis inclusions in the stated order.

We calculate $\mathcal Pv$ for an arbitrary column
$v=[v_0^{\T}\ v_p^{\T}]^{\T}$. The vector $z=\mathcal Pv$
is characterized by
\[
 Gz-\Gamma_L^{\T}\alpha=v,\qquad \Gamma_Rz=0.
\]
The second equation gives $z=Z_Rw$. Multiplying the first by
$Z_L^{\T}$ yields
$(Z_L^{\T}GZ_R)w=Z_L^{\T}v$, whose unique solution gives
$z=\mathcal Pv$. Conversely, for this solution,
$Z_L^{\T}(Gz-v)=0$. The vectors annihilated by $Z_L^{\T}$
are exactly the columns of $\Gamma_L^{\T}$; hence there is a
unique $\alpha$ with $Gz-v=\Gamma_L^{\T}\alpha$.
This proves both implications of the characterization.

Introduce $u=H_\nu z$ and eliminate the original polynomial
coefficients as before. The first block now gives
$z_0=v_0-C_0z_p-F_0u+X_0^{\T}\alpha$. The remaining pairing, mass and image equations have right-hand
sides $v_p-D_0v_0$, $-H_0v_0$ and $-Y_0v_0$, respectively.
Stacking these gives
\[
 \mathbb T\begin{bmatrix}z_p\\u\\\alpha\end{bmatrix}
 =\begin{bmatrix}v_p-D_0v_0\\-H_0v_0\\-Y_0v_0\end{bmatrix}
 =U_*v.
\]
Here $U_*$ has size $t\times(N+s_L)$ and $V_*$ has size
$(N+s_R)\times t$, so their product through $\mathbb T^{-1}$
has the same size as $\mathcal P$.
The original contribution is $Jv=[v_0^{\T}\ 0]^{\T}$, while
the remaining contribution is $V_*\mathbb T^{-1}U_*v$.
Since this holds for every $v$, it proves
$\mathcal P=J+V_*\mathbb T^{-1}U_*$. Multiplication by the two
ambient basis matrices proves \eqref{eq:kernel} and
\eqref{eq:correction}.

The ambient basis functions are linearly independent over $\CC$.
Consequently the smallest number of separated products equals the
rank of their coefficient matrix. This proves the correction-rank
assertion. In the expression defining $\mathcal P$, the middle
matrix is invertible and both inclusion matrices have rank $r$.
Thus $\mathcal P$ has rank $r$, with column and row spaces exactly
those of $Z_R$ and $Z_L^{\T}$, respectively.
\end{proof}

Set $\tau_0\coloneq \det\mathbb T$,
$\mathcal U(P)\coloneq \cA^+(P)V_*$ and
$\mathcal V(Q)\coloneq U_*\cB^+(Q)$. For every matrix entry,
\begin{equation}\label{eq:kernelborder}
 \cF_{i,j}(P,Q)
 =\frac{\det\begin{bmatrix}
       \mathbb T&\mathcal V_{:,j}(Q)\\
       -\mathcal U_{i,:}(P)&(K_E)_{i,j}(P,Q)
       \end{bmatrix}}{\tau_0}.
\end{equation}
The numerator determinant has order $t+1$. Add
$\mathcal U_{i,:}\mathbb T^{-1}$ times its upper block row to its
last row. The lower-left block becomes zero, and the determinant is
\[
 \det\mathbb T\,
 \bigl((K_E)_{i,j}+\mathcal U_{i,:}\mathbb T^{-1}
                                      \mathcal V_{:,j}\bigr).
\]
The minus sign in the lower-left block therefore produces the plus
sign in \eqref{eq:kernel}. Replacing the lower-right entry by zero
gives $(\cF-K_E)_{i,j}$.
An entrywise determinant for $\widehat K$ itself is obtained in
\eqref{eq:kernelborder} by replacing
$\mathcal U\mapsto S_R^{-1}\mathcal U$,
$\mathcal V\mapsto\mathcal V S_L^{-1}$ and
$K_E\mapsto S_R^{-1}K_ES_L^{-1}$.
Theorem~\ref{thm:kernel} places the factor spaces in $S_RY$ and
$XS_L$; multiplying by the inverse factors therefore cancels all poles and gives a tensor
in $Y\otimes_\CC X$.

 Empty matrices are permitted;
the determinant of the zero-order square matrix is one.

\subsubsection{A complete bilateral calculation for mixed-type orthogonality on a curve}
\label{sec:mixedcurveexample}

We now construct every input block and recover the transformed
polynomial on the curve
$A=\CC[x,y]/(x^2+y^2-1)$. Rows have two components and columns
have three. Take
\[
 L_1=\begin{bmatrix}x&0\\0&1\end{bmatrix},\quad
 L_0=\begin{bmatrix}1&y\\1&1\end{bmatrix},\qquad
 R_1=\diag(x,1,1),\quad
 R_0=\begin{bmatrix}1&1&0\\0&1&y\\1&0&1\end{bmatrix}.
\]
Their determinants are $x,1-y,x,1+y$, respectively, all nonzero
in this integral algebra. The products do not commute:
$(L_1L_0-L_0L_1)_{1,2}=(x-1)y$ and
$(R_1R_0-R_0R_1)_{1,2}=x-1$.
Keeping the specified order gives
\begin{equation}\label{eq:circlemaps}
 S_L=\begin{bmatrix}x^{-1}&y/x\\1&1\end{bmatrix},\qquad
 S_R=\begin{bmatrix}x^{-1}&1&0\\0&1&y\\x^{-1}&0&1\end{bmatrix}.
\end{equation}

For the original pairing use the two points $P_\pm=(\pm1,0)$ and
the rectangular matrix $W=\bigl[\begin{smallmatrix}1&0&0\\0&1&1\end{smallmatrix}\bigr]$:
\[
 \beta_\mu(f,g)\coloneq\tfrac12\sum_{P\in\{P_+,P_-\}}f(P)Wg(P).
\]
Evaluation at the same points defines its rational extension, since
$x$ does not vanish there. Choose the polynomial spaces with ordered bases
\[
 X=\Span\{[1\ 0],[x\ 0],[0\ 1]\},\qquad
 Y=\Span\{[1\ 0\ 0]^{\T},[x\ 0\ 0]^{\T},
                  [0\ 1\ 0]^{\T},[0\ 0\ 1]^{\T}\}.
\]
The normalization $\lambda(g)=g_3$ prescribes the constant third
component to be one. For the original dual pair take $E=X$ and
$H$ spanned by the first three displayed columns of $Y$.
The averages of $x$ and $x^2$ at $P_\pm$ are zero and one;
direct multiplication by $W$ therefore gives $G_0=I_3$.
The original pairing has global bilinear rank four: evaluation at
$P_+$ and $P_-$ can be prescribed independently, and each copy of
$W$ has rank two. Its normality in this calculation is the
nondegeneracy of the specified three-dimensional pair $E\times H$.
The target is likewise the single mixed-type orthogonality problem on $X\times Y$
with the displayed normalization.

At $Q=(0,1)$ the tangent derivative is $\partial_x$: it kills
$x^2+y^2-1$ at $Q$ and is well defined on classes. For a parameter
$\gamma\in\CC$, choose the homogeneous mass
\[
 \nu_\gamma(f,g)\coloneq
 \gamma\,\partial_x(f_1g_1)(Q),\qquad
 \widehat\beta(f,g)\coloneq
 \widetilde\beta_\mu(fS_L,S_Rg)+\nu_\gamma(f,g).
\]
The first components of $fL_1$ and $R_1g$ are $xf_1$ and $xg_1$.
Their product has a double zero at $Q$, so $L_1\nu_\gamma R_1=0$.
For $\gamma\ne0$, either unilateral multiplication instead leaves
$\gamma f_1(Q)g_1(Q)$, which is nonzero for constant first
components. Thus this mass is allowed by the bilateral equation
and is lost if unilateral vanishing is imposed separately.
It cannot even be a sum of balanced unilateral masses: the
$(1,1)$ entry of either unilateral term is annihilated by $x$,
whereas $x(\nu_\gamma)_{1,1}=\gamma\delta_Q\ne0$.
This is the obstruction of Theorem~\ref{thm:tormasses} in the
present curve calculation.

The finite rank of the whole transformed pairing can also be read
from these data. Both rational factors are invertible at $P_\pm$,
so the rational part still has rank four. In the coordinates given
by the value and the tangent derivative of the first component at
$Q$, the mass has matrix
\[
 \begin{bmatrix}0&\gamma\\\gamma&0\end{bmatrix}.
\]
The polynomial classes $x^2(1+x)/2$ and $x^2(1-x)/2$ independently
prescribe the two values at $P_\pm$ and have zero value and tangent
derivative at $Q$. Conversely, $y$ and $xy$ vanish at $P_\pm$;
their respective value--derivative pairs at $Q$ are $[1\ 0]$ and
$[0\ 1]$. Hence the two-point evaluation data and the two jet
coordinates are independent. For $\gamma\ne0$ the global rank of
$\widehat\beta$ is therefore $4+2=6$. Normality of the finite spaces
verified below is compatible with this rank bound; an exhaustive
sequence of normal square pairings would require unbounded rank.

\noindent\textbf{Rational vectors and image conditions.}
Append to the original bases the rows and columns
\[
 \phi_1=[x^{-1}\ y/x],\quad \phi_2=[0\ y],\qquad
 \psi_1=[x^{-1}\ 0\ x^{-1}]^{\T},\quad
 \psi_2=[0\ 0\ 1]^{\T},\quad \psi_3=[0\ y\ 0]^{\T}.
\]
These are independent additions in the coordinate algebra, giving
$N=r=3$, $s_L=2$, $s_R=3$ and $k=1$. Multiplication of the
three basis rows of $X$ by $S_L$ gives
$\phi_1,b_1+\phi_2,b_1+b_3$.
Multiplication of the four basis columns of $Y$ by $S_R$ gives
$\psi_1,a_1+\psi_2,a_1+a_3,\psi_2+\psi_3$.
Consequently the complete image constraints are
\begin{equation}\label{eq:circleconstraints}
 \Gamma_L=\begin{bmatrix}0&1&0&0&0\\1&0&-1&0&-1\end{bmatrix},\qquad
 \Gamma_R=\begin{bmatrix}0&1&0&0&0&0\\1&0&-1&0&-1&1\end{bmatrix}.
\end{equation}
For example, left coordinates $t$ must satisfy
$t_2=0$ and $t_1=t_3+t_5$; the resulting vector is
$t_4\phi_1+t_5(b_1+\phi_2)+t_3(b_1+b_3)$.
Its polynomial preimage is $[t_4+xt_5\ t_3]$.
On the right, coordinates $v$ satisfy $v_2=0$ and
$v_1=v_3+v_5-v_6$. This gives the recovery maps
\[
 \pi_L(t)=[t_4+xt_5\ t_3],\qquad
 \pi_R(v)=[v_4+x(v_5-v_6)\ v_3\ v_6]^{\T}.
\]
They enforce the stated component bounds as well as polynomiality.
Using the same expressions off the two kernels specifies their
linear extensions and gives $\eta=[0\ 0\ 0\ 0\ 0\ 1]$.

\noindent\textbf{Moments, mass and the reduced matrix.}
All rational moments are averages at $P_\pm$. Since $y=0$ and
$x^{-1}=x$ there, their blocks are
\[
 C_0=\begin{bmatrix}0&0&0\\1&0&0\\0&1&0\end{bmatrix},\quad
 D_0=\begin{bmatrix}0&1&0\\0&0&0\end{bmatrix},\quad
 E_0=\begin{bmatrix}1&0&0\\0&0&0\end{bmatrix}.
\]
For instance, the first column of $C_0$ is the pairing of
$b_1,b_2,b_3$ with $\psi_1$, whose values are the averages of
$x^{-1},1,x^{-1}$. The other entries follow from the same rule.
For the recovered polynomials the mass is
\[
 \nu_\gamma(\pi_Lt,\pi_Rv)
       =\gamma\{t_5v_4+t_4(v_5-v_6)\}.
\]
For $\gamma\ne0$ the mass has rank two. Taking the left mass
coordinates to be $[t_4\ t_5]$ gives
$h=2$, $F_0=0$, $F_p=I_2$, $H_0=0$ and
$H_p=\gamma\bigl[\begin{smallmatrix}0&1&-1\\1&0&0\end{smallmatrix}\bigr]$.
Here $h=2$ is the chosen factorization size; at $\gamma=0$ the
mass vanishes and its two coordinates can be omitted.
Substitution in \eqref{eq:small}, with column order
$(z_1,z_2,z_3,u_1,u_2,\alpha_1,\alpha_2)$, gives
\begin{equation}\label{eq:circleresidual}
 \mathbb T_1=
 \begin{bmatrix}
 0&0&0&1&0&1&0\\
 0&0&0&0&1&0&1\\
 0&\gamma&-\gamma&-1&0&0&0\\
 \gamma&0&0&0&-1&0&0\\
 -1&0&0&0&0&1&0\\
 0&0&1&0&0&0&2\\
 0&0&1&0&0&0&0
 \end{bmatrix},\qquad \det\mathbb T_1=-2\gamma^2.
\end{equation}
Its first two rows are the remaining pairing equations; the next
two define the mass coordinates. The last three impose the two
right image conditions and the normalization.
For an explicit determinant calculation, expand along its last row,
then along the resulting row with the sole entry $2$ in the
$\alpha_2$ column. Expanding next along the $u_2$ entry in the
second row and the remaining $\gamma$ entry in the $z_1$ column gives
\[
 \det\mathbb T_1
 =2\gamma\det\begin{bmatrix}0&1&1\\\gamma&-1&0\\0&0&1\end{bmatrix}
 =-2\gamma^2.
\]

\noindent\textbf{Reconstruction and direct verification.}
For $\gamma\ne0$, solve $\mathbb T_1w=e_7$. The last two rows give
$z_3=1$ and $\alpha_2=-1/2$. The second and fourth then give
$u_2=1/2$ and $z_1=1/(2\gamma)$. The fifth and first give
$\alpha_1=z_1$ and $u_1=-z_1$; the third gives
$z_2=1-1/(2\gamma^2)$.
The eliminated original coordinates are
$z_0=-C_0z_p+X_0^{\T}\alpha$.
Applying $\pi_R$ yields the polynomial column
\begin{equation}\label{eq:circleanswer}
 \widehat a=
 \begin{bmatrix}
 (\gamma-x)/(2\gamma^2)\\
 (1-\gamma^2)/(2\gamma^2)\\1
 \end{bmatrix}.
\end{equation}
It has exactly the required component bounds and normalization.
As an independent check, evaluation of the transformed pairing
on the prescribed bases gives
\[
 [\widehat\beta(f_i,g_j)]_{i,j}
 =\begin{bmatrix}1&\gamma&0&0\\\gamma&1&1&0\\0&2&2&1\end{bmatrix}.
\]
Its product with the coefficient column of \eqref{eq:circleanswer}
is zero. To verify the entire table in one calculation, write
\[
 f=[a+bx\ c],\qquad g=[u+vx\ w\ t]^{\T}.
\]
At $P_\pm$, the products in \eqref{eq:circlemaps} give
\[
 fS_L=[a/x+b+c\ c],\qquad
 S_Rg=[u/x+v+w\ w\ u/x+v+t]^{\T}.
\]
The averages of $1/x$ and $1/x^2$ are zero and one. Multiplying
by $W$ and averaging therefore gives
$au+bv+bw+2cv+2cw+ct$. The mass adds $\gamma(bu+av)$,
which proves every entry of the displayed Gram matrix.
For the column in \eqref{eq:circleanswer}, its coefficients are
\[
 u=\frac1{2\gamma},\quad v=-\frac1{2\gamma^2},\quad
 w=\frac{1-\gamma^2}{2\gamma^2},\quad t=1.
\]
The three orthogonality expressions are respectively
$u+\gamma v=0$, $\gamma u+v+w=0$ and $2v+2w+t=0$.
Appending $[0\ 0\ 0\ 1]$ to the Gram matrix gives determinant
$1(2-2)-\gamma(2\gamma)=-2\gamma^2$, exactly the reduced
determinant. At $\gamma=0$
the first two equations force the first coefficient to vanish
and the second to be the negative of the third; the last pairing
equation then contradicts the prescribed fourth coefficient one.
Thus the mass restores solvability in this finite mixed-type orthogonality problem.

The square problem using the first three columns of $Y$ has Gram
matrix $\widehat G$ given by the first three columns of the last
display. Its determinant is also $-2\gamma^2$. The additional
right image condition is $v_6=0$: replacing the
normalization by this image condition leaves the same reduced
matrix \eqref{eq:circleresidual}. With the original
lists $\cA,\cB$, its kernel is
$\widehat K(P,Q)=\cA(P)\widehat G^{-1}\cB(Q)$, where
\[
 \widehat G^{-1}=\begin{bmatrix}
 0&\gamma^{-1}&-(2\gamma)^{-1}\\
 \gamma^{-1}&-\gamma^{-2}&(2\gamma^2)^{-1}\\
 -\gamma^{-1}&\gamma^{-2}&(\gamma^2-1)/(2\gamma^2)
 \end{bmatrix}.
\]
This also gives an explicit check of the inverse pairing constructed
in Theorem~\ref{thm:kernel}. The calculation uses the curve relation,
independent matrix actions and a strictly bilateral mass throughout.

\subsection{Determinant orders: degrees, components and support}
\label{sec:orders}

The order of a Christoffel determinant is part of the constructive
result. It specifies how many original polynomial and rational data
must be assembled. This subsection gives two explicit counts for the
four-factor transformation: a formula obtained directly from degrees,
and a smaller formula obtained by polynomial row reduction. The order
is the number of rows of the square matrix whose determinant is taken.
It is distinct from the number of components of a polynomial and from
the degree of that polynomial. Both counts refer
to specified constructions. A determinant can be enlarged by an identity
block, or shortened after an invertible block has been eliminated;
there is therefore no representation-independent determinant order.

For example, the rational scalar perturbation in
Section~\ref{sec:scalarresidual} has numerator and denominator of degree
two. The usual formula retains as unknowns the coefficients of a finite
list of consecutive original polynomials. For this scalar example that
system has order four, whereas the determinant
displayed there has order three. The latter uses sums of products of
the original polynomials and rational moments. Its entries already
contain the coefficients eliminated from the larger system.

\subsubsection{Counting polynomial functions on the support}

Let $F_jA$ be the classes represented by polynomials of total degree
at most $j$, and set $F_jA=0$ for $j<0$. Define
\[
 H_A(j)\coloneq\dim_{\CC}F_jA,\qquad
 \rho_A(j)\coloneq H_A(j)-H_A(j-1).
\]
Thus $H_A(j)$ counts all available polynomial functions through degree
$j$, while $\rho_A(j)$ counts the new functions at degree $j$. The
relations of the support must be imposed before counting. For example,
on $y=x^2$, the two expressions $y$ and $x^2$ represent the same
function, but the inherited degree of that function is one.

The following formulas make the dependence on the support explicit.
Binomial terms with negative degree arguments are taken as zero.
\begin{center}\small
\begin{tabular}{@{}ll@{}}
\toprule
Support and degree convention & $H_A(j)$, $j\ge0$\\
\midrule
One independent variable & $j+1$\\[2pt]
$D$ independent variables & $\binom{j+D}{D}$\\[2pt]
Plane curve of degree $d$ & $\binom{j+2}{2}-\binom{j-d+2}{2}$\\[2pt]
Hypersurface of degree $d$ in $\CC^M$
 & $\binom{j+M}{M}-\binom{j-d+M}{M}$\\[2pt]
Circle, Laurent length & $2j+1$\\[2pt]
Two-dimensional torus, Laurent length & $2j^2+2j+1$\\[2pt]
Embedded torus of Section~\ref{part:classical} & $3j(j+1)+1$\\
\bottomrule
\end{tabular}
\end{center}
For a hypersurface, multiplication by its defining polynomial identifies
the relations of degree at most $j$ with the polynomials of degree at
most $j-d$. Subtracting their dimensions proves the displayed formula;
smoothness is not required. For the torus, count the Laurent monomials
$z^\alpha$ with $\sum_i|\alpha_i|\le j$. In $D$ Laurent variables the
same count is
\[
 H_A(j)=\sum_{i=0}^{D}2^i\binom{D}{i}\binom{j}{i}.
\]
Indeed choose the $i$ nonzero coordinates, their signs, and their
positive absolute values whose sum is at most $j$. For fixed coordinates
and signs, subtract one from each positive absolute value. The resulting
$i$ nonnegative integers have sum at most $j-i$, giving
$\binom ji$ possibilities. The case $i=0$ gives the constant monomial.

For any affine support, a Gr\"obner basis with an order compatible
with total degree provides a basis of standard monomials (those not
divisible by a leading monomial of a relation) and hence
$H_A$. The eventual polynomial count can also be obtained by
homogenizing the support equations. Write $A=\CC[x_1,\ldots,x_M]/I$
and let $I^{\mathrm h}$ be the ideal generated by the homogenizations
of all elements of $I$, using one new variable $x_0$. Put
\[
 A^{\mathrm h}\coloneq\CC[x_0,x_1,\ldots,x_M]/I^{\mathrm h}.
\]
Setting $x_0=1$ identifies the homogeneous degree-$j$ part of
$A^{\mathrm h}$ with $F_jA$. Surjectivity follows by homogenizing
each polynomial of degree at most $j$ to degree $j$. For injectivity,
a homogeneous polynomial whose dehomogenization belongs to $I$
is a power of $x_0$ times the homogenization of that element of $I$,
and therefore belongs to $I^{\mathrm h}$. Hence the two spaces
have the same dimension.

The Hilbert-polynomial theorem applied to $A^{\mathrm h}$ gives,
for sufficiently large $j$, a polynomial $P_A(j)$ of degree
$s=\dim A$. Its leading coefficient is $d/s!$, where $d$ is the
degree of the projective closure with its multiplicities; see
\cite[Section~1.9 and Chapter~12]{Eisenbud}.
Write $P_A(j)=dj^s/s!+\sum_{\ell=0}^{s-1}c_\ell j^\ell$.
For a fixed nonnegative integer $a$, subtraction gives
\[
 P_A(j+a)-P_A(j)
 =\frac{d}{s!}\bigl((j+a)^s-j^s\bigr)
   +\sum_{\ell=0}^{s-1}c_\ell\bigl((j+a)^\ell-j^\ell\bigr).
\]
The binomial expansion of the first difference starts with
$sa j^{s-1}$; each remaining difference has degree at most $s-2$.
Thus, in the polynomial range,
\begin{equation}\label{eq:hilbertdifference}
 H_A(j+a)-H_A(j)
 =\frac{da}{(s-1)!}j^{s-1}+\mathrm{O}(j^{s-2})
 \quad(j\longrightarrow\infty),\qquad s\ge2.
\end{equation}
For a curve, $P_A(j)=dj+c_0$, so this difference is exactly $da$
for sufficiently large $j$. The degree and the filtration are separate data: the section
spaces $H^0(C,F(nP))$ used in Section~\ref{sec:families} add one
dimension per step at a smooth point $P$ in their stable range.
They must not be counted as ambient total-degree blocks of size $d$.

\subsubsection{An explicit order for all four matrix factors}

Choose component bounds
$u_i,v_j\in\{0,1,\ldots\}$ and take
\[
 X=\bigoplus_{i=1}^q F_{u_i-1}A\,e_i^{\T},\qquad
 Y=\bigoplus_{j=1}^p F_{v_j-1}A\,e_j.
\]
The numbers of coefficients are
\begin{equation}\label{eq:orderdimensions}
 r\coloneq\sum_{i=1}^q H_A(u_i-1),\qquad
 N\coloneq\sum_{j=1}^p H_A(v_j-1),\qquad k=N-r\ge0.
\end{equation}
A zero bound gives a zero component. Choose once and for all a
degree-ordered basis of $A$, and order its component copies by degree.
Let $E_N,H_N$ be the spans of the first $N$ row and column vectors
in these orders. Their nested lists as $N$ increases will be called
the row and column \emph{flags}; thus a flag here simply adds one
specified independent vector at each step.
The only original normality assumption needed below is that
$\beta_\mu$ is nondegenerate on $E_N\times H_N$.
For $c$ components, the largest degree occurring among the first $N$
vectors is
\[
 \ell_c(N)\coloneq\min\{j:cH_A(j)\ge N\}.
\]
We take $N>0$; the zero-dimensional system is interpreted with the
usual empty determinant convention.

The four matrices enter the count through actual polynomial degrees:
\begin{equation}\label{eq:orderclearing}
 \begin{aligned}
 \delta_L&\coloneq\det L_1,& B_L&\coloneq\operatorname{adj}(L_1)L_0,
 &e_L&\coloneq\deg\delta_L,&\gamma_L&\coloneq\deg B_L,\\
 \delta_R&\coloneq\det R_1,& B_R&\coloneq R_0\operatorname{adj}(R_1),
 &e_R&\coloneq\deg\delta_R,&\gamma_R&\coloneq\deg B_R.
 \end{aligned}
\end{equation}
The degree of a matrix is the largest degree of its entries in $A$.
Thus cancellations caused by the support relations are included.
Write $U\coloneq\max_{i:u_i>0}(u_i-1)$ and
$V\coloneq\max_{j:v_j>0}(v_j-1)$. If $X=0$, omit the term containing
$U$ in the next maximum. Set
\begin{equation}\label{eq:orderK}
 K_L\coloneq\max\{\ell_q(N)+e_L,U+\gamma_L\},\qquad
 K_R\coloneq\max\{\ell_p(N)+e_R,V+\gamma_R\}.
\end{equation}
Finally, let $h$ be the rank of the matrix of the specified mass
$\nu$ on $X\times Y$. This rank is calculated from the mass data;
it is zero when the mass is absent.

\begin{theorem}[Degree count for bilateral determinants in mixed-type orthogonality]
\label{thm:orders}
Assume the four determinants are regular and the original pair
$E_N\times H_N$ is normal. There is a reduced determinant in
Theorem~\ref{thm:taucriterion} whose exact order is
\begin{equation}\label{eq:orderfull}
 t_{\mathrm{deg}}
 =qH_A(K_L)+pH_A(K_R)-2r-k+h.
\end{equation}
Its nonvanishing is equivalent to existence and uniqueness of the
prescribed normalized block. Each scalar polynomial entry has the
bordered determinant formula of order $t_{\mathrm{deg}}+1$.

The same data give a smaller construction by polynomial rank
computations. In bases of $(F_{K_L}A)^{1\times q}$ and
$(F_{K_R}A)^{p\times1}$, form the coefficient matrices
\begin{equation}\label{eq:orderrankmatrices}
 C_L\coloneq[\,\delta_LE_N\ \ XB_L\,],\qquad
 C_R\coloneq[\,\delta_RH_N\ \ B_RY\,].
\end{equation}
Here a displayed space means the list obtained from any of its
bases; these are ordinary complex coefficient matrices, with one
column for each vector on the list. The exact order using the
spans of these lists is
\begin{equation}\label{eq:orderrank}
 t_{\mathrm{span}}
 =\rank C_L+\rank C_R-2r-k+h
 \le t_{\mathrm{deg}}.
\end{equation}
This is the smallest ambient-space count for the fixed original
pair $E_N,H_N$ and a mass factorization of rank $h$.
\end{theorem}

\begin{proof}
First choose the explicit rational spaces
\[
 W_L=\delta_L^{-1}(F_{K_L}A)^{1\times q},\qquad
 W_R=\delta_R^{-1}(F_{K_R}A)^{p\times1}.
\]
Regularity of the denominators makes multiplication by their inverses
injective, so these spaces have dimensions $qH_A(K_L)$ and
$pH_A(K_R)$. They contain the original spaces: for example, an
element of $E_N$ has degree at most $\ell_q(N)$, and multiplying it
by $\delta_L$ raises degree by at most $e_L$. They also contain the
transformed spaces, because
$XS_L=\delta_L^{-1}XB_L$ and
$S_RY=\delta_R^{-1}B_RY$. The two maxima in
\eqref{eq:orderK} are exactly the bounds needed for these inclusions.

Thus $s_L=qH_A(K_L)-N$ and $s_R=pH_A(K_R)-N$.
Substitution in the order $N-r+s_L+s_R+h$ of
\eqref{eq:small} gives
\[
 N-r+(qH_A(K_L)-N)+(pH_A(K_R)-N)+h
 =qH_A(K_L)+pH_A(K_R)-N-r+h.
\]
Using $N=r+k$ proves \eqref{eq:orderfull}. A rank factorization of the mass matrix has
exactly $h$ terms, as required in that construction. The normality
criterion and the bordered formulas follow from
Theorem~\ref{thm:taucriterion}.

For the smaller spaces, clearing the regular denominator identifies
$E_N+XS_L$ with the span of $\delta_LE_N$ and $XB_L$.
Its dimension is therefore $\rank C_L$. The same argument gives
$\rank C_R$ on the right. Thus the numbers of added vectors are
$\rank C_L-N$ and $\rank C_R-N$, and their substitution gives
\[
 N-r+(\rank C_L-N)+(\rank C_R-N)+h
 =\rank C_L+\rank C_R-2r-k+h,
\]
which is \eqref{eq:orderrank}. Any ambient space containing the required two
lists contains their span, which proves the stated minimality. This
does not assert minimality among different original pairs or among
determinants obtained by further elimination.
\end{proof}

All quantities in \eqref{eq:orderfull} are now specified by polynomial
degrees, the component bounds, the Hilbert function of the support
and the given mass matrix. Formula~\eqref{eq:orderrank} replaces
the potentially redundant degree bounds by two ordinary rank
calculations. Numerator factors equal to the identity give pure
division; denominator factors equal to the identity give bilateral
Christoffel multiplication. A factor equal to the identity on each
side gives the corresponding unilateral specialization. No separate
order theorem is needed for these cases.

If a different original normal pair has dimension $N_0$, use its
actual largest row and column degrees in place of
$\ell_q(N),\ell_p(N)$ in \eqref{eq:orderK}. The same proof gives
\[
 t_{\mathrm{deg}}=qH_A(K_L)+pH_A(K_R)-N_0-r+h,
 \qquad
 t_{\mathrm{span}}=\rank C_L+\rank C_R-N_0-r+h.
\]
Thus an original flag that fails to be normal at $r+k$ does not
invalidate a construction based on another normal pair. The
worked mixed-type example on a curve uses precisely this option, with $N_0=r=3$.

For row and column sections of modules, replace $H_A$ on the left
and right by the respective dimensions of the section spaces, and
compute the coefficient matrices in bases of those spaces. The
proof is unchanged when multiplication by a degree-$a$ function
maps a section of level $j$ to level at most $j+a$ and the
denominator actions are injective. Use actual filtration shifts,
such as allowed pole orders, in place of ambient polynomial degrees.
In particular, for line bundles
filtered by allowed poles along a divisor $D_\infty$, each stable
increment has size $\deg D_\infty$, not necessarily the degree of
an ambient equation of the curve. In the family of
Section~\ref{sec:families}, the stated bases have dimension $N$
on both sides; their increments are one.

For a zero-dimensional coordinate algebra, $H_A(j)$ eventually
equals $\dim_{\CC}A$, and the same finite formulas apply; there is
no infinite sequence of independent higher-degree blocks. If a
factor is nonregular, the inverse-based formulas are not invoked:
Proposition~\ref{prop:compatibility} first decides existence of the
transformed functional, and its direct normalized Gram system,
when specified, has order $N$.

\subsubsection{Mass ranks and their contribution to the order}

The number of scalar parameters describing all possible masses is
different from the rank of one chosen mass. Only the latter is added
to \eqref{eq:orderfull}. It can be computed without a choice of global
mass coordinates. Define
\[
 \begin{aligned}
 d_X&\coloneq\dim X-\dim\bigl(X\cap A^{1\times q}L_1\bigr),\\
 d_Y&\coloneq\dim Y-\dim\bigl(Y\cap R_1A^{p\times1}\bigr).
 \end{aligned}
\]
These are polynomial module membership computations, over the same
coordinate algebra used to calculate $H_A$.

\begin{proposition}[Finite mass rank bound]\label{prop:ordermass}
For every admissible balanced or bivariate mass,
\begin{equation}\label{eq:ordermassbound}
 h\le\min\{r,N,d_X+d_Y\}.
\end{equation}
If $a_1\coloneq\deg L_1$ and $b_1\coloneq\deg R_1$, put
$u_{\min}\coloneq\min_i u_i$ and
$v_{\min}\coloneq\min_j v_j$. Then the entirely numerical bounds
\begin{equation}\label{eq:ordermasshilbert}
 d_X\le r-qH_A(u_{\min}-1-a_1),\qquad
 d_Y\le N-pH_A(v_{\min}-1-b_1)
\end{equation}
hold. On a curve, the global quotient dimensions in
Theorem~\ref{thm:masses} satisfy
\begin{equation}\label{eq:ordercurvemass}
 h\le d_L+d_R\le d(qa_1+pb_1).
\end{equation}
\end{proposition}

\begin{proof}
Choose complements of $X\cap A^{1\times q}L_1$ in $X$ and of
$Y\cap R_1A^{p\times1}$ in $Y$. The mass vanishes when both arguments
belong to the intersections. Its matrix, in these decompositions,
has, with the complements listed first, the form
\[
 M_\nu=\begin{bmatrix}A&B\\C&0\end{bmatrix}
      =\begin{bmatrix}A&B\\0&0\end{bmatrix}
       +\begin{bmatrix}0&0\\C&0\end{bmatrix}.
\]
The first summand has only $d_X$ possibly nonzero rows and the
second only $d_Y$ possibly nonzero columns. Subadditivity of rank
therefore gives $\rank M_\nu\le d_X+d_Y$. Its dimensions also
give $\rank M_\nu\le\min\{r,N\}$, proving
\eqref{eq:ordermassbound}.

Multiplication by $L_1$ maps
$(F_{u_{\min}-1-a_1}A)^{1\times q}$ injectively into
$X\cap A^{1\times q}L_1$. Counting its dimension proves the first
bound in \eqref{eq:ordermasshilbert}; the right bound follows in
the same way. On a curve apply this argument to the full spaces
$(F_jA)^q$. Let $Q_L\coloneq A^{1\times q}/A^{1\times q}L_1$.
Its finite dimension means that representatives of a basis all lie
in some $F_jA$. Thus, for all sufficiently large $j$, the map from
$(F_jA)^q$ to $Q_L$ is onto. Its kernel contains the injective image
of $(F_{j-a_1}A)^q$ under multiplication by $L_1$. Rank-nullity gives
\[
 d_L\le q\bigl(H_A(j)-H_A(j-a_1)\bigr)=qda_1.
\]
The same proof gives $d_R\le pdb_1$. This argument uses regularity
and the Hilbert polynomial, and includes singular, reducible and
nonreduced curves whenever the four-factor regularity assumptions hold.
It also allows isolated zero-dimensional components: their eventual
constant contribution to $H_A$ cancels in $H_A(j)-H_A(j-a_1)$.
Equivalently, a regular scalar on an isolated finite-dimensional
local component is a unit, so its determinant contributes no
remainder there. No assumption that all components have dimension
one is needed for this argument.
\end{proof}

For a purely additive perturbation one can use $E=X$, $H=Y$ in
the square kernel problem, with no extra rational vectors. Its
determinant has order $h$, the usual finite-rank update order.
For a normalized problem with $\dim Y-\dim X=k$, taking nested
original spaces $E_N,H_N$ containing $X,Y$ gives order $k+h$
in \eqref{eq:small}. When the unperturbed normalized problem is
already solved, eliminating that normalization block gives the
usual order-$h$ update instead. Masses arising from division can
already be included in spectral data; they must not be added a
second time to a count based on consecutive original degrees.

\subsubsection{One variable and curves: exact stable counts}

For a degree-ordered flag for mixed-type orthogonality, take $X=E_r$ and $Y=H_N$, with
$N=r+k$. On a curve of degree $d$ write
$H_A(j)=dj+\chi$ in the stable range. The integer $\chi$ is the
constant term of this Hilbert polynomial. Define the two endpoint
counts and two fixed degree shifts by
\[
 \begin{aligned}
 \Delta_q&\coloneq qH_A(\ell_q(N))-N,&
 \Delta_p&\coloneq pH_A(\ell_p(N))-N,\\
 \sigma_L&\coloneq\max\{e_L,\gamma_L+\ell_q(r)-\ell_q(N)\},&
 \sigma_R&\coloneq\max\{e_R,\gamma_R\}.
 \end{aligned}
\]
The endpoint counts measure the unused places in the last degree
block, and satisfy $0\le\Delta_q<qd$ and $0\le\Delta_p<pd$.
For flags the formulas for $K_L,K_R$ use the largest degrees
$\ell_q(r),\ell_p(N)$ directly; they do not require a partial block
to be a full component degree space.

\begin{corollary}[Bilateral orders on a curve]\label{cor:ordercurve}
Once all the degrees entering the count are in the linear Hilbert
range, the degree construction has exact order
\begin{equation}\label{eq:ordercurveexact}
 t_{\mathrm{deg}}
 =k+\Delta_q+\Delta_p+d(q\sigma_L+p\sigma_R)+h.
\end{equation}
In particular,
\begin{equation}\label{eq:ordercurvebound}
 t_{\mathrm{span}}\le t_{\mathrm{deg}}
 \le k+(q+p)d-2
 +d\bigl(q\max\{e_L,\gamma_L\}
          +p\max\{e_R,\gamma_R\}+qa_1+pb_1\bigr).
\end{equation}
For fixed factors and fixed $k$, this is independent of the level.
For a fixed original mass the rank $h$ stabilizes along exhaustive
nested flags. The degree-construction order is eventually periodic
in $r$, with a period dividing $\operatorname{lcm}(qd,pd)$.
\end{corollary}

\begin{proof}
Write $K_L=\ell_q(N)+\sigma_L$ and
$K_R=\ell_p(N)+\sigma_R$. The linear Hilbert formula gives
$qH_A(K_L)=N+\Delta_q+qd\sigma_L$ and its right counterpart.
Substitution in \eqref{eq:orderfull} proves
\eqref{eq:ordercurveexact}. Since $r\le N$, the monotonicity of
$\ell_q$ bounds $\sigma_L$ by $\max\{e_L,\gamma_L\}$.
The endpoint and mass bounds give \eqref{eq:ordercurvebound}.
For stabilization of $h$, let the global mass rank be $h_\infty$.
Choose $h_\infty$ rows and columns whose mass matrix is invertible;
such a minor exists by the definition of rank. Exhaustiveness puts
all of them in the two flags at some finite level. Every subsequent
restriction has rank at least $h_\infty$, and no restriction can
exceed the global rank. Hence $h=h_\infty$ from that level onward.
Finally,
\[
 \ell_c(N)=\left\lceil\frac{N-c\chi}{cd}\right\rceil
\]
in the stable range. Put $L\coloneq\operatorname{lcm}(qd,pd)$.
Increasing $r$ and $N$ by $L$ increases $\ell_q$ by $L/(qd)$
and $\ell_p$ by $L/(pd)$. Thus $\Delta_q,\Delta_p$ and
$\ell_q(r)-\ell_q(N)$ stay unchanged. Every term in
\eqref{eq:ordercurveexact} then repeats, proving eventual periodicity.
\end{proof}

This gives an explicit answer for the full bilateral mixed-type orthogonality problem
on a curve, including all four matrices and any admissible mass.
The order may oscillate with the position inside a degree block;
it need not be the same integer at consecutive scalar levels.
The rank computation \eqref{eq:orderrank} detects additional
reductions, including cancellations between matrix factors.

On the line $d=1$, $\chi=1$ and
$\ell_c(N)=\lceil N/c\rceil-1$. Thus an exact formula at every
level $r\ge1$ is
\begin{equation}\label{eq:orderlineexact}
 \begin{aligned}
 t_{\mathrm{deg}}={}&
 q\max\!\left\{\left\lceil\frac{r+k}{q}\right\rceil+e_L,
                  \left\lceil\frac rq\right\rceil+\gamma_L\right\}
 +p\left\lceil\frac{r+k}{p}\right\rceil
 +p\max\{e_R,\gamma_R\}-2r-k+h.
 \end{aligned}
\end{equation}
At $r=0$, the left target space is zero. Omit the second term
inside the first maximum, as prescribed in \eqref{eq:orderK};
the numerator factor then creates no left rational vectors.
For example, with $q=p=k=1$, $L_0=z^a$ and all other factors
equal to one, the construction has order one for every $a\ge0$.
The normalization alone determines the constant candidate.
For $r\ge1$, $q=p=1$ and $k=1$, if the four scalar factors have degrees
$a_0,a_1,b_0,b_1$, this becomes
\[
 t_{\mathrm{deg}}=
 \max\{1+a_1,a_0\}+\max\{b_1,b_0\}+h.
\]
In particular, four linear factors with no mass give order three,
in agreement with the worked example. This count belongs to the
reduced formula, and is different from the order
$a_0+b_0+a_1+b_1$ when the scalar numerator and denominator have
those degrees and the standard spectral formula applies.

For scalar orthogonal blocks of total degree $n\ge1$ on any support,
take $X=F_{n-1}A$ and $Y=F_nA$. Thus $r=H_A(n-1)$ and
$k=\rho_A(n)$. If the four scalar factors have filtered degrees
$a_0,a_1,b_0,b_1$, define
\[
 u\coloneq\max\{a_1,a_0-1\},\qquad
 v\coloneq\max\{b_1,b_0\}.
\]
Theorem~\ref{thm:orders} gives the exact bilateral order
\begin{equation}\label{eq:orderscalarblock}
 t_{\mathrm{deg}}(n)=
 H_A(n+u)+H_A(n+v)-H_A(n)-H_A(n-1)+h_n.
\end{equation}
This applies whenever the indicated original Gram is normal and
the four factors are regular; $h_n$ is the rank of the actual
mass on $F_{n-1}A\times F_nA$. On a curve in its stable range,
\[
 t_{\mathrm{deg}}(n)
 =d\bigl(\max\{1+a_1,a_0\}+\max\{b_1,b_0\}\bigr)+h_n.
\]
For two independent variables the exact formula is
\[
 t_{\mathrm{deg}}(n)
 =(u+v+1)n+\frac{u(u+3)+(v+1)(v+2)}2+h_n.
\]
In particular, with four linear factors and zero mass the orders
are $3$ on the line, $3n+5$ in two independent variables, and
$3d$ on a plane curve once $n\ge\max\{1,d-1\}$. These are exact orders
of the reduced construction; the consecutive-degree construction has
its separate count below.
Setting $a_1=b_0=b_1=0$ and $a_0=a\ge1$ in
\eqref{eq:orderscalarblock} recovers the usual Christoffel order
$H_A(n+a-1)-H_A(n-1)$.

The four entry degrees alone give a convenient bound even before
forming the adjugates:
\[
 \begin{aligned}
 e_L&\le qa_1,&\gamma_L&\le(q-1)a_1+a_0,\\
 e_R&\le pb_1,&\gamma_R&\le(p-1)b_1+b_0.
 \end{aligned}
\]
Inserting these inequalities in \eqref{eq:ordercurvebound} gives
a bound involving only $d,q,p,k$ and the four degrees. This bound
refers to the degree-based construction; polynomial rank reduction
and further elimination of invertible blocks can lower its order.
For example, for a $2\times3$ pairing on a conic,
four degree-one matrices give $t_{\mathrm{deg}}\le k+44$;
the exact count uses \eqref{eq:ordercurveexact}, and the reduced
count uses \eqref{eq:orderrank}.

\subsubsection{Several variables, surfaces and general component degree bounds}

For two independent variables, \eqref{eq:orderfull} reads
\begin{equation}\label{eq:orderplane}
 t_{\mathrm{deg}}=
 \frac q2(K_L+1)(K_L+2)
 +\frac p2(K_R+1)(K_R+2)-2r-k+h.
\end{equation}
Here
\[
 \ell_c(N)=
 \left\lceil\frac{\sqrt{1+8N/c}-3}{2}\right\rceil,
\]
so the formula is numerical once the factors and component bounds
are specified. On a surface in $\CC^3$ defined by one equation of
degree $d$, replace each binomial $H_A(j)$ by
\[
 H_A(j)=\binom{j+3}{3}-\binom{j-d+3}{3}.
\]
The resulting formula is again exact at every degree, including
degrees below the eventual Hilbert polynomial range. Other surfaces
and higher-dimensional supports use their own $H_A$ in
\eqref{eq:orderfull}; the polynomial rank formula applies unchanged.

\begin{proposition}[Growth with the dimension of the support]
\label{prop:ordergrowth}
Fix the four factors. Assume the original normal flags exist at
the levels under consideration. Let $X,Y$ be degree-ordered flags,
or component degree spaces for which the spread of the bounds among
all components on each side stays bounded. Let $n$ be the larger of the two largest
component degrees. If the support has dimension $s\ge1$ and
$k=\mathrm{O}(n^{s-1})$ as $n\longrightarrow\infty$, then
\begin{equation}\label{eq:ordergrowth}
 t_{\mathrm{span}}\le t_{\mathrm{deg}}
 =\mathrm{O}(n^{s-1})\qquad(n\longrightarrow\infty).
\end{equation}
This includes every admissible balanced or bivariate mass for the
fixed denominator factors. Thus the orders are bounded on curves,
grow at most linearly on surfaces, and at most quadratically on
three-dimensional supports. These upper bounds allow smaller orders
for special factors or mass data.
\end{proposition}

\begin{proof}
Write $U,V$ for the largest degrees occurring on the two sides,
so $n=\max\{U,V\}$, and retain $N=r+k$. First compare each
target space with its full degree space. For component spaces,
choose fixed nonnegative integers $D_L,D_R$ bounding the two
spreads. Then
\[
 \begin{aligned}
 qH_A(U-D_L)&\le r\le qH_A(U),\\
 pH_A(V-D_R)&\le N\le pH_A(V).
 \end{aligned}
\]
Indeed every left component contains $F_{U-D_L}A$ and is contained
in $F_UA$; the right argument is the same. Define the missing
dimensions by
\[
 \varepsilon_X\coloneq qH_A(U)-r,\qquad
 \varepsilon_Y\coloneq pH_A(V)-N.
\]
The preceding inequalities and the fixed-shift Hilbert difference give
\[
 \begin{aligned}
 0\le\varepsilon_X
   &\le q\bigl(H_A(U)-H_A(U-D_L)\bigr),\\
 0\le\varepsilon_Y
   &\le p\bigl(H_A(V)-H_A(V-D_R)\bigr).
 \end{aligned}
\]
Both bounds are $\mathrm{O}(n^{s-1})$ as $n\longrightarrow\infty$.
If $X,Y$ are flags, all degree blocks below the last one are full,
so the corresponding bounds are instead
\[
 0\le\varepsilon_X<q\rho_A(U),\qquad
 0\le\varepsilon_Y<p\rho_A(V).
\]
They satisfy the same upper bound. This also proves the assertion when a
flag is used on one side and a component space on the other.

These comparisons control the degrees of the original spaces as
well. The upper bound is $N\le pH_A(n)$. If $n=U$, the left
lower bound gives $N\ge r$ of order $n^s$; if $n=V$, the right
lower bound does so directly. The flag lower bounds use the full
block immediately below the last one and give the same conclusion.
Thus there are positive constants $C_1,C_2$ such that
\[
 C_1n^s\le N\le C_2n^s
 \qquad\text{for all sufficiently large }n.
\]
For $c\in\{q,p\}$ put $b_c\coloneq\ell_c(N)$.
Its defining minimality gives
\[
 cH_A(b_c-1)<N\le cH_A(b_c).
\]
Since $H_A(j)$ has leading term $dj^s/s!$, these inequalities imply
\[
 b_c\sim\left(\frac{s!N}{cd}\right)^{1/s}
 \qquad(N\longrightarrow\infty).
\]
In particular, each $b_c$ is bounded above and below by positive
constant multiples of $n$. The endpoint remainder
$\eta_c\coloneq cH_A(b_c)-N$ therefore satisfies
\[
 0\le\eta_c<c\rho_A(b_c)
       =\mathrm{O}(n^{s-1})\qquad(n\longrightarrow\infty).
\]

To compare the two entries in each maximum in \eqref{eq:orderK}, put
$M_L\coloneq\max\{b_q,U\}$ and
$M_R\coloneq\max\{b_p,V\}$. Monotonicity of $H_A$ gives the exact
identities
\[
 \begin{aligned}
 qH_A(M_L)
   &=\max\{N+\eta_q,r+\varepsilon_X\}
     =N+\max\{\eta_q,\varepsilon_X-k\},\\
 pH_A(M_R)
   &=\max\{N+\eta_p,N+\varepsilon_Y\}
     =N+\max\{\eta_p,\varepsilon_Y\}.
 \end{aligned}
\]
Both excesses over $N$ are nonnegative and bounded by a constant
times $n^{s-1}$ as $n\longrightarrow\infty$. The factor degrees are fixed,
and \eqref{eq:orderK} gives
\[
 \begin{aligned}
 M_L&\le K_L\le M_L+\max\{e_L,\gamma_L\},\\
 M_R&\le K_R\le M_R+\max\{e_R,\gamma_R\}.
 \end{aligned}
\]
The Hilbert difference for these fixed increments now proves
\[
 qH_A(K_L)=N+\mathrm{O}(n^{s-1}),\qquad
 pH_A(K_R)=N+\mathrm{O}(n^{s-1})
 \quad(n\longrightarrow\infty).
\]
The mass estimate is also a fixed-shift comparison. For component
spaces, $u_{\min}-1\ge U-D_L$, so
\eqref{eq:ordermasshilbert} implies
\[
 \begin{aligned}
 d_X&\le r-qH_A(U-D_L-a_1)
       \le q\bigl(H_A(U)-H_A(U-D_L-a_1)\bigr),\\
 d_Y&\le N-pH_A(V-D_R-b_1)
       \le p\bigl(H_A(V)-H_A(V-D_R-b_1)\bigr).
 \end{aligned}
\]
For a left flag, $(F_{U-1}A)^{1\times q}\subset X$, and multiplication
by $L_1$ maps $(F_{U-1-a_1}A)^{1\times q}$ injectively into
$X\cap A^{1\times q}L_1$. Hence
\[
 d_X\le q\bigl(H_A(U)-H_A(U-1-a_1)\bigr),
\]
with the analogous estimate on the right. In every case
\eqref{eq:ordermassbound} yields
$h\le d_X+d_Y=\mathrm{O}(n^{s-1})$ as $n\longrightarrow\infty$.
This argument uses the rank of the mass on the chosen finite spaces,
so it applies to all the balanced and bivariate masses in the statement.

Substitution in \eqref{eq:orderfull} gives
\[
 t_{\mathrm{deg}}
 =2N-2r-k+\mathrm{O}(n^{s-1})+h
 =k+\mathrm{O}(n^{s-1})+h
 =\mathrm{O}(n^{s-1})\quad(n\longrightarrow\infty),
\]
where the last step uses the assumed bound on $k$.
The inequality $t_{\mathrm{span}}\le t_{\mathrm{deg}}$ follows from
Theorem~\ref{thm:orders}.

The comparison is made separately on the two sides and therefore
allows $q\ne p$. For example, on $A=\CC[z]$ take $q=2$, $p=3$,
$u_1=u_2=3m$ and $v_1=v_2=v_3=2m$, with $m\ge1$. Then
$r=N=6m$, $k=0$, $U=b_q=3m-1$ and $V=b_p=2m-1$.
Each side has zero spread, although the two largest degrees differ
by $m$. The estimates above compare each of these degrees with its
own flag endpoint. With unrestricted component spreads the exact
order formulas still apply, but the bounds on
$\varepsilon_X,\varepsilon_Y,d_X,d_Y$ used here need not hold.
\end{proof}

\subsubsection{Consecutive degrees and conditions at infinity}

To compare with familiar scalar Christoffel formulas, suppose the
scalar relation is $Q\widehat\mu=P\mu$, where $P,Q$ have degrees
$a,b$, and use a spectral construction whose unknown coefficients
occupy the full original blocks of degrees from $n-b$ through
$n+a-1$. Its determinant order is
\begin{equation}\label{eq:orderspectral}
 t_{\mathrm{win}}=H_A(n+a-1)-H_A(n-b-1),\qquad n\ge b.
\end{equation}
This statement counts that specified interval of degrees. Its use as an invertible
spectral system requires the divisibility, degree and normality
hypotheses of the chosen spectral construction
\cite{AMRational,AMChristoffel}. A Hilbert-function calculation alone
does not prove those hypotheses on a new support. The general
bilateral construction above imposes the full image conditions
and uses the image spaces generated by the actual factors.

For one variable \eqref{eq:orderspectral} is $a+b$. For two
independent variables it is
\[
 \frac{(a+b)(2n+a-b+1)}2.
\]
On a plane curve it counts
\[
 \sum_{j=n-b}^{n+a-1}\min\{j+1,d\},
\]
and equals $d(a+b)$ when $n\ge b+d-1$. In $D$ independent
variables its exact value is
\[
 \binom{n+a+D-1}{D}-\binom{n-b+D-1}{D}.
\]
These familiar orders and \eqref{eq:orderfull} describe different
coefficient systems and should be reported separately.

In the one-variable spectral transformation for mixed-type orthogonality
$\widehat\mu R=L\mu$, the standard compatible step-line
construction has order $\deg\det L+\deg\det R$
\cite{MixedUvarov}. For invertible leading matrices of degrees
$a,b$, this is $qa+pb$. On a curve, simply multiplying this
integer by $d$ counts complete degree layers but does not by itself
construct that spectral determinant. The proved bilateral orders
on the curve are \eqref{eq:orderfull}, \eqref{eq:orderrank} and
\eqref{eq:ordercurveexact}.

Two examples explain why the extra degree information is needed.
On the parabola $A=\CC[x,y]/(y-x^2)$, the scalar polynomial $x$
has ambient degree one, but $\dim A/(x)=1$, whereas a stable
total-degree block has size two. The difference comes from a
condition on the highest allowed degree, in addition to divisibility
at the finite zero. Even on the line, the matrix
\[
 L(z)=\begin{bmatrix}1&z^a\\0&1\end{bmatrix},\qquad a>0,
\]
has determinant one and a polynomial inverse. It creates no
finite divisibility conditions, but changes component degree
bounds. Thus determinant degrees alone cannot cover arbitrary
matrix factors with fixed component degree bounds. The rank matrices in
\eqref{eq:orderrankmatrices} retain precisely these degree effects.
Finally, cancelling a common factor from a rational fraction does
not automatically cancel it from the functional equation: from
$z\widehat\mu=z\mu$ one may still have
$\widehat\mu=\mu+c\delta_0$. The remaining mass must be retained
when choosing the reduced representation and its order.

\subsubsection{How to report the size of a concrete formula}

For a given transformation, first state the degree filtration and
the component bounds, then compute $H_A$, the four quantities in
\eqref{eq:orderclearing}, and the mass matrix rank. Equations
\eqref{eq:orderK} and \eqref{eq:orderfull} give the exact order of
the degree construction, without any further dimension choices.
Reducing the two polynomial coefficient matrices gives
\eqref{eq:orderrank}. If a formula using consecutive degrees or a finite-rank update
is used instead, its own order is stated along with its hypotheses.
The normalized Gram system always has order $N=r+k$ and is an
additional option at low degrees. For the reduced construction,
the polynomial entry determinants have order $t+1$, and the same
order-$t$ inverse determines the kernel; no new large Gram
determinant is needed for that kernel calculation.

For reference, the concrete constructions already displayed in
this article have the following orders. Here the symbols $r,N$
in the last two rows retain their local meaning in those examples.
\begin{center}\small
\begin{tabular}{@{}p{.62\linewidth}p{.29\linewidth}@{}}
\toprule
Construction & Determinant order\\
\midrule
Scalar four-linear-factor reduced system, Section~\ref{sec:scalarresidual}
 & $3$; bordered entries $4$\\[2pt]
Mixed-type circle example, Section~\ref{sec:mixedcurveexample}
 & $7$ for $\gamma\ne0$; bordered entries $8$\\[2pt]
Same circle lists with the zero mass removed
 & $5$; determinant zero\\[2pt]
Additive correction of rank $h$, normal reference pairing
 & $h$; scalar bordered entries $h+1$\\[2pt]
Singular-curve family, Section~\ref{sec:families}
 & $4r$ for the update; full Gram order $2N$\\[2pt]
Distributed torus mass, Proposition~\ref{prop:torusmass}
 & $2N+1$ for the update; full Gram order $3N(N+1)+1$\\
\bottomrule
\end{tabular}
\end{center}
In the circle example the direct normalized Gram determinant has
order four, although the displayed reduced uses order seven.
At such a small level the direct system is smaller. In the
singular-curve family the update order $4r$ is independent of the
section level $N$; on the two-dimensional torus the update order
$2N+1$ grows linearly. These distinctions prevent a determinant
order from being confused with the rank of a weighted kernel
difference or with the number of free mass parameters.
Balanced and bivariate pairings use the same dimension counts;
their different mass freedoms enter through the specified rank $h$.
For the Hermitian and symmetric specializations, conjugation and
transposition do not change the sizes of the coefficient spaces.

\subsection{Finite updates of kernels and normalized vectors}
\label{res:sec:finitarytools}

The same finite pairings also admit a direct Gram description, useful
when a particular rationally transformed pairing has already been
specified. The following formulas retain that pairing as one block
and add the chosen mass to it. They will be applied to the residue
coordinates and spaces of sections in
Sections~\ref{res:sec:residueconstruction} and~\ref{sec:families}.

\subsubsection{Rational solutions with regular denominators}
The particular-solution argument in Theorem~\ref{thm:masses}(i)
uses only the inverses of the denominator factors. We record this
specialization for the balanced pairings used in the geometric
constructions, allowing arbitrary numerator factors. In this
subsubsection $A$ is integral and $\det L_1\det R_1\ne0$.
Extend $\mu$ complex linearly to the rational matrix tests that occur
below, and denote the extension by $\widetilde\mu$.
Such an extension exists by completing a vector-space basis of the
polynomial tests to one of the rational tests. Its additional values
are part of the data: the polynomial moments alone do not determine
them. Formula \eqref{eq:pairing}, with $\widetilde\mu$,
defines the corresponding rational pairing
$\widetilde\beta_\mu$. Set
\begin{equation}\label{res:geo:rational-factors}
 S_L\coloneq L_1^{-1}L_0,\qquad
 S_R\coloneq R_0R_1^{-1},\qquad
 \beta_{\mathrm{part}}(f,g)
    \coloneq\widetilde\beta_\mu(fS_L,S_Rg).
\end{equation}
The inverses exist over $K$ by the determinant assumption.

\begin{proposition}[Particular solution and homogeneous mass]
\label{res:geo:allsolutions}
Every balanced solution of \eqref{eq:four} is uniquely of the form
\begin{equation}\label{res:geo:pairing-allsolutions}
 \widehat\beta=\beta_{\mathrm{part}}+\nu,
 \qquad \nu(fL_1,R_1g)=0\quad\text{for every }f,g,
\end{equation}
where $\nu$ is a balanced bilinear pairing. Uniqueness refers to the
fixed rational extension in \eqref{res:geo:rational-factors}.
The pairing $\nu$ is called a homogeneous mass.
\end{proposition}

\begin{proof}
The specialization of Theorem~\ref{thm:masses}(i) follows by the
same substitution on the polynomial tests; it is also the calculation
needed for the section modules below. Only $L_1,R_1$ are inverted,
and the order of the factors gives
$(fL_1)S_L=fL_0$ and $S_R(R_1g)=R_0g$.
Thus $\beta_{\mathrm{part}}(fL_1,R_1g)
=\beta_\mu(fL_0,R_0g)$, because the extension agrees with
$\mu$ on polynomial tests. Subtracting this identity from
\eqref{eq:four} shows that
$\nu\coloneq\widehat\beta-\beta_{\mathrm{part}}$ vanishes
on the indicated pairs. Conversely, adding any such $\nu$ preserves
the equation. The difference of balanced pairings is balanced, and
the displayed subtraction also proves uniqueness.
\end{proof}

For rows and columns of sections of modules stable under
multiplication by $A$, the two cancellation identities in the proof
hold in their rational section spaces. Once a pairing on these
rational sections is specified and agrees with the original pairing
on regular tests, substitution and subtraction give the same
particular solution and homogeneous mass. Balancedness is preserved
because scalar multiplication commutes with the matrix actions.
In Subsection~\ref{res:bk:section} the two
modules are a line bundle and its dual with values in the canonical
differentials of the curve. Their products are differential tests.
Subsection~\ref{res:bq:sec:quotients} describes the finite quotient through
which the homogeneous masses factor, including at singular points.

\subsubsection{From additive masses to a finite reproducing kernel}
\label{res:bq:sec:finitekernel}

For scalar polynomial spaces in several variables, additions of
one Dirac mass and of finitely many point masses yield the polynomial
and kernel formulas of \cite{FPPX2010,Granada2010,Granada2016}.
The following finite-dimensional calculation gives the kernel and the
normalized vectors after an additive change of pairing. It applies to
polynomial or rational functions with any prescribed finite bounds.
Only the chosen spaces, their pairings and the rational factors
enter the calculation. We use the column--row convention of
Subsection~\ref{sec:kernelsetup}.

Let $E_L$ and $E_R$ be spaces of row and column functions of sizes
$q$ and $p$, with the same finite complex dimension $d$ and bases
$f_1,\ldots,f_d$ and $g_1,\ldots,g_d$. The original pairing
$\beta_\mu$ is extended to the rational functions occurring below;
this extension is denoted $\widetilde\beta_\mu$. Let $\nu$ be an
additive pairing on the target spaces.
For the rational changes of arguments, set $S_L\coloneq L_1^{-1}L_0$ and
$S_R\coloneq R_0R_1^{-1}$, with the factors of
Subsection~\ref{sec:algebra}, and write
\[
 \beta_{\mathrm{part}}(f,g)
   \coloneq\widetilde\beta_\mu(fS_L,S_Rg),\qquad
 \widehat\beta=\beta_{\mathrm{part}}+\nu.
\]
Choose bases $a_1,\ldots,a_s$ and $b_1,\ldots,b_t$ of finite
rational function spaces containing $E_LS_L$ and $S_RE_R$,
respectively. Define their coordinate matrices by
\begin{equation}\label{res:bq:eq:coordinates}
 f_iS_L=\sum_{\alpha=1}^s(T_L)_{i,\alpha}a_\alpha,
 \qquad
 S_Rg_j=\sum_{\beta=1}^tb_\beta (T_R)_{\beta,j},
\end{equation}
where $T_L$ is $d\times s$ and $T_R$ is $t\times d$. Finally put
\[
 \begin{aligned}
 G_{\mathrm{rat}}\coloneq
   {}&[\widetilde\beta_\mu(a_\alpha,b_\beta)]_{\alpha\in\{1,\ldots,s\},\,\beta\in\{1,\ldots,t\}},\\
 H_\nu\coloneq{}&[\nu(f_i,g_j)]_{i,j\in\{1,\ldots,d\}}.
 \end{aligned}
\]
Thus $G_{\mathrm{rat}}$ contains the rational pairings of the ambient
bases used to express $f_iS_L$ and $S_Rg_j$, and $H_\nu$ contains the mass evaluated on the target
basis pairs. When the mass factors through a finite quotient, its
entries can be computed from the remainders of these pairs and the
prescribed values of the mass on a quotient basis.

\begin{proposition}[Normality and the finite kernel]
\label{res:bq:prop:finitekernel}
With the rational extension, spaces and coordinates above,
the transformed Gram is
\begin{equation}\label{res:bq:eq:gram}
 \widehat G=T_LG_{\mathrm{rat}} T_R+H_\nu.
\end{equation}
The pairing is normal on $E_L,E_R$ if and only if
$\det\widehat G\ne0$. Under this condition its reproducing kernel is
\begin{equation}\label{res:bq:eq:finitekernel}
 \widehat K(P,Q)=\boldsymbol g(P)\widehat G^{-1}\boldsymbol f(Q),
 \qquad
 \boldsymbol g(P)\coloneq[g_1(P)\ \cdots\ g_d(P)],\quad
 \boldsymbol f(Q)\coloneq[f_i(Q)]_{i\in\{1,\ldots,d\}}.
\end{equation}
No invertibility of the possibly rectangular matrix $G_{\mathrm{rat}}$ is
assumed.

Suppose in addition that
$G_{\mathrm{part}}\coloneq T_LG_{\mathrm{rat}} T_R$ is invertible. Let
$\rho\coloneq\rank H_\nu$ and choose a rank factorization
$H_\nu=P_\nu Q_\nu$, where $P_\nu$ is $d\times\rho$ and
$Q_\nu$ is $\rho\times d$, both of rank $\rho$. Then
$\widehat G$ is invertible if and only if
\begin{equation}\label{res:bq:eq:woodburymatrix}
 W_\nu\coloneq I_\rho+Q_\nu G_{\mathrm{part}}^{-1}P_\nu
\end{equation}
is invertible. When these conditions hold,
\begin{align}
 \widehat K(P,Q)-K_{\mathrm{part}}(P,Q)
 ={}&-\bigl[\boldsymbol g(P)G_{\mathrm{part}}^{-1}P_\nu\bigr]
       W_\nu^{-1}
       \bigl[Q_\nu G_{\mathrm{part}}^{-1}\boldsymbol f(Q)\bigr],
       \label{res:bq:eq:woodburykernel}\\
 K_{\mathrm{part}}(P,Q)
 \coloneq{}&\boldsymbol g(P)G_{\mathrm{part}}^{-1}
                                      \boldsymbol f(Q).\notag
\end{align}
The separation rank of this correction is exactly $\rho$.
For $\rho=0$, the correction is zero.
\end{proposition}

\begin{proof}
Expanding both arguments by \eqref{res:bq:eq:coordinates} gives
$[\beta_{\mathrm{part}}(f_i,g_j)]=T_LG_{\mathrm{rat}} T_R$. Adding the mass
gives \eqref{res:bq:eq:gram}. Normality is precisely invertibility of
this finite Gram. If $g=\sum_jg_jc_j$, its pairings against the
rows $f_i$ form the column $\widehat Gc$; multiplication by
$\boldsymbol g(P)\widehat G^{-1}$ recovers $g(P)$. The corresponding
row calculation proves the other reproduction identity and yields
\eqref{res:bq:eq:finitekernel}.

For the additional assertion, form the square matrix of order
$d+\rho$
\[
 \mathcal D\coloneq
 \begin{bmatrix}G_{\mathrm{part}}&P_\nu\\-Q_\nu&I_\rho\end{bmatrix}.
\]
Adding $Q_\nu G_{\mathrm{part}}^{-1}$ times its upper block row
to its lower block row gives diagonal blocks $G_{\mathrm{part}}$
and $I_\rho+Q_\nu G_{\mathrm{part}}^{-1}P_\nu=W_\nu$.
Alternatively, subtracting $P_\nu$ times its lower block row from
its upper block row gives diagonal blocks
$G_{\mathrm{part}}+P_\nu Q_\nu=\widehat G$ and $I_\rho$.
Both operations preserve determinants, proving
\[
 \det\widehat G=\det\mathcal D
              =\det G_{\mathrm{part}}\det W_\nu.
\]
For the inverse formula put
$B\coloneq G_{\mathrm{part}}^{-1}P_\nu$ and
$C\coloneq Q_\nu G_{\mathrm{part}}^{-1}$.
Then
$\widehat G B=P_\nu(I_\rho+Q_\nu B)=P_\nu W_\nu$,
so
\[
 \widehat G\bigl(G_{\mathrm{part}}^{-1}-BW_\nu^{-1}C\bigr)
 =I_d+P_\nu C-P_\nu W_\nu W_\nu^{-1}C=I_d.
\]
This proves
\[
 \widehat G^{-1}=G_{\mathrm{part}}^{-1}
  -G_{\mathrm{part}}^{-1}P_\nu W_\nu^{-1}
                      Q_\nu G_{\mathrm{part}}^{-1}.
\]
Multiplying by the two basis matrices proves
\eqref{res:bq:eq:woodburykernel}. Moreover,
\[
 \widehat G^{-1}-G_{\mathrm{part}}^{-1}
       =-\widehat G^{-1}H_\nu G_{\mathrm{part}}^{-1}.
\]
Its rank is $\rho$, since both outside factors are invertible.
The two chosen bases are complex linearly independent, so this
matrix rank is the separation rank of their tensor, as in
Subsection~\ref{sec:kernelsetup}.
\end{proof}

A mass may make $\widehat G$ invertible even when
$G_{\mathrm{part}}$ is singular; then \eqref{res:bq:eq:finitekernel}
still applies. When both are normal, \eqref{res:bq:eq:woodburykernel}
compares two pairings on the same target spaces. The mass values,
the rational moments and the choice of finite spaces are the data
needed for this computation.

\subsubsection{The transformed polynomials and their normalization}
\label{res:pf:subsection}

The same Gram matrices determine the individual orthogonal vectors.
We retain the normalization information by prescribing their classes
modulo the lower spaces. For polynomial spaces these classes prescribe
the leading coefficients; the argument also applies whenever fixed
complements specify the desired new directions.

Write $E\coloneq E_L$ and $H\coloneq E_R$ for the finite spaces just
used, and retain their bases $f_1,\ldots,f_d$ and $g_1,\ldots,g_d$.
These are the spaces against which the next vectors must be orthogonal.
Let
\[
 \widehat\beta=\beta_{\mathrm{part}}+m,\qquad
 G\coloneq[\beta_{\mathrm{part}}(f_i,g_j)]_{i,j\in\{1,\ldots,d\}},
 \qquad \widehat G\coloneq[\widehat\beta(f_i,g_j)]_{i,j\in\{1,\ldots,d\}}.
\]
Here $m$ denotes the full additive mass pairing, including any scalar
coefficient multiplying it. Thus $m=\nu$ in
Proposition~\ref{res:bq:prop:finitekernel}, whereas a family written as
$\widehat\beta=\beta_{\mathrm{part}}+c\nu$ has $m=c\nu$.

\begin{proposition}[Transformation of normalized orthogonal vectors]
\label{res:pf:transformation}
Assume that $G$ and $\widehat G$ are invertible. Let $\Pi$ be a row
vector and $\Sigma$ a column vector satisfying
\[
 \beta_{\mathrm{part}}(\Pi,g_j)=0\quad(j\in\{1,\ldots,d\}),
 \qquad
 \beta_{\mathrm{part}}(f_i,\Sigma)=0\quad(i\in\{1,\ldots,d\}).
\]
Put
\[
 \begin{aligned}
 \boldsymbol f&\coloneq[f_i]_{i\in\{1,\ldots,d\}},
 &\qquad \boldsymbol g&\coloneq[g_1\ \cdots\ g_d],\\
 a_\Pi&\coloneq[m(\Pi,g_j)]_{j\in\{1,\ldots,d\}},
 &\qquad b_\Sigma&\coloneq[m(f_i,\Sigma)]_{i\in\{1,\ldots,d\}},
 \end{aligned}
\]
where $a_\Pi$ is a $1\times d$ row and $b_\Sigma$ a $d\times1$ column.
Then the unique transformed vectors with the same classes modulo
$E$ and $H$, respectively, are
\begin{equation}
 \widehat \Pi\coloneq \Pi-a_\Pi\widehat G^{-1}\boldsymbol f,
 \qquad
 \widehat \Sigma\coloneq \Sigma-\boldsymbol g\widehat G^{-1}b_\Sigma.
 \label{res:pf:eq:transformation}
\end{equation}
They satisfy
\[
 \widehat\beta(\widehat \Pi,H)=0,
 \qquad \widehat\beta(E,\widehat \Sigma)=0.
\]
In particular, any prescribed leading coefficients encoded by those
quotient classes are preserved.
\end{proposition}

\begin{proof}
All entries of the Gram matrices and of $a_\Pi,b_\Sigma$ are complex
scalars. In particular, the row correction in
\eqref{res:pf:eq:transformation} means
\[
 a_\Pi\widehat G^{-1}\boldsymbol f
   =\sum_{i,j=1}^d m(\Pi,g_j)(\widehat G^{-1})_{j,i}f_i.
\]
Pairing this identity with $g_k$ gives
\[
 \widehat\beta(\widehat \Pi,g_k)
 =m(\Pi,g_k)
  -\sum_{i,j=1}^d
     m(\Pi,g_j)(\widehat G^{-1})_{j,i}\widehat G_{i,k}=0.
\]
We used the original orthogonality of $\Pi$ in the first term and
$\widehat G^{-1}\widehat G=I_d$ in the second. For the column correction, its pairing with $f_k$ is
\[
 \widehat\beta(f_k,\widehat \Sigma)
 =m(f_k,\Sigma)-\sum_{i,j=1}^d
       \widehat G_{k,i}(\widehat G^{-1})_{i,j}m(f_j,\Sigma)=0,
\]
using $\widehat G\widehat G^{-1}=I_d$. Thus the order of the
inverse factors is verified separately for rows and columns.

Both corrections lie in the corresponding lower spaces. They therefore
preserve the prescribed classes. If two transformed rows had the same
class, their difference would be $\sum_i c_i f_i$; its pairings with
all $g_j$ would say $[c_1\ \cdots\ c_d]\widehat G=0$. Invertibility
forces every $c_i$ to vanish. For columns the same argument uses
$\widehat G[c_1\ \cdots\ c_d]^{\T}=0$. This proves uniqueness.
\end{proof}

To obtain normalized blocks, take larger spaces $E^+\supset E$ and
$H^+\supset H$ with equal dimension increments. Choose fixed
complements with row matrix $\Pi^0$ and column matrix $\Sigma^0$. If the
common increment is $a$, these matrices have sizes $a\times q$ and
$p\times a$. Define
\[
 \begin{aligned}
 \Pi&\coloneq \Pi^0-
     \beta_{\mathrm{part}}(\Pi^0,\boldsymbol g)G^{-1}\boldsymbol f,\\
 \Sigma&\coloneq \Sigma^0-
     \boldsymbol gG^{-1}\beta_{\mathrm{part}}(\boldsymbol f,\Sigma^0).
 \end{aligned}
\]
The pairing tables before $G^{-1}$ and after it have sizes
$a\times d$ and $d\times a$, respectively. Pairing the first
equation with $\boldsymbol g$ subtracts its own table multiplied by
$G^{-1}G=I_d$; pairing the second with $\boldsymbol f$ uses
$GG^{-1}=I_d$. Thus every row of $\Pi$ and column of $\Sigma$ satisfies
the orthogonality hypotheses of Proposition~\ref{res:pf:transformation},
while the subtracted terms lie in the lower spaces and preserve the
chosen classes.
Applying its formulas to each of them gives blocks
$\widehat{\boldsymbol \Pi}$ and $\widehat{\boldsymbol \Sigma}$ with the
chosen leading classes. If the transformed pairing is also normal
on $E^+\times H^+$, then
\[
 \widehat S\coloneq
 \widehat\beta(\widehat{\boldsymbol \Pi},\widehat{\boldsymbol \Sigma})
\]
is invertible. Indeed, in the completed bases the Gram matrix is
\[
 \begin{bmatrix}
 \widehat\beta(\boldsymbol f,\boldsymbol g)&
        \widehat\beta(\boldsymbol f,\widehat{\boldsymbol \Sigma})\\
 \widehat\beta(\widehat{\boldsymbol \Pi},\boldsymbol g)&
        \widehat\beta(\widehat{\boldsymbol \Pi},\widehat{\boldsymbol \Sigma})
 \end{bmatrix}
 =\begin{bmatrix}\widehat G&0\\0&\widehat S\end{bmatrix}.
\]
The two off-diagonal blocks vanish by the orthogonality just proved.
Its determinant is $\det\widehat G\det\widehat S$; normality on
the larger spaces and $\det\widehat G\ne0$ force
$\det\widehat S\ne0$. Keeping
$\widehat{\boldsymbol \Pi}$ and replacing
$\widehat{\boldsymbol \Sigma}$ by
$\widehat{\boldsymbol \Sigma}\widehat S^{-1}$ produces dual blocks whose
mutual pairing is the identity. This last normalization changes the
leading classes of the column block; it is a separate choice from
prescribing both leading blocks. Iterating along nested normal cuts
gives biorthogonal block families: every new row or column block
is orthogonal to all preceding blocks on the opposite side.

The computation uses the same small matrix as the kernel correction.
If $H_m\coloneq[m(f_i,g_j)]$ has a rank factorization
$H_m=P_mQ_m$ of rank $\rho$, substitution of the inverse formula in
the proof of Proposition~\ref{res:bq:prop:finitekernel} into
\eqref{res:pf:eq:transformation} uses
$I_\rho+Q_mG^{-1}P_m$, exactly as in
\eqref{res:bq:eq:woodburymatrix}. When the mass has finite rank on the
whole space, $a_\Pi$ and $b_\Sigma$ are obtained from its same finite
remainder data, even as the polynomial spaces increase. This does
not mean that only finitely many orthogonal polynomials change.

Finally, these formulas control parameter dependence. Suppose the
spaces, bases and prescribed leading classes are fixed, while the
moments of both pairings on the relevant enlarged spaces are
holomorphic in a parameter near a given value. If the needed Gram
determinants are nonzero at that value, their inverses are
holomorphic in a neighbourhood, by the adjugate formula. The
particular projections, the transformed vectors in
\eqref{res:pf:eq:transformation}, and the dual normalization by
$\widehat S^{-1}$ are consequently holomorphic there. The same
argument gives rational coefficients when the moments are rational
functions of the parameter. The curve family in
Section~\ref{res:bf:section} applies this observation at a collision
of the mass support with a singular point.

\section{Residues and matrix division on curves and surfaces}\label{res:sec:residueconstruction}

The finite formulas of Section~\ref{part:finite} require linear
conditions that recognize polynomial images and values of masses on
the remaining classes. This section constructs those data by residues.
A quotient records the remainder after division. A perfect pairing on
that quotient identifies every linear functional on its classes with
one dual class. Choosing bases then turns the pairing into a finite
matrix, whose entries are the residues given below. Inserting these
entries in the finite construction produces the normalized polynomials
as well as their kernel.

The curve construction first supplies its geometric Cauchy kernel and
then its matrix quotient pairing. The surface construction uses a
second equation to make the quotient finite and proves its residue
formula directly. The finite-dimensional update formulas apply to both.

The underlying finite linear algebra can be expressed in terms of
remainder classes and their dual functionals. A class
$[f]$ in a quotient $V/TV$ is the remainder information that survives
adding $Ta$ to $f$. A pairing $b:Q\times Q^D\to\CC$ is
\emph{perfect} when every linear functional on either space is
pairing with a unique element of the other. If
$G\coloneq[b(e_i,\eta_j)]_{i,j\in\{1,\ldots,m\}}$ is its matrix
in two bases, this means that $G$ is square and invertible. The
\emph{coevaluation tensor} is the corresponding inverse written as
\[
 \mathcal C\coloneq
 \sum_{i,j=1}^{m}(G^{-\T})_{i,j}e_i\otimes\eta_j.
\]
Indeed, contracting the second factor with $e_k$ gives
$\sum_{i,j}(G^{-\T})_{i,j}G_{k,j}e_i=e_k$.
This calculation fixes our index order: the first tensor factor
labels the rows, and the second labels the columns. The residue
formulas below construct $b$ and $\mathcal C$ without first
inverting $G$. Their value for the perturbation problem is that
they turn divisibility and the remaining mass values into finite
computations.
The quotient is chosen for the equation under consideration:
numerator-image conditions use the corresponding division map,
whereas homogeneous masses use the denominator test map. The
residue method applies separately to these choices.

The global rank of a mass pairing $\nu$ is the dimension
of the image of $g\mapsto[f\mapsto\nu(f,g)]$. Its kernel is the
right radical, consisting of the columns pairing to zero with every
row. Thus this global rank is the dimension of the column module
modulo that radical. Every finite matrix of mass pairings has rank at
most the global rank. The singular families below identify the finite
cuts on which the full rank is attained.

\subsection{A global Cauchy kernel on an integral curve}
\label{res:bk:section}\label{res:part:curves}

On a curve, division by a regular denominator leaves a finite-dimensional
space of remainder classes. Residue duality represents their linear
functionals and hence the values of a mass on the chosen polynomial or
section bases. Adding these values to the particular pairing gives the
Gram matrix used in Subsection~\ref{res:sec:finitarytools} to compute the
kernel and normalized polynomials. The construction below establishes
this duality also when the support of the mass meets a singular point.

The geometric kernel $S_F$ has a prescribed
diagonal pole. Multiplication by $T(P)-T(Q)$ produces a regular kernel
$\mathcal C_T$ whose class represents the inverse of a finite residue
pairing. The finite reproducing kernel $K_\beta$ then uses the full
orthogonality pairing and its Gram inverse. The first two objects
encode division and its remainder duality; the third reproduces the
chosen polynomial or section space. Subsection~\ref{res:bq:curvefiniteapplication}
assembles these data explicitly. We then give a noncommutative family
whose normality, mass rank and normalized polynomial blocks persist
through confluence at a singular point, including non-Gorenstein
examples.

The global construction of $S_F$ ensures that multiplying it by a
function vanishing on the diagonal produces a regular tensor even at
a singular pair. This property will justify the matrix divided
difference in Theorem~\ref{res:bq:thm:quotient}; its residue normalization
will give the reconstruction identity. A finite projection then supplies
an expression in trace-dual frames and identifies the restriction on
rational factors arising from a change of geometric data.

For concrete calculations, polynomial reduction provides representatives
of the quotient classes; it may be carried out by module Gr\"obner bases
or by a division adapted to the equations. The residue pairing provides
the dual coordinates of those classes, and $S_F$ gives their inverse
tensor through a matrix divided difference. These descriptions can be
used together. In Section~\ref{sec:families}, coefficient extraction
computes the moments and Grams directly, while the explicit geometric
kernel identifies the tensor reproducing the mass quotient.

\subsubsection{The geometric kernel and the finite reproducing kernel}
\label{res:bk:twokernels}

The finite reproducing kernel
of Subsection~\ref{sec:kernelsetup} is determined by a pairing on two
finite-dimensional spaces. Its coefficients contain the inverse Gram
matrix. Changing the weight or the masses therefore changes this kernel.
Where the basis sections are regular, the finite kernel is regular
on the diagonal as well.

The geometric Szeg\H{o} kernel has different defining data: a curve,
a line bundle and a prescribed diagonal residue. It provides the
analogue of the Cauchy kernel used to form transforms and residues.
Once these geometric data are fixed, changing the orthogonality
functional does not change this kernel. Bertola uses both objects:
his Szeg\H{o} kernel is normalized in \cite[(2.22)]{Bertola}, whereas
the correspondence in \cite[Section~3.1]{Bertola} concerns finite
Christoffel--Darboux kernels. The classical theta-function description
on smooth complex curves appears in Fay~\cite{Fay}. A bundle
formulation, including bundles of higher rank, and its residue
interpretation are given in
\cite[Section~3.4 and Proposition~5.1.1]{BenZviBiswas}.
For smooth projective curves over a finite field with a fixed point
at infinity, Hu and Yu~\cite[Theorems~6.2 and~6.3]{HuYuResidue}
identify the finite residual Casimir tensor with the intrinsic
principal part of the normalized Szeg\H{o} kernel and derive
residue reproduction, after extension to their coefficient field.
Thus the scalar smooth residue construction is an established
antecedent of the result below. The present construction extends
this quotient formula to integral singular curves, retains the
dualizing module at non-Gorenstein points, and specifies the ordered
matrix operations that enter the bilateral transformation of mixed-type orthogonality.

\subsubsection{A kernel defined along the entire diagonal}
\label{res:bk:globalconstruction}

Let $C$ be an integral projective curve over $\CC$: a compact
algebraic curve with one irreducible component and no nilpotent
functions in its structure sheaf $\OO_C$. We allow singular points.
For an open set $U$, $\OO_C(U)$ denotes the functions regular on $U$;
the local ring $\OO_{C,s}$ consists of their germs near $s$.
Its dualizing sheaf $\omega_C$ replaces the sheaf of regular
differentials in Serre duality. It agrees with regular differentials
on the smooth part of $C$, and it may fail to be a line bundle at a
singular point. The \emph{normalization} is the smooth curve mapping
finitely to $C$ that separates the local branches at its singularities.
Locally the dualizing sheaf can be described there:
a meromorphic differential $\eta$ represents a regular dualizing
differential at $s$ if
$\sum_{a\mapsto s}\operatorname{Res}_a(g\eta)=0$
for every $g\in\OO_{C,s}$. The sum runs over the branches above $s$.
This is Rosenlicht's definition of regular differentials on a
singular curve; see \cite[Chapter~IV]{SerreACF}.
Thus regularity includes the relations among the branches. Its
identification with the dualizing sheaf is part of the local and
global duality framework in
\cite[Introduction and Sections~5, 7--10]{Lipman1984}.
A residue on a branch means the coefficient of $t^{-1}\dd t$ in
a local parameter $t$. A differential that is regular in this
dualizing sense may have poles on individual branches. For example,
at a node the pair $(\dd u/u,-\dd v/v)$ has this property: a regular
function has the same value at the origins of the two branches, so
the two residues of its product with this pair cancel. The point is
called \emph{Gorenstein} when the dualizing module has one local
generator; the argument below also treats points where more than
one generator is needed.

The global residue theorem in this convention follows from the
usual theorem on the normalization. Indeed, let
$\varpi:\widetilde C\to C$ be the normalization and let $\eta$ be a
meromorphic differential, regarded also as a differential on the
smooth projective curve $\widetilde C$. By the definition of the
residue at a point of $C$,
\[
 \sum_{s\in C}\operatorname{Res}_s\eta
 =\sum_{s\in C}\sum_{a\in\varpi^{-1}(s)}
       \operatorname{Res}_a\eta
 =\sum_{a\in\widetilde C}\operatorname{Res}_a\eta=0.
\]
Only finitely many terms can be nonzero. At a point where $\eta$
is regular in Rosenlicht's sense, its contribution vanishes by
taking $g=1$ in the local condition. This is the global residue
identity used in Theorem~\ref{res:bq:thm:quotient}.

A line bundle has one local generator on each sufficiently small
open set; a section is a regular function times that generator, with
the prescribed changes of generator on overlaps. For a line bundle
$F$, write
\[
 F^{\mathrm D}\coloneq
 \mathcal{H}\!om_{\OO_C}(F,\omega_C)
 =F^{-1}\otimes\omega_C.
\]
The product of a section of $F$ and a section of $F^{\mathrm D}$
is thus a dualizing differential. We assume that $F$ is
\emph{acyclic}, meaning
\begin{equation}\label{res:bk:acyclic}
 H^0(C,F)=H^1(C,F)=0.
\end{equation}
The first vanishing excludes nonzero global regular sections of $F$;
the second will remove the obstruction to prescribing the residue.
Serre duality identifies $H^0(C,F^{\mathrm D})$ with $H^1(C,F)^*$
and $H^1(C,F^{\mathrm D})$ with $H^0(C,F)^*$. Thus both groups
for $F^{\mathrm D}$ vanish as well.

Let $p_1,p_2:C\times C\to C$ be the two projections and let
$i:C\to C\times C$ be the diagonal map, $i(P)=(P,P)$. Set
\[
 J\coloneq\II_\Delta,\qquad
 \mathscr M_F\coloneq F\boxtimes F^{\mathrm D}
       =p_1^*F\otimes p_2^*F^{\mathrm D}.
\]
The diagonal ideal $J$ consists of functions vanishing when the two
arguments coincide. On an affine product it is generated by the
differences $a(P)-a(Q)$ of regular functions on the curve.

On a smooth curve, a Cauchy kernel is a section of $\mathscr M_F$
with a simple pole on $\Delta$. At a singular point the diagonal
ideal need not have one generator, so we use the sheaf
$\mathcal{H}\!om(J,\mathscr M_F)$. Its sections are locally linear
maps from the diagonal ideal to $\mathscr M_F$. Equivalently, one
may view such a section as a meromorphic kernel $S$ for which
$jS$ is regular for every $j\in J$. Away from the diagonal,
$J=\OO_{C\times C}$, so this condition makes $S$ regular there.
On the smooth part it is precisely the usual simple-pole condition.
The next exact sequence describes this pole condition in elementary
terms: its kernel consists of regular kernels, and its last term
records the one freely prescribed residue along the diagonal.

\begin{proposition}[The diagonal residue sequence]
\label{res:bk:residuesequence}
There is an exact sequence
\begin{equation}\label{res:bk:residueexact}
 0\longrightarrow\mathscr M_F
 \longrightarrow\mathcal{H}\!om(J,\mathscr M_F)
 \xrightarrow{\ \operatorname{res}_\Delta\ }\OO_\Delta
 \longrightarrow0.
\end{equation}
Under \eqref{res:bk:acyclic}, there is a unique global section $S_F$ with
$\operatorname{res}_\Delta S_F=1$.
The sign of the residue map is chosen so that, at a smooth point,
in a local coordinate $z$ and a frame $e$ of $F$,
\begin{equation}\label{res:bk:localnormalization}
 S_F(P,Q)=
 \frac{e(P)e(Q)^{-1}\dd z(Q)}{z(P)-z(Q)}
       +\text{a section regular near the diagonal}.
\end{equation}
Thus $\operatorname{Res}_{Q=P}S_F(P,Q)f(Q)=-f(P)$ for a local
regular section $f$ of $F$.
\end{proposition}

\begin{proof}
An integral noetherian curve is Cohen--Macaulay: its local rings have
dimension at most one, and a nonzero element in the maximal ideal
of a one-dimensional local domain is a nonzerodivisor. This permits
the use of its dualizing sheaf also at a non-Gorenstein point.
We first compute the local quotient that records the diagonal
residue. We then check that its identification with functions on
the diagonal is independent of the coordinates. The final
cohomology calculation constructs the global kernel.

Choose an affine neighbourhood $U=\operatorname{Spec}A$ on which
$F$ is trivial; $A$ is its algebra of regular functions.
A Noether normalization chooses a nonconstant function $x$ such
that $A$ is a finite module over $\CC[x]$. Geometrically this is
a finite map from $U$ to the affine line. It is torsion free over this polynomial
ring, since $A$ is a domain, and hence it is a finite free module.
Denote the image of $x$ in $A$ by $a$, and put
\[
 R\coloneq A\otimes_\CC A,\qquad
 B\coloneq A\otimes_\CC\CC[x]=A[x].
\]
The second factor $\CC[x]$ is included in the second copy of $A$.
The diagonal is the $R$-module $A$ obtained by multiplying the two
factors. As a $B$-module, it is $B/(x-a)$. Consequently
\begin{equation}\label{res:bk:auxresolution}
 0\longrightarrow B\xrightarrow{\ x-a\ }B
   \longrightarrow A\longrightarrow0
\end{equation}
is exact. Injectivity follows by comparing the highest coefficient
of a polynomial multiplied by the monic polynomial $x-a$.

Here is the finite-duality identification used in the computation,
including its normalization. Let $K_A$ be the fraction field of $A$.
A meromorphic differential can be written $u\,\dd x$ with $u\in K_A$.
The description of regular dualizing differentials by the trace
condition is
\[
 \omega_A=
 \{u\,\dd x:\operatorname{Tr}_{K_A/\CC(x)}(bu)\in\CC[x]
                         \text{ for every }b\in A\}.
\]
It identifies the differential $u\,\dd x$ with the functional
$b\mapsto\operatorname{Tr}_{K_A/\CC(x)}(bu)\,\dd x$ in
$\operatorname{Hom}_{\CC[x]}(A,\CC[x])\dd x$;
see \cite[Definition~2.1, Example~2.1.2 and Lemma~2.2]{Lipman1984}.
The trace pairing on the finite separable field extension is
nondegenerate. Therefore every functional on $A$ extends to the
fraction fields and has a unique expression of this form. Its
values lie in $\CC[x]\dd x$ precisely when the displayed condition
holds. This also fixes the identification on the smooth part:
for an unramified local branch the differential trace is the
ordinary trace of the branch coefficients.

After using $\dd x$ as a basis of the differentials on the base,
finite freeness gives
\[
 A\otimes_\CC\omega_A
 \simeq \operatorname{Hom}_{B}(R,B).
\]
For clarity, the scalar-extension map sends $c\otimes\ell$ to
the $B$-linear map
$c'\otimes b\mapsto cc'\ell(b)$; it is an isomorphism because
a basis of the second copy of $A$ over $\CC[x]$ and its dual
basis remain dual after extending coefficients to the first copy.

The module $R$ is finite free over $B$. For an $R$-module $P$,
the adjunction needed here is the explicit isomorphism
\[
 \operatorname{Hom}_R(P,\operatorname{Hom}_B(R,B))
       \longrightarrow\operatorname{Hom}_B(P,B),\qquad
       \Phi\longmapsto[p\mapsto\Phi(p)(1)].
\]
Its inverse sends $\lambda$ to the map
$p\mapsto[r\mapsto\lambda(rp)]$. The $R$-action on
$\operatorname{Hom}_B(R,B)$ is
$(r\ell)(s)=\ell(rs)$, so both formulas are module maps and
their two compositions are the identity. A projective resolution
over $R$ is also projective over $B$: each term is a direct summand
of a free $R$-module, and $R$ is free over $B$. Applying these
inverse maps to the terms and differentials of the resolution
therefore identifies the corresponding Ext groups. Applied to the
diagonal module and \eqref{res:bk:auxresolution}, it yields
\begin{equation}\label{res:bk:localext}
 \operatorname{Ext}_R^1(A,A\otimes_\CC\omega_A)\simeq A,
 \qquad
 \operatorname{Ext}_R^j(A,A\otimes_\CC\omega_A)=0
       \quad(j\ne1).
\end{equation}
More explicitly, applying $\operatorname{Hom}_B(-,B)$ to the
two free terms of \eqref{res:bk:auxresolution} gives the complex
$B\xrightarrow{x-a}B$. Its degree-zero kernel is zero and its
degree-one quotient is $B/(x-a)=A$; there are no further terms.
Thus $\operatorname{Ext}^1$ records exactly one copy of the
functions on the diagonal. In the residue sequence it measures
the failure of a map on the diagonal ideal to extend to a map on
the whole ring. The element $x-a$ belongs to the auxiliary ring
$B$; this calculation does not require the diagonal ideal in
$R$ to have one generator.

We next verify compatibility on overlapping neighbourhoods.
Use the preceding trace normalization and fix the sign of the
Ext identification so that, wherever $x$ is a local coordinate,
the class of $\dd x(Q)/(x(P)-x(Q))$ has value $1$.
This assertion can be checked directly on the resolution.
At a smooth point where the finite map is unramified, the
completed local rings of the branch and the line are both
$\CC[[x-x_0]]$. In the completion along the corresponding
pair of points the trace on that branch is the identity, so the
adjunction above reduces to the two-term calculation for the
principal diagonal equation $x(Q)-x(P)$. Its connecting map
takes the class of $b(P,Q)\dd x(Q)/(x(Q)-x(P))$ to $b(P,P)$.
Taking the negative of this map gives precisely the stated
normalization with denominator $x(P)-x(Q)$.

This normalized value is unchanged if another local coordinate
$z$ is used. Indeed,
\[
 x(P)-x(Q)=u(P,Q)\bigl(z(P)-z(Q)\bigr),\qquad
 u(P,P)=\frac{\dd x}{\dd z}(P)\ne0.
\]
The numerator $\dd x(Q)$ contributes the same factor
$\dd x/\dd z$ on the diagonal, so it cancels the factor
$u(P,P)$ in the denominator. Changing a frame $e$ of $F$ to
$ve$ likewise multiplies the kernel by $v(P)/v(Q)$, whose
restriction to the diagonal is $1$.

For two choices of finite coordinate, these calculations apply
on a common dense open subset of their overlap: in characteristic
zero each finite map is unramified away from finitely many points,
and the smooth points form a dense open subset. The local Ext
calculation makes the degree-one Ext sheaf locally free of rank
one as a sheaf on the diagonal. The two normalized maps from
that sheaf to $\OO_\Delta$ agree on this dense open subset.
After a local trivialization their difference is multiplication
by a regular function vanishing on a dense open subset of the
integral curve. That function is zero. Hence the maps agree on
the whole overlap, including its singular points, and they glue.

The glued local calculations can also be recorded as
\[
 \mathbf R\mathcal{H}\!om_{C\times C}
       (\OO_\Delta,p_2^*\omega_C)
       \simeq i_*\OO_C[-1].
\]
Here derived Hom records all the Ext sheaves at once, and the
shift $[-1]$ places its sole nonzero sheaf in degree one.
Tensoring by the line bundle $p_1^*F\otimes p_2^*F^{-1}$
preserves this formula: its restriction to the diagonal is
canonically trivial. We obtain, in particular,
\[
 \mathcal{H}\!om(\OO_\Delta,\mathscr M_F)=0,\qquad
 \mathcal{E}\!xt^1(\OO_\Delta,\mathscr M_F)\simeq\OO_\Delta.
\]
Thus the residue identification is intrinsic. The coordinate
calculation and its compatibility establish here the diagonal
case of relative duality directly.

Apply $\mathcal{H}\!om(-,\mathscr M_F)$ to
$0\to J\to\OO_{C\times C}\to\OO_\Delta\to0$.
The two displayed identifications, together with
$\mathcal{E}\!xt^1(\OO_{C\times C},\mathscr M_F)=0$,
give \eqref{res:bk:residueexact}. Its first map sends a regular
section $m$ to the multiplication map $j\mapsto jm$.
On a smooth coordinate neighbourhood, $J$ is generated by
$z(P)-z(Q)$. Our choice of sign gives
\[
 \operatorname{res}_\Delta
 \left(\frac{b(P,Q)e(P)e(Q)^{-1}\dd z(Q)}
                  {z(P)-z(Q)}\right)=b(P,P).
\]
This is the normalization in \eqref{res:bk:localnormalization}.

Finally, K\"unneth and \eqref{res:bk:acyclic} imply
$H^0(C\times C,\mathscr M_F)=H^1(C\times C,\mathscr M_F)=0$:
the first group is $H^0(F)\otimes H^0(F^{\mathrm D})$, and the
second is
$[H^1(F)\otimes H^0(F^{\mathrm D})]\oplus
 [H^0(F)\otimes H^1(F^{\mathrm D})]$.
Every summand is zero. The cohomology sequence of
\eqref{res:bk:residueexact} now makes the residue map an isomorphism
\[
 H^0(C\times C,\mathcal{H}\!om(J,\mathscr M_F))
       \longrightarrow H^0(C,\OO_C)=\CC.
\]
The preimage of $1$ is the required section $S_F$.
\end{proof}

The proposition constructs $S_F$ also at $(s,s)$ when $s$ is
singular. Its residue there is the quotient map in
\eqref{res:bk:residueexact}. In particular, the construction retains
$\omega_C$ and the full ideal $J$ at that point. On a smooth curve
it recovers the usual Szeg\H{o} kernel by uniqueness. The proof is
the residue construction of that kernel, with relative duality
providing its formulation on the singular curve.

\subsubsection{An expression through a finite projection}
\label{res:bk:projection}

Let $\pi:C\to\mathbb P^1$ be finite of degree $b$. The direct image
$E\coloneq\pi_*F$ assigns to an open set $V$ of the base the
sections $F(\pi^{-1}(V))$. It is a vector bundle of rank $b$: locally on the
smooth base it is a finite torsion-free module over a discrete
valuation ring and is therefore free. Finite duality gives
\begin{equation}\label{res:bk:finiteduality}
 \pi_*F^{\mathrm D}
 \simeq\mathcal{H}\!om_{\mathbb P^1}(E,\omega_{\mathbb P^1}).
\end{equation}
The evaluation associated with this duality is denoted
$\operatorname{Tr}_\pi:\pi_*\omega_C\to\omega_{\mathbb P^1}$.
It is the trace on dualizing differentials.
In local module terms it sends a functional
$\ell\in\operatorname{Hom}_{\CC[x]}(A,\CC[x])\dd x$ to
$\ell(1)$. This formulation remains meaningful when several
points of a fibre meet or a fibre is ramified.

Since $\pi$ is finite, the cohomology of $E$ is the cohomology of
$F$. Write its splitting on $\mathbb P^1$ as
$E\simeq\bigoplus_{\alpha=1}^b\OO(a_\alpha)$.
The vanishing of $H^0(E)$ forces $a_\alpha\leq-1$, whereas the
vanishing of $H^1(E)$ forces $a_\alpha\geq-1$. Hence
\begin{equation}\label{res:bk:splitting}
 E\simeq\OO(-1)^b,\qquad \pi_*F^{\mathrm D}\simeq\OO(-1)^b.
\end{equation}
Here $\OO(a)$ is the line bundle of degree $a$ on the projective
line. The numerical conditions follow from
$\dim H^0(\OO(a))=\max\{a+1,0\}$ and
$\dim H^1(\OO(a))=\max\{-a-1,0\}$.
Fix such a splitting. On the affine base with coordinate $x$,
it gives a frame $e_1,\ldots,e_b$ of $E$ and a trace-dual frame
$\eta_1,\ldots,\eta_b$ of $\pi_*F^{\mathrm D}$ characterized by
\begin{equation}\label{res:bk:traceframes}
 \operatorname{Tr}_\pi(e_\alpha\eta_\beta)
       =\delta_{\alpha,\beta}\dd x.
\end{equation}
Thus the $e_\alpha$ are sections of $F$ above this affine chart,
while the $\eta_\alpha$ are sections of $F^{\mathrm D}$.

\begin{proposition}[Projection formula for the geometric kernel]
\label{res:bk:projectionformula}
With frames compatible with \eqref{res:bk:splitting}, the kernel of
Proposition~\ref{res:bk:residuesequence} is
\begin{equation}\label{res:bk:szegoprojection}
 S_F(P,Q)=
 \frac{\displaystyle\sum_{\alpha=1}^b
                 e_\alpha(P)\eta_\alpha(Q)}
      {x(P)-x(Q)}.
\end{equation}
The expression defines a section of
$\mathcal{H}\!om(J,\mathscr M_F)$ through ramified fibres and at
singular points. It has no poles away from the diagonal.
\end{proposition}

\begin{proof}
Write $d\coloneq x(P)-x(Q)$ and
$\mathcal C_\pi(P,Q)\coloneq
\sum_\alpha e_\alpha(P)\eta_\alpha(Q)$.
Modulo $d$, this numerator is the coevaluation tensor for
\eqref{res:bk:traceframes}: it represents the identity map of the
finite free module $E$ over $\CC[x]$.
For every regular function $a$ on $\pi^{-1}(\mathbb A^1)$,
multiplication by $a$ is self-adjoint for the evaluation pairing
between $E$ and its dual. The identity tensor consequently satisfies
\[
 (a(P)-a(Q))\mathcal C_\pi(P,Q)
       \in d\,\Gamma(\pi^{-1}(\mathbb A^1)^2,\mathscr M_F).
\]
One can verify this statement directly: applying multiplication
by $a$ to the first tensor factor represents the map $m\mapsto am$;
applying its dual to the second represents the same map.
In coordinates, if $ae_\alpha=\sum_\beta A_{\beta,\alpha}e_\beta$,
trace duality gives
$a\eta_\beta=\sum_\alpha A_{\beta,\alpha}\eta_\alpha$.
Thus both actions on $\sum_\alpha e_\alpha\otimes\eta_\alpha$
give $\sum_{\alpha,\beta}A_{\beta,\alpha}
e_\beta\otimes\eta_\alpha$ after identifying the two base
coordinates. This is the asserted divisibility by $d$.

The differences $a(P)-a(Q)$ generate $J$. It follows that
$j\mathcal C_\pi/d$ is regular for every $j\in J$.
Multiplication by $d$ is injective on $\mathscr M_F$ on this
affine product. For example, its module is finite free over
$\CC[x(P),x(Q)]$ by \eqref{res:bk:splitting}. Thus the quotient
determines an unambiguous map $J\to\mathscr M_F$.
This proves the required membership in the Hom sheaf.
It also explains the cancellation at distinct points with
$x(P)=x(Q)$: off the diagonal, $J$ is the unit ideal, so the
kernel itself is regular. The argument uses the whole finite
module, and does not require a fibre to consist of distinct points.

On a smooth unramified fibre, trace duality is ordinary duality
on the direct sum of the one-dimensional fibres. It gives
$\mathcal C_\pi(P,P)=\dd x(P)$, with the $F$ factors evaluated.
The quotient in \eqref{res:bk:szegoprojection} therefore has diagonal
residue $1$. Its residue is a regular function on the integral
diagonal, so this equality extends to the ramified and singular
points as well.

It remains to check the points above infinity. Write $e_{0,\alpha},\eta_{0,\alpha}$ for the affine
frames already used, and set $w=1/x$.
For compatible frames on the two charts, the transition of
$\OO(-1)$ and the equality $\dd x=-x^2\dd w$ give
\[
 e_{0,\alpha}=x e_{\infty,\alpha},\qquad
 \eta_{0,\alpha}=-x\eta_{\infty,\alpha}.
\]
The numerator and denominator in \eqref{res:bk:szegoprojection} both
acquire the factor $-x(P)x(Q)$. Their quotient is therefore
$\mathcal C_{\pi,\infty}(P,Q)/(w(P)-w(Q))$ on the other chart.
On the two crossed charts it takes the forms
\[
 \frac{\sum_\alpha e_{0,\alpha}(P)\eta_{\infty,\alpha}(Q)}
      {1-x(P)w(Q)},\qquad
 \frac{\sum_\alpha e_{\infty,\alpha}(P)\eta_{0,\alpha}(Q)}
      {1-w(P)x(Q)}.
\]
These expressions are regular near the corresponding points at
infinity. When both points lie above infinity, the same ideal
argument applies with $w$ in place of $x$. We have obtained a
global section with residue $1$, which equals $S_F$ by uniqueness.
\end{proof}

The trace-dual frames used here come from the global splitting
\eqref{res:bk:splitting}. Changes between such frames are constant
matrices, and the numerator in \eqref{res:bk:szegoprojection} is
invariant under their dual changes. An arbitrary nonconstant
change of affine frame also changes the behaviour at infinity
and cannot be substituted into this assertion. Independence of
the chosen projection follows from the global uniqueness of $S_F$.

\subsubsection{Geometric changes and bilateral changes of the weight}
\label{res:bk:bertolafactors}

\begin{remark}[The trace restriction on Bertola's geometric factors]
\label{res:bk:bertolarestriction}
On a smooth curve, let $F=\mathcal X\otimes\sqrt{\omega_C}$ be
acyclic and put $\Xi\coloneq\pi^*(\infty)$. Choose bases $\varphi_a$ of
$H^0(C,F(\Xi))$ and $\varphi_a^\vee$ of
$H^0(C,F^{\mathrm D}(\Xi))$. Above an unramified value $x$, form
\[
 \Phi(x)\coloneq[\varphi_a(x^{(c)})]_{a,c\in\{1,\ldots,b\}},
 \qquad
 \Phi^\vee(x)\coloneq[\varphi_c^\vee(x^{(a)})]_{a,c\in\{1,\ldots,b\}}.
\]
The trace product and its normalization are those of
\cite[Proposition~3.9 and Remark~3.10]{Bertola}:
\begin{equation}\label{res:bk:traceconstant}
 \Phi(x)\Phi^\vee(x)=C_0\dd x,
 \qquad C_0\in\mathrm{GL}_b(\CC).
\end{equation}
Indeed, these bases frame the two trivial bundles $E(1)$ and
$(E^*\otimes\omega_{\mathbb P^1})(1)$; their perfect pairing takes
values in $\omega_{\mathbb P^1}(2\infty)\simeq\OO$. Its matrix is
therefore constant and invertible.

For a second acyclic bundle $\widetilde F$, let $D$ be the diagonal
matrix of values of a meromorphic comparison section of
$\mathcal{H}\!om(F,\widetilde F)$. Bertola's weight factors are
\cite[Proposition~3.8]{Bertola}
\begin{equation}\label{res:bk:bertolalr}
 L=\widetilde\Phi D^{-1}\Phi^{-1},\qquad
 R=(\Phi^\vee)^{-1}D\widetilde\Phi^\vee.
\end{equation}
Substitution into $L(\Phi\Phi^\vee)R=
\widetilde\Phi\widetilde\Phi^\vee$ gives
\begin{equation}\label{res:bk:bertolaconstraint}
 L C_0 R=\widetilde C_0.
\end{equation}
In trace-dual bases both constant matrices are identities, so
$R=L^{-1}$. This is the restriction on factors obtained from that
change of geometric data. The independent factors of
\eqref{eq:four} instead act on the orthogonality pairing.
The geometric restriction also follows from the diagonal residue:
$\operatorname{res}_\Delta(A(P)S_F(P,Q)B(Q))=A(P)B(P)$.
Preserving the identity residue requires $B=A^{-1}$; for a
meromorphic bundle isomorphism the local section modules must also
be transported. The finite reproducing kernel changes with the
pairing even while $S_F$ remains fixed.
\end{remark}

For the matrix division below, fix an affine open $U\subset C$ and
use the notation
\begin{equation}\label{res:bk:affinemodules}
 A\coloneq\OO_C(U),\qquad
 M\coloneq F(U),\qquad N\coloneq F^{\mathrm D}(U).
\end{equation}
Here $M$ is a projective module of rank one and $N$ is its dual
with values in $\omega_A$. Their evaluation takes values in
$\omega_A$, not necessarily in $A$. The kernel $S_F$ remains
fixed. Multiplying it by differences of matrix entries will
produce regular tensors for finite quotient pairings, while the
orthogonality functional determines the separate finite
reproducing kernel.

\subsection{Matrix quotient duality and bilateral masses}\label{res:bq:sec:quotients}

Retain the curve, line bundle and affine modules of
\eqref{res:bk:affinemodules}, and write $\omega\coloneq\omega_C(U)$.
We use the normalization \eqref{res:bk:localnormalization}, so that
$\operatorname{Res}_{Q=P}S_F(P,Q)f(Q)=-f(P)$ for a local section
$f$ of $F$ at a smooth point. Recall that
$N=\operatorname{Hom}_A(M,\omega)$; its elements already carry
their canonical differential factor.

For a scalar denominator $h$, multiplying $S_F(P,Q)$ by
$h(P)-h(Q)$ removes its diagonal singularity. The resulting tensor
encodes the remainder classes modulo $h$. The same construction
works for a matrix denominator, provided that the inverse in the
residue pairing and the transpose in the dual quotient are placed
in the correct order.

\begin{theorem}[Residue kernel of a matrix quotient]\label{res:bq:thm:quotient}
\label{res:bq:matrixquotient}
Let $T\in\Mat_n(A)$ have nonzero determinant, and let $Z\subset U$
be the finite set of zeros of $\det T$. Define
\begin{equation}\label{res:bq:eq:quotients}
 Q_T\coloneq M^n/TM^n,\qquad
 Q_T^D\coloneq N^n/T^{\T}N^n.
\end{equation}
These are finite-dimensional complex vector spaces. The formula
\begin{equation}\label{res:bq:eq:residuepairing}
 b_T([f],[\eta])\coloneq
 \sum_{a\in Z}\operatorname{Res}_a
       \bigl(\eta^{\T}T^{-1}f\bigr)
\end{equation}
defines a perfect bilinear pairing between them. At a singular point,
$\operatorname{Res}_a$ denotes the sum of the residues on the branches
of the normalization above $a$.

The matrix of sections
\begin{equation}\label{res:bq:eq:divideddifference}
 \mathcal C_T(P,Q)\coloneq
       \bigl[T(P)-T(Q)\bigr]S_F(P,Q)
\end{equation}
is regular on $U\times U$. Its row index belongs to $M^n$ and its
column index to $N^n$. Its class in $Q_T\otimes_{\CC}Q_T^D$ is the
coevaluation tensor of \eqref{res:bq:eq:residuepairing}: contracting its
second factor against any $[f]\in Q_T$ returns $[f]$.
\end{theorem}

\begin{proof}
The adjugate identity shows that $\det T$ annihilates both quotients.
Since $A$ is the coordinate ring of an affine integral curve and
$\det T\ne0$, their support is finite and their complex dimensions
are finite. Replacing $f$ by $f+Ta$ changes the differential in
\eqref{res:bq:eq:residuepairing} by $\eta^{\T}a$. Replacing $\eta$ by
$\eta+T^{\T}b$ changes it by $b^{\T}f$. Both are regular dualizing
differentials near $Z$, so their residues there vanish. At a singular
point it is the sum over its branches that vanishes; individual
branch residues need not vanish. Thus the formula descends to the
two stated quotients.

We verify perfection by reducing the displayed residue pairing to
polynomial division in one variable. Choose a Noether normalization
$B=\CC[x]\subset A$. The module $A$ is finite and torsion free over
$B$, hence free. Since $M$ is a finite projective $A$-module, it is
a direct summand of a finite free $A$-module and therefore a finite
free $B$-module as well. Write its $B$-rank as $b$.
The trace description of the canonical module gives
\begin{equation}\label{res:bq:finitebaseduality}
 \omega_A\simeq\operatorname{Hom}_B(A,B)\dd x,
 \qquad
 N=\operatorname{Hom}_A(M,\omega_A)
     \simeq\operatorname{Hom}_B(M,B)\dd x;
\end{equation}
the first identification is the finite-map case of
\cite[Theorem~5.1 and Remark~5.2]{Lipman1984}. In this description,
a differential $\eta$ corresponds to the rule
$a\mapsto\operatorname{Tr}_\pi(a\eta)$, with values in $B\dd x$.
Thus the trace of $\eta$ itself is the value of this rule at $1$.

The second identification in \eqref{res:bq:finitebaseduality} can be
checked directly. Suppress the common factor $\dd x$ for this
calculation, and give $\operatorname{Hom}_B(A,B)$ the $A$-action
$(a\ell)(c)\coloneq\ell(ac)$. An $A$-linear map
$\psi:M\to\operatorname{Hom}_B(A,B)$ determines the $B$-linear
functional $\lambda_\psi(m)\coloneq\psi(m)(1)$.
Conversely, a $B$-linear functional $\lambda:M\to B$ determines
\[
 \psi_\lambda(m)(a)\coloneq\lambda(am).
\]
This rule is $A$-linear in $m$, since
$\psi_\lambda(cm)(a)=\lambda(acm)
 =(c\psi_\lambda(m))(a)$.
The two constructions are inverse: evaluation of
$\psi_\lambda(m)$ at $1$ returns $\lambda(m)$, while
$\lambda_\psi(am)=\psi(am)(1)=\psi(m)(a)$.
This proves the asserted adjunction with its normalization.
In particular, $N$ is free over $B$, even when it is not free over
$A$, because it is the $B$-dual of the finite free module $M$.
Choose a $B$-basis of $M^n$ and its trace-dual basis of $N^n$;
the trace of the evaluation of a basis element against its dual is
$\dd x$, and against any other dual basis element it is zero.
The adjugate identity makes $T$ injective on the torsion-free module
$M^n$. Thus its matrix $D(x)\in\Mat_{nb}(B)$ in these bases has
nonzero determinant, and its dual action $T^{\T}$ is represented
by $D(x)^{\T}$. For the coordinate columns $v,w$ of $f,\eta$,
respectively, finite duality gives the identity of meromorphic
differentials on the base
\[
 \operatorname{Tr}_\pi(\eta^{\T}T^{-1}f)
       =w^{\T}D(x)^{-1}v\,\dd x,
 \qquad \pi:\operatorname{Spec}A\longrightarrow\mathbb A^1_x.
\]

The trace preserves the sum of residues over each fibre:
\begin{equation}\label{res:bq:traceresidues}
 \operatorname{Res}_{x=c}\operatorname{Tr}_\pi\xi
      =\sum_{a\mapsto c}\operatorname{Res}_a\xi.
\end{equation}
Here residues at singular points mean the sums on their normalization
branches. To check the formula, pass to the normalization of the curve.
In the completed local ring of a branch above $c$, choose a parameter
so that
$x-c=z^e$. Since $\dd x=e z^{e-1}\dd z$, summing over
the $e$ conjugate parameters gives
\[
 \operatorname{Tr}(z^k\dd z)=
 \begin{cases}
  (x-c)^{(k+1)/e-1}\dd x,&e\mid k+1,\\
  0,&e\nmid k+1.
 \end{cases}
\]
For a Laurent differential $\sum_k a_kz^k\dd z$, the coefficient
of $\dd x/(x-c)$ is therefore $a_{-1}$. Summing this
calculation over the branches proves \eqref{res:bq:traceresidues}.
Away from $Z$, $\eta^{\T}T^{-1}f$ is a regular dualizing
differential and has zero branch-sum residue. Thus
\eqref{res:bq:eq:residuepairing} is exactly the sum of the finite
residues of $w^{\T}D^{-1}v\,\dd x$.

It remains to check this pairing for the polynomial matrix $D$.
Take a Smith form
$U D V=\operatorname{diag}(d_1,\ldots,d_{nb})$, where $U,V$
are invertible over $B$ and every $d_j$ is nonzero. The changes
$v\mapsto Uv$ and $w\mapsto V^{\T}w$ identify the two quotient
modules and the pairing with the direct sum of
\[
 B/(d_j)\times B/(d_j)\longrightarrow\CC,
 \qquad ([a],[b])\longmapsto
       \sum_{d_j(c)=0}\operatorname{Res}_{x=c}
                       \frac{a(x)b(x)}{d_j(x)}\,\dd x.
\]
Each scalar pairing is perfect. Indeed, near a zero $c$ of order
$k$, write $z\coloneq x-c$ and $d_j=z^k u(z)$ with $u(0)\ne0$.
The residue is the coefficient of $z^{k-1}$ in $ab/u$ modulo
$z^k$. Multiplication by the unit $u^{-1}$ is invertible, and
the coefficient pairing on $\CC[z]/(z^k)$ has an invertible
anti-diagonal Gram in the monomial basis. Explicitly, the classes
$u(z)z^{k-1-j}$ are dual to $z^j$, for
$j\in\{0,\ldots,k-1\}$: the coefficient of $z^{k-1}$ in
$z^i u(z)z^{k-1-j}/u(z)$ is $\delta_{i,j}$.
The Chinese remainder
decomposition combines these local pairings. A constant $d_j$
contributes a zero quotient. This proves perfection of the original
pairing and derives its residue normalization explicitly. The dual
quotient contains $T^{\T}$ because it is the adjoint action for the
evaluation into $\omega_A$, not because $\omega_A$ is identified
with $A$.

Each entry $T_{i,j}(P)-T_{i,j}(Q)$ belongs to the diagonal ideal.
The defining property of
$S_F\in\operatorname{Hom}(\II_\Delta,F\boxtimes F^D)$ therefore
makes \eqref{res:bq:eq:divideddifference} regular, including at pairs of
singular points. Since $U$ is affine, its entries are finite sums
of tensors from $M\otimes_{\CC}N$.

It remains to verify that this regular tensor is the inverse of the
residue pairing. For $f\in M^n$, its contraction is
\begin{equation}\label{res:bq:eq:contraction}
 I(P)\coloneq\sum_{a\in Z}\operatorname{Res}_{Q=a}
       \bigl(\mathcal C_T(P,Q)T(Q)^{-1}f(Q)\bigr).
\end{equation}
The inverse is on the right of $\mathcal C_T$. Keeping that order,
the differential inside the sum is
\begin{equation}\label{res:bq:eq:residueorder}
 T(P)S_F(P,Q)T(Q)^{-1}f(Q)-S_F(P,Q)f(Q).
\end{equation}
For a general smooth $P\in U\setminus Z$, the two terms together
have no pole at $Q=P$: this also follows directly from the regularity
of $\mathcal C_T$. Their remaining poles lie in $Z$ or in the finite
complement $C\setminus U$. The global residue theorem gives
\[
 I(P)=-T(P)A_f(P)+B_f(P),
\]
where
\begin{align*}
 A_f(P)&\coloneq\sum_{a\in C\setminus U}
       \operatorname{Res}_{Q=a}
          \bigl(S_F(P,Q)T(Q)^{-1}f(Q)\bigr),\\
 B_f(P)&\coloneq\sum_{a\in C\setminus U}
       \operatorname{Res}_{Q=a}\bigl(S_F(P,Q)f(Q)\bigr).
\end{align*}
These are regular sections on $U$. Indeed, the residues are taken
at fixed points outside $U$, so they do not meet the diagonal as
$P$ varies in $U$; taking their finite Laurent coefficients preserves
regularity in $P$. In particular $A_f\in M^n$.

Apply the global residue theorem once more to $S_F(P,Q)f(Q)$.
Inside $U$ its only nonzero residue is the diagonal residue $-f(P)$.
Consequently $B_f(P)=f(P)$, and hence
\[
 I(P)=f(P)-T(P)A_f(P).
\]
Both sides are regular on $U$, so the identity obtained for general
$P$ holds throughout $U$. It gives $[I]=[f]$ in $Q_T$.
Because the pairing is perfect, this reproduction identity determines
its coevaluation tensor uniquely and proves the theorem.
\end{proof}

The theorem supplies image-membership tests for any chosen matrix
$T$. If $[e_1],\ldots,[e_m]$ and $[e_1^D],\ldots,[e_m^D]$
are residue-dual quotient bases, put
\[
 \ell_j(f)\coloneq b_T([f],[e_j^D]),
 \qquad j\in\{1,\ldots,m\}.
\]
Then $[f]=\sum_j\ell_j(f)[e_j]$, so $f\in TM^n$ if and only if
all $\ell_j(f)$ vanish. Prescribed mass values $c_j$ on $[e_j]$
give the functional
$\nu[f]=\sum_jc_j\ell_j(f)=b_T([f],\sum_jc_j[e_j^D])$.
Thus the residue calculation tests division by $T$ and describes
every functional vanishing on its image. In a bilateral problem,
apply it separately to the division map being tested and to the
denominator map defining the homogeneous mass.

A small singular quotient makes the differential module concrete.
Take $A=\CC[t^3,t^4,t^5]$, whose only singular point is $t=0$,
and use the scalar denominator $T=t^3$. The allowed exponents of
functions in $A$ are $0,3,4,5,\ldots$. A Laurent differential
$t^k\dd t$ belongs to $\omega_A$ precisely when multiplying by
any such function gives zero residue at $0$. Thus its exponent
cannot be $-1,-4,-5,\ldots$, and
\[
 \omega_A=A\,t^{-3}\dd t+A\,t^{-2}\dd t.
\]
Both generators are needed at $0$. Division by $t^3$ leaves the
following bases of the two three-dimensional quotients:
\[
 A/(t^3):\quad 1,t^4,t^5,
 \qquad
 \omega_A/t^3\omega_A:\quad
 t^{-3}\dd t,t^{-2}\dd t,t^2\dd t.
\]
For example, the first list consists of the allowed exponents that
are not $3$ plus another allowed exponent. The same subtraction
gives the second list from $-3,-2,0,1,\ldots$.
The pairing is $\operatorname{Res}_{t=0}(f\eta/t^3)$, so its
entry is $1$ exactly when the two displayed exponents sum to $2$.
Its Gram matrix is
\[
 \begin{bmatrix}0&0&1\\0&1&0\\1&0&0\end{bmatrix}.
\]
This verifies perfection at a non-Gorenstein point and shows
explicitly why the dual remainder classes are differentials.

The construction has an ordered product rule. If $T_a,T_b$ are two
matrices of the same size with entries in $A$, then
\begin{equation}\label{res:bq:eq:productrule}
 \mathcal C_{T_aT_b}(P,Q)
   =T_a(P)\mathcal C_{T_b}(P,Q)
       +\mathcal C_{T_a}(P,Q)T_b(Q).
\end{equation}
Indeed, after removing the common scalar kernel, the right-hand side is
\[
 T_a(P)\bigl[T_b(P)-T_b(Q)\bigr]
 +\bigl[T_a(P)-T_a(Q)\bigr]T_b(Q)
 =T_a(P)T_b(P)-T_a(Q)T_b(Q).
\]
This proof retains the order of the factors and requires no
commutativity between the two matrices.

\begin{remark}[Row--column convention]\label{res:bq:rem:roworientation}
The arguments of $b_T$ are columns. If instead $f$ is a row with
entries in $M$, $g$ is a column with entries in $N$, and
\[
 \nu(f,g)=c\sum_{a\in Z}\operatorname{Res}_a
                   \bigl(fT^{-1}g\bigr),\qquad c\ne0,
\]
then $\nu(f,g)=c\,b_{T^{\T}}(f^{\T},g)$. Thus its quotient pair is
$M^n/T^{\T}M^n$ and $N^n/TN^n$, and its coevaluation in that order
is $c^{-1}\mathcal C_{T^{\T}}$. This transpose changes the matrix
indices, without exchanging the two points. The column--row kernel
convention of \eqref{eq:kerneldefinition} subsequently exchanges the
two tensor factors.
\end{remark}

\subsubsection{The bilateral test quotient}\label{res:bq:sec:tests}

We now connect the quotient theorem with the four factors of
Subsection~\ref{sec:algebra}. For a balanced $q\times p$ test matrix,
write
\[
 D_i(H)\coloneq L_i^{\T}HR_i^{\T},\qquad
 T_i\coloneq R_i\otimes L_i^{\T},\qquad i\in\{0,1\}.
\]
Vectorization by columns gives
\begin{equation}\label{res:bq:eq:vectorization}
 \operatorname{vec}(D_i(H))=T_i\operatorname{vec}(H).
\end{equation}
For clarity, the entry in row $a$ and column $b$ of $D_i(H)$ is
$\sum_{c,d}(L_i)_{c,a}H_{c,d}(R_i)_{b,d}$.
Placing $H_{c,d}$ in position $c+q(d-1)$ gives precisely the
Kronecker-product entry
$(R_i)_{b,d}(L_i)_{c,a}$ in $T_i$.
Here $T_i$ has size $qp$, regardless of whether a particular mass
also admits a representation using a smaller matrix quotient.
Assume $\det L_1\det R_1\ne0$. The numerator factors are independent
of the denominators and no matrix commutation is imposed.

Let $\mathcal E$ be the module of scalar test entries: it is $A$ for
function tests and $\omega$ for tests obtained by evaluating
$F\otimes F^D\to\omega_C$. Put
\[
 V_{\mathcal E}\coloneq\mathcal E^{qp},\qquad
 Q_{\mathcal E}\coloneq V_{\mathcal E}/T_1V_{\mathcal E}.
\]
Choose a complex linear extension $\widetilde\mu$ of $\mu$ to the
rational tests needed below. Applying the subtraction argument of
Proposition~\ref{res:geo:allsolutions} to these test modules gives
every solution in the form
\begin{equation}\label{res:bq:eq:allsolutions}
 \widehat\mu[v]=\widetilde\mu[T_0T_1^{-1}v]+\nu[v],
 \qquad \nu\in Q_{\mathcal E}^{*}.
\end{equation}
Indeed, evaluating the rational term at $T_1v$ returns $\mu[T_0v]$,
and the difference of any two solutions vanishes on $T_1V_{\mathcal E}$.
The ordered product is
\begin{equation}\label{res:bq:eq:rationalorder}
 T_0T_1^{-1}
  =(R_0R_1^{-1})\otimes(L_1^{-1}L_0)^{\T}.
\end{equation}
Reversing $T_0$ and $T_1^{-1}$ would generally fail to solve the
functional equation. The choice of rational extension is part of
the data; changing it changes the particular solution and the
remaining mass in \eqref{res:bq:eq:allsolutions}.

The two test modules have different canonical dual quotients:
\begin{align}
 A^{qp}/T_1A^{qp}\quad&\text{pairs with}\quad
       \omega^{qp}/T_1^{\T}\omega^{qp},\label{res:bq:eq:functiontests}\\
 \omega^{qp}/T_1\omega^{qp}\quad&\text{pairs with}\quad
       A^{qp}/T_1^{\T}A^{qp}.\label{res:bq:eq:differentialtests}
\end{align}
In both lines the pairing is the sum of residues of
$\eta^{\T}T_1^{-1}v$. The representative $\eta$ carries the
differential in the first line; $v$ carries it in the second.

Here is the precise link with Theorem~\ref{res:bq:thm:quotient}.
The quotients are supported on finitely many points, and the line
$F$ has a frame $e$ on a neighbourhood of that finite set. For
\eqref{res:bq:eq:functiontests}, use $v\mapsto ev$ and
$\eta\mapsto e^{-1}\eta$. These identify the two quotients with
$M^{qp}/T_1M^{qp}$ and $N^{qp}/T_1^{\T}N^{qp}$ near their support,
which suffices to identify the quotients themselves. For
\eqref{res:bq:eq:differentialtests}, use $v\mapsto e^{-1}v$ and
$\eta\mapsto e\eta$. The resulting pairing is the one for
$T_1^{\T}$ with its arguments exchanged, since
\[
 (e^{-1}v)^{\T}(T_1^{\T})^{-1}(e\eta)
       =\eta^{\T}T_1^{-1}v.
\]
Its coevaluation is therefore obtained by exchanging the tensor
factors of $\mathcal C_{T_1^{\T}}$ and transporting them by the
chosen frame. A change of frame changes these representatives and
their dual coordinates together, leaving the functional unchanged.
This local choice does not require $F$ to be trivial on all of $U$.

This distinction matters when $C$ is not Gorenstein: the differential
module $\omega$ need not be free of rank one. Writing a differential
in fractional coordinates does not turn its coefficient module into
$A$. Likewise, the present balanced test quotient retains its
$A$-module relations; it does not impose those relations on a
bivariate pairing that originally had independent arguments.

\subsubsection{Computing the finite pairing from the quotient}
\label{res:bq:curvefiniteapplication}

The residue coordinates enter the finite construction in two places.
For polynomial division, under the regularity assumptions of
Subsection~\ref{sec:finite}, take the lifts $P_L,P_R$ of
\eqref{eq:polynomiallifts}. A row candidate $b$ satisfies polynomial
divisibility precisely when the class of $(a_Lb)^{\T}$ modulo
$P_L^{\T}A^q$ vanishes; for a column candidate $c$, the corresponding
condition is that the class of $a_Rc$ modulo $P_RA^p$ vanishes.
Theorem~\ref{res:bq:thm:quotient}, with a local frame at the finite
support as above, supplies residue functionals detecting these classes.
Evaluating the functionals on the chosen ambient basis vectors and
selecting independent rows gives the divisibility rows of $\Gamma_L$
and $\Gamma_R$. The remaining
rows test membership of the divided vectors in the prescribed finite
spaces, as in Proposition~\ref{prop:lifteddivision}. Thus residue tests
and the degree or section bounds together give the full image conditions
used in Theorems~\ref{thm:finite} and~\ref{thm:kernel}.

For masses, the relevant quotient is instead the one defined by the
denominator action $D_1$. Its residue coordinates give the mass entries
in those same finite formulas. They also give a direct Gram computation
by Subsection~\ref{res:sec:finitarytools}, as follows.
Choose $d$ row sections $f_i\in M^{1\times q}$ and $d$ column
sections $g_j\in N^{p\times1}$ spanning the prescribed spaces.
Their products give differential tests
$H(f_i,g_j)\coloneq[(f_i)_a(g_j)_b]_{a\in\{1,\ldots,q\},\,b\in\{1,\ldots,p\}}$, so the homogeneous
mass is an element of the dual of
$\Mat_{q\times p}(\omega)/D_1\Mat_{q\times p}(\omega)$.
Choose a basis $\zeta_1,\ldots,\zeta_m$ of this quotient and reduce
each test $H(f_i,g_j)$ modulo $D_1$:
\[
 [H(f_i,g_j)]=\sum_{a=1}^m c_{i,j,a}\zeta_a,
 \qquad
 (H_\nu)_{i,j}=\sum_{a=1}^m c_{i,j,a}\nu[\zeta_a].
\]
If $\zeta_a^D$ denotes the residue-dual basis, the reduction coefficients
are the pairings of $[H(f_i,g_j)]$, after column vectorization, with
$\zeta_a^D$ under the duality
\eqref{res:bq:eq:differentialtests}. This computes every $c_{i,j,a}$
by a residue, and the displayed finite sum then inserts the prescribed
mass parameters $\nu[\zeta_a]$. Reduction takes place in the actual
differential module, including its relations at singular points.

Next express $f_iS_L$ and $S_Rg_j$ in finite rational bases as in
\eqref{res:bq:eq:coordinates}. Their rational pairings form $G_{\mathrm{rat}}$,
and their coordinate matrices are $T_L,T_R$. The complete target Gram is
$\widehat G=T_LG_{\mathrm{rat}}T_R+H_\nu$. If it is invertible,
\eqref{res:bq:eq:finitekernel} gives the reproducing kernel, and
Proposition~\ref{res:pf:transformation} gives the next orthogonal rows
and columns with their prescribed leading classes. This is the
calculation to be performed after the geometric quotient is known.

The three pairings have separate roles. The original pairing is
$\beta_\mu$. The chosen rational extension gives
$\beta_{\mathrm{part}}(f,g)=
\widetilde\beta_\mu(fS_L,S_Rg)$, and the target pairing is
$\widehat\beta=\beta_{\mathrm{part}}+\nu$. Formula
\eqref{res:bq:eq:woodburykernel} compares the target and particular kernels
on the same target spaces when both Grams are invertible. These
finite Grams belong to the orthogonality problem; the geometric
kernel $S_F$ is determined by the curve and the line bundle.

\subsection{Finite residue corrections on surfaces}\label{res:sf:section}

On a curve, a regular denominator has a finite zero set, and the
quotients in Theorem~\ref{res:bq:matrixquotient} are finite dimensional.
On a surface, its zero set is usually a curve. For example,
$\CC[x,y]/(x)=\CC[y]$ has infinite complex dimension. To obtain a
finite residue problem, we impose a second equation: a scalar
polynomial $h$ that meets the determinant divisor in finitely many
points. The resulting quotient retains the multiplicities at these
points. Its dual describes those masses that satisfy both the matrix
equation and the additional scalar equation.

We first work in the plane, where polynomial divided differences give
an explicit tensor reproducing this quotient. Scalar
complete-intersection residue duality and its B\'ezoutian tensor are
classical; we use the formulation in
\cite[Proposition~4.3 and Lemma~4.4]{BMPBezoutians}. Residue methods
for polynomial matrix division are developed in
\cite{AnderssonWulcan}. The calculation below specifies the quotient,
the transposed dual and the order of the matrix factors needed in
the present perturbation problem. Its output is a finite residue
tensor. A global geometric Szeg\H{o} kernel on a surface would require
a different construction.

\subsubsection{A matrix equation and a scalar cut}\label{res:sf:setup}

Let $A\coloneq\CC[x,y]$, let $T\in\Mat_n(A)$ and put
$d\coloneq\det T$. Assume that $d\ne0$ and that $h\in A$ satisfies
\begin{equation}\label{res:sf:eq:finiteintersection}
 0<\dim_{\CC}A/(d,h)<\infty.
\end{equation}
Thus $d$ and $h$ have no common curve component, and their common
zero set $Z$ is finite. In particular, $h$ is a nonzerodivisor in
$A$, and multiplication by $d$ is injective in $A/(h)$. Indeed,
$d$ and $h$ have no common irreducible factor, so $h\mid da$
implies $h\mid a$. These are the regularity properties used below.
They do not require the points of $Z$ to be simple.

Fix the volume form $\dd x\wedge\dd y$ and the order $(d,h)$ of
the two equations. Define the scalar residue functional by
\begin{equation}\label{res:sf:eq:scalarresidue}
 \lambda_{d,h}(a)\coloneq
   \sum_{\xi\in Z}\operatorname{Res}_{\xi}
       \frac{a(x,y)\,\dd x\wedge\dd y}{d(x,y)h(x,y)}.
\end{equation}
Here $\operatorname{Res}_{\xi}$ is the local Grothendieck residue
with these ordered denominators. Analytically it is the integral
\[
 \frac1{(2\pi\mathrm{i})^2}
 \int_{\Gamma_\xi}
 \frac{a(x,y)\,\dd x\wedge\dd y}{d(x,y)h(x,y)},
 \qquad
 \Gamma_\xi:\quad |d|=r_1,\quad |h|=r_2,
\]
on a sufficiently small cycle enclosing only $\xi$, with orientation
$\dd\arg d\wedge\dd\arg h>0$. The algebraic residue is the same
functional and is computed by finite polynomial reduction below.
For denominators $x^r,y^s$ at the origin it extracts
$[x^{r-1}y^{s-1}]a$, by applying the one-variable residue formula
twice. Thus repeated zeros retain derivative information.
This convention gives
$\lambda_{x,y}(1)=1$. If every intersection point is simple, the
same functional is
\[
 \lambda_{d,h}(a)=
 \sum_{\xi\in Z}
   \frac{a(\xi)}{(d_xh_y-d_yh_x)(\xi)}.
\]
The residue definition in \eqref{res:sf:eq:scalarresidue} also applies
when this Jacobian vanishes. It then retains the corresponding
local multiplicities. It annihilates $(d,h)$, and
$([a],[b])\mapsto\lambda_{d,h}(ab)$ is a perfect pairing on
$A/(d,h)$. We give the scalar reconstruction calculation that
proves this perfection immediately before
Theorem~\ref{res:sf:coevaluation}.

We use columns for both arguments of the matrix residue pairing.
Define
\begin{equation}\label{res:sf:eq:quotients}
 Q_{T,h}\coloneq\frac{A^n}{TA^n+hA^n},\qquad
 Q^D_{T,h}\coloneq\frac{A^n}{T^{\T}A^n+hA^n}.
\end{equation}
The superscript $D$ records the dual role of the second quotient.
Its entries are understood to multiply the fixed volume form.
This identification is available here because the canonical module
of the plane is freely generated by $\dd x\wedge\dd y$.

\begin{theorem}[Duality for a finite matrix quotient]\label{res:sf:duality}
Under \eqref{res:sf:eq:finiteintersection}, both quotients in
\eqref{res:sf:eq:quotients} are finite dimensional, and
\begin{equation}\label{res:sf:eq:matrixresidue}
 b_{T,h}([f],[\eta])\coloneq
     \lambda_{d,h}\bigl(\eta^{\T}\operatorname{adj}(T)f\bigr)
\end{equation}
is a perfect bilinear pairing between them. Consequently every
complex linear functional on $Q_{T,h}$ has a unique expression
$[f]\mapsto b_{T,h}([f],[\eta])$ with $[\eta]\in Q^D_{T,h}$.
\end{theorem}

\begin{proof}
Write $J\coloneq\operatorname{adj}(T)$. The identities
$TJ=JT=dI_n$ show that $d$ annihilates both quotients. Since $h$
also annihilates them, they are quotients of $(A/(d,h))^n$ and are
finite dimensional.

First check that the formula is independent of the representatives.
Replacing $f$ by $f+Ta$ adds $d\eta^{\T}a$ to its numerator.
Replacing $\eta$ by $\eta+T^{\T}b$ adds $db^{\T}f$.
Adding a multiple of $h$ to either argument adds a numerator
divisible by $h$. All these terms have zero residue. In particular,
the transpose in the second quotient is exactly what makes this
second substitution valid.

We next show that a class pairing to zero with every opposite
class must itself be zero. Suppose first that
$b_{T,h}([f],[\eta])=0$ for every $\eta\in A^n$.
Taking $\eta=ae_i$, with arbitrary $a\in A$, shows that each
component of $Jf$ pairs to zero with all of $A/(d,h)$ under the
scalar residue pairing. Its perfection implies
\[
 Jf\in(d,h)A^n.
\]
Work now in $A/(h)$, denoting classes by bars. The last inclusion
means that $\overline J\,\overline f=\overline d\,\overline a$
for some $\overline a\in(A/(h))^n$. Multiplication by
$\overline T$ gives
\[
 \overline d\,\overline f
       =\overline d\,\overline T\,\overline a.
\]
Multiplication by $\overline d$ is injective, so cancellation
yields $\overline f=\overline T\,\overline a$.
Lifting this equality to $A$ gives $f\in TA^n+hA^n$.
Thus $[f]=0$ in $Q_{T,h}$.

For the other argument, write the same scalar as
$\lambda_{d,h}(f^{\T}J^{\T}\eta)$. If it vanishes for every
$f$, scalar duality gives $J^{\T}\eta\in(d,h)A^n$.
Apply the preceding cancellation argument with $T^{\T}$ and
$J^{\T}$ to obtain $\eta\in T^{\T}A^n+hA^n$.
Hence the pairing has zero radical on either side. Since the two
spaces are finite dimensional, the induced maps into their complex
duals are injective in both directions and therefore are
isomorphisms. This proves perfection and the asserted uniqueness.

The sign agrees with our normalization: for $T=xI_n$ and $h=y$,
the formula uses $d=x^n$ and $J=x^{n-1}I_n$ and gives the ordinary
dot product on $\CC^n$.
\end{proof}

\begin{remark}[Canonical modules on a singular surface]
\label{res:sf:canonicalmodule}
The duality argument has the same form on a Cohen--Macaulay affine
surface with coordinate ring $R$ and canonical module $\omega_R$.
Take $T\in\Mat_n(R)$ and write $d\coloneq\det T$.
If $h$ is regular on $R$, $d$ is regular modulo $h$, and $(d,h)$
has finite support, the two quotients are
\[
 \frac{R^n}{TR^n+hR^n},\qquad
 \frac{\omega_R^n}{T^{\T}\omega_R^n+h\omega_R^n}.
\]
The second quotient uses the canonical module, which need not have
one generator at a singular point.

Put $R_{d,h}\coloneq R/(d,h)$ and
$\Omega_{d,h}\coloneq\omega_R/(d,h)\omega_R$.
The scalar duality needed for the matrix argument is the
$R_{d,h}$-linear identification
\[
 \Theta:\Omega_{d,h}\xrightarrow{\ \simeq\ }
                \operatorname{Hom}_\CC(R_{d,h},\CC),
 \qquad (a\ell)(b)\coloneq\ell(ab).
\]
Here is a local description of this identification and of the
regularity needed for cancellation. At a point of the finite
support, complete the local ring and denote the result by
$\widehat R$. The ordered pair $(d,h)$ is a system of parameters.
In a Cohen--Macaulay local ring every system of parameters is a
regular sequence, so the order can be reversed from $(h,d)$ to
$(d,h)$ here.
The power-series ring $S\coloneq\CC[[u,v]]$ acts by $u\mapsto d$
and $v\mapsto h$. To see that $\widehat R$ is finite over $S$,
lift a complex basis of $\widehat R/(d,h)\widehat R$.
Reduction modulo $(d,h)$ writes every element as a linear
combination of these lifts plus terms divisible by $d$ or $h$.
Repeating this reduction on the remaining coefficients gives
power series in $d,h$ multiplying the same finite set of lifts.
These series converge in $\widehat R$, because the $(d,h)$-adic
topology agrees with the maximal-ideal topology when the quotient
has finite length. The lifts therefore generate $\widehat R$
over $S$.
It has depth two as an $S$-module, because the two parameters form
a regular sequence. The Auslander--Buchsbaum formula therefore
gives projective dimension zero over the regular local ring $S$,
so $\widehat R$ is finite free over $S$.
The canonical-module construction for a finite algebra over a
regular ring identifies the completed canonical module with
\[
 \widehat\omega_R\simeq
       \operatorname{Hom}_S(\widehat R,S)\,
                         \dd u\wedge\dd v;
\]
This is \cite[Theorem~21.15]{Eisenbud} in equal dimensions;
Corollaries~21.17 and~21.18 there give compatibility with
localization and completion. Choosing an $S$-basis of $\widehat R$ and reducing its
dual basis modulo $(u,v)$ gives
\[
 \frac{\widehat\omega_R}{(d,h)\widehat\omega_R}
 \simeq
 \operatorname{Hom}_\CC
       \bigl(\widehat R/(d,h)\widehat R,\CC\bigr).
\]
Both sides are unchanged by completion, since they have finite
length. Combining these identifications over the finitely many
support points gives $\Theta$. The same free-module description
shows that $h,d$ is a regular sequence on $\omega_R$ near those
points: multiplication by $v$, and then by $u$ modulo $v$, is
injective on the free $S$-module
$\operatorname{Hom}_S(\widehat R,S)$.

The residue trace is now a concrete linear functional:
\[
 \lambda_{d,h}([\eta])\coloneq\Theta([\eta])(1),
 \qquad
 \lambda_{d,h}(a[\eta])=\Theta([\eta])(a).
\]
Consequently the pairing
$(a,[\eta])\mapsto\lambda_{d,h}(a[\eta])$ is perfect: every
functional on $R_{d,h}$ is $\Theta([\eta])$ for a unique class
$[\eta]$, and a nonzero element of $R_{d,h}$ is detected by a linear
functional. With the ordered parameters $(d,h)$, this trace is the
local-duality residue of the generalized fraction with numerator
$\eta$ and denominator $(d,h)$. In the completed description it
is the constant coefficient of the $S$-linear functional
corresponding to $\eta$, evaluated at $1$; it is summed over the
support points. This is the coefficient realization of the
residue trace; compare \cite[Sections~7 and~9]{Lipman1984}.

To obtain the matrix pairing, put $J\coloneq\operatorname{adj}(T)$
and use $\lambda_{d,h}(\eta^{\T}Jf)$.
The identities $TJ=JT=dI_n$ show that this expression vanishes on
each of the submodules divided out above. If it vanishes for every
$\eta$, scalar perfection gives $Jf\in(d,h)R^n$.
Modulo $h$, write $\overline J\,\overline f=
\overline d\,\overline a$ and multiply by $\overline T$.
Cancellation of the regular element $\overline d$ gives
$\overline f=\overline T\,\overline a$.
For the other radical, scalar perfection gives
$J^{\T}\eta\in(d,h)\omega_R^n$.
Multiplication by $T^{\T}$ and cancellation of $d$ on
$\omega_R/h\omega_R$ give
$\eta\in T^{\T}\omega_R^n+h\omega_R^n$.
These arguments apply at every support point; outside that set the
quotients are zero. Thus both radicals vanish and the finite
matrix pairing is perfect. The explicit divided differences in
the next subsection are stated in the planar coordinates $x,y$.
\end{remark}

\subsubsection{The polynomial tensor that reproduces the quotient}
\label{res:sf:divideddifferences}

For a perfect pairing of two finite-dimensional spaces, its
\emph{coevaluation tensor} is the inverse pairing written without
choosing coordinates. More explicitly, if $e_1,\ldots,e_s$ is a
basis of the first space and $e_1^D,\ldots,e_s^D$ is its dual basis
in the second, the tensor is $\sum_i e_i\otimes e_i^D$. Contracting
the second factor against $f$ gives
$\sum_i e_i\,b(f,e_i^D)=f$. We now construct this tensor for
$b_{T,h}$ by polynomial operations on $T$ and $h$.

Use two points $P=(x,y)$ and $Q=(u,v)$, with the intermediate point
$U=(u,y)$. For a scalar polynomial or a polynomial matrix $F$, define
the entrywise divided differences
\begin{equation}\label{res:sf:eq:differences}
 \delta_xF\coloneq\frac{F(P)-F(U)}{x-u},\qquad
 \delta_yF\coloneq\frac{F(U)-F(Q)}{y-v}.
\end{equation}
Both are polynomial: the respective numerator vanishes when
$x=u$ or $y=v$. Thus the displayed quotients introduce no poles,
including on the diagonal $P=Q$.

For the scalar complete intersection, put
\begin{equation}\label{res:sf:eq:scalarbezout}
 B_{d,h}(P,Q)\coloneq
     \delta_xd\,\delta_yh-\delta_yd\,\delta_xh.
\end{equation}
The scalar B\'ezoutian identity states that
\begin{equation}\label{res:sf:eq:scalarreproduction}
 \lambda_{d,h}^{(Q)}
       \bigl(B_{d,h}(P,Q)a(Q)\bigr)
       \equiv a(P)\pmod{(d(P),h(P))}
 \qquad(a\in A).
\end{equation}
The superscript means that the residue acts in $u,v$, leaving
$x,y$ as polynomial coefficients. This is the coevaluation identity
for the scalar residue pairing, valid for the entire quotient,
including its nilpotent elements
\cite[Proposition~4.3 and Lemma~4.4]{BMPBezoutians}.
For $d=x$ and $h=y$, both the B\'ezoutian and the residue of $1$
are $1$, fixing the same orientation as
\eqref{res:sf:eq:scalarresidue}.

Here is a derivation that also explains the multiplicities. The
residue change-of-equations rule says that, if ordered denominator
columns satisfy $g=Af$, then
\[
 \operatorname{Res}_{g}(\det(A)\,a)
       =\operatorname{Res}_{f}(a)
\]
whenever the corresponding quotients are finite and the sequences
are regular; it applies with finite coefficient algebras as well
\cite[Section~2, especially Proposition~2.1]{BerensteinYgerResidues}.
Work over the coefficient algebra
$R\coloneq\CC[x,y]/(d(P),h(P))$. The divided differences give
in $R[u,v]$ the exact identity
\[
 \begin{bmatrix}d(Q)\\h(Q)\end{bmatrix}
 =\begin{bmatrix}\delta_xd&\delta_yd\\
                  \delta_xh&\delta_yh\end{bmatrix}
   \begin{bmatrix}u-x\\v-y\end{bmatrix}.
\]
The determinant is $B_{d,h}$. Both denominator sequences are
regular over $R$, and their quotients are finite free $R$-modules:
respectively $R\otimes_\CC A/(d,h)$ and $R$.
The change-of-equations rule therefore gives
\[
 \lambda_{d,h}^{(Q)}(B_{d,h}a(Q))
   =\operatorname{Res}_{u-x,v-y}(a(u,v))=a(x,y)
   \quad\hbox{in }R.
\]
The last equality is the two successive Cauchy coefficient
formulas. It remains valid for nilpotent $x,y$: on a local component
write $x=c+n$ with $n^N=0$, and expand
$1/(u-x)=\sum_{j=0}^{N-1}n^j/(u-c)^{j+1}$; the residue is the
finite Taylor expansion at $c$. Do the same for $v-y$.
This proves \eqref{res:sf:eq:scalarreproduction} on all remainder
classes. Finally, if $\lambda_{d,h}(ab)=0$ for every $b$, expand
$B_{d,h}=\sum_\ell p_\ell(P)q_\ell(Q)$ in finite tensors. The
reconstruction identity gives
$[a]=\sum_\ell[p_\ell]\lambda_{d,h}(q_\ell a)=0$.
Thus the scalar pairing has zero radical and is perfect.

\begin{theorem}[Matrix coevaluation for a general scalar cut]
\label{res:sf:coevaluation}
Under \eqref{res:sf:eq:finiteintersection}, define
\begin{equation}\label{res:sf:eq:matrixbezout}
 \mathcal C_{T,h}(P,Q)\coloneq
       \delta_xT\,\delta_yh-\delta_yT\,\delta_xh.
\end{equation}
Its class in $Q_{T,h}\otimes_{\CC}Q^D_{T,h}$ is the coevaluation
tensor of \eqref{res:sf:eq:matrixresidue}. The row index of
$\mathcal C_{T,h}$ belongs to the first quotient and its column
index to the second. Equivalently, for every $f\in A^n$,
\begin{equation}\label{res:sf:eq:reproduction}
 \lambda_{d,h}^{(Q)}\!\left(
     \mathcal C_{T,h}(P,Q)\operatorname{adj}(T(Q))f(Q)\right)
       \equiv f(P)\pmod{T(P)A^n+h(P)A^n}.
\end{equation}
The residue on a vector is taken entrywise. This statement permits
arbitrary intersection multiplicities and imposes no restriction
that $h$ depend on only one coordinate.
\end{theorem}

\begin{proof}
Put $J\coloneq\operatorname{adj}(T)$ and form
$\mathcal C_{J,h}$ by the same definition
\eqref{res:sf:eq:matrixbezout}. We first establish an identity of
polynomial matrices before taking any residue or quotient:
\begin{align}
 B_{d,h}(P,Q)I_n
   ={}&\mathcal C_{T,h}(P,Q)J(Q)
          +T(P)\mathcal C_{J,h}(P,Q)\notag\\
     &+\bigl[h(P)-h(Q)\bigr]\delta_xT\,\delta_yJ.
     \label{res:sf:eq:adjugateidentity}
\end{align}
To derive it, take the two consecutive differences of $TJ=dI_n$.
The product rule, with the factors kept in their original order,
gives
\begin{align*}
 \delta_xd\,I_n&=\delta_xT\,J(U)+T(P)\delta_xJ,\\
 \delta_yd\,I_n&=\delta_yT\,J(Q)+T(U)\delta_yJ.
\end{align*}
Multiply the first equality by $\delta_yh$ and subtract the
second multiplied by $\delta_xh$. After subtracting
$\mathcal C_{T,h}(P,Q)J(Q)$, the remaining terms are
\begin{align*}
 &\delta_xT\bigl[J(U)-J(Q)\bigr]\delta_yh\\
 &\hspace{1em}+T(P)\delta_xJ\,\delta_yh
                 -T(U)\delta_yJ\,\delta_xh.
\end{align*}
Use $J(U)-J(Q)=(y-v)\delta_yJ$ and the scalar identity
\[
 (y-v)\delta_yh=h(P)-h(Q)-(x-u)\delta_xh.
\]
Since $(x-u)\delta_xT=T(P)-T(U)$, the expression becomes
\begin{align*}
 &\bigl[h(P)-h(Q)\bigr]\delta_xT\,\delta_yJ\\
 &\quad-\bigl[T(P)-T(U)\bigr]\delta_yJ\,\delta_xh
       +T(P)\delta_xJ\,\delta_yh-T(U)\delta_yJ\,\delta_xh.
\end{align*}
The two terms containing $T(U)$ cancel. The other terms containing
$T(P)$ give $T(P)\mathcal C_{J,h}$, proving
\eqref{res:sf:eq:adjugateidentity}. Only the divided differences of the
scalar $h$ have been moved past matrices.

Now multiply \eqref{res:sf:eq:adjugateidentity} on the right by $f(Q)$
and apply $\lambda_{d,h}^{(Q)}$. By
\eqref{res:sf:eq:scalarreproduction}, the left-hand side is $f(P)$
modulo $(d(P),h(P))A^n$. The second term on the right belongs to
$T(P)A^n$. In the last term, the contribution with $h(P)$ belongs
to $h(P)A^n$, while the contribution with $h(Q)$ has zero residue.
Finally,
\[
 d(P)A^n\subseteq T(P)A^n,
 \qquad d(P)I_n=T(P)J(P).
\]
Passing to $Q_{T,h}$ therefore leaves precisely
\eqref{res:sf:eq:reproduction}.

The contraction in that formula pairs the second tensor factor
with $[f]$ by \eqref{res:sf:eq:matrixresidue}; this explains the order
$\mathcal C_{T,h}J(Q)f(Q)$. Theorem~\ref{res:sf:duality} makes this
pairing perfect. Hence a tensor whose contraction is the identity
is its unique coevaluation tensor, as asserted.
\end{proof}

For a cut $h=k(y)$, the second summand in
\eqref{res:sf:eq:matrixbezout} vanishes, and the formula reduces to
\begin{equation}\label{res:sf:eq:coordinatecut}
 \mathcal C_{T,k}(P,Q)=
 \frac{T(x,y)-T(u,y)}{x-u}\,
 \frac{k(y)-k(v)}{y-v}.
\end{equation}
Repeated roots of $k$ are allowed. For a general cut, both
summands are needed, as the following small example shows. Take
\[
 T=\begin{bmatrix}x&y\\0&x\end{bmatrix},\qquad h=y^2-x.
\]
Then $d=x^2$ and $A/(d,h)\simeq\CC[y]/(y^4)$ by substituting
$x=y^2$. The scalar residue is the coefficient of $y^3$ in this
remainder. Its normalization agrees with
$\lambda_{d,h}(d_xh_y-d_yh_x)=\lambda_{d,h}(4xy)=4$.
The four divided differences give
\[
 \delta_xT=I_2,\qquad
 \delta_yT=\begin{bmatrix}0&1\\0&0\end{bmatrix},\qquad
 \delta_xh=-1,\qquad\delta_yh=y+v.
\]
Consequently
\begin{equation}\label{res:sf:eq:nonverticalexample}
 \mathcal C_{T,h}(P,Q)=
       \begin{bmatrix}y+v&1\\0&y+v\end{bmatrix}.
\end{equation}
The off-diagonal entry comes from the second summand of
\eqref{res:sf:eq:matrixbezout}. In the matrix quotient the relations
are $y^2e_1=0$ and $ye_1+y^2e_2=0$, after eliminating $x$.
They give the basis $e_1,e_2,ye_2,y^2e_2$, on which
we now check the claimed reproduction. In the residue variable,
put $u=v^2$ and take the coefficient of $v^3$. Multiplication by
the adjugate gives
\[
 \begin{bmatrix}y+v&1\\0&y+v\end{bmatrix}
 \begin{bmatrix}v^2&-v\\0&v^2\end{bmatrix}
 =\begin{bmatrix}yv^2+v^3&-yv\\0&yv^2+v^3\end{bmatrix}.
\]
Consequently contraction sends the four basis vectors to
$e_1,e_2,ye_2,-ye_1$, respectively. The last image equals
$y^2e_2$ in the quotient by its defining relation, so all four
classes are reproduced. The basis also follows directly from the
relations: $y^2e_1=0$ and $ye_1=-y^2e_2$ imply $y^3e_2=0$;
the four remaining classes are independent, as the residue Gram
against $e_2,e_1,ye_1,y^2e_1$ in the transposed quotient is
\[
 \begin{bmatrix}
 0&0&1&0\\0&0&0&-1\\1&0&-1&0\\0&-1&0&0
 \end{bmatrix},
 \qquad \det=1.
\]

\subsubsection{From a finite residue to a bilateral correction}
\label{res:sf:corrections}

The constructions above allow $T$ to be an ordered product of
noncommuting polynomial matrices. The second equation is the
scalar $h$, so it commutes with all the matrix operations in the
proofs.
For row--column pairings, the orientation is the same as in
Remark~\ref{res:bq:rem:roworientation}. Namely, if $f$ is a polynomial
row, $g$ a polynomial column and $c\ne0$, then
\[
 \nu(f,g)\coloneq c\,\lambda_{d,h}
        \bigl(f\operatorname{adj}(T)g\bigr)
       =c\,b_{T^{\T},h}(f^{\T},g).
\]
Its quotient pair is therefore $Q_{T^{\T},h},Q_{T,h}$, and its
coevaluation in this order is $c^{-1}\mathcal C_{T^{\T},h}$.
The transpose changes matrix indices without interchanging $P$
and $Q$. Passing to the column--row convention of the finite
reproducing kernel then exchanges the tensor factors.

More generally, return to the polynomial test matrix $H$ of the
bilateral functional equation. As in
\eqref{res:bq:eq:vectorization}, column vectorization gives
\[
 \operatorname{vec}(L_1^{\T}HR_1^{\T})
       =T_1\operatorname{vec}(H),\qquad
 T_1\coloneq R_1\otimes L_1^{\T}.
\]
Assume that $T_1$ and a chosen scalar $h$ satisfy
\eqref{res:sf:eq:finiteintersection}. The masses which annihilate both
the image of $T_1$ and multiplication by $h$ are exactly the
functionals on $Q_{T_1,h}$. Theorem~\ref{res:sf:duality} writes every
such mass uniquely as
\begin{equation}\label{res:sf:eq:testmass}
 \nu[H]=\lambda_{\det T_1,h}
       \bigl(\eta^{\T}\operatorname{adj}(T_1)
                                      \operatorname{vec}(H)\bigr),
 \qquad [\eta]\in Q^D_{T_1,h}.
\end{equation}
The coevaluation of Theorem~\ref{res:sf:coevaluation} recovers these
finite remainder classes. The scalar equation specifies the finite
sector: the full dual of $\coker T_1$ may also contain masses that
do not annihilate $h$.

To use such a mass in an orthogonality problem, choose the prescribed
finite row and column spaces and evaluate its entries
$H_\nu=[\nu(f_i,g_j)]_{i,j}$. The finite-coordinate calculation of
Proposition~\ref{res:bq:prop:finitekernel} then gives
$\widehat G=T_LG_{\mathrm{rat}}T_R+H_\nu$. Its invertibility is exactly the
normality condition on those spaces. If both $\widehat G$ and
$G_{\mathrm{part}}=T_LG_{\mathrm{rat}}T_R$ are invertible,
\eqref{res:bq:eq:woodburykernel} gives the mass correction to their
finite reproducing kernels, of separation rank $\rank H_\nu$.
This step uses the residue data together with the rational moments
and the chosen polynomial spaces; perfection of the residue
quotient alone does not imply normality of the complete problem.

\subsubsection{A polynomial residue tensor on the singular cone}
\label{res:sf:cone}

The finite residue construction also gives an explicit reproducing
tensor on a singular surface. In this example the calculation reduces
to the plane by writing each function on the surface in two
components. The resulting tensor reproduces classes in a finite
polynomial quotient. Its construction requires no orthogonality
measure or choice of polynomial degree bounds.

Take the complex coordinate algebras
\[
 A\coloneq\CC[x,y,z]/(z^2-x^2-y^2),\qquad B\coloneq\CC[x,y].
\]
The polynomial defining the surface is irreducible: its quadratic
form has rank three and therefore cannot be a product of two linear
forms. The surface has its unique singular point at the origin,
where the three partial derivatives vanish. Every class in $A$ has
a unique representative $p_0(x,y)+zp_1(x,y)$, because division by
the monic polynomial $z^2-x^2-y^2$ leaves degree at most one in $z$.
Thus $A$ is a free $B$-module with basis $1,z$. Define
\begin{equation}\label{res:sf:cone-phi}
 \varphi[p_0+zp_1]\coloneq p_1.
\end{equation}
This map takes the coefficient of $z$ after using the surface
equation to reduce higher powers. The matrix of the $B$-valued
pairing $(p,q)\mapsto\varphi[pq]$ in the basis $1,z$ is
\[
 \begin{bmatrix}0&1\\1&0\end{bmatrix}.
\]
It is invertible over $B$, including above the singular point.
Consequently every $B$-linear functional on $A$ is uniquely of
the form $p\mapsto\varphi[cp]$ for some $c\in A$. In the
canonical-module description associated with the finite projection
to the $(x,y)$-plane, this says that
$\varphi\,\dd x\wedge\dd y$ is a generator. The reproduction
calculation below uses the coefficient map \eqref{res:sf:cone-phi}
directly.

Fix a real parameter $\varepsilon$ and put
\begin{equation}\label{res:sf:cone-scalars}
 a\coloneq1+\varepsilon^2,\qquad
 h\coloneq ax^2+\varepsilon^2y^2,\qquad
 k\coloneq y^2+\varepsilon^2.
\end{equation}
The parameter is real, while $x,y,z$ are coordinate classes in
the complex algebra $A$. In particular, the finite quotients below
are taken over $\CC$; their supporting points need not be real.
The matrix needed for the residue calculation is
\begin{equation}\label{res:sf:cone-product}
 T\coloneq\begin{bmatrix}h&2xy\\0&2h\end{bmatrix},\qquad
 \operatorname{adj}T=\begin{bmatrix}2h&-2xy\\0&h\end{bmatrix},
 \qquad\det T=2h^2.
\end{equation}
This matrix also has a bilateral denominator interpretation. Put
$\ell\coloneq x-\mathrm{i}\varepsilon z$ and
$r\coloneq x+\mathrm{i}\varepsilon z$, and choose
\[
 L_1=\begin{bmatrix}\ell&y\\0&\ell\end{bmatrix},\qquad
 R_1=\begin{bmatrix}r&y\\0&2r\end{bmatrix}.
\]
The surface equation gives $\ell r=h$ and $\ell+r=2x$.
Consequently $R_1L_1=T$, although the factors themselves involve
the surface coordinate $z$ and do not commute in general.
All entries of their product belong to $B$.
The polynomial $h$ is monic up to
the nonzero factor $a$ as a polynomial in $x$, and $k$ is monic
in $y$. Their simultaneous quotients can therefore be computed by
successive polynomial division, including at $\varepsilon=0$.

For $j\in\{1,2\}$ and $p\in B$, divide first in $x$ by $(h/a)^j$ and
then divide the coefficients of the remainder in $y$ by $k$.
Denote the final remainder by $\operatorname{rem}_j(p)$. It has
degree less than $2j$ in $x$ and less than two in $y$. Define
\begin{equation}\label{res:sf:cone-residues}
 \rho_j[p]\coloneq
   a^{-j}[x^{2j-1}y]\operatorname{rem}_j(p),\qquad
 \Lambda_j\coloneq\rho_j\circ\varphi.
\end{equation}
Here $[x^uy^v]$ extracts the indicated coefficient. The factor
$a^{-j}$ compensates for the leading coefficient $a^j$ of $h^j$
in $x$. With the residue convention of
\eqref{res:sf:eq:scalarresidue}, $\rho_j$ is the scalar residue for
the ordered denominators $(h^j,k)$ and differential
$\dd x\wedge\dd y$.
Indeed, summing the residues in $x$ of a reduced numerator over
$h^j$ picks its $x^{2j-1}$ coefficient divided by $a^j$;
terms of lower $x$ degree have no $x^{-1}$ term at infinity.
The subsequent sum of the $y$ residues over the monic quadratic
$k$ picks the coefficient of $y$. This gives exactly
\eqref{res:sf:cone-residues}, including when either polynomial has
multiple roots.

These functionals give perfect finite pairings. To see this without
assuming that the supporting points are simple, first work over
$R\coloneq\CC[y]/(k)$. Division by the monic polynomial $(h/a)^j$
gives the free $R$-basis $1,x,\ldots,x^{2j-1}$. The pairing that
takes the coefficient of $x^{2j-1}$ in a product has zeros whenever
the sum of the two basis exponents is less than $2j-1$, and has
ones when that sum equals $2j-1$. Reversing the columns of its
matrix therefore gives a triangular matrix with diagonal entries
one. Its determinant is a unit in $R$.

Next, the coefficient of $y$ gives an invertible pairing on $R$:
its matrix in the basis $1,y$ is
$\begin{bmatrix}0&1\\1&0\end{bmatrix}$, since $y^2=-\varepsilon^2$.
Composing the two coefficient pairings and multiplying by
$a^{-j}$ proves that $([p],[q])\mapsto\rho_j[pq]$ is perfect on
$B/(h^j,k)$. Composing once more with \eqref{res:sf:cone-phi} proves
that $([p],[q])\mapsto\Lambda_j[pq]$ is perfect on $A/(h^j,k)$.
This last quotient has the $8j$ basis classes
$x^uy^vz^w$, with $u\in\{0,\ldots,2j-1\}$ and
$v,w\in\{0,1\}$.

Polynomial multiplication of a functional means
$(b\Lambda)[p]\coloneq\Lambda[bp]$. The definitions imply
\begin{equation}\label{res:sf:cone-residue-relations}
 h\Lambda_2=\Lambda_1,\qquad h\Lambda_1=0,\qquad
 h^2\Lambda_2=0,\qquad k\Lambda_1=k\Lambda_2=0.
\end{equation}
For the first identity, divide a polynomial in $B$ by $h$ and
write its remainder as $c_0(y)+xc_1(y)$. Modulo $h^2$, its
product with $h$ has $x^3$ coefficient $ac_1(y)$. The factors
$a^{-2}$ and $a^{-1}$ in \eqref{res:sf:cone-residues} give
$h\rho_2=\rho_1$. Since $h$ belongs to $B$, it commutes with
$\varphi$, giving the stated identity for $\Lambda_j$. The
remaining identities follow because the remainders annihilate
the respective defining ideals.

Define the matrix functional
\begin{equation}\label{res:sf:cone-functionals}
 \nu\coloneq(\operatorname{adj}T)\Lambda_2
 =\begin{bmatrix}2\Lambda_1&-2xy\Lambda_2\\0&\Lambda_1\end{bmatrix}.
\end{equation}
For a row $f$ and a column $g$, it acts as
$\nu(f,g)=\Lambda_2[f\operatorname{adj}(T)g]$.
It is a homogeneous mass for the denominator pair just defined.
Indeed, $T=R_1L_1$ and the adjugate product rule give
\begin{align*}
 L_1\operatorname{adj}(T)R_1
 &=L_1\operatorname{adj}(L_1)\operatorname{adj}(R_1)R_1\\
 &=\det(L_1)\det(R_1)I_2=2h^2I_2.
\end{align*}
Therefore
\[
 \nu(fL_1,R_1g)=2\Lambda_2[h^2fg]=0,
 \qquad k\nu=0,
\]
by \eqref{res:sf:cone-residue-relations}. Thus both defining
annihilation equations hold without choosing an underlying measure.

The adjugate identities and \eqref{res:sf:cone-residue-relations} show that this
pairing depends only on the classes of $f^{\T}$ and $g$ in
\[
 Q_L\coloneq A^2/(T^{\T}A^2+kA^2),\qquad
 Q_R\coloneq A^2/(TA^2+kA^2),
\]
respectively. For example, adding $Tb$ to $g$ contributes
$2\Lambda_2[h^2fb]=0$. Adding a multiple of $k$ also contributes
zero. These quotients are finite dimensional, because both $h^2$
and $k$ annihilate them. More precisely, each has the sixteen basis
classes
\[
 x^uy^vz^we_i,\qquad u,v,w\in\{0,1\},\quad i\in\{1,2\}.
\]
To obtain this basis, first separate the coefficients of $1,z$.
Over $\CC[y]/(k)$, both $T$ and $T^{\T}$ are quadratic matrix
polynomials in $x$ with invertible leading matrix
$\operatorname{diag}(a,2a)$. Polynomial division therefore leaves
a unique vector remainder of degree less than two in $x$.
The two coefficients in $y$ and the two coefficients in $z$ give
the displayed basis. This also proves its independence.

The pairing is perfect on these quotients. There is a direct
determinant check: order the eight scalar monomials $x^uy^vz^w$
lexicographically by $(u,v,w)$. The matrix $H$ of $\Lambda_1[pq]$
is $a^{-1}$ times the reversal matrix. Only complementary exponent
triples contribute the coefficient $xyz$. Thus $\det H=a^{-8}$.
The matrix of $\nu$, listing component $1$ before component $2$,
has the block form
$\begin{bmatrix}2H&*\\0&H\end{bmatrix}$, and its determinant is
\[
 2^8(\det H)^2=\frac{256}{(1+\varepsilon^2)^{16}}.
\]
It remains nonzero at $\varepsilon=0$, where the support reaches
the vertex. A \emph{coevaluation tensor} for this
pairing reproduces a quotient class by contraction with the
opposite argument, as in the following explicit identity.

\begin{proposition}[Polynomial coevaluation on the cone]
\label{res:sf:cone-coevaluation}
For the independent arguments $P=(x,y,z)$ and $Q=(u,v,w)$, the
coevaluation of \eqref{res:sf:cone-functionals} has the representative
\begin{equation}\label{res:sf:cone-coevaluation-formula}
 \mathcal C_\nu(P,Q)\coloneq
 \frac{(z+w)(y+v)}2
 \begin{bmatrix}a(x+u)&0\\2y&2a(x+u)\end{bmatrix}.
\end{equation}
It represents a tensor in $Q_L\otimes_\CC Q_R$, with the first
factor indexed by its rows and the second by its columns. For
every column $f\in A^2$,
\begin{equation}\label{res:sf:cone-reproduction}
 \Lambda_{2,Q}\!\left[
   \mathcal C_\nu(P,Q)\operatorname{adj}(T(Q))^{\T}f(Q)
 \right]
 =f(P)\pmod{T(P)^{\T}A^2+k(P)A^2}.
\end{equation}
Here $\Lambda_{2,Q}$ acts entry by entry in the $Q$ variables.
The transposed identity, with $P,Q$ exchanged, reproduces the
classes of $Q_R$. The displayed representative is polynomial in
$\varepsilon$ and remains defined at the singular vertex.
\end{proposition}

\begin{proof}
First apply Theorem~\ref{res:sf:coevaluation} on the plane to
$T^{\T}$ and $k$. Their divided differences are
\[
 \frac{T(x,y)^{\T}-T(u,y)^{\T}}{x-u}
 =\begin{bmatrix}a(x+u)&0\\2y&2a(x+u)\end{bmatrix},
 \qquad
 \frac{k(y)-k(v)}{y-v}=y+v.
\]
Let $B_T(x,y;u,v)$ be the product of these two displayed
expressions. The theorem pairs this tensor with the scalar
residue for the denominators $(\det T,k)=(2h^2,k)$. That scalar
residue is $\rho_2/2$. Consequently its reproduction identity
for a polynomial column $b\in B^2$ reads
\[
 \frac12\rho_{2,Q}\!\left[
 B_T(x,y;u,v)\operatorname{adj}(T(u,v))^{\T}b(u,v)
 \right]
 =b(x,y)\pmod{T(x,y)^{\T}B^2+k(x,y)B^2}.
\]
Thus $B_T/2$ is the tensor appropriate to the functional
$\rho_2$ used in \eqref{res:sf:cone-residues}. This explains the
factor $1/2$ in \eqref{res:sf:cone-coevaluation-formula}.

It remains to pass from the plane to the surface. Write
$f(Q)=f_0(u,v)+wf_1(u,v)$ with $f_0,f_1\in B^2$. Reduction
by $w^2=u^2+v^2$ gives
\[
 (z+w)(f_0+wf_1)
   =zf_0+(u^2+v^2)f_1+w(f_0+zf_1).
\]
Taking its $w$ coefficient therefore returns $f_0+zf_1$.
The other matrices in the reproduction identity have entries in
$B$, so they commute with this coefficient extraction. Applying
the planar identity separately to $f_0$ and $f_1$ now gives
$f_0(x,y)+zf_1(x,y)$ modulo $T^{\T}A^2+kA^2$, which is
\eqref{res:sf:cone-reproduction}. The identity on the other quotient
follows from the corresponding transposed planar identity and
the same coefficient extraction. This proof uses no division by
$z$ or by a branch discriminant, so it also applies at the vertex.
\end{proof}

This polynomial residue tensor reproduces the finite quotient
at every parameter value, including at the cone vertex.

\section{An explicit bilateral family on singular curves}\label{sec:families}\label{res:bf:section}

The family constructed in this section has explicit geometric and
algebraic data: the curve, the two spaces of sections, the four
polynomial matrix factors, the added mass and the resulting polynomial
blocks. Normality is established for every real parameter, including
the value at which the denominator zeros meet the singular point.
This collision of zeros is the \emph{confluence} considered here;
the underlying curve remains fixed.

The construction uses different left and right section modules whose
products are canonical differentials. Their scalar Gram matrices are
shown to be invertible, providing the initial normal pairings.
Residue duality on the corresponding quotients supplies the masses,
with the quotient pairings remaining nondegenerate when their support
reaches the singular point. The resulting matrix mass has global rank
$4r$, and an explicit constant minor shows that this rank is attained
on every section space with $N\ge2r$, including at confluence.
These scalar pairings determine upper triangular $2\times2$
matrix pairings. The common triangular structure makes their Gram
determinants products of the two diagonal block determinants. These
blocks control normality,
meaning invertibility of the Gram matrix, while the nonzero
off-diagonal entries retain the order dependence of the matrix factors.
An explicit complex coefficient is chosen so that the same mass
coupling applies to every real parameter and every finite space in
the family. The proof also shows that almost every complex coupling
has this simultaneous property. Once these nonvanishing determinants
are established, the Gram formulas give the normalized polynomials
and their regular limits by evaluation at the confluent parameter.

The finite spaces have the same scalar bound $N\ge r$ in both
components. Explicit bases and the nonvanishing of all required
determinants establish the existence of normalized polynomials at
precisely these bounds. The matrix character of the denominators is
verified separately from triangularity: their entry ideals prevent
either denominator from becoming a scalar matrix under invertible
polynomial changes of rows and columns.

\subsection{The curve and two explicit modules of sections}

An \emph{acyclic} line bundle, with no global sections and no first
cohomology, provides the geometric data for the construction. These
vanishing conditions give matching dimensions for the initial section
spaces and determine the geometric residue kernel. Pairing with the
canonical dual then supplies the differentials needed for the residue
calculation.

For an integer $r\ge2$, consider the coordinate algebra
\begin{equation}\label{res:bf:curve-algebra}
 A_r\coloneq\CC+t^r\CC[t]
       =\CC[t^r,t^{r+1},\ldots,t^{2r-1}].
\end{equation}
Let $C_r$ be the curve obtained by gluing $\operatorname{Spec}A_r$
and $\operatorname{Spec}\CC[t^{-1}]$ along
$\operatorname{Spec}\CC[t,t^{-1}]$. Its normalization is
$\mathbb P^1_t$, and the point $\infty$ is smooth. For any integer
$b\ge r$, the function $t^b$ defines a finite map
$\pi:C_r\to\mathbb P^1$, so $C_r$ is projective. The quotient
$\CC[t]/A_r$ has basis $t,\ldots,t^{r-1}$. Consequently the
arithmetic genus is $g=r-1$. The canonical module at $t=0$ will have $r-1$ minimal generators.
For $r\ge3$ it is therefore not free of rank one: this is the
non-Gorenstein singularity used below.

The affine section modules and the coefficients of their canonical differentials
are described by the following polynomial spaces:
\begin{align}
 \ell_0&\coloneq1+t^{r-1},
 &M_r&\coloneq\CC\ell_0+t^r\CC[t],\label{res:bf:modules-left}\\
 K_r&\coloneq\Span_{\CC}\{1,t,\ldots,t^{r-2}\}+t^r\CC[t],
 &N_r&\coloneq
 \{k\in\CC[t]:[t^{r-1}]k=-[t^0]k\}.
 \label{res:bf:modules-right}
\end{align}
Here $[t^j]k$ denotes the coefficient of $t^j$ in $k$. In
particular,
\[
 N_r=\CC(1-t^{r-1})
       +\Span_{\CC}\{t,\ldots,t^{r-2}\}+t^r\CC[t],
\]
where the middle span is zero if $r=2$. Each of these spaces is an
$A_r$-module. The polynomial factor of a differential in the
canonical module belongs to $K_r$:
\begin{equation}\label{res:bf:canonical-module}
 \omega_{A_r}=K_r\frac{\dd t}{t^r}.
\end{equation}
Indeed, the residue condition
$\operatorname{Res}_0(a\eta)=0$ for every $a\in A_r$ excludes
$t^{-1}\dd t$ and every $t^j\dd t$ with $j<-r$, and allows
$t^{-r}\dd t,\ldots,t^{-2}\dd t$ and all polynomial
differentials. These are exactly the terms in
\eqref{res:bf:canonical-module}. More explicitly, if
$\eta=\sum_j c_jt^j\dd t$, the test $a=1$ forces $c_{-1}=0$,
and the tests $a=t^m$, $m\ge r$, force $c_{-m-1}=0$.
No other negative coefficient is constrained. The polynomial
factors $1,t,\ldots,t^{r-2}$ generate $K_r$ over $A_r$; their
classes are independent modulo the maximal ideal times $K_r$.
Indeed every product by a nonconstant element of $A_r$ has degree
at least $r$. This also verifies the stated number of generators.

\begin{proposition}[The acyclic line bundle]\label{res:bf:curve}
The module $M_r$, glued to $\OO$ at infinity by its inclusion in
$\CC[t,t^{-1}]$, defines a line bundle $\mathcal L_r$ of degree
zero on $C_r$. The line bundle
\[
 F_r\coloneq\mathcal L_r((r-2)\infty)
\]
has $H^0(C_r,F_r)=H^1(C_r,F_r)=0$. Its canonical dual
$F_r^D\coloneq\operatorname{Hom}(F_r,\omega_{C_r})$ has affine
module $N_r\dd t/t^r$, and
$M_rN_r\subset K_r$.
\end{proposition}

\begin{proof}
Away from $t=0$, the element $t^r$ is a unit and $M_r$ is the
structure sheaf, whereas at the origin $M_r$ is generated by $\ell_0$.
The latter description follows in the local ring of the singular curve
from the identity
\[
 \frac{t^r}{1+t^{r-1}}
   =\frac{t^r-t^{2r-1}}{1-t^{2r-2}}\in(A_r)_0.
\]
The denominator has constant term one and belongs to $A_r$.
Multiplying the numerator by an arbitrary polynomial still leaves
all its exponents at least $r$. Thus every element of
$t^r\CC[t]$ is a local multiple of $\ell_0$, and the module is
locally free of rank one.

The quotient $\CC[t]/M_r$ has dimension $r-1$, the same as
$\CC[t]/A_r$. The normalization sequences for these two modules
therefore give $\chi(\mathcal L_r)=\chi(\OO_{C_r})$, hence
$\deg\mathcal L_r=0$. A global section of $F_r$ is an element
of $M_r$ whose degree is at most $r-2$. Such an element has the
form $c(1+t^{r-1})+t^rp(t)$. Its degree bound first forces $c=0$
and then $p=0$. Hence $H^0(C_r,F_r)=0$. Since
$\deg F_r=r-2=g-1$, Riemann--Roch gives $\chi(F_r)=0$, and
$H^1(C_r,F_r)=0$ follows.

The canonical dual is determined by those $k\in\CC(t)$ satisfying
$kM_r\subset K_r$. For such an element, both $k\ell_0$ and $kt^r$
are polynomials. Since $\ell_0$ and $t^r$ are relatively prime,
$k$ must itself be a polynomial. Its product with $t^r\CC[t]$
already belongs to $K_r$, so the condition on the other generator
is exactly
\[
 [t^{r-1}](k\ell_0)=[t^{r-1}]k+[t^0]k=0.
\]
Thus $(K_r:M_r)=N_r$, proving the assertion about $F_r^D$ and
the inclusion $M_rN_r\subset K_r$.
\end{proof}

For $N\ge r$, the complete spaces of sections
\[
 U_N\coloneq H^0(C_r,F_r(N\infty)),\qquad
 V_N\coloneq H^0(C_r,F_r^D(N\infty))
\]
have the following ordered bases:
\begin{align}
 \boldsymbol f_N
   &\coloneq[1+t^{r-1},t^r,t^{r+1},\ldots,t^{N+r-2}],
   \label{res:bf:left-basis}\\
 \boldsymbol g_N
   &\coloneq[1-t^{r-1},t,\ldots,t^{r-2},t^r,t^{r+1},\ldots,t^N].
   \label{res:bf:right-basis}
\end{align}
The entries of the second list are multiplied by $\dd t/t^r$
when regarded as sections; products are expressed through their
polynomial parts, with this differential factor understood.
The degree bound $N+r-2$ in the first list comes from the twist
of the original line bundle. The bound in the second list is
determined by its global dual,
\[
 F_r^D=\mathcal L_r^{-1}\otimes\omega_{C_r}(-(r-2)\infty).
\]
The order at infinity of $k\dd t/t^r$ is $r-2-\deg k$.
Membership in $F_r^D(N\infty)$ requires that order to be at
least $r-2-N$, giving the bound $\deg k\le N$. The lists have
$1+(N-1)=N$ and $1+(r-2)+(N-r+1)=N$ elements, respectively.
They are therefore bases of two spaces of the same dimension.

A finite projection relates these section spaces to vector polynomials
on the line. For $\pi=t^b$ with $b\ge r$, acyclicity and the
splitting of vector bundles on $\mathbb P^1$ give
\[
 \pi_*F_r=\OO(-1)^b,\qquad \pi_*F_r^D=\OO(-1)^b.
\]
The pushforwards are torsion-free of rank $b$ over the smooth
base $\mathbb P^1$, so they are vector bundles. A finite map
has no higher direct images, and hence preserves the cohomology
of these sheaves. In the splitting into line bundles on
$\mathbb P^1$, a summand $\OO(a)$ has both cohomology groups
zero exactly when $a=-1$. Serre duality gives the same cohomology
vanishing for $F_r^D$. Thus the bounds $N=bk$ correspond to
polynomial components of degree at most $k-1$ in both projected
bundles. The modules on the curve remain $M_r$ and
$N_r\dd t/t^r$; projection does not replace the latter by a
free $A_r$-module.

The same module relations determine an explicit geometric kernel:
\begin{equation}\label{res:bf:geometric-kernel}
 S_r(t,s)=
 \left[\frac{t^{r-1}}{s^{r-1}(t-s)}-\frac1{s^r}\right]\dd s
 =\left[\frac{s\,t^{r-1}}{t-s}-1\right]\frac{\dd s}{s^r}.
\end{equation}
Its normalization and compatibility with the two modules follow
from its local expansions. At a fixed smooth point $t$, the factor
$t^{r-1}/s^{r-1}$ is one at $s=t$, so the polar term is $\dd s/(t-s)$,
whose residue in the second variable is $-1$, the convention used
for the geometric kernel. In the first variable,
the coefficients of $t^0$ and $t^{r-1}$ are both $-s^{-r}$ and
the intervening coefficients are zero, as required by $M_r$.
These first-variable coefficients are obtained by expanding
$1/(t-s)=-\sum_{j\ge0}t^js^{-j-1}$ at $t=0$.
The corresponding expansion in the second variable is
$s t^{r-1}/(t-s)=\sum_{j\ge1}s^jt^{r-1-j}$ at $s=0$.
After extracting $\dd s/s^r$, the constant coefficient is $-1$
and its coefficient of $s^{r-1}$ is $1$, as required by $N_r$.
The behavior at both points at infinity is expressed in the coordinates
$w=1/t$ and $v=1/s$, in which
\[
 S_r=w^{-(r-2)}v^{r-2}
       \left[\frac1{w-v}+w^{r-2}\right]\dd v.
\]
The factors $w^{-(r-2)}$ and $v^{r-2}\dd v$ are the local
frames of $F_r$ and $F_r^D$ at infinity. The remaining expression
is regular off the smooth diagonal and has its prescribed simple
pole there.

At the singular pair, extension of the kernel is determined by
multiplication by the ideal of the diagonal. On the affine product
this ideal is generated by
$t^m-s^m$, $r\le m\le2r-1$. The elementary identity
\[
 (t^m-s^m)\left[\frac{s t^{r-1}}{t-s}-1\right]
 =\sum_{j=0}^{m-1}t^{m+r-2-j}s^{j+1}-t^m+s^m
\]
shows that multiplication of \eqref{res:bf:geometric-kernel} by
such a generator gives a polynomial times $\dd s/s^r$. Its coefficients of $t^0$ and $t^{r-1}$ are
both $s^m$, and its intervening $t$-coefficients vanish. Its
coefficients of $s^0$ and $s^{r-1}$ are $-t^m$ and $t^m$.
Consequently it lies in $M_r\otimes_{\CC}N_r\dd s/s^r$.
The acyclic residue construction therefore identifies
\eqref{res:bf:geometric-kernel} with the normalized global section
of $\operatorname{Hom}(I_\Delta,F_r\boxtimes F_r^D)$.

\subsection{A scalar pairing with invertible Gram matrices at every stated bound}

The scalar pairing used in the diagonal matrix blocks is obtained
from the following integral functional on $K_r$:
\begin{equation}\label{res:bf:gamma}
 \Gamma[p]\coloneq\int_{1/4}^{1/2}p(t)t^{-r}\,\dd t.
\end{equation}
It pairs $f\in M_r$ with $g\dd t/t^r\in F_r^D$ by
$\Gamma[fg]$. The different bases in
\eqref{res:bf:left-basis}--\eqref{res:bf:right-basis} need not produce a
symmetric Gram matrix. Its invertibility follows instead from the
signs of the two evaluation determinants and their integral relation
to the Gram determinant.

\begin{lemma}[Invertibility of the scalar Gram matrices]\label{res:bf:scalar-normality}
Let $N\ge r$ and let $w$ be a continuous positive function on
$[1/4,1/2]$. Then
\[
 \det\left[\int_{1/4}^{1/2}
       f_i(t)g_j(t)w(t)t^{-r}\,\dd t\right]_{i,j\in\{1,\ldots,N\}}>0,
\]
where $f_i,g_j$ are the lists
\eqref{res:bf:left-basis}--\eqref{res:bf:right-basis}. For polynomial
weights in $A_r$ with rational coefficients, all these entries
are rational numbers.
\end{lemma}

\begin{proof}
Take $1/4<t_1<\cdots<t_N<1/2$ and put
$V(t_1,\ldots,t_N)\coloneq\prod_{i<j}(t_j-t_i)>0$.
The determinant $\det[f_i(t_j)]$ is the sum of the two
alternants whose exponent lists are
\[
 (0,r,r+1,\ldots,N+r-2),\qquad
 (r-1,r,r+1,\ldots,N+r-2).
\]
Both lists are strictly increasing. For exponents
$0\le a_1<\cdots<a_N$, the alternant formula gives
$\det[t_j^{a_i}]=V s_\lambda(t)$, with
$\lambda_i=a_{N-i+1}-(N-i)$. Here $\lambda$ is a decreasing
list of nonnegative integers and $s_\lambda$ is its Schur
polynomial. Its tableau expansion is a nonempty sum of monomials
with nonnegative integer coefficients. Consequently it is strictly
positive when every $t_j>0$. Applying this formula to each of the
two displayed exponent lists proves positivity of the left
evaluation determinant.

For the right determinant, the constant term in its first row
gives the exponent list $0,1,\ldots,r-2,r,\ldots,N$.
The corresponding alternant is
$V e_{N-r+1}(t_1,\ldots,t_N)$, where $e_j$ denotes the
elementary symmetric polynomial of degree $j$. This identity
follows by expanding $V(t_1,\ldots,t_N,z)=V(t)\prod_j(z-t_j)$
along its last column and comparing coefficients of $z^{r-1}$.
To put the other
term, $-t^{r-1}$, in increasing exponent order requires $r-2$
row interchanges. Its contribution is therefore
$(-1)^{r-1}V\prod_jt_j$. The resulting decomposition is
\begin{align}
 \det[g_i(t_j)]
 &=V\left[e_{N-r+1}(t)+(-1)^{r-1}\prod_jt_j\right]
       \notag\\
 &=V\prod_jt_j
       \left[e_{r-1}(t_1^{-1},\ldots,t_N^{-1})
                     +(-1)^{r-1}\right]>0.
 \label{res:bf:right-alternant}
\end{align}
The last inequality follows because each $t_j^{-1}>2$ and
$e_{r-1}(t_1^{-1},\ldots,t_N^{-1})>1$.

The signs of the evaluation determinants pass to the Gram determinant
through an integral identity. With $\dd\sigma(t)\coloneq w(t)t^{-r}\dd t$,
expansion of $\det[\int f_i g_j\,\dd\sigma]$ by permutations
expresses each product of $N$ integrals as an integral on the product
interval. Antisymmetrization of the $f$ factors then gives
\[
 \det G_w=\frac1{N!}\int_{[1/4,1/2]^N}
       \det[f_i(t_j)]\det[g_i(t_j)]
                      \prod_{j=1}^N\dd\sigma(t_j).
\]
A permutation of the $t_j$ changes the sign of each determinant,
so their product is invariant. The diagonals have measure zero,
and the $N!$ ordered regions have equal integrals. Canceling the
factor $1/N!$ gives
\[
 \det G_w=
 \int_{1/4<t_1<\cdots<t_N<1/2}
   \det[f_i(t_j)]\det[g_i(t_j)]
   \prod_{j=1}^Nw(t_j)t_j^{-r}\,\dd t_1\cdots\dd t_N.
\]
Every factor is positive in the interior of this simplex, so the
integral is strictly positive.

The rationality of the moments follows from the missing coefficient
in $f_ig_j\in K_r$: its coefficient of $t^{r-1}$ is zero, and
multiplication by a polynomial in $A_r$ preserves this property.
Division by $t^r$ therefore produces no $t^{-1}$ term, so every
monomial integral has the form
\[
 \int_{1/4}^{1/2}t^{m-r}\,\dd t
 =\frac{(1/2)^{m-r+1}-(1/4)^{m-r+1}}{m-r+1},
 \qquad m\ne r-1,
\]
which is rational. This proves the assertion about moments.
\end{proof}

\subsection{Two residue functionals and their perfect quotients}

For the residue masses, the coordinates $x\coloneq t^r$,
$y\coloneq t^{r+1}$ give the denominator family
\[
 h_u\coloneq x^2+u=t^{2r}+u.
\]
For a polynomial $p\in K_r$, define
\begin{align}
 \rho_{1,u}[p]&\coloneq
   \sum_{j\ge0}(-u)^j[t^{3r-1+2rj}]p,
       \label{res:bf:rho1}\\
 \rho_{2,u}[p]&\coloneq
   \sum_{j\ge0}(j+1)(-u)^j[t^{5r-1+2rj}]p.
       \label{res:bf:rho2}
\end{align}
Only finitely many terms occur for each test. These are the sums
of the finite residues of
$p\dd t/(t^rh_u)$ and $p\dd t/(t^rh_u^2)$, respectively.
The coefficient formulas also define them at $u=0$, when all the
zeros of $h_u$ meet the singular point. Their behavior under
multiplication by the denominator is expressed by
\begin{equation}\label{res:bf:rho-identities}
 \rho_{2,u}[h_up]=\rho_{1,u}[p],\qquad
 \rho_{1,u}[h_up]=0,\qquad
 \rho_{2,u}[h_u^2p]=0.
\end{equation}

\begin{lemma}[Perfect residue pairings]\label{res:bf:quotient-duality}
For every $u\in\CC$, the two multiplication pairings induced by
$\rho_{1,u}$ and $\rho_{2,u}$ are perfect, meaning nondegenerate in both arguments, on
\begin{equation}\label{res:bf:quotient-pairs}
 \frac{M_r}{h_uM_r}\times\frac{N_r}{h_uN_r},\qquad
 \frac{M_r}{h_u^2M_r}\times\frac{N_r}{h_u^2N_r}.
\end{equation}
The dimensions of each pair of arguments are $2r$ and $4r$,
respectively. On polynomial representatives with rational
coefficients, their Gram entries are polynomials in $u$ with
rational coefficients.
\end{lemma}

\begin{proof}
Set $v=t^{2r}$ and $B=\CC[v]$. The map
\[
 \Psi[p]\coloneq
    \sum_{j\ge0}[t^{3r-1+2rj}]p\,v^j,
 \qquad p\in K_r,
\]
is $B$-linear. In fact,
$\Psi[vp]=v\Psi[p]+[t^{r-1}]p=v\Psi[p]$.
The missing coefficient is thus exactly the condition needed for
linearity over this polynomial ring; it is supplied by the canonical
module used in the pairing.

Bases adapted to this polynomial-valued pairing are indexed by the
following set of $2r$ exponents:
\[
 E_r\coloneq\{0\}\cup\{r,\ldots,2r-1\}
                       \cup\{2r+1,\ldots,3r-1\}.
\]
For $e\in E_r$, let $m_e=t^e$, except that $m_0=1+t^{r-1}$.
Let $n_e=t^{3r-1-e}$, except that
$n_{3r-1}=1-t^{r-1}$. These are $B$-bases of $M_r$ and $N_r$.
Their spanning property includes the polynomial $t^{2r}$ absent
from the monomial lists, since
\[
 v(1+t^{r-1})=t^{2r}+t^{3r-1},\qquad
 v(1-t^{r-1})=t^{2r}-t^{3r-1}.
\]
The other monomials are generated by powers of $v$. For
independence, the constant residue class modulo $2r$ first forces
the coefficient of the modified generator to vanish; the remaining
exponents occupy distinct residue classes.

These two bases are dual for $\Psi$:
\begin{equation}\label{res:bf:B-duality}
 \Psi[m_e n_{e'}]=\delta_{e,e'}.
\end{equation}
The set $E_r$ has exactly one representative of each residue class
modulo $2r$: its three parts give the classes $0$,
$r,\ldots,2r-1$, and $1,\ldots,r-1$, respectively. The map
$\Psi$ retains the exponents $3r-1+2rj$ with $j\ge0$.
The products can therefore be checked in four cases.

If $e\ne0$ and $e'\ne3r-1$, both factors are monomials and
\[
 m_en_{e'}=t^{3r-1+e-e'}.
\]
A retained exponent requires $e-e'=2rj$. Distinct elements of
$E_r$ have distinct residue classes, so this forces $e=e'$ and
$j=0$. Its contribution to $\Psi$ is $1$.

If $e=0$ and $e'\ne3r-1$, then
\[
 m_0n_{e'}=t^{3r-1-e'}+t^{4r-2-e'}.
\]
The first exponent is retained precisely when $e'=0$; the equation
$-e'=2rj$ has no other solution with $e'\ge0$ and $j\ge0$.
For the second exponent, the required equation is
$e'=r-1-2rj$. When $j=0$, the value $r-1$ is absent from
$E_r$, and for $j\ge1$ the value is negative. Thus this product
contributes $\delta_{0,e'}$.

If $e\ne0$ and $e'=3r-1$, then
\[
 m_en_{3r-1}=t^e-t^{e+r-1}.
\]
Since $e\le3r-1$, the first term contributes only at $e=3r-1$,
with value $1$. The second exponent would require
$e=2r+2rj$. The value $2r$ is absent from $E_r$, and all
values with $j\ge1$ exceed its largest element. This case
therefore contributes $\delta_{e,3r-1}$.

Finally, if $e=0$ and $e'=3r-1$, both generators are modified:
\[
 m_0n_{3r-1}=(1+t^{r-1})(1-t^{r-1})=1-t^{2r-2}.
\]
Both exponents are smaller than $3r-1$, so the contribution is
zero, as required by $\delta_{0,3r-1}=0$. These four computations
prove \eqref{res:bf:B-duality} for every pair of basis elements.

Evaluation and differentiation of this polynomial-valued pairing
recover the two residue functionals:
\[
 \rho_{1,u}[p]=\Psi[p](-u),\qquad
 \rho_{2,u}[p]=\left.\frac{\dd}{\dd v}\Psi[p](v)
                                         \right|_{v=-u}.
\]
In particular, for $b,c\in B$, $B$-linearity and
\eqref{res:bf:B-duality} give
$\rho_{1,u}[b m_e c n_{e'}]=\delta_{e,e'}b(-u)c(-u)$.
Thus the $2r\times2r$ Gram on the quotient by $h_u=v+u$
is the identity in these bases. For the quotient by $h_u^2$,
write $b=b_0+b_1(v+u)$ and $c=c_0+c_1(v+u)$. Then
$(bc)'(-u)=b_0c_1+b_1c_0$. The scalar Gram in the basis
$1,v+u$ is
\[
 \begin{bmatrix}0&1\\1&0\end{bmatrix}
\]
and has determinant $-1$. The full $4r\times4r$ Gram is a
direct sum of $2r$ copies of this matrix. This proves both
nondegeneracy assertions and both dimensions without requiring
that the roots of $h_u$ be distinct. The same evaluation and
derivative formulas show that
$((v+u)b)'|_{v=-u}=b(-u)$, $((v+u)b)|_{v=-u}=0$ and
$((v+u)^2b)'|_{v=-u}=0$ give
\eqref{res:bf:rho-identities}. Polynomial dependence on $u$ follows
directly from \eqref{res:bf:rho1}--\eqref{res:bf:rho2}.
\end{proof}

The perfect pairings in \eqref{res:bf:quotient-pairs} relate the
section module to its canonical dual. Their persistence at $u=0$
is therefore compatible with the non-Gorenstein nature of $A_r$
for $r\ge3$: the second argument retains the distinct canonical
module, rather than replacing it by a second copy of $A_r/(h_u)$.
Perfection is established for this dual pair of modules.

\subsection{A fixed coefficient giving normality for every real parameter}

The matrix construction below reduces normality to the nonvanishing
of determinants that are polynomial in the parameter and in the mass
coefficient. A consequence of algebraic independence supplies one
coefficient for which all such determinants remain nonzero. This fixes
the coupling for the whole family before any finite bound is chosen.

\begin{lemma}[A uniform coefficient]\label{res:bf:uniform-coupling}
Let $P(u,z)\in\mathbb Q[u,z]$ and suppose that $P(u,0)\ne0$
for every $u\in\mathbb R$. Then
\[
 P(u,\kappa)\ne0\quad\hbox{for every }u\in\mathbb R,
 \qquad
 \kappa\coloneq\mathrm{e}^{-1}
                   +\mathrm{i}\,\mathrm{e}^{-\sqrt2}.
\]
The same $\kappa$ works for any collection of polynomials satisfying
these hypotheses.
\end{lemma}

\begin{proof}
Write $a=\mathrm{e}^{-1}$ and $b=\mathrm{e}^{-\sqrt2}$. These
two real numbers are algebraically independent over $\mathbb Q$.
Indeed, a polynomial relation would be a linear relation over
$\mathbb Q$ among the numbers
$\exp(-j-k\sqrt2)$ for distinct pairs of nonnegative integers
$(j,k)$. Their exponents are distinct algebraic numbers.
The Lindemann--Weierstrass theorem states that such exponentials
are linearly independent over the algebraic numbers; see
\cite[Chapter~1, Section~3]{BakerTranscendental} and
\cite{BezivinRobba}.

Suppose now that $P(u,a+\mathrm{i} b)=0$ for a real $u$.
The specialization $P(u,z)$ is a nonzero polynomial, because
$P(u,0)\ne0$. It follows that $a+\mathrm{i} b$ is algebraic
over $\mathbb Q(u)$. Its conjugate is a root of the same
polynomial, since the coefficients are real. Therefore $a$ and
$b$ both belong to an algebraic extension of $\mathbb Q(u)$.
Adjoining algebraic elements does not increase transcendence
degree, whereas a field generated by the single real number $u$
has transcendence degree at most one over $\mathbb Q$.
The field $\mathbb Q(a,b)$ would consequently have transcendence
degree at most one. Their algebraic independence gives degree two,
a contradiction. Since this argument is independent of the degree
of $P$ and requires no choice of exceptional values of $u$, the
single coefficient $\kappa$ has the stated simultaneous property.
\end{proof}

Almost every complex coupling has the same uniform property.
More precisely, there is one set $\mathscr E\subset\CC$ of planar
Lebesgue measure zero such that, for every $\kappa\notin\mathscr E$,
\[
 P(u,\kappa)\ne0\qquad(u\in\mathbb R)
\]
simultaneously for every $P\in\mathbb Q[u,z]$ satisfying
$P(u,0)\ne0$ on the real line. To prove this assertion, define
\[
 \mathscr E\coloneq
 \bigcup_{0\ne Q\in\mathbb Q[X,Y]}
       \{a+\mathrm{i}b:a,b\in\mathbb R,\ Q(a,b)=0\}.
\]
Suppose that $P(u_0,\kappa)=0$ for some real $u_0$, and write
$\kappa=a+\mathrm{i}b$. The condition $P(u_0,0)\ne0$ makes
$P(u_0,z)$ a nonzero polynomial. Hence $\kappa$ is algebraic
over $\mathbb Q(u_0)$, and so is its conjugate, since this
polynomial has real coefficients. The equalities
\[
 a=\frac{\kappa+\overline\kappa}{2},\qquad
 b=\frac{\kappa-\overline\kappa}{2\mathrm{i}}
\]
show that both numbers are algebraic over
$\mathbb Q(u_0,\mathrm{i})$. Consequently
$\operatorname{trdeg}_{\mathbb Q}\mathbb Q(a,b)\le1$.
They satisfy a nonzero polynomial relation over $\mathbb Q$,
so $\kappa\in\mathscr E$. This proves the assertion for every
coupling outside the single set $\mathscr E$.

There are countably many rational polynomials $Q$. Each nonzero
one has a zero set of planar measure zero: write
$Q(X,Y)=\sum_{j=0}^m q_j(X)Y^j$ and choose one nonzero coefficient
$q_j$. Except at its finitely many real roots, a vertical section
of the zero set consists of finitely many roots of the nonzero
polynomial $Q(a,Y)$. The exceptional vertical lines also have
planar measure zero. Fubini's theorem on each bounded rectangle,
followed by a countable union of rectangles, therefore gives
measure zero for the whole zero set. The countable union defining
$\mathscr E$ has measure zero as well. This argument controls
all real parameters together. The explicit coefficient in
Lemma~\ref{res:bf:uniform-coupling} belongs to its complement by
the algebraic independence already proved.

\subsection{The four factors and the normal matrix pairings}

For $\varepsilon\in\mathbb R$, set
\[
 u=\varepsilon^2,\qquad
 \ell_\varepsilon\coloneq x-\mathrm{i}\varepsilon,
 \qquad a_\varepsilon\coloneq x+\mathrm{i}\varepsilon,
 \qquad h\coloneq x^2+\varepsilon^2.
\]
The four factors are
\begin{equation}\label{res:bf:four-factors}
 L_1=\begin{bmatrix}\ell_\varepsilon&y\\0&\ell_\varepsilon\end{bmatrix},
 \quad
 R_1=\begin{bmatrix}a_\varepsilon&y\\0&2a_\varepsilon\end{bmatrix},
 \quad
 L_0=\begin{bmatrix}1&0\\0&x\end{bmatrix},
 \quad
 R_0=\begin{bmatrix}x&0\\0&1\end{bmatrix}.
\end{equation}
Their determinants are respectively $\ell_\varepsilon^2$,
$2a_\varepsilon^2$, $x$ and $x$, all nonzero in the domain
$A_r$. The order of the two denominator matrices matters:
\begin{align}
 T\coloneq R_1L_1
     &=\begin{bmatrix}h&2xy\\0&2h\end{bmatrix},
 &L_1R_1
     &=\begin{bmatrix}h&(3x+\mathrm{i}\varepsilon)y\\0&2h\end{bmatrix},
       \label{res:bf:ordered-products}\\
 \operatorname{adj}(T)
     &=\operatorname{adj}(L_1)\operatorname{adj}(R_1)
       =\begin{bmatrix}2h&-2xy\\0&h\end{bmatrix},
 &\det T&=2h^2.
       \label{res:bf:adjugate}
\end{align}

The integral and residue parts have distinct supports. The functional
$\Gamma$ integrates over a real interval, while for real
$\varepsilon\ne0$ the zeros of $h=t^{2r}+\varepsilon^2$ are complex
points of the curve outside that interval. These zeros support the
residue masses. The resulting complex matrix pairings are normal
on the stated finite spaces, as proved below.

For a matrix $W$ with entries in $A_r$, the notation $\Gamma W$ means the
matrix of functionals with entries
$(\Gamma W)_{a,b}[p]=\Gamma[W_{a,b}p]$. With this convention,
the original pairing, the particular solution and the mass correction
are specified by
\begin{align}
 \mu_\varepsilon
   &\coloneq\Gamma
        \begin{bmatrix}h&x(3x+\mathrm{i}\varepsilon)y\\0&2h\end{bmatrix},
 &\mu_{\mathrm{part}}&\coloneq x\Gamma I_2,
       \label{res:bf:original-particular}\\
 \nu_u
   &\coloneq\operatorname{adj}(T)\rho_{2,u}
       =\begin{bmatrix}
            2\rho_{1,u}&-2xy\rho_{2,u}\\
            0&\rho_{1,u}
         \end{bmatrix},
 &\widehat\mu_\varepsilon
   &\coloneq\mu_{\mathrm{part}}+\kappa\nu_u.
       \label{res:bf:mass-target}
\end{align}
For example, $(xy\rho_{2,u})[p]=\rho_{2,u}[xyp]$.
For a row $\boldsymbol f\in M_r^{1\times2}$ and a column
$\boldsymbol g\in N_r^{2\times1}$, a matrix of functionals
defines the pairing
$\beta_\mu(\boldsymbol f,\boldsymbol g)
=\sum_{a,b}\mu_{a,b}[f_ag_b]$. The differential in the second
argument is understood. With this convention,
\[
 \beta_{L\mu R}(\boldsymbol f,\boldsymbol g)
       =\beta_\mu(\boldsymbol fL,R\boldsymbol g).
\]

\begin{theorem}[A normal bilateral family]\label{res:bf:family}
For every integer $r\ge2$ and every $\varepsilon\in\mathbb R$,
the functionals \eqref{res:bf:original-particular}--\eqref{res:bf:mass-target}
satisfy
\begin{equation}\label{res:bf:bilateral-equation}
 L_1\widehat\mu_\varepsilon R_1=L_0\mu_\varepsilon R_0.
\end{equation}
Both pairings are normal on
$U_N^{1\times2}\times V_N^{2\times1}$ for every $N\ge r$,
using the same scalar bound in each component. The homogeneous mass
$\nu_u$ has global rank $4r$. Its restriction to these finite spaces has
rank $4r$ if and only if $N\ge2r$, for every value of the
parameter, including $\varepsilon=0$.
\end{theorem}

\begin{proof}
The functional equation follows from the particular solution and
the homogeneous mass separately. The polynomial densities in \eqref{res:bf:original-particular} give
\[
 L_1(xI_2)R_1
 =x\begin{bmatrix}h&(3x+\mathrm{i}\varepsilon)y\\0&2h\end{bmatrix}
 =L_0\begin{bmatrix}h&x(3x+\mathrm{i}\varepsilon)y\\0&2h\end{bmatrix}R_0.
\]
Thus $\mu_{\mathrm{part}}$ is a particular solution. For the
mass, preserving the adjugate order gives
\begin{align*}
 L_1\nu_uR_1
 &=L_1\operatorname{adj}(L_1)
       \operatorname{adj}(R_1)R_1\rho_{2,u}\\
 &=2h^2\rho_{2,u}I_2=0
\end{align*}
by \eqref{res:bf:rho-identities}. This proves
\eqref{res:bf:bilateral-equation}.

When all $N$ basis vectors of the first component precede those of
the second, the Gram matrices inherit the triangular form of the
matrix pairings. Their blocks are expressed through
\[
 \begin{aligned}
 G_w&\coloneq[\Gamma[w f_i g_j]]_{i,j\in\{1,\ldots,N\}},\\
 H_1(u)&\coloneq[\rho_{1,u}[f_i g_j]]_{i,j\in\{1,\ldots,N\}},\\
 H_{12}(u)&\coloneq[-2\rho_{2,u}[xyf_i g_j]]_{i,j\in\{1,\ldots,N\}}.
 \end{aligned}
\]
The original Gram has the block form
\[
 G_{\mathrm{orig}}
   =\begin{bmatrix}G_h&G_{x(3x+\mathrm{i}\varepsilon)y}\\0&2G_h\end{bmatrix},
 \qquad
 \det G_{\mathrm{orig}}=2^N(\det G_h)^2.
\]
On the integration interval, $h(t)>0$ for every real
$\varepsilon$. Lemma~\ref{res:bf:scalar-normality} therefore proves
that this determinant is nonzero. The determinant depends only on the diagonal blocks, so positivity
of the scalar pairing is sufficient even when the matrix density
is not positive.

The same ordering gives the target Gram in triangular form:
\begin{equation}\label{res:bf:target-gram}
 \widehat G_N(u)=
 \begin{bmatrix}
   G_x+2\kappa H_1(u)&\kappa H_{12}(u)\\
   0&G_x+\kappa H_1(u)
 \end{bmatrix}.
\end{equation}
The entries of $G_x$ are rational and $\det G_x>0$, because
$x(t)>0$ on the same interval. The entries of $H_1(u)$ are
polynomials in $u$ with rational coefficients. The two determinant
polynomials relevant to Lemma~\ref{res:bf:uniform-coupling} are
\[
 \det(G_x+2zH_1(u)),\qquad\det(G_x+zH_1(u)).
\]
They belong to $\mathbb Q[u,z]$, because each determinant is
a finite signed sum of products of the specified entries.
Their constant term in $z$ is $\det G_x$, independent of $u$
and nonzero. Both diagonal blocks in \eqref{res:bf:target-gram} are
therefore invertible for every real $u$, simultaneously for all
$r$ and $N$. This proves target normality for $u=\varepsilon^2$.

The global mass rank is determined by the perfect pairing $\rho_{2,u}$
from Lemma~\ref{res:bf:quotient-duality}, applied in each component.
Composition in its second argument with $\operatorname{adj}(T)$
gives the mass. Perfection therefore characterizes its right radical,
for $v\in N_r^2$, by
\[
 \operatorname{adj}(T)v\in h^2N_r^2
       \quad\text{if and only if}\quad v\in TN_r^2.
\]
One implication follows from
$\operatorname{adj}(T)T=2h^2I_2$. Conversely, if
$\operatorname{adj}(T)v=h^2w$, multiplication by $T$ gives
$2h^2v=h^2Tw$. The module $N_r$ is torsion-free and $h\ne0$,
so cancellation gives $v=T(w/2)$. Hence
\begin{equation}\label{res:bf:rank-cokernel}
 \rank\nu_u=\dim_{\CC}(N_r^2/TN_r^2).
\end{equation}
The triangular matrix in \eqref{res:bf:ordered-products} gives an
exact sequence
\[
 0\longrightarrow N_r/hN_r
   \longrightarrow N_r^2/TN_r^2
   \longrightarrow N_r/hN_r\longrightarrow0.
\]
The maps send $[a]$ to $[[a\ 0]^{\T}]$ and $[[a\ b]^{\T}]$
to $[b]$. If $[a\ 0]^{\T}=T[c\ d]^{\T}$, then $2hd=0$,
so $d=0$ and $a=hc$; this proves injectivity. If $b$ is a
multiple of $h$, subtracting a suitable $T[0\ d]^{\T}$ leaves
a representative in the first component, proving exactness in
the middle. Both end terms have dimension $2r$, so
\eqref{res:bf:rank-cokernel} gives $\rank\nu_u=4r$.

On the finite section spaces, this global rank is attained at
$N=2r$, where the lists
\eqref{res:bf:left-basis}--\eqref{res:bf:right-basis} represent bases
of $M_r/hM_r$ and $N_r/hN_r$. Relative to the dual bases used
in Lemma~\ref{res:bf:quotient-duality}, their only replacement is
the high monomial $t^{3r-1}$ by $t^{2r}$. The relations are
\[
 t^{2r}=-u(1+t^{r-1})-t^{3r-1}\pmod{hM_r},\qquad
 t^{2r}=t^{3r-1}-u(1-t^{r-1})\pmod{hN_r}.
\]
In each replacement the coefficient of the removed monomial is
$-1$ or $1$; the other term is a multiple of the retained first
basis vector. Therefore each change of basis has determinant
$+1$ or $-1$, independently of $u$. The quotient pairing had
identity Gram in the dual $B$-bases, so its Gram in the displayed
section bases satisfies $\det H_1(u)=\pm1$ at this bound, and the
matrix of the full mass has determinant
\begin{equation}\label{res:bf:constant-minor}
 \det\begin{bmatrix}2H_1(u)&H_{12}(u)\\0&H_1(u)\end{bmatrix}
     =2^{2r}(\det H_1(u))^2=2^{2r}.
\end{equation}
Every larger section space contains this minor. Together with the global
upper bound $4r$, this proves the claimed rank for $N\ge2r$.
If $N<2r$, the mass Gram has size $2N<4r$ and cannot attain
that rank.
\end{proof}

The target Gram $\widehat G_N(u)$ remains nonsingular simultaneously
for every $r\ge2$, $N\ge r$ and real $u$ for any coupling
$\kappa\notin\mathscr E$.
Indeed, each diagonal determinant is one of the rational-coefficient
polynomials covered above, with nonzero constant term $\det G_x$.
The set $\mathscr E$ is independent of these indices and of $u$,
so almost every complex coupling works for the entire collection of
target Grams. At $u=\varepsilon^2$ this gives target normality for
every real parameter of the bilateral family.
Such a coupling is nonzero, because $0\in\mathscr E$; multiplication
of the mass by it preserves the rank $4r$ on the stated cuts.
The fixed explicit coefficient will continue to be used in the
formulas below.

Theorem~\ref{res:bq:matrixquotient} gives a geometric description
of this particular mass through the matrix $T=R_1L_1$ in
\eqref{res:bf:ordered-products}. This $2\times2$ product defines
the quotient relevant to the chosen mass, whereas the $4\times4$
vectorized denominator operator $R_1\otimes L_1^{\T}$ acts on
all matrix tests. For the mass constructed here, the residue formula is
\begin{equation}\label{res:bf:geometric-mass}
 \nu_u(\boldsymbol f,\boldsymbol g)
 =2\sum_{a\in\{h=0\}}\operatorname{Res}_a
       \left(\boldsymbol f T^{-1}\boldsymbol g\,
                         \frac{\dd t}{t^r}\right)
 =2b_{T^{\T}}(\boldsymbol f^{\T},\boldsymbol g).
\end{equation}
At a singular support point the residues are interpreted through
the canonical module, or equivalently by
\eqref{res:bf:rho1}--\eqref{res:bf:rho2}. The transpose in the last
expression is necessary because the first argument of
$b_{T^{\T}}$ is a column. Consequently the two quotients
are $M_r^2/T^{\T}M_r^2$ and $N_r^2/TN_r^2$, and the tensor
of dual basis elements (the coevaluation tensor) of $\nu_u$ is
\begin{equation}\label{res:bf:mass-coevaluation}
 \frac12\mathcal C_{T^{\T}}(t,s)
 =\frac12\bigl(T^{\T}(t)-T^{\T}(s)\bigr)S_r(t,s)
\end{equation}
in their tensor product. Both the transpose and the factor $1/2$
follow from the pairing convention: $\det T=2h^2$ and
$\nu_u=2b_{T^{\T}}$. Formula
\eqref{res:bf:mass-coevaluation} links the global kernel with the
finite mass, including the confluence value $u=0$.

\subsection{The finite kernel correction and the noncommuting factors}

The constant full-rank minor obtained above gives a polynomial
factorization of the mass Gram that persists through confluence.
This factorization reduces the kernel correction to a matrix whose
size is the global mass rank, independently of the section bound.

\begin{corollary}[A uniform correction of size $4r$]\label{res:bf:woodbury}
Let $N\ge2r$. Relative to the particular pairing
$\mu_{\mathrm{part}}=x\Gamma I_2$, the reproducing kernel of
$\widehat\mu_\varepsilon$ has a correction of rank exactly
$4r$. It can be computed by inverting a $4r\times4r$ matrix
that is nonsingular for every real $\varepsilon$. Its
factorization can be chosen polynomially in $u=\varepsilon^2$.
\end{corollary}

\begin{proof}
Let $H_N(u)$ be the $2N\times2N$ mass Gram and put
$G_{\mathrm{part}}\coloneq\diag(G_x,G_x)$. Choose row and
column index sets $I,J$ corresponding to the bound $N=2r$ in
the two components. By \eqref{res:bf:constant-minor},
$H_*(u)\coloneq H_N(u)[I,J]$ has constant determinant
$2^{2r}$. Let $E_I,E_J$ have as their columns the standard
coordinate vectors indexed by $I,J$, respectively, so that
$H_*=E_I^{\T}H_NE_J$. Define
\[
 P_N(u)\coloneq H_N(u)E_J,\qquad
 Q_N(u)\coloneq H_*(u)^{-1}E_I^{\T}H_N(u).
\]
Both matrices have polynomial entries, since the denominator of
$H_*^{-1}$ is constant. The identity $H_N=P_NQ_N$ follows by
moving the selected rows and columns to the front, where the block
decomposition is $H_N=[H_*\ A;\ B\ D]$. The first $4r$ columns are independent,
and every remaining column is their linear combination because
$\rank H_N=4r$. Its first $4r$ entries force the coefficients
to be $H_*^{-1}A$, so $D=BH_*^{-1}A$. Thus
\[
 H_N=\begin{bmatrix}H_*\\ B\end{bmatrix}
             H_*^{-1}\begin{bmatrix}H_*&A\end{bmatrix}
\]
in the permuted coordinates. Returning to the original order
gives exactly $P_NQ_N$. The equality first holds over $\CC(u)$;
since all entries are polynomials, it is a polynomial identity. Thus the
factorization is fixed uniformly through the confluence.

The matrix governing the inverse correction is
\[
 Z_N(u)\coloneq I_{4r}
           +\kappa Q_N(u)G_{\mathrm{part}}^{-1}P_N(u).
\]
The determinant identity
\[
 \det Z_N(u)=\frac{\det\widehat G_N(u)}{\det G_{\mathrm{part}}}
\]
shows that $Z_N(u)$ is invertible for every real $u$.
The relation to the full inverse follows by multiplication. For
$G=G_{\mathrm{part}}$, $U=\kappa P_N$ and $V=Q_N$, the defining
equality $Z_N=I+VG^{-1}U$ gives
$Z_N^{-1}+Z_N^{-1}VG^{-1}U=I$, which cancels the additional terms
in $(G+UV)(G^{-1}-G^{-1}UZ_N^{-1}VG^{-1})$.
The resulting inverse identity is the Woodbury formula
\begin{equation}\label{res:bf:woodbury-inverse}
 \widehat G_N(u)^{-1}=G_{\mathrm{part}}^{-1}
 -\kappa G_{\mathrm{part}}^{-1}P_N(u)Z_N(u)^{-1}
                         Q_N(u)G_{\mathrm{part}}^{-1}.
\end{equation}
This factorization introduces no parameter poles in $P_N,Q_N$.

To pass from the Gram inverse to the reproducing kernel, let
$\mathsf F_N(Q)$ be the $2N\times2$
matrix whose rows are the evaluations of the row-section basis,
and let $\mathsf G_N(P)$ be the $2\times2N$ matrix whose
columns are the evaluations of the column-section basis. In the
column--row convention of the manuscript, the reproducing kernel is
\[
 \widehat K_N(P,Q)
   =\mathsf G_N(P)\widehat G_N(u)^{-1}\mathsf F_N(Q).
\]
Substitution of \eqref{res:bf:woodbury-inverse} gives
\begin{align}
 \widehat K_N(P,Q)-K_{\mathrm{part},N}(P,Q)
 ={}&-\kappa
   [\mathsf G_N(P)G_{\mathrm{part}}^{-1}P_N(u)]
       Z_N(u)^{-1}\notag\\
 &\hspace{1.5em}\cdot
   [Q_N(u)G_{\mathrm{part}}^{-1}\mathsf F_N(Q)].
 \label{res:bf:kernel-correction}
\end{align}
Proposition~\ref{res:bq:prop:finitekernel} gives separation rank
$\rank(\kappa H_N)=4r$. This compares the target kernel with the
particular kernel on the same section spaces. The original functional
$\mu_\varepsilon$ has a different Gram and is not the comparison
pairing in \eqref{res:bf:kernel-correction}.
\end{proof}

The same Gram inverse determines normalized polynomial sections.
Their leading classes are fixed by the passage to the next section
space, and orthogonality determines the remaining coefficients.
The construction requires normality of consecutive section spaces
and therefore applies from the initial bound of the family. The
additional full-rank condition is needed only for the fixed-size
factorization of the correction.

\begin{corollary}[Normalized polynomial sections through confluence]
\label{res:pf:curve-polynomials}
Fix $r\ge2$ and $N\ge r$. Use the lower spaces
$U_N^{1\times2}$ and $V_N^{2\times1}$, and the upper spaces with
$N$ replaced by $N+1$, in the bases
\eqref{res:bf:left-basis}--\eqref{res:bf:right-basis}.
There are unique row and column blocks for
$\widehat\mu_\varepsilon$ orthogonal to the respective lower
spaces and with prescribed leading classes
\begin{equation}
 \widehat{\mathcal P}_N(t)
       \equiv t^{N+r-1}I_2\pmod{U_N^{1\times2}},\qquad
 \widehat{\mathcal Q}_N(t)
       \equiv t^{N+1}I_2\pmod{V_N^{2\times1}}.
 \label{res:pf:eq:curve-leading}
\end{equation}
The congruences are interpreted row by row and column by column.
The rows are polynomial representatives in $M_r$, and the columns
represent sections in $N_r\dd t/t^r$, with the differential
understood. Their mutual $2\times2$ pairing matrix is invertible.

Let $\mathcal P_N,\mathcal Q_N$ be the blocks with these same
leading classes for the particular pairing $x\Gamma I_2$.
With the basis matrices used in \eqref{res:bf:kernel-correction},
the transformed blocks are explicitly
\begin{equation}
 \begin{aligned}
 \widehat{\mathcal P}_N
   &=\mathcal P_N-
     \kappa\nu_u(\mathcal P_N,\mathsf G_N)
                \widehat G_N(u)^{-1}\mathsf F_N,\\
 \widehat{\mathcal Q}_N
   &=\mathcal Q_N-
     \kappa\mathsf G_N\widehat G_N(u)^{-1}
                \nu_u(\mathsf F_N,\mathcal Q_N).
 \end{aligned}
 \label{res:pf:eq:curve-polynomials}
\end{equation}
Pairings of blocks are taken entry by entry. Their coefficients
are rational functions of $u=\varepsilon^2$, with no poles for
$u\ge0$. As functions of $\varepsilon$, the blocks are
holomorphic near every real parameter, and their limit at zero
is the normalized block computed directly from
$\widehat\mu_0$. If $N\ge2r$, their correction can use the same
$4r\times4r$ matrix $Z_N(u)$ as the kernel correction.
\end{corollary}

\begin{proof}
The new scalar basis elements at level $N+1$ are exactly
$t^{N+r-1}$ on the left and $t^{N+1}$ on the right. These
leading monomials determine the particular normalized blocks after
subtraction of their projections onto the lower spaces. For
$A_N(t)\coloneq t^{N+r-1}I_2$ and
$B_N(t)\coloneq t^{N+1}I_2$, with $\beta_{\mathrm{part}}$ denoting
the particular pairing, this gives
\[
 \begin{aligned}
 \mathcal P_N&=A_N-
    \beta_{\mathrm{part}}(A_N,\mathsf G_N)
                      G_{\mathrm{part}}^{-1}\mathsf F_N,\\
 \mathcal Q_N&=B_N-
    \mathsf G_NG_{\mathrm{part}}^{-1}
                       \beta_{\mathrm{part}}(\mathsf F_N,B_N).
 \end{aligned}
\]
The two moment blocks have sizes $2\times2N$ and $2N\times2$,
respectively, so both resulting blocks are $2\times2$. The leading
monomials specify fixed complements to the lower spaces; subtraction
of their particular-pairing projections constructs
$\mathcal P_N,\mathcal Q_N$ while preserving these leading classes.
The particular lower Gram is invertible by
Lemma~\ref{res:bf:scalar-normality}; Theorem~\ref{res:bf:family} gives
invertibility of the transformed Grams at both $N$ and $N+1$.
For the row correction, the target pairing of $\mathcal P_N$ with
the lower column basis is $\kappa\nu_u(\mathcal P_N,\mathsf G_N)$,
since its particular pairing with that basis vanishes.
A correction $C\mathsf F_N$ has target pairing $C\widehat G_N$.
The equation for orthogonality is therefore
$C\widehat G_N=\kappa\nu_u(\mathcal P_N,\mathsf G_N)$,
whose unique solution is the first correction in
\eqref{res:pf:eq:curve-polynomials}. The column equation is
$\widehat G_ND=\kappa\nu_u(\mathsf F_N,\mathcal Q_N)$ and
gives the second correction. This is the explicit specialization
of Proposition~\ref{res:pf:transformation} with $m=\kappa\nu_u$.
In the completed bases, these subtractions are invertible unit
triangular changes. The upper Gram becomes block diagonal with
the lower Gram and the mutual block pairing. Its invertibility proves that
the latter is nondegenerate.
The fixed bases at $N=r$ provide an initial block of $2r$ rows
and columns with invertible Gram $\widehat G_r$. Appending the
successive blocks above gives biorthogonal block families
exhausting the two section modules: each later block is orthogonal
to the lower spaces containing all earlier blocks. The families are
therefore constructed entirely from the section spaces $N\ge r$.

The bases, complements and leading classes are independent of the
parameter. The particular moments are constant, and the mass moments
are polynomial in $u$. Writing $\widehat G_N^{-1}=\operatorname{adj}(\widehat G_N)/
\det\widehat G_N$ makes each resulting coefficient a rational
function of $u$, with the fixed number $\kappa$ in its coefficient
field. Its denominator is nonzero for every real $u$ by the
normality proof, in particular for $u\ge0$. For each real
$\varepsilon_0$, the denominator composed with $u=\varepsilon^2$
is nonzero at $\varepsilon_0$ and hence in a complex neighborhood
of that point. The coefficients are holomorphic there.
At zero the formulas can consequently be evaluated by direct
substitution. Their leading classes and orthogonality remain the
same, so uniqueness identifies the limit with the normalized
blocks for $\widehat\mu_0$. This also gives the parameter conclusion
of Subsubsection~\ref{res:pf:subsection} in explicit coordinates. For
$N\ge2r$, \eqref{res:bf:woodbury-inverse} expresses the same computation
through $Z_N^{-1}$.
The analogous blocks for the original pairing are rational in
$\varepsilon$ and regular near the real axis; its density in
\eqref{res:bf:original-particular} also contains a term linear in
$\varepsilon$. The rational dependence on $u$ asserted above
concerns the particular and transformed blocks.
\end{proof}

The Gram calculation establishes normality and confluence, while
the factor identities describe the algebraic character of the
transformation. If $E_{1,2}$ denotes the elementary
matrix with its only nonzero entry in position $(1,2)$, then
\begin{equation}\label{res:bf:commutators}
 [L_0,L_1]=(1-x)yE_{1,2},\qquad
 [R_0,R_1]=(x-1)yE_{1,2},\qquad
 [L_1,R_1]=a_\varepsilon yE_{1,2}.
\end{equation}
Thus the order matters between the numerator and denominator
on each side, and also in the two matrix products displayed in
\eqref{res:bf:ordered-products}. The actions on separate arguments
of a pairing still commute as left and right actions; this is
distinct from the matrix commutators in
\eqref{res:bf:commutators}.

The entry ideals give a further obstruction to reducing the
denominators to scalar matrices by unimodular left and right
multiplications over $A_r$. Such
multiplications preserve the ideal generated by all matrix
entries. For $L_1$ this ideal is $(\ell_\varepsilon,y)$,
whereas the ideal of $\ell_\varepsilon I_2$ is
$(\ell_\varepsilon)$. At $\varepsilon=0$, their inequality
follows from $y/x=t\notin A_r$. At $\varepsilon\ne0$,
$y$ is not divisible by $\ell_\varepsilon$ even in
$\CC[t]$, since the roots of $\ell_\varepsilon$ are nonzero
and $y$ does not vanish there. Comparing determinants shows that
any scalar matrix equivalent to $L_1$ would have to be a
nonzero constant multiple of $\ell_\varepsilon I_2$:
the only units of $A_r$ are constants. This excludes every
scalar possibility. The same argument applies to $R_1$ with
$a_\varepsilon$. The numerators have determinant $x$, also
a nonunit, and are not polynomial changes of basis.

The associated rational changes of arguments make their relation
to dual changes of basis explicit:
\[
 S_L=L_1^{-1}L_0
   =\begin{bmatrix}
       1/\ell_\varepsilon&-xy/\ell_\varepsilon^2\\
       0&x/\ell_\varepsilon
     \end{bmatrix},\qquad
 S_R=R_0R_1^{-1}
   =\begin{bmatrix}
       x/a_\varepsilon&-xy/(2a_\varepsilon^2)\\
       0&1/(2a_\varepsilon)
     \end{bmatrix}.
\]
Their determinants satisfy
$\det S_L\det S_R=x^2/(2h^2)\ne1$, which shows that the two
rational transformations do not form inverse changes of dual bases.
The common triangular filtration preserved by all four matrices
provides the block structure that makes the diagonal Gram calculations
sufficient for normality. The entry-ideal obstruction preserves the
distinction between these denominators and scalar matrices over
$A_r$. The resulting normality statement applies to the complete
section spaces
\eqref{res:bf:left-basis}--\eqref{res:bf:right-basis}, repeated in the
two components with the same scalar bound; it does not extend
the assertion to arbitrary component degree bounds.

\section{Standard and Laurent specializations}\label{part:classical}

The general finite system acquires more explicit forms when the
coordinate algebra and the transformation factors are specified.
On the line, partial fractions describe the possible poles, while
Taylor coefficients determine whether the inverse polynomial matrix
product has a pole at any of these points. Once these finite conditions
hold, the coefficients at infinity impose the prescribed degree bounds.
Together, the two sets of conditions determine the constraint matrices
associated with the four factors.

The Christoffel and Geronimus cases simplify different parts of this
construction. Christoffel multiplication introduces neither denominator
poles nor homogeneous masses, so only the polynomial-image and
normalization equations remain. Geronimus division instead fixes a Gram
block on the denominator images, which can be eliminated before the
remaining equations are solved. The Laurent specialization retains both
features but also requires the adjoint convention and the degree
filtration to be specified. The concluding torus example calculates
their effect on determinants, polynomial blocks and kernel ranks.

\subsection{Rational bilateral mixed-type orthogonality in one variable}
\label{sec:univariate}

In one variable, principal parts and Jordan chains describe the
divisibility conditions, and Cauchy transforms give the rational moments
needed in the same finite system. Let $A=\CC[z]$ and let the
nonnegative multiindices
$\boldsymbol m\in\NN^q$, $\boldsymbol n\in\NN^p$ be arbitrary.
The polynomial row and column spaces are
$X_{\boldsymbol m}\coloneq\{f\in A^{1\times q}:\deg f_i<m_i\}$ and
$Y_{\boldsymbol n}\coloneq\{g\in A^{p\times1}:\deg g_j<n_j\}$, with
$|\boldsymbol m|\coloneq\sum_i m_i=r$ and $|\boldsymbol n|\coloneq\sum_j n_j=r+k$.
A zero bound makes the corresponding component zero. The four factors
are regular matrix polynomials, so their rational inverses exist, while
their order remains relevant because factors on the same side need not
commute. The construction also includes shared poles and confluent masses.

Let $d_L=\deg\det L_1$, $d_R=\deg\det R_1$. For every zero $\lambda$
of $\det L_1$, of multiplicity $a_\lambda$, use the row pole functions
$\frac{e_i^{\T}}{(z-\lambda)^a}$, $i\in\{1,\ldots,q\}$, $a\in\{1,\ldots,a_\lambda\}$.
Similarly use $e_j/(z-\xi)^b$ at the zeros of $\det R_1$ with their
determinant multiplicities. The sum of the root multiplicities is $d_L$ on the left and
$d_R$ on the right, so these families have sizes $s_L=qd_L$ and
$s_R=pd_R$. The determinant multiplicities bound the pole order of every entry of
the inverse matrix. The displayed families therefore contain all required
pole terms, although individual entries may require fewer of them.
Let an original normal pair of polynomial spaces $E,H$ of common
dimension $N$ containing the polynomial parts of $X_{\boldsymbol m}S_L$
and $S_RY_{\boldsymbol n}$ be fixed. Adjoining the pole families
to its dual bases gives the required ambient bases.
Retaining only the pole directions needed by the rational images can
reduce the order of the remaining system.

\Needspace{4\baselineskip}
A rational function splits into a polynomial part and pole terms.
For example, $(z^2+1)/(z-1)=z+1+2/(z-1)$: the last term is its
\emph{principal part} at $1$. In general this is the sum of negative
Laurent powers at the chosen point. The symbol $[z^a]_\infty r$
extracts the coefficient of $z^a$ from the polynomial part.

The distinction between the principal part and the polynomial part
separates the two requirements on an inverse product. Vanishing of
all finite principal parts makes it polynomial; vanishing of its
forbidden coefficients at infinity then places it in the prescribed
polynomial space. Both requirements are linear in the coordinates of the
chosen ambient basis. The following theorem identifies their common
kernel and applies the general finite construction to it.

\begin{theorem}[One-variable bilateral rational formulas]
\label{thm:univariate}
Consider on these ambient spaces the linear conditions
\begin{align}
 &\text{every finite principal part of }bL_0^{-1}L_1\text{ vanishes},
 & [z^a]_\infty(bL_0^{-1}L_1)_i&=0\quad(a\ge m_i),\label{eq:unileft}\\
 &\text{every finite principal part of }R_1R_0^{-1}c\text{ vanishes},
 & [z^b]_\infty(R_1R_0^{-1}c)_j&=0\quad(b\ge n_j).
 \label{eq:uniright}
\end{align}
A choice of independent rows forms $\Gamma_L,\Gamma_R$. Their kernels
are exactly the two rational images, and their row counts are
$N+qd_L-r$ and $N+pd_R-r-k$.

The allowed masses are those of Theorem~\ref{thm:masses}, under
either convention. The balanced mass space satisfies
\begin{equation}\label{eq:unimassdim}
 \dim\{\nu\text{ balanced}:L_1\nu R_1=0\}=pd_L+qd_R.
\end{equation}
The bivariate masses are \eqref{eq:bivariantmasses}; each mass has rank
$h\le d_L+d_R$ in either case. Using \eqref{eq:secondkindblocks}
and the complete constraints in \eqref{eq:small}, the remaining
square system, after eliminating the original dual block, has order
\begin{equation}\label{eq:uniorder}
                  t=N-r+qd_L+pd_R+h.
\end{equation}
Its determinant tests final normality and its inverse gives the normalized right block
\eqref{eq:normalizedanswer}, its left counterpart, and the kernel
\eqref{eq:kernel}. All finite principal parts in the recovered polynomial
vectors cancel, and the prescribed component degree bounds hold.
Neither a denominator-only nor a numerator-only intermediate pairing is
required to be normal.
\end{theorem}

\begin{proof}
The adjugate expression for $L_1^{-1}$ shows that each component
of $fS_L$ has denominator dividing $\det L_1$. Partial fractions
therefore express it as a polynomial plus the stated row pole
functions. On the right, $S_Rg=R_0\operatorname{adj}(R_1)g/\det R_1$.
Thus its scalar denominators divide $\det R_1$, and partial
fractions give the stated column pole functions, with the same
bound on each pole order.
Since the original spaces contain the polynomial parts, adjoining
these pole functions expresses every required rational product in
the enlarged bases.

The image conditions follow from the uniqueness of the inverse
rational product: an ambient row $b$ has only the possible preimage
$f=bL_0^{-1}L_1$. A rational function of one variable is
polynomial if and only if all its finite principal parts vanish.
After this condition is satisfied, the coefficient tests at
infinity impose exactly $\deg f_i<m_i$ in each component.
Thus the common kernel of the left tests is $X_{\boldsymbol m}S_L$.
The right inverse product is $R_1R_0^{-1}c$, in that order,
and the same argument proves the right kernel assertion.
Because all factors are invertible over $\CC(z)$, multiplying
by $S_L$ or $S_R$ is injective. The two image dimensions are
therefore $r$ and $r+k$. Polynomial and proper pole terms have
zero intersection, so the ambient dimensions are exactly
$N+qd_L$ and $N+pd_R$. The ranks of the full lists of tests
are consequently $N+qd_L-r$ and $N+pd_R-r-k$. Selecting
independent rows preserves their common kernels and gives the
asserted row counts. The lists are finite: each ambient basis
vector has finitely many poles of finite order and a polynomial
part of finite degree.

The same decomposition supplies explicit extensions of $\pi_L,\pi_R$ on the whole
ambient space: each extension takes the allowed component powers
of the polynomial part of the inverse rational product. On a
constrained image the inverse product is already an allowed
polynomial, so the extension agrees there with the actual inverse.
These are the recovery maps required by Theorems~\ref{thm:finite} and
\ref{thm:kernel}. In the reduced system of Theorems~\ref{thm:finite} and
\ref{thm:kernel}, the column blocks have sizes
$pd_R$, $h$ and $N+qd_L-r$, in the order of \eqref{eq:small}. Their sum is
$N-r+qd_L+pd_R+h$, proving \eqref{eq:uniorder}.
The right image conditions, remaining pole equations, mass
equations and $k$ normalization equations give the same total
number of rows. The resulting square system is therefore precisely
the specialization of those theorems, with their determinant
criterion and polynomial formulas.

The mass dimensions follow independently from Smith normal form
over $\CC[z]$, which splits each polynomial
matrix quotient into scalar quotients by its invariant factors.
If an invariant factor is $d(z)$ of degree $a$, division by
$d$ gives the basis $1,z,\ldots,z^{a-1}$ of $\CC[z]/(d)$;
a constant factor gives the zero quotient. Hence its dimension
is $a$. The product of the invariant factors is a nonzero
constant times the determinant, so their degrees add to the
determinant degree. Thus the quotient dimensions
are $d_L,d_R$, and Theorem~\ref{thm:masses} gives both
\eqref{eq:unimassdim} and the mass-rank bound.
\end{proof}

\subsubsection{Canonical Jordan chains and polynomial quotients}
\label{sec:jordandivision}

Divisibility by $(z-\lambda)^k$ requires the first $k$ Taylor
coefficients to vanish. A matrix division has the same structure after
an invertible local change of components, but those component directions
can depend on $z$. A \emph{Jordan chain} here is the finite list of
Taylor coefficients of one such direction. The \emph{partial
multiplicity} is the power of $z-\lambda$ by which the corresponding
component must be divisible.
For the polynomial lifts $P_L,P_R$ from \eqref{eq:polynomiallifts},
these directions are obtained at each zero $\lambda$ from a local
Smith normal form,
$P_L=U_\lambda\diag((z-\lambda)^{\theta^L_{\lambda,i}})V_\lambda$.
For every positive partial multiplicity, the Taylor coefficients of
$V_\lambda^{-1}e_i$ through order $\theta^L_{\lambda,i}-1$ give the
complete right chain $[v_{\lambda,i,\ell}]_{\ell\in\{0,\ldots,\theta^L_{\lambda,i}-1\}}$. For $P_R$, the rows of the
inverse left Smith factor similarly give complete left chains
$[w_{\xi,j,\ell}]_{\ell\in\{0,\ldots,\theta^R_{\xi,j}-1\}}$.
Their equations are
\begin{equation}\label{eq:jordanequations}
 \sum_{a=0}^{\ell}\frac{P_L^{(a)}(\lambda)}{a!}
 v_{\lambda,i,\ell-a}=0,\qquad
 \sum_{a=0}^{\ell}w_{\xi,j,\ell-a}
 \frac{P_R^{(a)}(\xi)}{a!}=0,
\end{equation}
in the respective chain ranges. The chain equations follow by expanding
$P_L(V_\lambda^{-1}e_i)=(z-\lambda)^{\theta^L_{\lambda,i}}
U_\lambda e_i$, and its right counterpart, below the indicated
power. One chain of full length for each positive Smith exponent records
every component of the local quotient. This completeness is needed
for division, since a shorter collection of solutions
of \eqref{eq:jordanequations} can leave a pole undetected.
These classical Smith and maximal root-polynomial constructions
\cite{DopicoNoferini} enter matrix Christoffel quasideterminants in
\cite[Theorem~2]{MatrixChristoffelIMRN} and rational transformations in
\cite{Matrix2016}.

\Needspace{4\baselineskip}
\begin{proposition}[The Jordan entries of the constraint matrices]
\label{prop:jordanconstraints}
Put $p_L\coloneq a_Lb$ and $p_R\coloneq a_Rc$. Complete finite division tests are
\begin{equation}\label{eq:jordanentries}
\begin{aligned}
 \gamma_{\lambda,i,\ell}(b)&\coloneq
 \sum_{a=0}^{\ell}\frac{p_L^{(a)}(\lambda)}{a!}
                  v_{\lambda,i,\ell-a},
 &&\ell\in\{0,\ldots,\theta^L_{\lambda,i}-1\},\\
 \delta_{\xi,j,\ell}(c)&\coloneq
 \sum_{a=0}^{\ell}w_{\xi,j,\ell-a}
                  \frac{p_R^{(a)}(\xi)}{a!},
 &&\ell\in\{0,\ldots,\theta^R_{\xi,j}-1\}.
\end{aligned}
\end{equation}
Their vanishing is equivalent to $p_LP_L^{-1}$ and $P_R^{-1}p_R$ being
polynomial. The component bounds impose the additional conditions
\begin{equation}\label{eq:jordantails}
 [z^a]_\infty(p_LP_L^{-1})_i=0\ (a\ge m_i),\qquad
 [z^b]_\infty(P_R^{-1}p_R)_j=0\ (b\ge n_j).
\end{equation}
Evaluating these tests on the ambient basis and selecting independent
rows gives exactly the matrices \eqref{eq:constraints}.
\end{proposition}

\begin{proof}
Near $\lambda$, write the right Smith direction as
$V_\lambda^{-1}e_i=\sum_{j\ge0}v_{\lambda,i,j}(z-\lambda)^j$.
Multiplying this Taylor expansion by that of $p_L$ shows that
the coefficient of $(z-\lambda)^\ell$ in
$(p_LV_\lambda^{-1})_i$ is exactly the first expression in
\eqref{eq:jordanentries}.

The inverse product in Smith coordinates is
\[
 p_LP_L^{-1}
 =p_LV_\lambda^{-1}
    \diag((z-\lambda)^{-\theta^L_{\lambda,i}})
       U_\lambda^{-1}.
\]
Both $U_\lambda$ and its inverse are regular at the point.
Thus the product has no pole there if and only if the $i$th
component of $p_LV_\lambda^{-1}$ vanishes to order at least
$\theta^L_{\lambda,i}$ for every $i$. These are precisely
the listed Taylor-coefficient conditions. Applying them at all
roots removes every finite pole and proves polynomiality.
On the right write $P_R=U_\xi D_\xi V_\xi$. Then
$P_R^{-1}p_R=V_\xi^{-1}D_\xi^{-1}U_\xi^{-1}p_R$.
The $j$th row of $U_\xi^{-1}$ has Taylor coefficients
$w_{\xi,j,\ell}$. Multiplying it by the Taylor series of $p_R$
gives the second sum in \eqref{eq:jordanentries}. Its first
$\theta^R_{\xi,j}$ coefficients must vanish before multiplication
by $(z-\xi)^{-\theta^R_{\xi,j}}$. The factor $V_\xi^{-1}$
is invertible and regular, so these conditions are also sufficient.

Once the finite poles have been removed, the quotient bounds remain
to be imposed. The polynomial quotients are $bL_0^{-1}L_1$ and
$R_1R_0^{-1}c$. The coefficient conditions in
\eqref{eq:jordantails} therefore select exactly the prescribed
polynomial spaces. Only finitely many such coefficients can occur:
each quotient of an ambient basis vector has a finite polynomial
part, and there are finitely many basis vectors. Evaluating all
tests on that basis and selecting independent rows consequently
gives finite complete constraint matrices.
\end{proof}

For $P_L=\bigl[\begin{smallmatrix}z^2&0\\-z&1\end{smallmatrix}\bigr]$
and a row $p_L=[p_1\ p_2]$, the chain $v_0=e_1$, $v_1=e_2$ gives
$p_1(0)=0$ and $p_1'(0)+p_2(0)=0$.
The second condition couples the derivative of the first component
to the value of the second because the chain direction itself
varies with the variable. The inverse product exhibits this coupling:
\[
 P_L^{-1}=\begin{bmatrix}z^{-2}&0\\z^{-1}&1\end{bmatrix},\qquad
 p_LP_L^{-1}=\left[\frac{p_1+zp_2}{z^2}\ \ p_2\right].
\]
The numerator in the first component is divisible by $z^2$
exactly when its constant and linear coefficients vanish, giving
the two tests. The row $[-z\ 1]$ has quotient $[0\ 1]$,
whereas $[z\ 0]$ has quotient $[z^{-1}\ 0]$: it passes the
point-value test but retains a pole. The degree conditions have an
independent role, as is already apparent for a unimodular factor.
Its finite chain list can be empty while its quotient violates
the prescribed component bounds: for $P_L=\bigl[\begin{smallmatrix}1&z\\0&1\end{smallmatrix}\bigr]$,
the row $[1\ 0]$ has quotient $[1\ {-z}]$, outside the space of constant rows.

\noindent\textbf{Cauchy data and confluence.}
With the divisibility conditions determined, the remaining input
consists of the rational moments. For either convention,
\eqref{eq:secondkindblocks} gives the polynomial--pole
entries $\widetilde\beta_\mu(b_\ell,e_j/(z-\xi)^b)$ and
$\widetilde\beta_\mu(e_i^{\T}/(z-\lambda)^a,a_\ell)$ in the original
dual bases $b_\ell,a_\ell$.
For a finite bivariate measure on a compact support avoiding both poles, put
$\mathfrak F_{i,j}(\lambda,\xi)
\coloneq\iint\frac{\dd\mu_{i,j}(x,y)}{(x-\lambda)(y-\xi)}$. Then
\begin{equation}\label{eq:bivariatecauchy}
 (E_0)_{(i,\lambda,a),(j,\xi,b)}
   =\frac{\partial_\lambda^{a-1}\partial_\xi^{b-1}
                   \mathfrak F_{i,j}(\lambda,\xi)}{(a-1)!(b-1)!}.
\end{equation}
Because the two integrations use independent arguments, that
independence persists when $\lambda=\xi$.
For a balanced measure, $F_{i,j}(w)\coloneq\int(z-w)^{-1}\,\dd\mu_{i,j}(z)$ gives
$\mathfrak F_{i,j}(\lambda,\xi)
=\frac{F_{i,j}(\lambda)-F_{i,j}(\xi)}{\lambda-\xi}$.
The sign follows from
\[
 \frac1{(z-\lambda)(z-\xi)}
 =\frac1{\lambda-\xi}
       \left(\frac1{z-\lambda}-\frac1{z-\xi}\right).
\]
Differentiation under the integral is legitimate because the
compact support stays a positive distance from the poles.
Since $\partial_w^k(z-w)^{-1}=k!(z-w)^{-k-1}$,
when the poles coincide the product is $(z-\lambda)^{-a-b}$
and the entry is
$F_{i,j}^{(a+b-1)}(\lambda)/(a+b-1)!$.
For independent arguments the same differentiation applies to each
denominator separately, giving \eqref{eq:bivariatecauchy} with the same positive sign.
When no measure representation is specified, the rational moments
\eqref{eq:secondkindblocks} provide the same blocks as algebraic data.

These Cauchy entries are the ones used in
Subsection~\ref{sec:scalarresidual}, where they determine the
three-by-three reduced matrix from the original polynomials.
Subsubsection~\ref{sec:mixedcurveexample} extends that calculation
to matrix factors and a bilateral mass on a curve. In the present
one-variable setting, the corresponding image constraints are the
Jordan tests of Proposition~\ref{prop:jordanconstraints}, so the
moments and constraints together specify the reduced system.

\subsubsection{Bilateral masses and their parameters}

The distinction between the rank of a mass and the dimension of
its parameter space already occurs for
$q=2$, $p=3$, $L_1=\diag(z,1)$ and $R_1=\diag(z,1,1)$.
For these factors, the balanced mass $\nu(f,g)=f_1'(0)g_1(0)+f_1(0)g_1'(0)$
is the first derivative of the product $f_1g_1$ at zero.
On the two value--derivative coordinates its Gram is
$[0\ 1;\ 1\ 0]$, so its rank is two. If both arguments are
multiplied by their denominator factors, their first components
contain $z$ and the product contains $z^2$; its derivative at
zero vanishes. This proves $L_1\nu R_1=0$.
The unilateral actions do not annihilate this mass: with the left
denominator alone,
$\nu([z\ 0],e_1)=1$, while the right denominator alone gives
$\nu(e_1^{\T},[z\ 0\ 0]^{\T})=1$.
Neither unilateral equation therefore holds. The determinant degrees
are $d_L=d_R=1$, giving balanced mass dimension
$pd_L+qd_R=3+2=5$. The bivariate space is
infinite dimensional: the rank-one masses
$\nu_j(f,g)\coloneq f_1(0)[z^j]g_1$, $j\in\NN$, are homogeneous and independent,
since $\nu_j(e_1^{\T},z^ke_1)=\delta_{j,k}$. Thus finite rank does not
mean finitely many scalar mass parameters.

\subsubsection{An explicit order bound for the one-variable reduced system}
\label{sec:stepbound}

The pole counts also yield a numerical bound on the reduced order
once the original spaces have been chosen by degree. For a nonzero
matrix polynomial $M$, let $\deg M$ be the largest degree of an
entry, and put
\[
 \begin{aligned}
 b_L&\coloneq\max\{0,\deg\operatorname{adj}L_1+\deg L_0-d_L\},\\
 b_R&\coloneq\max\{0,\deg R_0+\deg\operatorname{adj}R_1-d_R\}.
 \end{aligned}
\]
The adjugate formulas show that the polynomial parts of $fS_L$
and $S_Rg$ have degrees below $\max_i m_i+b_L$ and
$\max_j n_j+b_R$, respectively. For example, on the left the
numerator degree is at most
$\max_i m_i-1+\deg\operatorname{adj}L_1+\deg L_0$;
polynomial division subtracts $d_L$ from this bound whenever
the polynomial part is nonzero.

For $c$ components, order the column monomials as
$e_1,\ldots,e_c,ze_1,\ldots,ze_c,z^2e_1,\ldots$, and let
$V_{c,N}$ be the span of the first $N$ terms. Assume the original
pairing normal on every
$V_{q,N}^{\T}\times V_{p,N}$. A sufficient common dimension is
\[
 N=\max\{q(\max_i m_i+b_L),\ p(\max_j n_j+b_R)\}.
\]
Indeed the first $qM$ row monomials contain all components of
degree below $M$, and the same assertion holds on the column side.
Thus these original spaces contain both required polynomial parts.

If the components of each multiindex differ by at most one, then
$q\max_i m_i\le r+q-1$ and
$p\max_j n_j\le r+k+p-1$. Subtracting $r$ from the displayed
choice of $N$ and using \eqref{eq:uniorder} gives
\begin{equation}\label{eq:unibound}
 t\le\max\{qb_L+q-1,\ pb_R+p-1+k\}+qd_L+pd_R+h.
\end{equation}
For fixed factors and fixed $k$, every term on the right is bounded
independently of the level, since $h\le d_L+d_R$. More generally,
a fixed bound on the spread of each multiindex gives the same
conclusion with the corresponding larger constants. For unrestricted
multiindices the exact order remains \eqref{eq:uniorder}.
This choice can be compared with the construction in
Subsection~\ref{sec:orders}, which fixes its ambient spaces directly
from the four factors and component bounds. Its one-variable order is \eqref{eq:orderlineexact};
\eqref{eq:orderrank} calculates a smaller order by removing
dependent polynomial coefficient columns. The different counts reflect the ambient spaces used in the same
finite construction.

\subsection{Christoffel and Geronimus reductions}\label{sec:special}

The roles of the numerator and denominator factors become separate
in the Christoffel and Geronimus cases. In the former, the transformed
arguments are polynomial images, so the construction reduces to their
membership conditions. In the latter, the functional equation fixes
the pairing on denominator images, providing a nondegenerate block
that can be eliminated. Both reductions retain independent factors
and arbitrary finite polynomial row and column spaces.

\subsubsection{Bilateral Christoffel multiplication}

For identity denominators set $L_1=I_q$, $R_1=I_p$, $L\coloneq L_0$, $R\coloneq R_0$,
and assume that $\det L$ and $\det R$ are regular.
The homogeneous mass vanishes and $\widehat\beta(f,g)=\beta_\mu(fL,Rg)$.
Choose finite target spaces $X\subset A^{1\times q}$ and
$Y\subset A^{p\times1}$, and an original normal pair $E,H$ of
dimension $N$ containing $XL$ and $RY$. In its dual bases
$\cB,\cA$, choose matrices $\Gamma_L,\Gamma_R$ with independent
rows whose kernels are the coefficient spaces of $XL$ and $RY$.
Thus their equations test whether a polynomial equals $fL$ or $Rg$
with $f\in X$, $g\in Y$. For the dimensions in the proposition,
$\Gamma_L$ has size $(N-r)\times N$ and $\Gamma_R$ has size
$(N-r-k)\times N$.

\Needspace{4\baselineskip}
Since these polynomial spaces already contain the images, no
additional rational vectors are required. Orthogonality restricts
the unknown coefficients to the directions described by the left
constraints. The resulting system determines the multiplied polynomial $R\widehat a$;
the image conditions then guarantee that division by $R$ gives
a polynomial in $Y$.

\begin{proposition}[Bilateral Christoffel formulas]\label{prop:Christoffel}
Suppose $\dim X=r$, $\dim Y=r+k$ and $\lambda:Y\to\CC^k$ is the
prescribed right normalization. Extend $Rg\mapsto\lambda(g)$ from
$RY$ to $H$, with matrix $\eta\in\CC^{k\times N}$ in the basis $\cA$. Put
\begin{equation}\label{eq:ChristoffelT}
 T_{C,k}\coloneq\begin{bmatrix}\Gamma_R\Gamma_L^{\T}\\
                         \eta\Gamma_L^{\T}\end{bmatrix},
 \qquad J_k=\begin{bmatrix}0_{(N-r-k)\times k}\\I_k\end{bmatrix}.
\end{equation}
The matrix has order $N-r$. For $k>0$ the normalized problem has a
unique solution if and only if $T_{C,k}$ is invertible, and then
\begin{equation}\label{eq:ChristoffelA}
 \widehat a=R^{-1}\cA\Gamma_L^{\T}T_{C,k}^{-1}J_k.
\end{equation}
For a square pair, let $M_C\coloneq\Gamma_R\Gamma_L^{\T}$. Normality is
equivalent to $\det M_C\ne0$, and
\begin{equation}\label{eq:ChristoffelK}
 R(P)\widehat K(P,Q)L(Q)
 =\cA(P)\bigl[I_N-\Gamma_L^{\T}M_C^{-1}\Gamma_R\bigr]\cB(Q).
\end{equation}
The correction in \eqref{eq:ChristoffelK} has tensor rank exactly $N-r$.
\end{proposition}

\begin{proof}
Original duality identifies the pairing with the ordinary product
of coefficient rows and columns. Orthogonality to $XL$ therefore
means that the column $z$ annihilates $\ker\Gamma_L$, and consequently
$z=\Gamma_L^{\T}u$, uniquely, because $\Gamma_L$ has independent
rows. Membership in $RY$ and normalization become respectively
$\Gamma_R\Gamma_L^{\T}u=0$ and
$\eta\Gamma_L^{\T}u=I_k$. This is the square system
$T_{C,k}u=J_k$, proving the criterion and the formula. Its image
condition also ensures that division by $R$ gives a polynomial in $Y$.

For a square pair, a coefficient column $z\in\ker\Gamma_R$
is in the right radical precisely when every coefficient row of
$XL$ annihilates it. The latter condition is
$z\in\im\Gamma_L^{\T}$. Thus normality is equivalent to
$\ker\Gamma_R\cap\im\Gamma_L^{\T}=0$. This is equivalent to invertibility of $M_C$.
Then $P_C\coloneq I_N-\Gamma_L^{\T}M_C^{-1}\Gamma_R$ is the projection
onto $\ker\Gamma_R$ along $\im\Gamma_L^{\T}$: it is the
identity on the first space and zero on the second. In coordinates, $\Gamma_RP_C=0$ and $P_Cz=z$ whenever
$\Gamma_Rz=0$. Also $x^{\T}P_C=x^{\T}$ for every
$x\in\ker\Gamma_L$, because $x^{\T}\Gamma_L^{\T}=0$.
Consequently $\cA(P)P_C\cB(Q)$ has the required column and
row images and the reproducing property, giving
\eqref{eq:ChristoffelK}. Its difference from
$\cA(P)\cB(Q)$ has coefficient matrix
$-\Gamma_L^{\T}M_C^{-1}\Gamma_R$.
The first and last factors have full column and row rank $N-r$,
respectively, and the middle one is invertible. The coefficient
matrix therefore has rank $N-r$. Independence of the ambient
basis functions identifies this matrix rank with the tensor rank
of the kernel correction.
\end{proof}

If $\dim X=r+k$, $\dim Y=r$ and $\zeta$ extends the left normalization
through $f\mapsto fL$, put
$T^L_{C,k}\coloneq[\Gamma_L\Gamma_R^{\T};\ \zeta\Gamma_R^{\T}]$.
Here $T^L_{C,k}$ has order $N-r$ and the lower identity in
$J_k$ has $k$ rows, with the zero block of height $N-r-k$.
Indeed left orthogonality writes the row coefficients as
$u^{\T}\Gamma_R$; left image membership and normalization are
$\Gamma_L\Gamma_R^{\T}u=0$ and
$\zeta\Gamma_R^{\T}u=I_k$. Thus the normalized row block is
$\widehat b=((T^L_{C,k})^{-1}J_k)^{\T}\Gamma_R\cB L^{-1}$.
The image equations ensure that its displayed division is polynomial.
Both constraint matrices include quotient bounds at infinity,
which finite-zero conditions alone do not impose.

For the matrix specialization of \cite[Theorem~2]{MatrixChristoffelIMRN},
take $A=\CC[z]$, $q=p$, $L=W$ monic of degree $d$, and $R=I_p$.
Write $P_j$ for the original monic left matrix polynomials and $\pi_j$
for their complete Jordan-jet blocks, as in Definition~13 of that
reference. Put
$\Pi_{n,d}\coloneq[\pi_n^{\T}\ \cdots\ \pi_{n+d-1}^{\T}]^{\T}$ and
$\mathcal P_{n,d}\coloneq[P_n^{\T}\ \cdots\ P_{n+d-1}^{\T}]^{\T}$.
The complete finite conditions have $pd$ scalar columns, because
$\deg\det W=pd$. Thus each $\pi_j$ has size $p\times pd$,
$\Pi_{n,d}$ has size $pd\times pd$, and
$\mathcal P_{n,d}$ has size $pd\times p$.
The relation with the matrix formula follows by expanding $\widehat P_nW$ in the monic original
basis. Its leading coefficient and degree give the term $P_{n+d}$.
Pairing with the original right blocks of degrees below $n$
forces all coefficients at $P_j$, $j<n$, to vanish. The remaining
coefficient row is $C\in\CC^{p\times pd}$, so
$\widehat P_nW=P_{n+d}+C\mathcal P_{n,d}$.
Applying the complete division tests to this identity gives
$\pi_{n+d}+C\Pi_{n,d}=0$. Solving for $C$ yields
\[
 \widehat P_nW=P_{n+d}
               -\pi_{n+d}\Pi_{n,d}^{-1}\mathcal P_{n,d},
\]
the Schur complement in formula~(10) of that theorem. Conversely,
complete jets give divisibility and a monic quotient, so
$\det\Pi_{n,d}\ne0$ is precisely the invertibility condition for the remaining finite system.
When the leading coefficient is singular, the same distinction
between finite divisibility and component degrees persists: the
quotient bounds supply conditions not detected by the finite jets.

\begin{remark}[Normality with vanishing auxiliary adjacent minors]
\label{rem:adjacentminors}
On $A=\CC[x,y]/(x^2+y^2-1)$ put
$\mu_k[f]\coloneq\tfrac12\int_{1/4}^{1/2}x^k
[f(x,\sqrt{1-x^2})+f(x,-\sqrt{1-x^2})]\dd x$, $k\ge0$.
For every nonzero real class $f=a(x)+yb(x)$,
$\mu_k[f^2]=\int_{1/4}^{1/2}x^k
[a(x)^2+(1-x^2)b(x)^2]\dd x>0$.
All principal Gram matrices in $1,x,y,x^2,xy,\ldots$ are therefore
positive definite. For $n,k\ge1$, the auxiliary minor
$\Theta_{2n}^{(k)}$ of \cite[Section~2.2, Remark~2.6]{WuHu}
has rows $x^{k+j},x^{k+j-1}y$, $j\in\{0,\ldots,n-1\}$, and the first $2n$
columns of that basis: $n$ rows of each $y$-parity, but $n+1$ even
and $n-1$ odd columns. The even--odd cross moments vanish because the two branches are
averaged with equal weights. After grouping rows and columns by
parity, the only nonzero rectangles have sizes $n\times(n+1)$
and $n\times(n-1)$. Their total rank is at most
$n+(n-1)=2n-1$, so $\Theta_{2n}^{(k)}=0$. These minors occur as denominators in
\cite[(2.12), (2.13) and the coefficients of (2.14)]{WuHu}, while Proposition~\ref{prop:Christoffel}
applies to the normal transformation $\mu_{k+1}=x\mu_k$.
\end{remark}

\subsubsection{Polynomial elimination and Geronimus division}
\label{sec:polynomialelimination}

For factors satisfying \eqref{eq:compatibility}, fix a solution $\widehat\beta$.
Its moments on $U\coloneq X\cap A^{1\times q}L_1$ and
$V\coloneq Y\cap R_1A^{p\times1}$ are fixed by
$\widehat\beta(fL_1,R_1g)=\beta_\mu(fL_0,R_0g)$, independently of preimages.
A nonsingular order-$d$ minor gives families $P_i\coloneq f_iL_1$, $Q_i\coloneq R_1g_i$
with $\beta_\mu(f_iL_0,R_0g_j)=\delta_{i,j}$ after changing bases.
Collect them as rows of $P$ and columns of $Q$, and complete to bases
$[P;\Phi]$, $[Q\ \Psi]$ of $X,Y$.
Put $C\coloneq\widehat\beta(\Phi,Q)$, $D\coloneq\widehat\beta(P,\Psi)$,
$E\coloneq\widehat\beta(\Phi,\Psi)$, $Z\coloneq\Phi-CP$, $W\coloneq\Psi-QD$ and $S_G\coloneq E-CD$.
These subtractions make the remaining families orthogonal to the
known block $P,Q$. Their Gram matrix is therefore the reduced
matrix $S_G$, the \emph{Schur complement} of that block.

The known block comes from moments fixed by the functional equation,
so its construction does not require rational inversion. Subtracting
the pairings with this block leaves polynomial families orthogonal
to it, and the remaining calculation involves only finite numerical
inverses. This is the reason that the reduction also applies to
compatible nonregular factors.

\begin{proposition}[Bilateral formulas without polynomial division]\label{prop:Geronimus}
The target pairing has rank $d+\rank S_G$. For $\dim X=r$,
$\dim Y=r+k$, let $\widehat M$ be its Gram matrix in these bases
with the $k$ normalization rows appended. Then
\begin{equation}\label{eq:GeronimusT}
 T_{G,k}\coloneq\begin{bmatrix}S_G\\\lambda(W)\end{bmatrix},\qquad
 \det\widehat M=\det T_{G,k},\qquad
 \widehat a=W T_{G,k}^{-1}
                 \begin{bmatrix}0_{(r-d)\times k}\\I_k\end{bmatrix}.
\end{equation}
For $k>0$, the normalized problem has a unique solution exactly when
$\det T_{G,k}\ne0$. For a square pair, normality is equivalent to
$\det S_G\ne0$, and
\begin{equation}\label{eq:GeronimusK}
 \widehat K(P_0,Q_0)
 =\sum_{i=1}^d Q_i(P_0)P_i(Q_0)
                  +W(P_0)S_G^{-1}Z(Q_0).
\end{equation}
The last term has tensor rank $r-d$.
\end{proposition}

\begin{proof}
With $r$ rows and $r+k$ columns, the block sizes are
$C:(r-d)\times d$, $D:d\times(r+k-d)$ and
$E:(r-d)\times(r+k-d)$. Direct multiplication gives
\[
 \begin{bmatrix}I_d&0\\-C&I_{r-d}\end{bmatrix}
 \begin{bmatrix}I_d&D\\C&E\end{bmatrix}
 \begin{bmatrix}I_d&-D\\0&I_{r+k-d}\end{bmatrix}
 =\begin{bmatrix}I_d&0\\0&S_G\end{bmatrix}.
\]
The two invertible changes of basis have determinant one when
square and implement the subtractions $Z=\Phi-CP$ and $W=\Psi-QD$.
The diagonal form consequently gives the rank $d+\rank S_G$ of
the original Gram. Orthogonality to $P$ forces a
candidate into $Wc$; its remaining equations are
$S_Gc=0$, $\lambda(W)c=I_k$. In the augmented Gram the column change replaces the normalization
block $[\lambda(Q)\ \lambda(\Psi)]$ by
$[\lambda(Q)\ \lambda(W)]$. Subtracting $\lambda(Q)$ times
the first $d$ rows then makes its first block zero. These
unit triangular operations give a block diagonal square matrix
with blocks $I_d$ and
$T_{G,k}=[S_G;\lambda(W)]$, of respective orders $d$ and
$r+k-d$. Hence $\det\widehat M=\det T_{G,k}$.
The system $T_{G,k}c=[0;I_k]$ proves both the formula and the
uniqueness criterion in \eqref{eq:GeronimusT}.
For a square pair, the inverse diagonal Gram gives
\eqref{eq:GeronimusK}. The independent families $Z,W$ have size
$r-d$, and $S_G^{-1}$ is invertible, so the last tensor has rank
$r-d$.
\end{proof}

When the dimension excess occurs on the left, transposition of the
construction gives
$[S_G^{\T};\lambda_L(Z)]$. It allows $d=0$ and inverts only numerical matrices.

For Geronimus division set $L_0=I_q$, $R_0=I_p$, $L\coloneq L_1$, $R\coloneq R_1$
with regular determinants. Formula~\eqref{eq:pairing} computes $C,D,E$;
for example,
$C_{\alpha,i}=\widetilde\beta_\mu(\phi_\alpha L^{-1},g_i)
+\nu(\phi_\alpha,Rg_i)$.
The first term of \eqref{eq:GeronimusK} is
$R(P_0)K_d^\mu(P_0,Q_0)L(Q_0)$, with $K_d^\mu\coloneq\sum_i g_i f_i$.
Similarly
$D_{i,\alpha}=\widetilde\beta_\mu(f_i,R^{-1}\psi_\alpha)
+\nu(f_iL,\psi_\alpha)$ and
$E_{\alpha,\beta}=\widetilde\beta_\mu(\phi_\alpha L^{-1},
R^{-1}\psi_\beta)+\nu(\phi_\alpha,\psi_\beta)$.
The mass terms remain in these three Gram blocks because
$L\nu R=0$ removes the mass only when both arguments belong
to their denominator images. Consequently the terms in $C$ and $D$
are part of the Geronimus data even though the known block has no
mass contribution.

\subsubsection{Additive and scalar specializations}

For $\widehat\beta=\beta_\mu+\nu$, any bilinear mass has finite rank
$h$ on prescribed finite spaces. Setting both rational factors equal to the identity in
Subsection~\ref{sec:finite} gives, for an original normal square pair,
$\mathbb T=-I_h-H_0F_0$ and
$\widehat K=\cA\bigl[I_N-F_0(I_h+H_0F_0)^{-1}H_0\bigr]\cB$.
Here $F_0\in\CC^{N\times h}$ and $H_0\in\CC^{h\times N}$ factor
the mass matrix in the original dual bases, as in \eqref{eq:blocks}.
For scalar polynomials in several variables and masses at finitely
many points, this factorization is obtained from the evaluation
matrices of the original polynomial bases. The displayed formula
then gives the multivariate Uvarov kernel update studied in
\cite{FPPX2010,Granada2010,Granada2016}, together with the corresponding
polynomial formulas and the invertibility criterion for the modified
moment functional. The finite-space construction here also permits
component degree bounds and general bilinear masses.

This Uvarov addition is the four-factor case
$L_1=L_0=aI_q$, $R_1=R_0=I_p$ when $a\in A$ is regular and
$\nu(af,g)=0$ for all $f,g$.
For a balanced mass at finitely many points of an integral
positive-dimensional support, a nonzero function vanishing at those
points exists. A sufficiently high power annihilates all its value
and derivative terms and remains regular, providing $a$ even for confluent masses.
Taking $q=p=1$, $L_1=b$, $L_0=a$, and $R_0=R_1=1$ gives
$b\widehat\mu=a\mu$. Proposition~\ref{prop:compatibility} and
Subsubsection~\ref{sec:polynomialelimination} also cover nonregular $a,b$.
For example, on $A=\CC[x,y]/(xy)$ take $a=x(x+2)$ and $b=x(x+1)$.
Both annihilators are $(y)$. Choose
$u(f)=\int_1^2(t+1)f(t,0)\,\dd t+\int_1^2f(0,t)\,\dd t$ and
$v(f)=\int_1^2(t+2)f(t,0)\,\dd t+3\int_1^2f(0,t)\,\dd t$.
Then $bv=au$, but cancellation of $x$ loses the free component mass:
$[(x+1)v-(x+2)u](f)=\int_1^2f(0,t)\,\dd t$.
Both functionals are positive definite: a class has the unique form
$c+xp(x)+yq(y)$, and vanishing on an interval of each axis forces
both one-variable restrictions, hence the class, to vanish.
The same example also exhibits the finite elimination: the
degree-two space has basis $1,x,y,x^2,y^2$ and dimension
five. Take $P=b=x(x+1)$ and $Q=4/55$. Since
\[
 v(b)=\int_1^2 t(t+1)(t+2)\,\dd t=\frac{55}{4},
\]
the known block is $v(PQ)=1$. Thus $d=1$ and the remaining
Gram in \eqref{eq:GeronimusK} has order $5-1=4$.

\subsection{Laurent supports and a bilateral perturbation on a torus}
\label{sec:torus}

For Laurent polynomials, the finite construction depends on the
degree filtration induced by the chosen embedding and on the
adjoint used in the second argument. Expressing conjugation in
the bilinear notation preserves the determinant formulas, while
the filtration determines the spaces to which they apply. A mass
distributed along a circle in a torus provides an explicit
calculation of both effects.

The product of circles $\mathbb T^r\coloneq\{z:|z_j|=1,
j\in\{1,\ldots,r\}\}$ has Laurent coordinate algebra
$A\coloneq\CC[z_1^{\pm1},\ldots,z_r^{\pm1}]/J$ after imposing relations $J$.
Indeed $z_j\coloneq x_j+\mathrm{i} y_j$ and $z_j^{-1}=x_j-\mathrm{i} y_j$
identify the full Laurent algebra with
$\CC[x_1,y_1,\ldots,x_r,y_r]/(x_j^2+y_j^2-1:j\in\{1,\ldots,r\})$.
Laurent polynomials allow positive and negative integer powers.
Real-coordinate degree becomes Laurent length $|\alpha|_1\coloneq\sum_j|\alpha_j|$.
The space $F_NA$ consists of classes admitting a representative with
all monomials of length at most $N$.
The full-torus Christoffel construction is developed in \cite{LaurentTorus}. Mixed-type Laurent orthogonality on the unit circle,
including diagonal Christoffel and Geronimus perturbations, is treated
in \cite{MixedLaurent2026}.

The Laurent algebra also encodes conjugation of values through the
operation
$(\sum_\alpha c_\alpha z^\alpha)^\star
\coloneq\sum_\alpha\overline{c_\alpha}z^{-\alpha}$.
Assume $J^\star=J$, so this operation is well defined on classes;
the spaces $F_NA$ are then $\star$-stable. For a matrix let
$B^\dagger\coloneq(B^\star)^{\T}$. For rows $f,h$ of lengths $q,p$ define
$\langle f,h\rangle_\mu\coloneq\beta_\mu(f,h^\dagger)$.
It is linear in $f$ and conjugate-linear in $h$; Hermitian symmetry
is an additional condition. The second transformed row is $hR_1^\dagger$, since
$(hR_1^\dagger)^\dagger=R_1h^\dagger$; the order follows
from reversing products under the adjoint. Thus the four-factor
equation becomes
\begin{equation}\label{eq:torusfour}
 \langle fL_1,hR_1^\dagger\rangle_{\widehat\mu}
 =\langle fL_0,hR_0^\dagger\rangle_\mu.
\end{equation}
For regular factors, the particular solution and mass decomposition
therefore take the form
$\langle f,h\rangle_{\widehat\mu}
=\widetilde\beta_\mu(fS_L,S_Rh^\dagger)+\nu(f,h^\dagger)$.
For arbitrary row spaces $X,Y$, including different component bounds, the Gram
matrix is the bilinear one for $X,Y^\dagger$, so its determinant criterion,
bordered formulas and correction ranks apply. Left rows are unchanged;
right rows and their normalization are adjoints of the bilinear columns
and normalization.
For $R_j=L_j^\dagger$, $j\in\{0,1\}$, one has
$S_R=L_0^\dagger(L_1^\dagger)^{-1}
=(L_1^{-1}L_0)^\dagger=S_L^\dagger$.
The rational part is then $\langle fS_L,hS_L\rangle_\mu$:
the same change acts on both row arguments.

For a balanced scalar functional $u:A\to\CC$ with all filtered
Gram matrices invertible, choose a basis column $\chi_{[n]}$ of a
$\star$-stable complement of $F_{n-1}A$ in $F_nA$, with $F_{-1}A\coloneq0$.
Write $\chi_{[n]}^\star=C_{[n]}\chi_{[n]}$, so
$\overline C_{[n]}C_{[n]}=I$. The monic block $P_{[n]}$ equals
$\chi_{[n]}$ up to lower-degree terms and is orthogonal to $F_{n-1}A$;
its norm is $H_{[n]}\coloneq u(P_{[n]}P_{[n]}^{\T})$, with $u$ acting entrywise.
The corresponding sesquilinear blocks are obtained without changing
the prescribed leading classes:
$\Phi_{[n]}\coloneq P_{[n]}$ and
$\Psi_{[n]}\coloneq C_{[n]}^{-1}P_{[n]}^\star$.
The identity $\overline C_{[n]}C_{[n]}=I$ gives
$\overline{C_{[n]}^{-1}}=C_{[n]}$, so
$\Psi_{[n]}^\star=C_{[n]}P_{[n]}$. Consequently
\[
 u\bigl(\Phi_{[n]}(\Psi_{[j]}^\star)^{\T}\bigr)
 =u(P_{[n]}P_{[j]}^{\T})C_{[j]}^{\T}
 =\delta_{n,j}H_{[n]}C_{[n]}^{\T}.
\]
The construction preserves the leading column in both families,
namely $\chi_{[n]}$, and the
displayed pairing determines their sesquilinear norm. The same
change of convention acts by a fixed basis matrix on each finite
Gram: in a finite basis column $\chi$ with $\chi^\star=C\chi$, put
$G^{\rm bil}\coloneq u(\chi\chi^{\T})$ and
$G^{\rm ses}\coloneq u(\chi(\chi^\star)^{\T})=G^{\rm bil}C^{\T}$.
The fixed right factor cancels in the ratio of the perturbed and
original Gram determinants when both use the same finite space
and basis. It also cancels between the column basis and inverse
Gram in the kernel:
\[
 (\chi^\star(P))^{\T}(G^{\rm ses})^{-1}\chi(Q)
 =\chi(P)^{\T}C^{\T}(C^{\T})^{-1}
                  (G^{\rm bil})^{-1}\chi(Q).
\]
The resulting kernel is therefore the one in the bilinear
convention. For Hermitian data the two polynomial families also
coincide: conjugating the orthogonality conditions and
using $u^\star=u$ shows that $C_{[n]}^{-1}P_{[n]}^\star$
has the same leading class and lower orthogonality as $P_{[n]}$.
Uniqueness gives equality of the two families. Conjugation acts
compatibly on the functional equation as well: conjugating $bv=au$
gives $b^\star v^\star=a^\star u^\star$, where
$v^\star(f)\coloneq\overline{v(f^\star)}$; the conjugate of an evaluation mass
at $c$ is at $c^\star\coloneq(1/\overline c_j)_j$.

The effect of the degree filtration can be calculated on the
subtorus $(s,t)\mapsto(s,t,st)$ of $\mathbb T^3$. Its algebra is
$A=\CC[s^{\pm1},t^{\pm1}]$, while the inherited degree counts
the three coordinates
$s,t,st$; thus $st$ has degree one. Writing $s^at^b=s^{a-c}t^{b-c}(st)^c$ shows that its
inherited degree is the minimum of $|a-c|+|b-c|+|c|$ over
integers $c$. A median of $0,a,b$ minimizes this sum, whose
value is the largest minus the smallest of these three numbers.
Equivalently,
$\ell(a,b)=\max(|a|,|b|,|a-b|)$.
Write $\mathcal H_N\coloneq\{(a,b)\in\mathbb Z^2:\ell(a,b)\le N\}$.
For fixed $b\in\{-N,\ldots,N\}$ the admissible $a$ form the
integer interval
\[
 \max\{-N,b-N\}\le a\le\min\{N,b+N\},
\]
of size $2N+1-|b|$. Summing these sizes gives
$(2N+1)^2-N(N+1)=3N(N+1)+1$. Subtracting the corresponding
count at $N-1$ gives $6N$ boundary monomials for $N\ge1$.
Let $w$ be uniform averaging of both angular variables over $[0,2\pi]$
and $\lambda$ uniform averaging on the circle $s=1$; these are the
normalized Haar functionals.
Put $d\coloneq s-1$, $q_G\coloneq dd^\star=2-s-s^{-1}$,
$u=q_Gw$ and $v_M\coloneq w+M\lambda$, $M\in\CC$.
Then $d v_M d^\star=u$: both denominator factors are nonconstant.
For $M\ge0$ both forms are positive definite, since $q_G>0$ off
the Haar-null circle and $v_M$ contains $w$.
For $b\in\{-N,\ldots,N\}$ set
$I_{N,b}\coloneq\{a\in\mathbb Z:(a,b)\in\mathcal H_N\}$,
$n_{N,b}\coloneq2N+1-|b|$ and
$A_{N,b}(s,t)\coloneq\sum_{a\in I_{N,b}}s^at^b$.
The two denominator factors thus express the scalar division as a
bilateral transformation. The associated kernels follow the
convention established above, with conjugation in the first variable.

This indexing separates the monomials into groups according to
their exponent $b$ of $t$.
Full-torus averaging distinguishes every pair of exponents, whereas
restriction to $s=1$ forgets the exponent of $s$. Consequently
the mass couples all monomials within one group and none between
different groups. Each group therefore contributes a rank-one mass matrix. The
determinant and inverse can then be computed group by group,
as in the following proposition.

\begin{proposition}[A distributed mass on a restricted torus]
\label{prop:torusmass}
Use the Laurent monomial basis $[s^at^b]_{(a,b)\in\mathcal H_N}$
of $F_NA$. Let $G_N(v)$ and $K_N^v$ be its sesquilinear Gram matrix
and kernel as defined above. Then
\begin{equation}\label{eq:torustau}
 \det G_N(v_M)=\prod_{b=-N}^N(1+Mn_{N,b}),\qquad
 \frac{\det G_N(v_M)}{\det G_N(u)}
 =\prod_{b=-N}^N\frac{1+Mn_{N,b}}{n_{N,b}+1}.
\end{equation}
Normality holds precisely when the displayed factors are nonzero.
For all $N\ge0$ simultaneously, this means
$M\notin\{-1/j:j\in\mathbb Z_{>0}\}$.
At every normal degree,
\begin{equation}\label{eq:toruskernel}
 K_N^{v_M}(P,Q)=K_N^w(P,Q)
 -\sum_{b=-N}^N\frac{M\,A_{N,b}^\star(P)A_{N,b}(Q)}{1+Mn_{N,b}}.
\end{equation}
For $M\ne0$ this correction has rank $2N+1$, whereas the kernel multiplied by both rational factors satisfies
\begin{equation}\label{eq:torusbilrank}
 \rank\left(\frac{K_N^{v_M}(P,Q)}{d^\star(P)d(Q)}-K_N^u(P,Q)\right)=4N+2.
\end{equation}
For $N\ge1$, assume also that $G_{N-1}(v_M)$ is invertible.
Then the monic blocks with prescribed leading classes are unique.
For $(a,b)\in\mathcal H_N\setminus\mathcal H_{N-1}$,
the left monic polynomial is
$Q_{N;a,b}=s^at^b-MA_{N-1,b}/(1+Mn_{N-1,b})$ if $|b|<N$,
and $s^at^b$ otherwise; the right one replaces $M$ by $\overline M$.
\end{proposition}

\begin{proof}
Haar orthogonality makes the Gram matrix of $w$ the identity.
On the circle $s=1$, every monomial in the group with exponent
$b$ restricts to $t^b$. Different groups remain orthogonal under
$\lambda$. Hence the block of $G_N(v_M)$ for this group is
$I+M\boldsymbol1\boldsymbol1^{\T}$, where $\boldsymbol1$
has length $n_{N,b}$.
This matrix is the identity on the subspace with coordinate sum
zero, and multiplies $\boldsymbol1$ by $1+Mn_{N,b}$.
Its determinant is that factor, and, when nonzero, its inverse is
\[
 I-\frac{M}{1+Mn_{N,b}}\boldsymbol1\boldsymbol1^{\T}.
\]
Since the groups are mutually orthogonal, their determinants
multiply and their inverse matrices contribute separately to the
kernel expansion. This proves the numerator in \eqref{eq:torustau}
and \eqref{eq:toruskernel}.

For $u=(2-s-s^{-1})w$, the block is tridiagonal with diagonal
$2$ and adjacent diagonals $-1$. Expanding its order-$j$
determinant gives the recurrence $d_j=2d_{j-1}-d_{j-2}$,
with $d_0=1,d_1=2$. Thus $d_j=j+1$, proving the denominator
in \eqref{eq:torustau}. For each $N$ the group sizes are the integers from $N+1$ to
$2N+1$. Their union over $N\ge0$ is all positive integers:
$1$ occurs at $N=0$, and any $j\ge2$ occurs at $N=j-1$
in an extreme group. Thus the excluded parameters are exactly
$M=-1/j$, $j\in\mathbb Z_{>0}$, proving the simultaneous
normality condition.

For $M\ne0$, each kernel correction block has rank one, and
the different groups have disjoint monomial supports. The group decomposition therefore makes their ranks add to $2N+1$. The same inverse matrices determine the monic polynomials:
for a leading monomial with $|b|<N$, its Haar pairing with the
lower group vanishes and its mass pairing is
$M\boldsymbol1^{\T}$. If $n=n_{N-1,b}$, multiplication by
the inverse lower Gram gives the coefficient row
\[
 M\boldsymbol1^{\T}
 \left[I-\frac{M}{1+Mn}\boldsymbol1\boldsymbol1^{\T}\right]
 =\frac{M}{1+Mn}\boldsymbol1^{\T}.
\]
Thus the projection to subtract is exactly
$MA_{N-1,b}/(1+Mn_{N-1,b})$. If $|b|=N$, there is no lower
group to subtract. The right orthogonality equations conjugate
$M$, giving the two monic formulas.

The bilateral comparison involves the rationally transformed
spaces, whose intersection with the original space controls the
additional rank. In this case take $E=X=F_NA$ and conjugate
right spaces. A function in this intersection has $f=dg$ with both $g$
and $dg$ in $F_NA$. Within a fixed group, multiplication by
$d=s-1$ raises the largest $s$ exponent by one unless its
coefficient vanishes. Precisely one condition is imposed, leaving
$n_{N,b}-1$ dimensions. Thus each side has excess $2N+1$.
Every such $f$ vanishes at $s=1$, so the restricted mass is
zero. This vanishing determines the rank of the bilateral kernel
difference through a finite linear-algebra calculation.
Write $D\coloneq\dim F_NA$ and $s\coloneq2N+1$. In minimal
ambient bases the two kernels have difference matrix
$B_RCB_L^{\T}$, where $B_L,B_R$ concatenate the original and
bases multiplied by the rational factors, and
$C=\diag(-G_u^{-1},G_{v_M}^{-1})$ is invertible.
Choose independent bases of the two sums of original and
rationally transformed functions. Each sum has dimension $D+s$, by the
overlap count, so $B_L,B_R$ have sizes $(D+s)\times2D$,
rank $D+s$, and kernels of dimension $D-s$. Each kernel describes the relations between the two basis lists;
a vector in it therefore represents a function common to the
original spaces and their rational images. The pairing
$u^{\T}C^{-1}v$ on these two kernels subtracts their original
pairing from their target pairing. The rational terms cancel;
the mass is zero because both candidates vanish at $s=1$.

Consequently, for every $u\in\ker B_L$ and $v\in\ker B_R$
the cancellation just proved reads $u^{\T}C^{-1}v=0$.
Hence $C^{-1}v$ lies in the annihilator of $\ker B_L$,
which is $\im B_L^{\T}$. The latter image is the annihilator:
it is contained there and both have dimension $D+s$.
Thus $\ker B_R\subset C\im B_L^{\T}$. The first map
$CB_L^{\T}$ has rank $D+s$, and applying $B_R$ kills exactly
its $D-s$ dimensional kernel. The difference matrix has rank
$(D+s)-(D-s)=2s=4N+2$. Independent ambient bases make this
its separation rank as well.
The cancellation used in this rank calculation is justified by
the rational extension
\[
 \widetilde\beta_u(f/d,g/d^\star)=w(fg),
\]
since the factor $dd^\star$ in $u$ cancels the two denominators.
The same rank is visible in the coefficient matrix of each group
of size $n$.
After removal of the common Laurent monomial factor, the remaining
powers are ordered from $0$ to $n-1$. Multiplying the kernel difference by
the two denominators gives the coefficient matrix
\[
 B_n\coloneq E_n(I_n+M\boldsymbol1_n\boldsymbol1_n^{\T})^{-1}E_n^{\T}
               -D_n(D_n^{\T}D_n)^{-1}D_n^{\T},
\]
where $E_n\coloneq[I_n;0]$ and the $j$th column of $D_n$ is
$e_{j+1}-e_j$, $j\in\{1,\ldots,n\}$. The matrix
$D_n^{\T}D_n$ is precisely the tridiagonal original Gram.
Its column image is the hyperplane of coordinate sum zero, so
$D_n(D_n^{\T}D_n)^{-1}D_n^{\T}
=I_{n+1}-\boldsymbol1_{n+1}\boldsymbol1_{n+1}^{\T}/(n+1)$.
Put $e\coloneq[\boldsymbol1_n;0]$ and $k\coloneq e_{n+1}$. Substituting the
already computed target inverse gives
\[
 B_n=\begin{bmatrix}e&k\end{bmatrix}
 \begin{bmatrix}
 \dfrac1{n+1}-\dfrac{M}{1+Mn}&\dfrac1{n+1}\\[3pt]
 \dfrac1{n+1}&-\dfrac n{n+1}
 \end{bmatrix}
 \begin{bmatrix}e^{\T}\\k^{\T}\end{bmatrix}.
\]
In this factorization, the middle determinant is $-1/((n+1)(1+Mn))$, nonzero at
every normal parameter, and $e,k$ are independent. Each group
therefore has rank two. Adding the $2N+1$ independent groups
again gives \eqref{eq:torusbilrank}.
\end{proof}

The distinction between the unrestricted quotient and the
prescribed degree bounds can also be described in this example.
Since $q_G=-(s-1)^2/s$, its quotient is
$A/(q_G)=\CC[t^{\pm1},\epsilon]/(\epsilon^2)$, where
$\epsilon\coloneq[s-1]$ and $s^at^b=(1+a\epsilon)t^b$.
The formula for negative $a$ follows as well from
$(1+\epsilon)^{-1}=1-\epsilon$ modulo $\epsilon^2$.
For $N\ge1$, each group contains at least two consecutive
$a$ values. Their classes $(1+a\epsilon)t^b$ and
$(1+(a+1)\epsilon)t^b$ span $t^b,\epsilon t^b$.
These are the value and first transverse coefficient at $s=1$;
the $2N+1$ different $b$ values consequently give $4N+2$
independent classes.
Restriction to that circle removes the $\epsilon$ terms; the remaining
$2N+1$ modes have nondegenerate Haar pairing.
By contrast, $\dim(F_NA/q_GF_{N-1}A)=6N$ removes only multiples with
bounded multipliers, whereas the image above removes $F_NA\cap q_GA$.
These can differ because multiplication by $q_G$ need not raise inherited
degree, as $q_Gst^2=2st^2-s^2t^2-t^2$ shows.
For instance $st^2$ has degree two, while each term of
$q_Gst^2$ still has degree at most two. Thus this multiple is
in $F_2A\cap q_GA$ although its unique multiplier is outside
$F_1A$. This example identifies the difference between the image in the
unrestricted quotient and the quotient by bounded multiples: the
former also removes multiples whose multipliers exceed the prescribed
bound. The quotient conditions must therefore be accompanied by the
degree constraints.

The filtration also affects Christoffel multiplication. Replacing
the division factors by the multiplier
$q_C\coloneq c-\tfrac12(s+s^{-1}+t+t^{-1})$, $c>2$,
permits an exact degree calculation. Each linear weight $\pm a,\pm b,\pm(a-b)$ has maximum one on the
exponents of $q_C$. The maximum-weight parts multiply to a nonzero
Laurent polynomial, so their weights add. Hence
$\deg_F(q_Cf)=1+\deg_F f$.
Also $q_C\ge c-2>0$. Proposition~\ref{prop:Christoffel}, with ambient
$F_NA$ and target $F_{N-1}A$, therefore gives $\dim F_NA-\dim F_{N-1}A=6N$ independent
conditions on each image and kernel correction rank $6N$. Under the usual Laurent-length filtration on $\mathbb T^2$, the
filtered space has dimension $2N(N+1)+1$ and the corresponding count is
$4N$ for the same $q_C$. Thus the filtration determines the number
of image conditions, which is reflected in the perturbation formula.

\section{Symmetric rational congruences}\label{sec:symmetric}

Symmetry and positivity impose additional structure on the finite
systems obtained above. A \emph{congruence} applies the same rational
change $f\mapsto fS$ to both arguments of a symmetric pairing, so
symmetry first restricts the admissible homogeneous masses. When the
rational pairing and the mass are positive, the part common to the
original and transformed spaces admits an orthogonal decomposition.
Eliminating this common part reduces the matrices needed to construct
the normalized polynomials and their kernels. The reduction depends
on an intersection of polynomial spaces, which can itself be determined
from their coefficient matrices.

Throughout this section the determinants are assumed regular, that is,
non-zero-divisors in the coordinate algebra; the corresponding matrix
changes are therefore invertible over its total quotient ring.
The congruence specialization is given by $q=p=m$ and $R_i=L_i^{\T}$.
In the bivariate case symmetry means
$\beta_\mu(f,g^{\T})=\beta_\mu(g,f^{\T})$, including interchange of
the two arguments; in the balanced case it is $\mu^{\T}=\mu$. Then
$L_1\widehat\mu L_1^{\T}=L_0\mu L_0^{\T}$,
$\widehat\mu=S\widetilde\mu S^{\T}+\nu$, with $S\coloneq L_1^{-1}L_0$.
With a symmetric rational extension, the rational term is symmetric,
so $\widehat\mu$ is symmetric exactly when $\nu$ is symmetric.
The two polynomial spaces may consequently be taken as transposes, while
their common component bounds remain arbitrary. These statements use
bilinear transposition; the assertions concerning positivity below
are made for real data and real vectors.
Massless adjugate congruences on the line appear in
\cite[Section~3.6.2]{Matrix2016}; matrix Geronimus--Sobolev forms
in \cite{SymGeronimus}.
The additional restriction imposed by positivity follows directly from
the homogeneous equation. For a positive semidefinite homogeneous
mass, $\nu(fL_1,(fL_1)^{\T})=0$ for every row $f$.
Its pairing with an arbitrary row $g$ therefore satisfies the
Cauchy--Schwarz bound
\[
 |\nu(fL_1,g^{\T})|^2
 \le\nu(fL_1,(fL_1)^{\T})\nu(g,g^{\T})=0.
\]
Consequently $L_1\nu=0=\nu L_1^{\T}$: vanishing on the product
of the two denominator images extends to vanishing whenever either
argument lies in its denominator image. The mass is therefore
unilateral on both sides. For balanced masses this means that its
class in Theorem~\ref{thm:tormasses} vanishes, in every dimension.

\Needspace{4\baselineskip}
This annihilation property separates symmetric masses from positive
ones. Symmetric masses may contain indefinite derivative terms,
whereas a zero-norm vector for a positive semidefinite form pairs
to zero with every vector. On a curve, applying this property to
the denominator image permits passage to a finite-dimensional
positive quotient.
For balanced pairings, multiplication by the coordinates is
self-adjoint on that quotient, excluding nontrivial nilpotent
multiplication and hence the corresponding derivative directions.
The resulting masses are described by point-evaluation matrices in
the following proposition.

\begin{proposition}[Symmetric and positive homogeneous masses]
\label{prop:symmetricmass}
Write $\operatorname{Sym}_m(A)\coloneq\{H\in\Mat_m(A):H^{\T}=H\}$ and
$d_L=\dim_\CC(A^{1\times m}/A^{1\times m}L_1)$.
For balanced pairings on a curve the symmetric homogeneous masses are the dual of
\begin{equation}\label{eq:symquotient}
 \QQ_{\rm sym}\coloneq\operatorname{Sym}_m(A)\big/
       \{L_1^{\T}HL_1:H\in\operatorname{Sym}_m(A)\}.
\end{equation}
On a smooth affine curve their dimension is $(m+1)d_L$. Locally, if
the denominator has Smith orders $\alpha_1,\ldots,\alpha_m$, the
independent symmetric entry $(i,j)$ has jets of orders less than
$\alpha_i+\alpha_j$, for $i\le j$.

For balanced pairings over a real curve coordinate algebra, the positive
semidefinite homogeneous masses are exactly the finite sums
\begin{equation}\label{eq:positivemass}
 \nu(f,g)=\sum_{P\in C^\circ(\mathbb R)}f(P)C_Pg(P),\qquad
 C_P=C_P^{\T}\ge0,\qquad L_1(P)C_P=0.
\end{equation}
Only zeros of $\det L_1$ contribute. Such a mass has rank
$\sum_P\rank C_P\le d_L$.
\end{proposition}

\begin{proof}
The decomposition of a matrix test $H$ into its symmetric and
skew-symmetric parts is
$(H+H^{\T})/2+(H-H^{\T})/2$. Since a symmetric functional
annihilates the second summand, only the symmetric part contributes
to its homogeneous freedom. This decomposition is preserved by
congruence, as follows from
$(L_1^{\T}HL_1)^{\T}=L_1^{\T}H^{\T}L_1$.
Moreover, if $L_1^{\T}HL_1$ is symmetric, multiplication by
the rational inverses on both sides gives $H=H^{\T}$.
Hence the symmetric part of the full homogeneous quotient is
exactly the symmetric tests modulo the congruence images of
symmetric tests. A functional on this quotient gives a unique
symmetric homogeneous mass, and every such mass factors through
it. This proves \eqref{eq:symquotient}.

The local dimension is determined by a Smith factorization
$L_1=UDV$ at a smooth point, with
$D=\diag(t^{\alpha_i})$. The change $H\mapsto U^{\T}HU$
is invertible on symmetric tests, and the outer factors
$V^{\T},V$ likewise induce an invertible change on the target.
The quotient therefore has the same dimension as that of the
diagonal action $H'\mapsto DH'D$. Since its independent
entry $(i,j)$ is multiplied by $t^{\alpha_i+\alpha_j}$,
exactly that many Taylor coefficients survive in the quotient.
The diagonal test $(i,i)$ contributes $2\alpha_i$ jets.
Each of the $m-1$ off-diagonal pairs containing $i$ contributes
one further $\alpha_i$. Thus each exponent occurs $m+1$ times,
and
\[
 \sum_{i\le j}(\alpha_i+\alpha_j)
   =(m+1)\sum_i\alpha_i.
\]
Summing over all denominator zeros proves the dimension formula.

For a positive mass, the quotient of the real row space by the
radical carries a positive inner product. The Cauchy--Schwarz
argument above places $A^{1\times m}L_1$ in that radical,
so the quotient has dimension at most $d_L$. The classification
can therefore be obtained from operators on a finite-dimensional
inner-product space. To define coordinate multiplication on this
space, the radical must be invariant under scalar multiplication.
Balance gives this invariance: if $f$ pairs to zero with every
$g$, then $af$ does as well, since
$\nu(af,g^{\T})=\nu(f,(ag)^{\T})$.
Thus each real coordinate defines an operator on the quotient,
and the same identity makes that operator self-adjoint.

Commutativity of the coordinate algebra makes these self-adjoint
operators commute. After the first operator is diagonalized, all
the others preserve its eigenspaces; the second can then be
diagonalized on each of them. Repetition through the finite list
of coordinate generators gives a common orthonormal eigenbasis.
The defining equations of the curve hold for the operators and
therefore for their joint eigenvalues, which are real points $P$
of the affine curve. On the joint eigenspace associated with $P$,
a polynomial consequently acts by its value at $P$. If $\Pi_P$
denotes the orthogonal projection onto that eigenspace at $P$,
the expansion of a row as $\sum_i f_i e_i^{\T}$ gives, by
this evaluation property at $P$,
\[
 \Pi_P[f]=\sum_i f_i(P)\Pi_P[e_i^{\T}].
\]
With $(C_P)_{i,j}$ defined as the inner product of
$\Pi_P[e_i^{\T}]$ and $\Pi_P[e_j^{\T}]$, the contribution of
each joint eigenspace is represented by a positive semidefinite
Gram matrix. Taking inner products in the displayed expansion and
summing over the mutually orthogonal eigenspaces then gives
\eqref{eq:positivemass}.

The homogeneity condition can now be imposed separately at each
joint eigenvalue. Their distinct maximal ideals are comaximal,
so the Chinese remainder theorem supplies a polynomial equal to
one at any chosen point and zero at the others. Its insertion
in the homogeneity identity isolates that point; varying the
constant row factors yields $L_1(P)C_PL_1(P)^{\T}=0$.
The positive factorization $C_P=C_P^{1/2}C_P^{1/2}$ then forces
$L_1(P)C_P^{1/2}=0$, hence $L_1(P)C_P=0$.
Only determinant zeros can therefore contribute. Conversely,
positive matrices with this property give a positive homogeneous
mass by direct substitution in \eqref{eq:positivemass}.
Polynomial evaluation separates the finitely many distinct support
points, so their contributions have additive ranks. Their sum is
bounded by the dimension of the quotient through which the mass
factors: Cauchy--Schwarz places $A^{1\times m}L_1$ in its
radical, and the quotient by this subspace has dimension $d_L$.
This proves the stated rank bound.
\end{proof}

For $L_1=\diag(z,1)$, define
$\mathsf c_1(f)\coloneq[z]f$. The symmetric mass
$\nu\coloneq E_{1,1}\mathsf c_1$ is homogeneous because its
only transformed entry applies $\mathsf c_1$ to a multiple of
$z^2$. On $e_1^{\T},ze_1^{\T}$ its Gram is
\[
 \begin{bmatrix}0&1\\1&0\end{bmatrix},
\]
whose eigenvalues $1$ and $-1$ show that the mass is indefinite.
The sign of the homogeneous remainder in a full pairing depends on
the particular solution chosen for its decomposition. For the positive form
$\Theta(f,g^{\T})\coloneq\int_1^2f(t)g(t)^{\T}\dd t$,
the choices $L_0=I_2$, $\mu=L_1\Theta L_1^{\T}$ give two
solutions: both $\Theta$ and $\Theta-\nu$ satisfy
$L_1\widehat\mu L_1^{\T}=\mu$.
Choosing $\Theta-\nu$ as the particular solution leaves the
indefinite homogeneous term $\nu$ in the positive target
$\Theta=(\Theta-\nu)+\nu$. Thus positivity of the full
pairing alone does not determine the sign of a chosen homogeneous
remainder; the positive reduction below assumes positivity of
its two terms separately.

For complex rows in real coordinates, the corresponding Hermitian
form is
$\langle f,h\rangle_\mu\coloneq\widetilde\beta_\mu(f,\overline h^{\T})$,
where the bar conjugates coefficients. When $S$ is real rational
and the rational extension is Hermitian, separate positivity is
likewise imposed on the rational space used in the congruence:
if $\langle fS,fS\rangle_\mu>0$ for $f\ne0$ and
$\nu(f,\overline f^{\T})\ge0$, then
$\langle f,f\rangle_{\widehat\mu}
=\langle fS,fS\rangle_\mu+\nu(f,\overline f^{\T})>0$.

\subsection{Positive congruences: smaller matrices and the orthogonal polynomials}
\label{sec:positiveprojection}

Separate positivity of the rational pairing and the homogeneous mass
allows the congruence to be reduced by orthogonal projection. Since
the same rational change acts on both arguments, the common part of
the original and transformed spaces can be removed together from the
two arguments. One positive matrix determines this projection; a
second is the Gram matrix on the remaining polynomial rows. The
latter matrix determines both the new orthogonal polynomials and the
reproducing kernel. Its entries depend only on the original
polynomials, the rational moments and the specified mass.

The construction is stated for real rows and symmetric bilinear
forms; replacing transposes by adjoints gives the Hermitian version.
Let $X$ be the space of lower polynomial rows against which
orthogonality is required, and let $X^+\supset X$ contain the
prescribed leading rows. The comparison space $E$ is an original
polynomial space of the same dimension
$N\coloneq\dim X=\dim E$; the component bounds in $E$ and $X$
remain arbitrary. The transformed form is written as
\[
 \langle f,g\rangle_{\widehat\mu}
 \coloneq\langle fS,gS\rangle_{\widetilde\mu}+\nu(f,g),
 \qquad S=L_1^{-1}L_0.
\]
Assume that the rational pairing is positive definite on
$E+X^+S$, that $S$ is invertible over the total quotient ring,
and that $\nu$ is positive semidefinite on $X^+$. The homogeneous
condition is imposed whenever this is a solution of the bilateral
functional equation. The projection argument itself applies to any
positive additive form.

\subsubsection{The common space and the contribution of the mass}

The intersection of the two spaces is measured before and after
imposing vanishing of the mass. These two requirements define
\begin{equation}\label{eq:positivecommon}
 \begin{aligned}
 Z&\coloneq\{f\in X:fS\in E\},&
 s&\coloneq N-\dim Z,\\
 Z_0&\coloneq\{f\in Z:\nu(f,f)=0\},&
 \eta&\coloneq\rank(\nu|_{Z\times Z}),\\
 F&\coloneq Z_0S\subset E,&
 t&\coloneq N-\dim Z_0.
 \end{aligned}
\end{equation}
Here $Z$ consists of target rows whose rational images already belong
to the original space. The number $s$ counts the dimensions lost from
$X$ when that requirement is imposed. Within $Z$, the nullspace of
the mass is $Z_0$, and $\eta$ counts the additional directions removed
by requiring the mass to vanish. Only the mass on this overlap enters $\eta$;
it can have smaller rank than the mass on all of $X$. Positivity
implies that $\nu(f,f)=0$ makes $f$ pair to zero with every element
of $X^+$, by Cauchy--Schwarz. Thus $Z_0$ is a vector space, and
the mass vanishes whenever one argument belongs to $Z_0$.

An orthonormal basis $p_1,\ldots,p_N$ of $E$ for the original
pairing gives the $N\times m$ polynomial matrix
$\boldsymbol p\coloneq[p_i]_{i\in\{1,\ldots,N\}}$ by stacking
its rows. Such a basis is obtained from the original orthogonal
polynomials and their norms, and a row in $E$ then has the
coordinate representation $c^{\T}\boldsymbol p$. In these
coordinates the common subspace is the kernel of a full-row-rank
$t\times N$ matrix $\Lambda$, characterized by
\[
 c^{\T}\boldsymbol p\in F\quad\hbox{if and only if}\quad\Lambda c=0.
\]
The orthogonality of the chosen coordinates expresses the projection
onto this kernel through the first reduced matrix:
\begin{equation}\label{eq:positivefirstmatrix}
 A_F\coloneq\Lambda\Lambda^{\T},\qquad
 H_F\coloneq I_N-\Lambda^{\T}A_F^{-1}\Lambda.
\end{equation}
For any rational row $u$ in the spaces under consideration, set
\[
 c(u)\coloneq[\langle u,p_i\rangle_{\widetilde\mu}]_{i\in\{1,\ldots,N\}},
 \qquad \Pi_Fu\coloneq c(u)^{\T}H_F\boldsymbol p.
\]
The vector $c(u)$ is an $N\times1$ column. It gives the
coordinates of the original projection onto $E$, since
$\langle u-c(u)^{\T}\boldsymbol p,p_j\rangle=0$ for every
$j$. The symmetric $N\times N$ matrix $H_F$ then projects
those coordinates onto $\ker\Lambda$, which represents $F$.
Thus $\Pi_Fu$ is the original orthogonal projection onto $F$.

For a complement $f_1,\ldots,f_t\in X$ to a basis of $Z_0$
in $X$, subtraction of this projection gives the rows
\begin{equation}\label{eq:positiveresidualrows}
 \rho_j\coloneq f_jS-\Pi_F(f_jS),\qquad
 w_j\coloneq\rho_jS^{-1}
       =f_j-\Pi_F(f_jS)S^{-1},\qquad j\in\{1,\ldots,t\}.
\end{equation}
The row $w_j$ remains polynomial and belongs to $X$, because
$\Pi_F(f_jS)\in F=Z_0S$ and its product with $S^{-1}$ therefore
lies in $Z_0$. With the notation
$\boldsymbol w\coloneq[w_j]_{j\in\{1,\ldots,t\}}$ and
$\boldsymbol\rho\coloneq[\rho_j]_{j\in\{1,\ldots,t\}}$,
both $\boldsymbol w$ and $\boldsymbol\rho$ have $t$ rows and
$m$ columns. They represent the same complementary directions
before and after the rational change $S$, with
$\boldsymbol wS=\boldsymbol\rho$. Their Gram matrix for the
transformed pairing is the second reduced matrix,
\begin{equation}\label{eq:positivesecondmatrix}
 T_\nu\coloneq
 [\langle\rho_i,\rho_j\rangle_{\widetilde\mu}+\nu(f_i,f_j)]_{i,j\in\{1,\ldots,t\}}.
\end{equation}
Both reduced matrices are thus determined before the perturbed
orthogonal polynomials are constructed. The rational pairings
appearing in their entries are precisely the second-kind data of
Subsubsection~\ref{sec:secondkind}; for measure pairings, these are
Cauchy data in one variable and Fantappi\`e data in several variables.

\begin{theorem}[Positive reduction and normalized polynomials]
\label{thm:positiveprojection}
Under the preceding assumptions, both $A_F$ and $T_\nu$ are positive
definite, and their common order is
\begin{equation}\label{eq:positiveorder}
 t=s+\eta.
\end{equation}
For a prescribed row $v\in X^+$, put
\[
 v_0\coloneq v-\Pi_F(vS)S^{-1},\qquad
 b(v)\coloneq
 [\langle vS,\rho_j\rangle_{\widetilde\mu}+\nu(v,f_j)]_{j\in\{1,\ldots,t\}}.
\]
Then the unique row with the class of $v$ modulo $X$ and orthogonal
to $X$ for $\widehat\mu$ is
\begin{equation}\label{eq:positivepolynomial}
 \widehat p=v_0-b(v)^{\T}T_\nu^{-1}\boldsymbol w.
\end{equation}
In particular, every scalar component has the bordered formula
\begin{equation}\label{eq:positiveborder}
 \widehat p_a=
 \frac{\det\begin{bmatrix}
 T_\nu&\boldsymbol w_{:,a}\\ b(v)^{\T}&(v_0)_a
 \end{bmatrix}}{\det T_\nu},
 \qquad a\in\{1,\ldots,m\}.
\end{equation}
Thus the denominator has order $t$ and is positive, and each bordered
numerator has order $t+1$. The preliminary projection uses the
additional order-$t$ matrix $A_F$.

Let $K_F(P,Q)\coloneq\boldsymbol p(P)^{\T}H_F\boldsymbol p(Q)$.
The kernel on $X$ for the transformed pairing is
\begin{equation}\label{eq:positivekernel}
 \widehat K_X(P,Q)=S(P)^{-\T}
 \bigl[K_F(P,Q)+\boldsymbol\rho(P)^{\T}T_\nu^{-1}
                         \boldsymbol\rho(Q)\bigr]S(Q)^{-1}.
\end{equation}
The apparent rational singularities cancel: this is a polynomial
kernel, interpreted by the polynomial rows in
\eqref{eq:positiveresidualrows} wherever an evaluated inverse is undefined.
\end{theorem}

\begin{proof}
Positivity identifies the kernel of the restricted mass on $Z$
with $Z_0$, so its rank is $\eta=\dim Z-\dim Z_0$.
The dimension definitions in \eqref{eq:positivecommon} then give
$t=s+\eta$. Since the change $f\mapsto fS$ is injective,
the common image has dimension $\dim F=N-t$.

The constraint matrix has full row rank, which determines the
positivity of the first reduced matrix. For every nonzero column
$a\in\mathbb R^t$, the equality
$a^{\T}A_Fa=\|\Lambda^{\T}a\|^2>0$ follows from injectivity
of $\Lambda^{\T}$. Hence $A_F$ is positive definite, and the
definition of $H_F$ gives the identities
\[
 H_F^{\T}=H_F,\qquad
 \Lambda H_F=0,\qquad
 H_Fc=c\quad(c\in\ker\Lambda)
\]
which follow by substituting $A_F=\Lambda\Lambda^{\T}$.
These identities show that $H_F$ is the Euclidean orthogonal projection
onto $\ker\Lambda$ and in particular $H_F^2=H_F$.
For $v=d^{\T}\boldsymbol p\in F$ one then has
$\langle u-\Pi_Fu,v\rangle
=c(u)^{\T}(I-H_F)d=0$.
This verifies the projection formula for every rational row $u$
in the stated finite space.

The form induced by $\widehat\mu$ on $X^+$ is positive definite:
the first term is positive on every nonzero row, and the second
term is nonnegative. On $Z_0$, the map $z\mapsto zS$ preserves
the inner product, since the mass vanishes there. For $z\in Z_0$,
\[
 \langle w_j,z\rangle_{\widehat\mu}
 =\langle\rho_j,zS\rangle_{\widetilde\mu}+\nu(w_j,z)=0.
\]
The two terms vanish for separate reasons: the definition of
$\Pi_F$ gives orthogonality of the rational term, and the mass
radical property gives vanishing of the mass term. Moreover,
$w_j$ and $f_j$ have the same classes modulo $Z_0$, so they
still span the complementary directions. It follows that $X$
is the orthogonal direct sum of $Z_0$ and the span of
$w_1,\ldots,w_t$. The latter rows are independent: a relation
$\sum_j a_jw_j=0$ would give $\sum_j a_jf_j\in Z_0$, and
the chosen complement forces $a_j=0$ for every $j$.
Their Gram matrix is determined without any further projection:
subtracting elements of $Z_0$ changes no mass pairing, while
$w_jS=\rho_j$, giving
\[
 \langle w_i,w_j\rangle_{\widehat\mu}
 =\langle\rho_i,\rho_j\rangle_{\widetilde\mu}
                         +\nu(f_i,f_j)=(T_\nu)_{i,j}.
\]
Independence and positivity now imply positivity of the second
reduced matrix. For a nonzero coefficient vector $a$, the row
$\sum_j a_jw_j$ is nonzero, so its squared norm for the
transformed pairing is $a^{\T}T_\nu a>0$. This proves that
$T_\nu$ is positive definite.

The orthogonal decomposition of the lower space also determines the
correction of a prescribed leading row. The subtraction defining
$v_0$ belongs to $Z_0$ and removes the pairings with that summand.
Its image $\Pi_F(vS)\in F$ is orthogonal to $\rho_j$, while
its pullback belongs to $Z_0$ and has zero mass pairing. The
pairings that remain after this first subtraction are therefore
\[
 \langle v_0,w_j\rangle_{\widehat\mu}
 =\langle vS,\rho_j\rangle_{\widetilde\mu}
                              +\nu(v,f_j)=b(v)_j.
\]
For a coefficient row $a^{\T}$, the polynomial
$v_0-a^{\T}\boldsymbol w$ is orthogonal to all $w_j$ exactly
when $a^{\T}T_\nu=b(v)^{\T}$. Its unique solution is
$a^{\T}=b(v)^{\T}T_\nu^{-1}$, proving
\eqref{eq:positivepolynomial}. Both corrections
belong to $X$, so the prescribed leading class is preserved.
Two rows with this class and orthogonal to $X$ differ by an
element of $X$ with zero norm; positivity makes their difference
zero. The component formula can then be written as a bordered
determinant: left multiplication of the matrix in
\eqref{eq:positiveborder} by
$[I_t\ 0;\ -b(v)^{\T}T_\nu^{-1}\ 1]$ preserves its determinant,
since the multiplier has determinant one. Its lower-left block
vanishes, and its lower-right entry is
$(v_0)_a-b(v)^{\T}T_\nu^{-1}\boldsymbol w_{:,a}$.
Taking the determinant of this block triangular matrix gives the
bordered formula, with the stated positive denominator and sizes.

The same orthogonal decomposition yields the kernel. An orthonormal
basis of $F$, pulled back by $S^{-1}$ to $Z_0$, combines with
the $w_j$ to give a basis with Gram matrix
$\diag(I_{N-t},T_\nu)$. The contribution of its first summand is
$S(P)^{-\T}K_F(P,Q)S(Q)^{-1}$, while inversion of the second
Gram block gives the contribution
$\boldsymbol w(P)^{\T}T_\nu^{-1}\boldsymbol w(Q)$.
Substituting $\boldsymbol w=\boldsymbol\rho S^{-1}$ and
factoring the two outer inverses gives
\eqref{eq:positivekernel}. Every basis row before and after this
orthogonal decomposition belongs to $X$, which proves the
polynomial cancellation assertion. Empty matrices have determinant
one; the formulas also include $t=0$.
\end{proof}

Applying \eqref{eq:positivepolynomial} to a basis of $X^+/X$ gives
the whole normalized block, using the same $T_\nu$ for every
leading row. A nonzero constant linear combination of these rows
has a nonzero class in $X^+/X$, hence is a nonzero polynomial
row. Its positive squared norm proves that the block norm matrix
is positive definite. Iterating along
nested component bounds constructs the complete orthogonal family.

The decomposition used in this proof depends on positivity of the
restricted form on the common subspace. Nondegeneracy of the full
symmetric form alone does not ensure it: for the form with matrix
$\diag(1,-1)$, the restriction to the line spanned by $[1\ 1]$
is zero, so that line cannot be a summand in an orthogonal direct
sum. Symmetric indefinite data are still governed by the determinant
of the general reduced system, but a reduction by the common space
requires a separate nondegeneracy assumption on its restriction.

\subsubsection{Computing the order from polynomial coefficients}

The order of the projection construction can be determined from the
polynomial spaces and factors, without first finding the new
orthogonal polynomials. With $\delta\coloneq\det L_1$ and
$B\coloneq\operatorname{adj}(L_1)L_0$, in the displayed order,
clearing the common denominator gives
\begin{equation}\label{eq:positivecoefficientcount}
 s=\rank[\,\delta E\ \ XB\,]-N.
\end{equation}
As in \eqref{eq:orderrankmatrices}, the displayed lists are
coefficient columns in a polynomial basis on the support. Their
rank computes the dimension of the sum of spaces because
multiplication by the regular scalar $\delta$ is injective on
rational rows and maps $E+XS$ onto $\delta E+XB$. The
coefficient matrix consequently has rank $\dim(E+XS)$.
The change $f\mapsto fS$ identifies $Z$ with $E\cap XS$;
since both $E$ and $XS$ have dimension $N$, the dimension
formula for a sum of subspaces gives
\[
 \dim(E+XS)=2N-\dim(E\cap XS)=2N-\dim Z
\]
and hence \eqref{eq:positivecoefficientcount}.
The same coefficient calculation recovers the common space itself:
the solutions of $fB=\delta e$ with $f\in X$, $e\in E$ give
a basis of $Z$. The matrix of the mass in that basis determines
its rank $\eta$ and nullspace $Z_0$, and consequently determines
both $F$ and $\Lambda$. Since this calculation uses the actual
coefficient spaces on the support, it also applies to unequal
component bounds, singular supports and noncommuting matrix factors.

Without a mass, $t=s$. With a positive mass of restricted rank $h$
on $X$, one has $t=s+\eta\le s+h$, and always $t\le N$.
For $S=I$ and $E=X$, this recovers the ordinary additive update:
$s=0$, $Z=X$ and $t=h$. A cancellation between rational factors
therefore does not discard a remaining mass. The two order-$t$
matrices above describe the orthogonal projection construction.
As in Subsection~\ref{sec:orders}, their size is determined by
the particular spaces eliminated in this representation.

\subsubsection{Exact degree counts under a coprimality condition}

For equal filtered spaces $E=X=(F_{n-1}A)^{1\times m}$, the
coefficient count reduces to a difference of Hilbert functions under
coprimality and degree assumptions. Here
$H_A(j)=\dim F_jA$ is as in Subsection~\ref{sec:orders}.
The factors $L_0,L_1$ are assumed to have filtered degrees
$a_0,a_1$ and to satisfy
\begin{equation}\label{eq:positivebezout}
 L_0U_0+L_1U_1=I_m
\end{equation}
for polynomial matrices $U_0,U_1$. This B\'ezout identity is
the coprimality condition used here and is compatible with
noncommuting $L_0,L_1$. The additional degree assumption is
\begin{equation}\label{eq:positivedegreestrict}
 \deg(hL_i)=\deg h+a_i
 \quad\hbox{for every nonzero polynomial row }h,
 \qquad i\in\{0,1\}.
\end{equation}
On the line, invertible leading coefficient matrices ensure this
property. On an algebraic support, multiplication can interact with
the defining relations, so the property is an independent assumption
on the inherited filtration rather than a consequence of matrix size.

\begin{corollary}[Orders for positive coprime congruences]
\label{cor:positiveorders}
Assume \eqref{eq:positivebezout}--\eqref{eq:positivedegreestrict}
and let $M\coloneq\max\{a_0,a_1\}$. For every positive homogeneous
mass the order in Theorem~\ref{thm:positiveprojection} is exactly
\begin{equation}\label{eq:positivehilbertorder}
 {\ t_n=m\bigl[H_A(n-1)-H_A(n-M-1)\bigr].\ }
\end{equation}
In particular, the mass does not increase this order. The same
formula, with $m=1$, holds for coprime scalar polynomials $P,Q$
in $A=\mathbb R[x_1,\ldots,x_D]$, with $S=P/Q$, even when a
polynomial B\'ezout identity is unavailable.
\end{corollary}

\begin{proof}
An element $f\in Z$ determines the rational row
$h=fL_1^{-1}$, whose products $hL_1=f$ and $hL_0=fS$ are
both polynomial. The B\'ezout identity then gives
$h=(hL_0)U_0+(hL_1)U_1$, proving that $h$ is itself
polynomial. The degree assumptions can therefore be
applied to it: for nonzero $h$, the bounds on $hL_1$ and
$hL_0$ are
$\deg h+a_1\le n-1$ and $\deg h+a_0\le n-1$.
Their conjunction is $\deg h\le n-M-1$.
Including the zero row, the two requirements are consequently
equivalent to $h\in(F_{n-M-1}A)^{1\times m}$, where a
negative filtered bound means the zero space. Conversely, every row in this
space gives $f=hL_1\in Z$. Injectivity of multiplication by $L_1$
therefore yields
\[
 Z=(F_{n-M-1}A)^{1\times m}L_1,
 \qquad \dim Z=mH_A(n-M-1).
\]
A positive homogeneous mass annihilates the entire denominator
image in either argument, as proved at the beginning of this
section. Hence $\nu|_{Z\times Z}=0$, so $\eta=0$ and $t=s$.
Subtracting dimensions proves \eqref{eq:positivehilbertorder}.

In the scalar case, the factorization of $Q$ into irreducibles
replaces the B\'ezout argument. Each irreducible factor occurring
in $Q$ has exponent zero in the coprime polynomial $P$.
Consequently $Q\mid fP$ forces each such power to divide
$f$, proving $Q\mid f$. With $f=Qh$, the remaining argument
again concerns the degrees of polynomial products.
The nonzero highest homogeneous parts of two polynomials have
nonzero product in $\mathbb R[x_1,\ldots,x_D]$; hence total
degrees add. The bounds on $Qh$ and $Ph$ give the same space
of $h$ as above. Finally positivity and homogeneity annihilate
$QA$, so the mass still vanishes on $Z$. This proves the scalar
claim without using a B\'ezout identity.
\end{proof}

On the line the order is $m\min\{n,M\}$, hence $mM$ for $n\ge M$.
For two independent variables it is
\[
 t_n=\frac{mM(2n-M+1)}2,\qquad n\ge M.
\]
In $D$ independent variables the exact expression is
\[
 t_n=m\left[\binom{n+D-1}{D}-\binom{n-M+D-1}{D}\right],
\]
where the second term is zero when $n-M-1<0$.
On a curve of degree $d$, once both arguments of $H_A$ are in
its linear range, the order is exactly $mdM$. On a support of
dimension $r\ge2$ and degree $d$,
\[
 t_n=\frac{mdM}{(r-1)!}n^{r-1}+\mathrm{O}(n^{r-2})
 \qquad(n\longrightarrow\infty).
\]
The first three counts follow by substituting
$H_A(j)=j+1$, $H_A(j)=(j+1)(j+2)/2$, and
$H_A(j)=\binom{j+D}{D}$, respectively. For a curve in its
linear Hilbert range, $H_A(j)=dj+c$, whose difference across
$M$ consecutive degrees is $dM$. In dimension $r\ge2$, the
leading Hilbert term is $dj^r/r!$; subtracting its translate
by $M$ gives $dMn^{r-1}/(r-1)!$, and the remaining terms have
order at most $n^{r-2}$. Thus all displayed counts retain the
assumptions of the corollary.
For arbitrary matrix factors or component bounds,
\eqref{eq:positivecoefficientcount} and $t=s+\eta$ give the
exact order instead.

\subsubsection{A noncommutative example with a positive derivative mass at a node}

The degree and coprimality assumptions above admit an explicit
realization on a nodal curve, where a derivative mass can also be
included. The coordinate algebra and parametrization are
\[
 A=\mathbb R[x,y]/(y^2-x^2(x+1)),\qquad
 P(t)=(t^2-1,t(t^2-1)),\qquad O=(0,0).
\]
Let
\begin{equation}\label{eq:positivenodalfactors}
 C=\begin{bmatrix}0&1\\0&1/4\end{bmatrix},\qquad
 D=\begin{bmatrix}2&0\\0&3\end{bmatrix},\qquad
 L_1=yI_2-C,\qquad L_0=yI_2-D.
\end{equation}
Direct multiplication gives
$[C,D]=E_{1,2}$ and $[L_0,L_1]=[D,C]=-E_{1,2}$.
Thus the two polynomial factors do not commute.
Put $\delta\coloneq y(y-1/4)$ and define the original balanced pairing by
\[
 \langle f,g\rangle_\mu
 \coloneq\int_2^3\delta(P(t))^2 f(P(t))g(P(t))^{\T}\,\dd t.
\]
Along the integration arc, the polynomial $y$ ranges from $6$
to $24$, keeping the denominator zeros away from the integration
support. The parametrization identifies the coordinate algebra
with a subalgebra of $\mathbb R[t]$, so a nonzero polynomial
row cannot vanish on the whole arc. The positive density therefore
gives a positive definite pairing. The same reasoning, after
clearing denominators, proves positivity of its integral extension
on all the rational spaces used in this example.

For $\alpha,\beta\ge0$, choose the homogeneous mass
\begin{equation}\label{eq:positivenodalmass}
 \nu_{\alpha,\beta}(f,g)
 \coloneq\alpha f_1(O)g_1(O)
       +\beta\,\partial_xf_1(O)\partial_xg_1(O).
\end{equation}
Both tests descend to the coordinate algebra. The defining polynomial
and its $x$ derivative vanish at $O$, so changing a representative
by a multiple of the relation changes neither its value nor the
derivative used in the mass. Nonnegative coefficients then make
this a positive form built from evaluation and a first derivative
at the singular point.

\begin{proposition}[Positive orthogonal families with component degree bounds on the nodal curve]
\label{prop:positivenodalfamily}
For the preceding data,
\[
 \langle f,g\rangle_{\widehat\mu}
 =\langle fS,gS\rangle_{\widetilde\mu}+\nu_{\alpha,\beta}(f,g),
 \qquad S=L_1^{-1}L_0,
\]
is a positive definite solution of the bilateral congruence equation.
It has uniquely normalized orthogonal blocks for every nested
choice of component bounds.

For equal bounds $E=X=(F_{n-1}A)^{1\times2}$, the two matrices
in Theorem~\ref{thm:positiveprojection} have orders
\begin{equation}\label{eq:positivenodalorders}
 t_1=2,\qquad t_2=4,\qquad t_n=6\quad(n\ge3).
\end{equation}
The orders are independent of $\alpha,\beta\ge0$.
For arbitrary positive bounds $u,v$ on the two components, put
$E=X=F_{u-1}A\oplus F_{v-1}A$. The exact order is
\begin{equation}\label{eq:positivenodalmixed}
 t_{u,v}=H_A(u-1)+H_A(v-1)
          -H_A(\min\{u-2,v-1\})-H_A(v-2),
\end{equation}
where $H_A(0)=1$, $H_A(j)=3j$ for $j\ge1$, and $H_A(j)=0$
for $j<0$. In particular, the component degree bounds $(u,v)=(n+3,n)$ give
order $12$ for $n\ge3$.
\end{proposition}

\begin{proof}
Homogeneity follows from the first component of $fL_1$, namely
$yf_1$. Its value at $O$ is zero, and
$\partial_x(yf_1)(O)=y(O)\partial_xf_1(O)=0$, since the
ambient coordinates $x,y$ are differentiated independently.
Both calculations descend to the coordinate algebra because the
defining relation and its derivative vanish at $O$. The mass
therefore annihilates the denominator image in either argument.
Together with $L_1S=L_0$, which gives the functional equation
for the rational part, this proves the congruence equation for
the sum. Its rational term is positive definite and its mass is
nonnegative, so the full form is positive definite. Existence and
uniqueness of the normalized blocks follow at every finite level.

The matrix $C-D$ is invertible, and
\[
 L_0(C-D)^{-1}-L_1(C-D)^{-1}=I_2.
\]
This verifies the B\'ezout hypothesis for the noncommuting factors.
The remaining hypothesis concerns degree, which is determined by
normal forms for the coordinate algebra. Writing the relation as
$x^3=y^2-x^2$, polynomial division by the monic cubic in $x$
gives the unique representative $a_0(y)+xa_1(y)+x^2a_2(y)$. Each reduction
replaces a degree-three factor by terms of degree two, so the
filtered basis consists precisely of
$x^iy^j$ with $i\in\{0,1,2\}$, $j\in\NN$ and $i+j\le k$.
For $k\ge2$ their count is $(k+1)+k+(k-1)=3k$;
for $k=0,1$ it is $1,3$. Multiplication by $y$ leaves these
normal forms unchanged except for increasing every $j$ by one,
so it increases the filtered degree by exactly one;
the leading part of each matrix $L_i$ is $yI_2$, so the matrix
degree property also holds. Corollary~\ref{cor:positiveorders}
with $m=2$ and $M=1$ gives \eqref{eq:positivenodalorders}.

The B\'ezout argument also identifies the common space for unequal
bounds: any row in that space has the form $f=hL_1$.
The diagonal numerator makes the condition $hL_0\in E$ equivalent to
$h_1\in F_{u-2}A$, $h_2\in F_{v-2}A$. The second component
of $hL_1$ is $-h_1+(y-1/4)h_2$. The second term already has
degree at most $v-1$, so requiring the same bound for their
sum is equivalent to $h_1\in F_{v-1}A$. Therefore
\[
 h_1\in F_{\min\{u-2,v-1\}}A,\qquad h_2\in F_{v-2}A.
\]
Counting these independent choices proves
\eqref{eq:positivenodalmixed}. Every such $f$ is in the
denominator image, so the mass vanishes on the overlap and
does not change the count.
\end{proof}

For $\beta>0$ this derivative mass is not balanced. Indeed, with
$e_1^{\T}=[1\ 0]$,
\[
 \nu_{\alpha,\beta}(x^2e_1^{\T},e_1^{\T})=0,
 \qquad
 \nu_{\alpha,\beta}(xe_1^{\T},xe_1^{\T})=\beta.
\]
The rational part of the transformed pairing is balanced, so the
same calculation shows that the full form is not balanced when
$\beta>0$. Its positive Sobolev term therefore lies outside the
point-evaluation classification in Proposition~\ref{prop:symmetricmass},
which assumes balance as well as positivity. Orthogonal projection
uses only positivity and the specified finite spaces, so
Theorem~\ref{thm:positiveprojection} applies for every
$\alpha,\beta\ge0$. The coordinate identities in
Theorem~\ref{thm:darboux} require both pairings to be balanced:
within this nodal family the evaluation case $\beta=0$ satisfies
that requirement, whereas the derivative case $\beta>0$ is governed
by the positive polynomial and kernel formulas proved here.

The order calculation determines the size of the formulas for every
normalized block and its kernel. From degree three onward, all
leading rows at a fixed level with equal component bounds use the
same $6\times6$ matrix $T_\nu$ and scalar bordered determinants
of order seven. Its dependence on the mass parameters can be
separated from the projection. Since the mass vanishes on $Z$,
one has $Z_0=Z$, and the spaces and matrices $F$ and $A_F$
are fixed. For a fixed complementary list $f_1,\ldots,f_t$, let
$a_j\coloneq f_{j,1}(O)$ and $b_j\coloneq\partial_xf_{j,1}(O)$. Then
\[
 T_\nu=T_0+\alpha aa^{\T}+\beta bb^{\T},
\]
where $T_0>0$ is the massless matrix. The two additions are
positive semidefinite and each has rank at most one; at levels
where the respective test is nonzero its rank is one.

\subsubsection{A consequence for extremal norms}

The positive dependence on the mass also controls the extremal norms
associated with the kernel. With the rational part and a finite
polynomial space $X$ fixed, increasing the positive additive mass
increases the Gram matrix in the positive semidefinite order.
The effect on the kernel follows from the reverse order of inverses:
for $G>0$ and $M\ge0$, one has
\[
 (G+M)^{-1}
 =G^{-1/2}(I+G^{-1/2}MG^{-1/2})^{-1}G^{-1/2}
 \le G^{-1}.
\]
Evaluation at a real point against the same polynomial basis therefore
decreases the kernel on the diagonal. For scalar real polynomials,
when evaluation at a real point $P$ is nonzero on $X$, the extremal identity is
\[
 \lambda_X(P)\coloneq
 \min\{\langle f,f\rangle:f\in X,\ f(P)=1\}
 =\frac1{K_X(P,P)}.
\]
The reproducing identity gives
$1=f(P)=\langle f,K_X(\cdot,P)\rangle$ for an admissible $f$.
Cauchy--Schwarz and
$\|K_X(\cdot,P)\|^2=K_X(P,P)$ give
$\|f\|^2\ge1/K_X(P,P)$. The assumption on evaluation makes
$K_X(P,P)>0$, and the polynomial
$K_X(\cdot,P)/K_X(P,P)$ attains equality. Thus increasing a positive
mass increases this minimum. The reduced kernel formula computes
the resulting change in the extremal norm as well as the
orthogonal polynomials.

\subsection{Matrix orthogonal polynomials and the Dur\'an--Gr\"unbaum weights}
\label{sec:matrixorthogonal}

When the component bounds coincide on the line, the normalized rows
of the preceding construction assemble into ordinary matrix
orthogonal polynomials. For $p=q=m$ and a real symmetric positive
matrix measure $\dd\Sigma$ on $\mathbb R$, the corresponding
matrix inner product is
\[
 \langle P,Q\rangle_\Sigma\coloneq
       \int_{\mathbb R}P(x)\,\dd\Sigma(x)\,Q(x)^{\T}.
\]
Its $(a,b)$ entry pairs row $a$ of $P$ with row $b$ of $Q$.
Assume that all polynomial moments exist and that every nonzero
polynomial row has positive squared norm. The monic matrix
orthogonal polynomials $P_n$ then have leading coefficient $I_m$;
write $H_n\coloneq\langle P_n,P_n\rangle_\Sigma>0$.

Let $L_0,L_1$ be regular matrix polynomials and $S=L_1^{-1}L_0$.
Assume that $S$ is defined $\Sigma$-almost everywhere, that the
required rational moments exist, and that the rational inner product
satisfies the positivity hypotheses of Theorem~\ref{thm:positiveprojection}
at every level. With positive semidefinite matrices $C_\xi$ at
finitely many real points, take
\begin{equation}\label{eq:matrixmeasurecongruence}
 \dd\widehat\Sigma(x)
   =S(x)\,\dd\Sigma(x)\,S(x)^{\T}
       +\sum_\xi C_\xi\,\delta_\xi,
 \qquad L_1(\xi)C_\xi=0.
\end{equation}
Multiplication on the two sides gives the measure identity
$L_1\dd\widehat\Sigma L_1^{\T}=L_0\dd\Sigma L_0^{\T}$.
The measure remains positive definite on polynomial rows, since its
rational part is positive definite and each added term is nonnegative.

\begin{corollary}[Monic matrix polynomials from the reduced matrix]
\label{cor:matrixorthogonal}
Under these assumptions, take
$E=X_n=(\mathbb R[x]_{<n})^{1\times m}$ and
$X_n^+=(\mathbb R[x]_{\le n})^{1\times m}$, for $n\ge1$.
Let the spaces and matrices be those of
Theorem~\ref{thm:positiveprojection}, with common reduced order
$t_n$. For $a\in\{1,\ldots,m\}$, the leading row is
$v_a=x^ne_a^{\T}$. The corresponding projected rows $(v_a)_0$
form the $m\times m$ polynomial matrix $V_{0,n}$, while
their pairings with the complementary rows form the
$m\times t_n$ matrix
\[
 B_n\coloneq[b(v_a)_j]_{a\in\{1,\ldots,m\},\ j\in\{1,\ldots,t_n\}}.
\]
Write $W_n\coloneq\boldsymbol w$ and $T_n\coloneq T_\nu$ for the
polynomial rows after projection and their Gram matrix. Then
\begin{equation}\label{eq:matrixmonicreduced}
 \widehat P_n=V_{0,n}-B_nT_n^{-1}W_n
\end{equation}
is the monic matrix orthogonal polynomial for $\widehat\Sigma$.
Its norm and the kernel for degrees less than $n$ are
\begin{align}
 \widehat H_n
   &=\langle V_{0,n},V_{0,n}\rangle_{\widehat\Sigma}
                      -B_nT_n^{-1}B_n^{\T}>0,
       \label{eq:matrixnormreduced}\\
 \widehat K_n(x,y)
   &=\sum_{j=0}^{n-1}\widehat P_j(x)^{\T}
                      \widehat H_j^{-1}\widehat P_j(y).
       \label{eq:matrixCDreduced}
\end{align}
The latter is also given directly by \eqref{eq:positivekernel}.
Every entry of $\widehat P_n$ has the bordered determinant
\eqref{eq:positiveborder}, with denominator of order $t_n$ and
numerator of order $t_n+1$. The preliminary projection uses the
additional order-$t_n$ matrix $A_F$.
If the B\'ezout and degree hypotheses of
Corollary~\ref{cor:positiveorders} hold, then
\[
 t_n=m\min\{n,M\},\qquad
 M\coloneq\max\{\deg L_0,\deg L_1\}.
\]
In particular, these orders are independent of $n$ for $n\ge M$.
\end{corollary}

\begin{proof}
The rows of $H_j^{-1/2}P_j$, with $j\in\{0,\ldots,n-1\}$,
form the original orthonormal basis of $E$ needed in the
projection construction. All its entries are consequently determined
by the original matrix polynomials and norms, together with the
rational moments and masses. Applied to each leading row $v_a$,
Theorem~\ref{thm:positiveprojection} subtracts only elements of
$X_n$, so the assembled matrix retains the leading term $x^nI_m$.
Its row orthogonality to $X_n$ means that its matrix inner product
with every lower-degree matrix polynomial vanishes. The result is
therefore a monic matrix orthogonal polynomial, whose uniqueness
follows from uniqueness of each normalized row.

The matrix of pairings of $V_{0,n}$ with $W_n$ is $B_n$, and the
Gram matrix of $W_n$ is $T_n$. Expanding the inner product of
\eqref{eq:matrixmonicreduced} with itself gives
\[
 \begin{aligned}
 \widehat H_n={}&\langle V_{0,n},V_{0,n}\rangle
 -B_nT_n^{-1}B_n^{\T}-B_nT_n^{-1}B_n^{\T}\\
 &\quad+B_nT_n^{-1}T_nT_n^{-1}B_n^{\T}
 =\langle V_{0,n},V_{0,n}\rangle-B_nT_n^{-1}B_n^{\T},
 \end{aligned}
\]
which is \eqref{eq:matrixnormreduced}. Its positivity follows because any
nonzero constant combination of the rows of $\widehat P_n$ has
a nonzero leading row. Orthogonality of the matrix blocks gives
\eqref{eq:matrixCDreduced}, whereas their rows form an orthogonal
basis of $X_n$, identifying it with \eqref{eq:positivekernel}.
The componentwise bordered formulas give the determinant statement.
For the order count, Corollary~\ref{cor:positiveorders} applies with
$H_A(j)=j+1$ for $j\ge0$, yielding the stated expression.
The construction begins with $\widehat P_0=I_m$.
\end{proof}

The argument extends to the Hermitian inner product
$\int P\,\dd\Sigma\,Q^*$ because complex orthogonal projection
has the same block decomposition. Positive Hermitian Gram matrices
therefore give the Schur complement \eqref{eq:matrixmonicreduced},
with $B_n^*$ as the last factor in \eqref{eq:matrixnormreduced}
and $\widehat P_j(x)^*$ as the first factor in
\eqref{eq:matrixCDreduced}. Here a star denotes conjugate transpose;
the measure congruence \eqref{eq:matrixmeasurecongruence} uses the
same adjoints and masses satisfying $C_\xi=C_\xi^*\ge0$.

\subsubsection{An explicit connection with the Hermite matrix weights}

The factorized matrix weights studied by Dur\'an and Gr\"unbaum
provide examples of this congruence construction. Their families of
matrix orthogonal polynomials satisfy second-order differential
equations \cite{DG2004,DG2007}. A weight $w(x)T(x)T(x)^*$
has the form \eqref{eq:matrixmeasurecongruence} whenever $T$ is
polynomial or rational and the stated moments exist, with
$\dd\Sigma=wI_m\dd x$, $S=T$ and zero mass. The broader
differential constructions also include nonrational $T$.

In particular, the two Hermite weights of
\cite{DGStructural} have the explicit factorization
\begin{equation}\label{eq:DGHermiteweights}
 \mathcal H_{a,r}(x)
   =\mathrm e^{-x^2}
     \begin{bmatrix}1+|a|^2x^{2r}&ax^r\\\overline a x^r&1\end{bmatrix}
   =T_r(x)\bigl(\mathrm e^{-x^2}I_2\bigr)T_r(x)^*,
 \quad T_r(x)\coloneq\begin{bmatrix}1&ax^r\\0&1\end{bmatrix},
 \quad r\in\{1,2\}.
\end{equation}
They correspond to $L_1=I_2$, $L_0=T_r$ and $\nu=0$, so
these established matrix families are specializations of the
bilateral polynomial transformation. Their triangular factors also
make it possible to compute the reduced order directly from the
common space. For $a\in\mathbb R\setminus\{0\}$ and
$E=X_n$ as above, the image of a row $f=[f_1\ f_2]$ is
\[
 fT_r=[f_1\ ax^rf_1+f_2].
\]
Since both components of $f$ have degree less than $n$, the second
component of $fT_r$ has degree less than $n$ exactly when
$\deg f_1<n-r$. A negative bound means the zero polynomial.
Consequently
\begin{equation}\label{eq:DGorder}
 \dim Z=n+\max\{n-r,0\},\qquad
 t_n=2n-\dim Z=\min\{n,r\}.
\end{equation}
The reduced matrices consequently have stable order one in the first
family and two in the second. These orders are obtained from the
common-space calculation because the highest-degree coefficient of
$T_r$ is singular. The count $mM$ in
Corollary~\ref{cor:matrixorthogonal} is subject to its additional
degree hypotheses.

For the first family, scalar Hermite orthogonality also gives a
closed expression for the monic matrix polynomial. Let $h_n$
denote the monic scalar Hermite polynomials, with normalization
\[
 \int_{\mathbb R}h_n(x)h_j(x)\mathrm e^{-x^2}\dd x
       =\gamma_n\delta_{n,j},\qquad
 \gamma_n\coloneq\frac{\sqrt\pi\,n!}{2^n},\qquad
 D_n\coloneq1+\frac{a^2n}{2},
\]
and set $h_j=0$ for $j<0$. The result is
\begin{equation}\label{eq:DGmonicHermite}
 \widehat P_n(x)=
 \begin{bmatrix}
 h_n(x)&-\dfrac{an}{2}h_{n-1}(x)\\[2pt]
 -\dfrac{an}{2D_n}h_{n-1}(x)&
 h_n(x)+\dfrac{a^2n(n-1)}{4D_n}h_{n-2}(x)
 \end{bmatrix},\qquad
 \widehat H_n=\gamma_n\diag(D_{n+1},D_n^{-1}).
\end{equation}
Multiplication by the triangular factor reduces verification of the
expression to scalar Hermite orthogonality. In fact, the recurrence
$xh_n=h_{n+1}+(n/2)h_{n-1}$ gives
\[
 \widehat P_nT_1=
 \begin{bmatrix}
 h_n&a h_{n+1}\\
 -\dfrac{an}{2D_n}h_{n-1}&D_n^{-1}h_n
 \end{bmatrix}.
\]
A lower test $[u\ v]$ has image $[u\ axu+v]$, which reduces
the matrix orthogonality conditions to scalar Hermite pairings.
For the first row these vanish by the bounds $\deg u<n$ and
$\deg(xu),\deg v\le n$. In the second row, the only possible
contribution comes from the coefficient of $h_{n-1}$ in $u$;
its two scalar pairings cancel because
$-an\gamma_{n-1}/(2D_n)+a\gamma_n/D_n=0$.
The same scalar inner product gives squared norm
$\gamma_n+a^2\gamma_{n+1}=\gamma_nD_{n+1}$ for the first
image row. For the second row it gives
\[
 \frac{a^2n^2\gamma_{n-1}}{4D_n^2}
       +\frac{\gamma_n}{D_n^2}
 =\frac{\gamma_n(1+a^2n/2)}{D_n^2}
 =\frac{\gamma_n}{D_n},
\]
using $\gamma_n=n\gamma_{n-1}/2$ when $n\ge1$.
Their cross pairing is zero: it pairs $h_n$ with $h_{n-1}$
in the first component and $h_{n+1}$ with $h_n$ in the second.
At $n=0$, the displayed image matrix gives the same norm without
using the recurrence for its previous index. The diagonal and cross
pairings therefore establish the entire norm matrix. Since the
leading coefficient in \eqref{eq:DGmonicHermite} is $I_2$,
uniqueness of the monic orthogonal polynomial identifies this
expression with \eqref{eq:matrixmonicreduced}.

\subsubsection{The additional differential structure}

The Hermite matrix families also possess the differential property
studied in \cite{DG2004,DGStructural,DG2007}: there are fixed
polynomial matrices $F_2,F_1,F_0$ satisfying
\[
 \widehat P_n''F_2+\widehat P_n'F_1+\widehat P_nF_0
       =\Lambda_n\widehat P_n
\]
for every $n$, with $\deg F_2\le2$, $\deg F_1\le1$ and
$F_0$ constant. Together with the matrix three-term recurrence,
this is the differential and discrete spectral structure of those
families. Its preservation under a further rational congruence
requires compatibility of the transformed weight with such an
operator. The rational congruence determines the orthogonal
polynomials, norms and kernel through the formulas above, while
the differential equation requires this additional compatibility.
Coordinate multiplication provides the corresponding discrete
relations: Section~\ref{sec:recurrence} expresses the change of
recurrence matrices through four-factor intertwining identities.

\subsection{A positive bivariate example on three nodes}
\label{sec:atomicbivariate}

The congruence argument also applies to positive bivariate pairings,
where different support points can interact. A finite example is
obtained from the nodes $x_1=-1$, $x_2=0$, $x_3=1$ and the
matrix blocks
\[
 U(x)\coloneq\begin{bmatrix}1&x\\0&1\end{bmatrix},\qquad
 K_{i,j}\coloneq\delta_{i,j}I_2+U(x_i)U(x_j)^{\T}.
\]
For complex polynomial rows $f,h$ of length two, define
\[
 \langle f,h\rangle\coloneq
 \sum_{i,j=1}^3 f(x_i)K_{i,j}\overline{h(x_j)}^{\T}.
\]
The $6\times6$ matrix of this pairing in evaluation coordinates is
\[
 \mathsf K\coloneq[K_{i,j}]_{i,j\in\{1,\ldots,3\}}=I_6+\mathsf U\mathsf U^{\T},
 \qquad
 \mathsf U\coloneq\begin{bmatrix}U(x_1)\\U(x_2)\\U(x_3)\end{bmatrix}.
\]
For a nonzero real or complex column $v$,
$v^*\mathsf Kv=\|v\|^2+\|\mathsf U^{\T}v\|^2>0$,
so it is positive definite. Evaluation at the three distinct nodes is
an isomorphism from rows of component degrees below three to
$\CC^6$, by scalar polynomial interpolation in each component.
This interpolation isomorphism makes the polynomial pairing positive
definite on the stated space. Its determinant is also explicit.
The identity $U(x_i)^{\T}U(x_i)=[1\ x_i;\ x_i\ 1+x_i^2]$,
summed over $x_i\in\{-1,0,1\}$, yields
$\mathsf U^{\T}\mathsf U=\diag(3,5)$; the determinant identity
$\det(I+AB)=\det(I+BA)$ then gives
$\det\mathsf K=4\cdot6=24$. This value is the determinant
in evaluation coordinates. If $\mathsf V$ has as columns the
evaluations of a real monomial basis, the change to that basis
gives Gram matrix $\mathsf V^{\T}\mathsf K\mathsf V$ and
hence determinant $(\det\mathsf V)^2\det\mathsf K$.

A rational congruence of this pairing is defined by
\[
 S(x)=\begin{bmatrix}(x+2)^{-1}&x\\0&x+3\end{bmatrix}.
\]
It is obtained from the polynomial factors
\[
 L_1=\begin{bmatrix}x+2&0\\0&1\end{bmatrix},\qquad
 L_0=\begin{bmatrix}1&x(x+2)\\0&x+3\end{bmatrix},
 \qquad R_i=L_i^{\T},
\]
since $S=L_1^{-1}L_0$ and $S_R=S^{\T}$. These factors make
the example a massless solution of the four-factor equation.
Regularity and invertibility at all three nodes ensure that the
transformed blocks $\widehat K_{i,j}=S(x_i)K_{i,j}S(x_j)^{\T}$
again define a positive pairing. In evaluation coordinates,
$\mathsf D\coloneq\diag(S(x_1),S(x_2),S(x_3))$ gives
$\widehat{\mathsf K}=\mathsf D\mathsf K\mathsf D^{\T}>0$.
The determinant transforms by the square of the determinant of
this change: $\det\mathsf D=2\cdot(3/2)\cdot(4/3)=4$
therefore gives $\det\widehat{\mathsf K}=4^2\cdot24=384$.
The two nonzero Gram determinants establish normality on the
specified polynomial space. At the same time, for $i\ne j$ the
blocks $K_{i,j}$ are nonzero and stay nonzero under the invertible
congruence. Positivity and normality in this example therefore
coexist with interactions between distinct nodes, as permitted by
the bivariate pairing convention.

\subsection{Positive systems on real curves}\label{sec:positivecurves}

A positive density on a real arc gives congruences that are normal
on every finite row space. Let the support be an integral real
curve with a nontrivial real smooth arc, and let a compact nontrivial
subarc avoid $\det L_0\det L_1=0$. For a continuous positive
definite symmetric density $W$ on that subarc, both $W$ and
$SWS^{\T}$, with any additional mass \eqref{eq:positivemass},
define positive definite pairings on every finite row space and
its transpose. In particular, for $f\ne0$,
$\int fSWS^{\T}f^{\T}\,\dd t+\sum_P f(P)C_Pf(P)^{\T}>0$.
Since $S$ is invertible over the function field, $fS\ne0$.
After clearing its denominators, at least one component is a
nonzero regular function on the integral curve. Such a function
has only finitely many zeros on this compact arc unless it is
identically zero: its zero locus is a proper zero-dimensional
closed subset of the curve. Hence $fS$ is nonzero on a
subinterval. The positive density makes the integral strictly
positive there, and the atomic term is nonnegative. The same
argument with $S=I$ proves the assertion for $W$.

For the elliptic curve $y^2=x^3-x+1$, an explicit instance is
given by
$L_1=\bigl[\begin{smallmatrix}x&1\\0&1\end{smallmatrix}\bigr]$ and
$L_0=\bigl[\begin{smallmatrix}1&y\\1&x+1\end{smallmatrix}\bigr]$.
These factors determine $S=L_1^{-1}L_0$. With $W=I_2$ on the
arc $3\le x\le4$, $y=\sqrt{x^3-x+1}$, the required
invertibility follows from $\det L_1=x$ and $\det L_0=x+1-y$:
the identity $y^2-(x+1)^2=x(x^2-x-3)>0$ excludes zeros of
the second determinant as well. The inverse product has the explicit
form
\[
 S=\begin{bmatrix}
 0&\dfrac{y-x-1}{x}\\[3pt]
 1&x+1
 \end{bmatrix}.
\]
The original measure has density $I_2\,\dd x$ on the selected
arc, and the transformed density $SS^{\T}\,\dd x$
is therefore positive and gives normality on every finite row space
paired with its transpose. The associated right factors are
$R_i=L_i^{\T}$, and either zero mass or any mass satisfying
\eqref{eq:positivemass} preserves this conclusion. The integral
and the added point masses are balanced, since scalar multiplication
can be moved between their two arguments. This supplies the
positive matrix measure used in Subsection~\ref{sec:ellipticdarboux}:
the factors and arc give positivity and normality directly, and
balance permits the two coordinate recurrence identities there.

\section{Coordinate multiplication and bilateral Darboux identities}
\label{sec:recurrence}

The polynomial formulas obtained above also determine recurrence
relations: multiplication by a coordinate produces a polynomial whose
expansion in the orthogonal basis defines the corresponding multiplication
matrix. Since the coordinates on an algebraic support satisfy polynomial
relations, their multiplication matrices satisfy those relations as well.
The bilateral transformation changes these matrices through the
four connections associated with $L_0,L_1,R_0,R_1$. The connections
are determined by the factors and the polynomial bases, so a single
set of them relates multiplication by every coordinate.

For scalar Jacobi matrices, Darboux transformations express spectral
modifications through triangular factorizations and reversal of the
factors; see \cite{BuenoMarcellan2004}. For scalar bilinear
functionals, Bueno and Marcell\'an~\cite{BuenoMarcellan2006} introduce
polynomial factors on the two arguments and relate the associated
Hessenberg matrices through LU and QR factorizations in the
one-sided and symmetric cases developed there. These results provide
the scalar functional precedent for changing the two arguments.
Recurrence formulations for matrix and mixed-type orthogonality are developed in
\cite{Matrix2016,MixedChristoffel,MixedGeronimus}.
Here the term \emph{bilateral Darboux identities} denotes the
four-factor intertwining relations and the equality between the
two products of connection matrices proved below. Their formulation
retains four independent matrix factors on the two modules over the
coordinate algebra, including numerator and denominator factors.
The same connections act simultaneously on all coordinate
multiplication matrices and preserve the equations of the support.
The theorem is stated in terms of finite
polynomial expansions; triangular factorization and
reversal impose additional conditions on the connection matrices.

Write $\beta\coloneq\beta_\mu$ and $\widehat\beta\coloneq\beta_{\widehat\mu}$.
Throughout this section both pairings are bilinear and balanced:
\[
 \beta(af,g)=\beta(f,ag),\qquad
 \widehat\beta(af,g)=\widehat\beta(f,ag),\qquad a\in A.
\]
Balance identifies the row and column matrices of scalar
multiplication and supplies the lower bound on their nonzero
blocks. It is also the condition that makes these matrices symmetric
for positive forms in orthonormal bases.
Take nested spaces $E_n\subset A^{1\times q}$ and
$H_n\subset A^{p\times1}$ of equal dimensions, exhausting the two
modules, with both pairings normal on every $E_n\times H_n$.
Put $E_{-1}=H_{-1}=0$ and
$r_n\coloneq\dim E_n-\dim E_{n-1}>0$.
Choose dual polynomial blocks $B_n,A_n$ spanning successive levels,
with $B_n$ of size $r_n\times q$, $A_n$ of size $p\times r_n$ and
$\beta(B_n,A_n)=I_{r_n}$ and $\beta(B_n,A_m)=0$ for $n\ne m$.
The analogous blocks for $\widehat\beta$ are denoted with hats.
The vertical arrays $\cB,\widehat\cB$ collect the row blocks, whereas
the horizontal arrays $\cA,\widehat\cA$ collect the column blocks.
Thus the rows of $B_0,B_1,\ldots$ occur consecutively in the vertical
array, and the columns of the $A_n$ occur consecutively in the horizontal
array. These arrays record polynomial bases constructed, one block at a
time, from finite Gram matrices. Exhaustion places every polynomial
in a finite level, so each expansion has only finitely many basis
elements. The arrays therefore act algebraically on the polynomial
modules; no bounded-operator interpretation on a completed function
space is required.

For each scalar polynomial class $a\in A$, assume finite nonnegative
shifts satisfying
\[
 aE_n\subset E_{n+d_L(a)},\qquad
 aH_n\subset H_{n+d_R(a)}\qquad(n\ge0).
\]
Define
\[
 J(a)\coloneq[\beta(aB_n,A_m)]_{n,m\ge0},\qquad
 \widehat J(a)\coloneq[\widehat\beta(a\widehat B_n,\widehat A_m)]_{n,m\ge0},
\]
and define the four connection matrices by the polynomial expansions
\begin{equation}\label{eq:fourconnections}
 \cB L_0=U_0\cB,\qquad \cB L_1=U_1\widehat\cB,
 \qquad R_0\cA=\cA V_0,\qquad R_1\cA=\widehat\cA V_1.
\end{equation}
Their blocks are computed directly from the dual polynomials:
\[
 \begin{aligned}
 U_0[n,m]&=\beta(B_nL_0,A_m),&
 U_1[n,m]&=\widehat\beta(B_nL_1,\widehat A_m),\\
 V_0[n,m]&=\beta(B_n,R_0A_m),&
 V_1[n,m]&=\widehat\beta(\widehat B_n,R_1A_m).
 \end{aligned}
\]
Indeed, the $n$th block row of $U_1$ is the coefficient row in
$B_nL_1=\sum_mU_1[n,m]\widehat B_m$, so pairing with $\widehat A_m$
recovers the indicated block by duality. Applying this coefficient
extraction to each of the other three expansions gives the remaining
formulas. Since each polynomial has a finite expansion, every block row
of $U_i$ and every block column of $V_i$ has finitely many nonzero blocks.
The choice of bases in these definitions expands the numerators in the
original system and the denominators in the transformed system, as
required by the two sides of the functional equation.

\begin{theorem}[Multiplication matrices and the four-factor identities]
\label{thm:darboux}\label{prop:recurrence}\label{sec:symmetricdarboux}
The dual arrays satisfy
\[
 a\cB=J(a)\cB,\qquad a\cA=\cA J(a),
\]
and the hatted counterparts. Their blocks vanish unless
$n-d_R(a)\le m\le n+d_L(a)$.
The maps $a\mapsto J(a)$ and $a\mapsto\widehat J(a)$ are faithful
unital algebra representations. Thus, for
$A=\CC[x_1,\ldots,x_d]/I$,
\[
 f(J(x_1),\ldots,J(x_d))=0\quad\text{if and only if}\quad f\in I,
\]
and likewise for the transformed matrices.

The functional equation \eqref{eq:four} holds if and only if
\begin{equation}\label{eq:factorization}
 U_1V_1=U_0V_0.
\end{equation}
For every $a\in A$ the connections satisfy
\begin{align}
 J(a)U_1&=U_1\widehat J(a),&
 \widehat J(a)V_1&=V_1J(a),\label{eq:intertwine}\\
 J(a)U_0&=U_0J(a),&
 J(a)V_0&=V_0J(a).\label{eq:commuteconnections}
\end{align}
Consequently, under the functional equation,
\begin{equation}\label{eq:darboux}
 U_1\widehat J(a)V_1=U_0J(a)V_0.
\end{equation}

For real factors with $q=p$, $R_i=L_i^{\T}$, $H_n=E_n^{\T}$
and real symmetric positive definite forms on real spaces at every
level, real orthonormal blocks
make $J(a),\widehat J(a)$ symmetric for real $a$ and give
$V_i=U_i^{\T}$. Under the functional equation,
\begin{equation}\label{eq:symmetricdarboux}
 U_1U_1^{\T}=U_0U_0^{\T},\qquad
 U_1\widehat J(a)U_1^{\T}=U_0J(a)U_0^{\T}.
\end{equation}
All products used here are finite entry by entry.
\end{theorem}

\begin{proof}
The dual blocks are constructed inductively by eliminating their
pairings with the preceding levels. Denote the preceding row blocks
by $B_{<n}$ and their dual columns by $A_{<n}$, so that
$\beta(B_{<n},A_{<n})=I$. If rows $b^+$ complete the row basis to $E_n$
and columns $a^+$ complete the column basis to $H_n$, the additional
Gram blocks and their Schur complement are
\[
 C\coloneq\beta(b^+,A_{<n}),\qquad D\coloneq\beta(B_{<n},a^+),
 \qquad G\coloneq\beta(b^+,a^+),\qquad S\coloneq G-CD.
\]
Subtracting the components detected by the preceding dual bases gives
\[
 B_n\coloneq b^+-CB_{<n},\qquad A_n^0\coloneq a^+-A_{<n}D.
\]
These subtractions ensure that $\beta(B_n,A_{<n})=C-C=0$ and
$\beta(B_{<n},A_n^0)=D-D=0$, while expansion of the four terms in
the pairing between the new blocks gives
\[
 \beta(B_n,A_n^0)=G-CD-CD+CD=G-CD=S.
\]
The same elimination applied to the full Gram matrix reads
\[
 \begin{bmatrix}I&0\\-C&I\end{bmatrix}
 \begin{bmatrix}I&D\\C&G\end{bmatrix}
 \begin{bmatrix}I&-D\\0&I\end{bmatrix}
 =\begin{bmatrix}I&0\\0&S\end{bmatrix}.
\]
Since the full Gram matrix is invertible by normality, its Schur
complement $S$ is invertible, and the normalization
$A_n\coloneq A_n^0S^{-1}$ makes the new blocks dual. At the initial
level $n=0$, the lower arrays are empty and $S=G$, so the same argument
starts the induction. Applying it to each pairing constructs the two
dual polynomial systems.

The multiplication matrix can now be read from these dual bases.
The expansion of $aB_n$ has no block above $n+d_L(a)$, and pairing
with $A_m$ gives its coefficient $J(a)[n,m]$. The other side of the
band follows from balance, since
$\beta(aB_n,A_m)=\beta(B_n,aA_m)=0$ when $m+d_R(a)<n$.
Moreover, in the column expansion $aA_m=\sum_nA_nC_{n,m}$, pairing
with $B_n$ gives $C_{n,m}=\beta(B_n,aA_m)=J(a)[n,m]$.
Thus the same matrix, with the two band bounds just established,
represents multiplication in both arrays. Its compatibility with
products follows from the successive expansions
\[
 ab\cB=a\bigl(J(b)\cB\bigr)=J(b)J(a)\cB.
\]
Interchanging $a,b$ and using uniqueness of polynomial expansions
therefore gives $J(ab)=J(a)J(b)=J(b)J(a)$; uniqueness also gives
additivity and $J(1)=I$. To prove faithfulness, suppose that $J(a)=0$.
Exhaustion then implies that multiplication by $a$ annihilates every
polynomial row, including each constant coordinate row, and hence
$a=0$ in $A$. This argument applies even when the algebra has zero
divisors, and identifies the polynomial relations among the coordinate
matrices with the defining relations of the support. The same reasoning
applies to the transformed bases.

For the denominator expansions, duality gives the explicit block
calculation
\begin{align*}
 \widehat\beta(B_nL_1,R_1A_m)
 &=\sum_{j,l}U_1[n,j]\widehat\beta(\widehat B_j,\widehat A_l)V_1[l,m]\\
 &=\sum_jU_1[n,j]V_1[j,m]=(U_1V_1)[n,m].
\end{align*}
Both sums are finite because both polynomial arguments have finite
expansions. With the original dual bases, the corresponding calculation
gives $(U_0V_0)[n,m]$, so equality of the two products is precisely
\[
 \widehat\beta(B_nL_1,R_1A_m)=\beta(B_nL_0,R_0A_m)
 \qquad(n,m\ge0).
\]
Since the bases exhaust the modules, this equality is equivalent to
\eqref{eq:pairfour}. By balance, testing that identity on constant
coordinate rows multiplied by arbitrary elements of $A$ recovers
\eqref{eq:four} on every matrix test, proving the asserted equivalence.

Each intertwining identity comes from two expansions of the same
polynomials. Keeping the order of all matrix factors, these expansions are
\begin{align*}
 a(\cB L_1)&=J(a)U_1\widehat\cB
             =U_1\widehat J(a)\widehat\cB,\\
 a(R_1\cA)&=\widehat\cA\widehat J(a)V_1
             =\widehat\cA V_1J(a),\\
 a(\cB L_0)&=J(a)U_0\cB=U_0J(a)\cB,\\
 a(R_0\cA)&=\cA J(a)V_0=\cA V_0J(a).
\end{align*}
These equalities use only the fact that $a$, being a scalar in $A$,
commutes with the polynomial matrix entries. Uniqueness of each row
or column expansion proves the corresponding intertwining or
commutation identity. Together with \eqref{eq:factorization}, they give
\[
 U_1\widehat J(a)V_1
 =J(a)U_1V_1=J(a)U_0V_0=U_0J(a)V_0.
\]
The finite rows of $U_i$, finite columns of $V_i$ and finite bands
of $J(a),\widehat J(a)$ justify these sums and regroupings entry by
entry, without taking the inverse of a semi-infinite matrix.

Finally, in the positive real case take $A_n=B_n^{\T}$ and
$\widehat A_n=\widehat B_n^{\T}$ using orthonormal blocks.
Balance and symmetry give
$\beta(aB_n,B_m^{\T})=\beta(aB_m,B_n^{\T})^{\T}$.
This proves symmetry of the multiplication matrices. In addition,
transposition turns the left expansions in \eqref{eq:fourconnections}
into the right expansions, whose uniqueness gives $V_i=U_i^{\T}$.
Substitution in the identities already proved yields
\eqref{eq:symmetricdarboux}.
\end{proof}

On the polynomial line $A=\CC[z]$, the band bounds have a direct
interpretation in terms of the numbers of components. Order the row
monomials as
\[
e_1^{\T},\ldots,e_q^{\T},ze_1^{\T},\ldots,ze_q^{\T},z^2e_1^{\T},\ldots,
\]
where $e_i$ are the coordinate columns of $\CC^q$, and order the
column monomials in the same way with the $p$ coordinate columns
of $\CC^p$. Let $E_n,H_n$ consist of the first $n+1$ monomials
in their respective lists. Assume, as above, that both pairings
are normal on every such pair of spaces. Multiplication by $z$
shifts the row list by $q$ positions and the column list by $p$
positions. Hence both recurrence matrices vanish outside
$n-p\le m\le n+q$.

\subsection{The two coordinate identities on an elliptic curve}
\label{sec:ellipticdarboux}

Two choices of spaces on the same curve exhibit the effects of
different row and column component counts and positivity. The rectangular calculation
assumes the hypotheses of Theorem~\ref{thm:darboux} and determines
the resulting bands and identities. The square calculation uses
the explicit positive matrix measure of
Subsection~\ref{sec:positivecurves}, whose positivity supplies
normality at every level.

The coordinate algebra is
$A=\CC[x,y]/(y^2-x^3+x-1)$. Its smooth projective completion
has one point $O=[0:1:0]$ at infinity in projective coordinates
$[X:Y:Z]$, where the equation is $Y^2Z=X^3-XZ^2+Z^3$.
In the chart $Y=1$, put $u\coloneq X/Y$ and $v\coloneq Z/Y$. Near $O$ the equation is
\[
 v(1+uv-v^2)=u^3.
\]
The factor in parentheses is invertible at $O$, where it equals one.
Moreover, the derivative of the left side minus $u^3$ with respect to
$v$ is one at $O$, so $u$ is a local parameter and the equation forces
$v$ to have order three. It follows that $x=u/v$ has a pole of order
two and $y=1/v$ a pole of order three. These orders determine the
spaces $L(mO)$ of rational functions regular away from $O$ and with a
pole of order at most $m$ there: a basis of $L(mO)$ consists of the
monomials $x^a$ with $2a\le m$ and $x^ay$ with $2a+3\le m$, where
$a\ge0$. Indeed, every element of $A$ has a unique expression
$p(x)+yq(x)$, so the monomials are independent; their distinct pole
orders then characterize the bound without cancellation.
For $m\ge3$, the two lists contain respectively
$\lfloor m/2\rfloor+1$ and $\lfloor(m-3)/2\rfloor+1$ monomials.
Their sum is $m$, as is seen by taking $m$ even and odd separately.
For $m=1,2$ the bases are $1$ and $1,x$. Thus
$\dim L(mO)=m$ for every $m\ge1$.

For the rectangular case, take two row components and three column
components with spaces
\[
 X_\ell=L(3\ell O)^{1\times2},\qquad
 Y_\ell=L(2\ell O)^{3\times1},\qquad \ell\ge1.
\]
The dimension calculation above shows that both spaces have dimension
$6\ell$; increasing the pole bound makes them nested and exhaustive
in the polynomial row and column modules. Set $E_n=X_{n+1}$ and
$H_n=Y_{n+1}$, so the block index is $n=\ell-1$ and every block,
including the initial one, has size six. Suppose the two balanced
pairings and the four polynomial matrix factors satisfy the hypotheses
and the functional equation of Theorem~\ref{thm:darboux} on these
spaces. In particular, both pairings are normal at every level.

The shifts in the recurrence follow by comparing the poles of the
coordinates with the increments in the chosen filtrations. Multiplication
by $x$ adds at most two to a pole order, whereas multiplication by $y$
adds at most three. One level on the left allows three additional
orders, and hence accommodates either coordinate. On the right one
level allows two additional orders, so $x$ requires one level and $y$
requires two. Consequently,
\[
 xX_\ell\subset X_{\ell+1},\qquad xY_\ell\subset Y_{\ell+1},
 \qquad
 yX_\ell\subset X_{\ell+1},\qquad yY_\ell\subset Y_{\ell+2}.
\]
The theorem now gives block ranges $n-1\le m\le n+1$ for
$J(x),\widehat J(x)$ and $n-2\le m\le n+1$ for
$J(y),\widehat J(y)$. Since these matrices represent multiplication
in $A$, the two coordinate matrices commute and satisfy
\[
 J(y)^2=J(x)^3-J(x)+I,\qquad
 \widehat J(y)^2=\widehat J(x)^3-\widehat J(x)+I.
\]
The four connections in \eqref{eq:fourconnections} are defined by the
matrix factors and the polynomial bases; their definition does not
involve a choice of coordinate. Thus the same $U_0,U_1,V_0,V_1$
satisfy all four intertwining and commutation identities
\eqref{eq:intertwine}--\eqref{eq:commuteconnections} for both $x$ and
$y$. In particular,
\[
 U_1\widehat J(x)V_1=U_0J(x)V_0,\qquad
 U_1\widehat J(y)V_1=U_0J(y)V_0.
\]
The bilateral transformation therefore acts simultaneously on the
two multiplication operators and respects their elliptic relation.

The positive realization is square: both arguments have two
components and are paired by the matrix measures in
Subsection~\ref{sec:positivecurves}. Use
$X_\ell=L(3\ell O)^{1\times2}$ and $Y_\ell=X_\ell^{\T}$.
The densities and point masses specified there are balanced, and
positivity guarantees normality of both pairings at every level. The
common pole bound increases by three on each side. Since either
coordinate increases the pole order by at most three, all four
multiplication matrices have block range $n-1\le m\le n+1$, again
with blocks of size six. Choosing real orthonormal blocks also makes
both coordinate matrices symmetric in each system, so the same
connections have $V_i=U_i^{\T}$ and satisfy
\[
 U_1U_1^{\T}=U_0U_0^{\T},\qquad
 U_1\widehat J(a)U_1^{\T}=U_0J(a)U_0^{\T}
 \quad(a\in\{x,y\}).
\]
These are the positive four-factor Darboux identities for the two coordinates
of the explicit elliptic weight, with any of the admissible positive
masses specified there.

\section{Conclusions}

The bilateral functional equation leads to a determinant construction
for mixed-type orthogonality on algebraic supports. Starting from the original
biorthogonal data, the ordered matrix factors, the rational moments
required by division and a chosen homogeneous mass, the construction
expresses the transformed orthogonality conditions as a finite linear
system. Its image conditions preserve the relations of the support and
the prescribed component bounds. Consequently, the determinant of the
reduced system characterizes normality, while its bordered determinants
give the normalized polynomials and their reproducing kernel.

The functional equation also determines the full freedom in the
transformed pairing. For balanced masses with regular denominator
determinants, its exact sequence distinguishes the obstruction to
unilateral decomposition from the ambiguity of the summands. Common
components of the denominator zero sets can support mass directions
that require the joint bilateral equation.

The size of these determinants reflects both the transformation and
the polynomial spaces on which it acts. In the degree construction,
the Hilbert function counts the independent polynomial classes on the
support, and the finite mass rank accounts for the contribution of the
chosen mass. Polynomial rank computations then remove redundant
directions and give the size of the smaller system. These counts explain
why fixed perturbations, under compatible growth of the component degrees,
have bounded determinant orders on curves and at most linear orders on
surfaces. The formulas based on consecutive scalar degrees, matrix
division and reduced linear systems eliminate different polynomial
spaces; their respective counts specify the sizes of those
representations under the stated hypotheses.

To form these systems explicitly, the quotient conditions must themselves
be expressed in computable coordinates. On integral curves, dualizing
differentials and the geometric Cauchy kernel supply dual coordinates
for matrix remainders, including at singular points. On surfaces, an
additional scalar equation selects a finite quotient, and ordered
matrix divided differences reconstruct each quotient class as a finite
sum of its residue pairings. Applied to the respective quotients, these
coordinates describe the masses and the divisibility tests required by
the polynomial formulas. The explicit curve family also shows how this
description behaves when the mass support reaches a singular point.
The residue quotients remain nondegenerately paired and an explicit
minor proves that the matrix mass keeps its rank. Normality of the
finite Gram matrices, proved for the stated component bounds, then
gives regular limits of the normalized polynomials and kernels.

For positive symmetric congruences, the intersection of the original
space with its rationally transformed counterpart permits a further
reduction. The resulting two positive matrices determine both the
normalized polynomials and the kernel; the coprimality count and the
nodal matrix family identify situations in which positive masses,
including bivariate derivative masses, leave the reduced order unchanged.
This reduction applies to both balanced forms and the nonbalanced
Sobolev forms exhibited at the node, because the orthogonal
decomposition uses the positive forms on the specified finite spaces.
With equal component bounds on the line, this construction assembles
the polynomial rows into monic matrix orthogonal polynomials and yields
their matrix norms and Christoffel--Darboux kernel. In particular, the
triangular Hermite weights of Dur\'an and Gr\"unbaum illustrate how the
common space determines the small orders of the corresponding formulas.

The one-variable and Laurent reductions connect the general formulas
with their standard specializations, while the description of positive
masses identifies the admissible additions in the symmetric setting.
When both the original and transformed pairings are balanced,
their polynomial bases also give four-factor intertwining identities
for coordinate multiplication. These bilateral Darboux identities
use the same connection matrices for every coordinate and preserve
the defining relations of the support. The explicit positive
elliptic matrix measure supplies a square realization in which
both coordinate matrices are symmetric. Thus, under the balance
and normality hypotheses of the recurrence construction, the
polynomial and kernel formulas determine the transformed recurrence
matrices with the original order of the four matrix factors.

\section*{Funding}
This work was supported by the Spanish Ministry of Science, Innovation
and Universities and the Spanish State Research Agency [grant number
PID2024-155133NB-I00], co-funded by the European Union through the
European Regional Development Fund.

\section*{Data availability}
The article gives the matrix factors, masses, polynomial spaces,
coefficient formulas and finite linear systems used in its constructions.
The accompanying verification package supplies exact symbolic checks of
the examples and of the matrix identities. The general arguments are
proved in the text; tests at individual degrees supplement those proofs.
No external datasets are used.

\section*{Declaration of AI use}

During the preparation of this manuscript, the author used ChatGPT
and Codex (OpenAI) and Claude (Anthropic) as aids in exposition,
mathematical exploration, and the discussion and review of mathematical
arguments. These tools also assisted with computational checks,
the development of formal verification code, and the review of
notation and references.

\end{document}